%% file: drt.tex
\def\draft{n}

\documentclass[12pt, reqno, a4paper]{amsart} 
\usepackage{amsfonts,amssymb,amsmath,amsthm,amscd}
\usepackage{ahsymb}
\usepackage{epic,eepic}
\usepackage{fullpage,enumerate}
\numberwithin{equation}{section}

\input macros.tex

\input defs.tex

\begin{document}
\newdimen\captionwidth\captionwidth=\hsize

\input{page1}
\input{intro}
\input{pre}

\input{Aarrow}
\input{qtrqH}

\input{EK}

\input{refs}
\end{document}

%% file: macros.tex
\theoremstyle{plain}

\newtheorem{theorem}{Theorem}

\newtheorem{proposition}{Proposition}[section]
\newtheorem{lemma}[proposition]{Lemma}
\newtheorem{corollary}[proposition]{Corollary}

\newtheorem{conjecture}{Conjecture}

\theoremstyle{definition}
\newtheorem{definition}[proposition]{Definition}

\newtheorem{question}{Question}

\theoremstyle{remark}
\newtheorem{example}[proposition]{Example}
\newtheorem{exercise}[proposition]{Exercise}

\newtheorem{remark}[proposition]{Remark}

\def\printname#1{
  \if\draft y
    \smash{\makebox[0pt]{\hspace{-0.5in}
      \raisebox{8pt}{\tt\tiny #1}}}
  \fi
}

\newcommand{\mathmode}[1]{$#1$}

\newlength{\standardunitlength}
\catcode`\@=11
\long\def\@makecaption#1#2{%
    \vskip 10pt
    \setbox\@tempboxa\hbox{
      \small\sf{\bfcaptionfont #1. }\ignorespaces #2}%
    \ifdim \wd\@tempboxa >\captionwidth {%
        \rightskip=\@captionmargin\leftskip=\@captionmargin
        \unhbox\@tempboxa\par}%
      \else
        \hbox to\hsize{\hfil\box\@tempboxa\hfil}%
    \fi}
\font\bfcaptionfont=cmssbx10 scaled \magstephalf
\newdimen\@captionmargin\@captionmargin=2\parindent
\catcode`\@=12

\def\lbl#1{\label{#1}\printname{#1}}

\newenvironment{myitemize}{
	\begin{list}{$\bullet$}{\setlength{\leftmargin}{16pt}
        \setlength{\labelwidth}{12pt}
        \setlength{\labelsep}{4pt}}
}{
	\end{list}
}

\newlength{\globalparindent}

%% file: defs.tex
\def\a{\mathfrak a}
\def\A{\mathcal A}
\def\Aar{\overset{\ora}{\A}}
\def\a{\mathfrak a}
\def\B{\mathcal B}
\def\Bar{\overset{\ora}{\B}}
\def\CSg{S(\g)^{\g}}
\def\CUg{U(\g)^{\g}}
\def\CUgh{{\CUg\hb}}
\def\CoUg{U(\g)_{\g}}
\def\CoUgh{{\CoUg\hb}}
\def\catFT{\text{\bf FT}}
\def\catFor{\text{\bf For}}

\def\catNAT{\text{\bf PFT}_{\mathbf{0}}}
\def\catPFT{\text{\bf PFT}}
\def\coinv{\text{coinv}}
\def\b{\mathfrak b}
\def\bt{\boxtimes}
\def\C{\mathcal C}
\def\D{\mathcal D}
\def\Daar{\overset{_{\ra}}{\D^a}}
\def\Dar{\overset{\ora}{\D}}
\def\Dcar{\overset{_{\ra}}{\D^c}}
\def\da{\downarrow}
\def\ek{\text{EK}}
\def\eps{\epsilon}
\def\ew{\varnothing}
\def\F{\mathcal F}
\def\G{\mathcal G}
\def\g{\mathfrak g}
\def\H{\mathcal H}
\def\HS{\text{\bf HS}} 

\def\h{\mathfrak h}
\def\hb{[[\hbar]]}
\def\kz{\text{KZ}}
\def\KZg{\mathbf{KZ}_{\g}}

\def\k{\Bbbk}
\def\kh{\k\hb}
\def\I{\mathcal I}
\def\Iar{\overset{\ora}{\I}}
\def\inv{\text{inv}}
\def\J{\mathcal J}
\def\lla{\longleftarrow}
\def\lra{\longrightarrow}
\def\LHC{\raisebox{4pt}{$
  \printname{LeftHalfCirc}
  \setlength{\unitlength}{0.25\standardunitlength}
  \begin{array}{c}  \hspace{-1.7mm}
    \raisebox{-8pt}{\input draws/LeftHalfCirc.tex }
      \hspace{-1.9mm}
  \end{array}
$}}
\def\LF{{\mathcal L}^F}

\def\longLRa{\longleftrightarrow}
\def\M{\mathbf M}
\def\m{\mathbf m}
\def\n{\mathfrak n}
\def\ot{\otimes}
\def\ora{_{\ \ra}}
\def\Par{\overset{\ora}{\mathcal P}}
\def\Pcar{\overset{_{\ra}}{{\mathcal P}^c}}
\def\qg{\text{DJ}}
\def\RHC{\raisebox{4pt}{$
  \printname{RightHalfCirc}
  \setlength{\unitlength}{0.25\standardunitlength}
  \begin{array}{c}  \hspace{-1.7mm}
    \raisebox{-8pt}{\input draws/RightHalfCirc.tex }
      \hspace{-1.9mm}
  \end{array}
$}}
\def\RT{\mathbf{RT}}
\def\RTg{\RT_{\g}}
\def\ra{\to}

\def\T{\mathcal T}
\def\TP{\raisebox{1pt}{$
  \printname{Tadpole}
  \setlength{\unitlength}{0.25\standardunitlength}
  \begin{array}{c}  \hspace{-1.7mm}
    \raisebox{-8pt}{\input draws/Tadpole.tex }
      \hspace{-1.9mm}
  \end{array}
$}}
\def\TPD{\raisebox{1pt}{$
  \printname{TadpoleDown}
  \setlength{\unitlength}{0.25\standardunitlength}
  \begin{array}{c}  \hspace{-1.7mm}
    \raisebox{-8pt}{\input draws/TadpoleDown.tex }
      \hspace{-1.9mm}
  \end{array}
$}}
\def\TPU{\raisebox{1pt}{$
  \printname{TadpoleUp}
  \setlength{\unitlength}{0.25\standardunitlength}
  \begin{array}{c}  \hspace{-1.7mm}
    \raisebox{-8pt}{\input draws/TadpoleUp.tex }
      \hspace{-1.9mm}
  \end{array}
$}}
\def\Tar{\overset{\ora}{\T}}
\def\Targ{{\Tar_{\hspace{-3.5pt}\g}}}
\def\Targh{{\Tar_{\hspace{-3.5pt}\g}^{\raisebox{-4pt}{$_\hbar$}}}}
\def\Tarr{{\Tar_{\hspace{-3.5pt} r}}}
\def\Tartg{{\Tar_{\hspace{-3.5pt}\tg}}}
\def\Tartgh{{\Tar_{\hspace{-3.5pt}\tg}^{\raisebox{-4pt}{$_\hbar$}}}}
\def\Tgh{{\T_{\g}^{\hbar}}}

\def\tg{\tilde{\g}}
\def\twoheadra{\twoheadrightarrow}
\def\Ug{U(\g)}
\def\Ugh{{\Ug\hb}}
\def\Uhg{U_{\hbar}(\g)}
\def\Utg{U(\tg)}
\def\Utgh{{\Utg\hb}}
\def\ua{\uparrow}
\def\upchi{\raisebox{2pt}{$\chi$}}
\def\VM{\mathcal M}
\def\Vect{\text{Vect}}
\def\W{\mathcal W}
\def\wtJ{\widetilde{J}}
\def\Z{\mathbf Z}
\def\Zhg{\Z_{\hbar}(\g)}
\def\ZK{\Z_{\text{K}}}

\def\hs{\hspace}
\def\NATcp{\hs{-3pt}
  \setlength{\unitlength}{0.25\standardunitlength}
  \begin{array}{c}  \hspace{-1.7mm}
    \raisebox{-2pt}{\input draws/natCP.tex }
    \hspace{-1.9mm}
  \end{array}
\hs{-5pt}} 
\def\NATcn{\hs{-3pt}
  \setlength{\unitlength}{0.25\standardunitlength}
  \begin{array}{c}  \hspace{-1.7mm}
    \raisebox{-2pt}{\input draws/natCN.tex }
    \hspace{-1.9mm}
  \end{array}
\hs{-5pt}} 
\def\NATap{\hs{-3pt}
  \setlength{\unitlength}{0.25\standardunitlength}
  \begin{array}{c}  \hspace{-1.7mm}
    \raisebox{-2pt}{\input draws/natAP.tex }
    \hspace{-1.9mm}
  \end{array}
\hs{-5pt}} 
\def\NATan{\hs{-3pt}
  \setlength{\unitlength}{0.25\standardunitlength}
  \begin{array}{c}  \hspace{-1.7mm}
    \raisebox{-2pt}{\input draws/natAN.tex }
    \hspace{-1.9mm}
  \end{array}
\hs{-5pt}} 
\def\NATov{\hs{-3pt}
  \setlength{\unitlength}{0.25\standardunitlength}
  \begin{array}{c}  \hspace{-1.7mm}
    \raisebox{-2pt}{\input draws/natOV.tex }
    \hspace{-1.9mm}
  \end{array}
} 
\def\NATun{\hs{-3pt}
  \setlength{\unitlength}{0.25\standardunitlength}
  \begin{array}{c}  \hspace{-1.7mm}
    \raisebox{-2pt}{\input draws/natUN.tex }
    \hspace{-1.9mm}
  \end{array}
} 
\def\NATra{\hs{-3pt}
  \setlength{\unitlength}{0.25\standardunitlength}
  \begin{array}{c}  \hspace{-1.7mm}
    \raisebox{-2pt}{\input draws/natRA.tex }
    \hspace{-1.9mm}
  \end{array}
} 
\def\NATla{\hs{-3pt}
  \setlength{\unitlength}{0.25\standardunitlength}
  \begin{array}{c}  \hspace{-1.7mm}
    \raisebox{-2pt}{\input draws/natLA.tex }
    \hspace{-1.9mm}
  \end{array}
} 

\newcommand{\bbinom}[2]{
	\begin{bmatrix} #1 \\ #2 \end{bmatrix}
}

\newcommand{\lowupper}[2]{
	{#1}^{\raisebox{-4pt}{$_{#2}$}} }

\newcommand{\ad}{\operatorname{ad}}

\newcommand{\End}{\operatorname{End}}
\newcommand{\Hom}{\operatorname{Hnd}}
\newcommand{\img}{\operatorname{Im}}
\newcommand{\Ind}{\operatorname{Ind}}
\newcommand{\id}{\operatorname{id}}
\newcommand{\mor}{\operatorname{Mor}}
\newcommand{\op}{\operatorname{op}}
\newcommand{\R}{\operatorname{R}}
\newcommand{\Span}{\operatorname{span}}

\newcommand{\sym}{\operatorname{sym}}
\newcommand{\tr}{\operatorname{Tr}}
\newcommand{\udot}{{\mathaccent\cdot\cup}}

%% file: draws/LeftHalfCirc.tex
%
\begingroup\makeatletter\ifx\SetFigFont\undefined%
\gdef\SetFigFont#1#2#3#4#5{%
  \reset@font\fontsize{#1}{#2pt}%
  \fontfamily{#3}\fontseries{#4}\fontshape{#5}%
  \selectfont}%
\fi\endgroup%
{\renewcommand{\dashlinestretch}{30}
\begin{picture}(2235,819)(0,-10)
\put(1034.500,93.500){\arc{1401.000}{3.1409}{6.2839}}
\path(453,484)(661,578)
\path(452,486)(527,698)
\thicklines
\path(33,93)(2202,93)
\blacken\path(2022.000,33.000)(2202.000,93.000)(2022.000,153.000)(2022.000,33.000)
\end{picture}
}

%% file: draws/RightHalfCirc.tex
%
\begingroup\makeatletter\ifx\SetFigFont\undefined%
\gdef\SetFigFont#1#2#3#4#5{%
  \reset@font\fontsize{#1}{#2pt}%
  \fontfamily{#3}\fontseries{#4}\fontshape{#5}%
  \selectfont}%
\fi\endgroup%
{\renewcommand{\dashlinestretch}{30}
\begin{picture}(2235,820)(0,-10)
\put(1032.500,94.500){\arc{1401.000}{3.1409}{6.2839}}
\thicklines
\path(33,93)(2202,93)
\blacken\path(2022.000,33.000)(2202.000,93.000)(2022.000,153.000)(2022.000,33.000)
\thinlines
\path(1614,485)(1406,579)
\path(1615,487)(1540,699)
\end{picture}
}

%% file: draws/Tadpole.tex
%
\begingroup\makeatletter\ifx\SetFigFont\undefined%
\gdef\SetFigFont#1#2#3#4#5{%
  \reset@font\fontsize{#1}{#2pt}%
  \fontfamily{#3}\fontseries{#4}\fontshape{#5}%
  \selectfont}%
\fi\endgroup%
{\renewcommand{\dashlinestretch}{30}
\begin{picture}(1251,1553)(0,-10)
\path(638,1526)(450,1426)(652,1359)
\path(559,92)(559,692)
\put(560,1068){\ellipse{750}{750}}
\thicklines
\path(33,93)(1218,93)
\blacken\path(1038.000,33.000)(1218.000,93.000)(1038.000,153.000)(1038.000,33.000)
\end{picture}
}

%% file: draws/TadpoleDown.tex
%
\begingroup\makeatletter\ifx\SetFigFont\undefined%
\gdef\SetFigFont#1#2#3#4#5{%
  \reset@font\fontsize{#1}{#2pt}%
  \fontfamily{#3}\fontseries{#4}\fontshape{#5}%
  \selectfont}%
\fi\endgroup%
{\renewcommand{\dashlinestretch}{30}
\begin{picture}(1251,1553)(0,-10)
\path(638,1526)(450,1426)(652,1359)
\path(471,506)(560,310)(641,507)
\path(559,92)(559,692)
\put(560,1068){\ellipse{750}{750}}
\thicklines
\path(33,93)(1218,93)
\blacken\path(1038.000,33.000)(1218.000,93.000)(1038.000,153.000)(1038.000,33.000)
\end{picture}
}

%% file: draws/TadpoleUp.tex
%
\begingroup\makeatletter\ifx\SetFigFont\undefined%
\gdef\SetFigFont#1#2#3#4#5{%
  \reset@font\fontsize{#1}{#2pt}%
  \fontfamily{#3}\fontseries{#4}\fontshape{#5}%
  \selectfont}%
\fi\endgroup%
{\renewcommand{\dashlinestretch}{30}
\begin{picture}(1251,1553)(0,-10)
\path(638,1526)(450,1426)(652,1359)
\path(471,313)(560,509)(641,312)
\path(559,92)(559,692)
\put(560,1068){\ellipse{750}{750}}
\thicklines
\path(33,93)(1218,93)
\blacken\path(1038.000,33.000)(1218.000,93.000)(1038.000,153.000)(1038.000,33.000)
\end{picture}
}

%% file: draws/natOV.tex
%
\begingroup\makeatletter\ifx\SetFigFont\undefined%
\gdef\SetFigFont#1#2#3#4#5{%
  \reset@font\fontsize{#1}{#2pt}%
  \fontfamily{#3}\fontseries{#4}\fontshape{#5}%
  \selectfont}%
\fi\endgroup%
{\renewcommand{\dashlinestretch}{30}
\begin{picture}(324,339)(0,-10)
\path(12,12)(312,312)
\path(12,312)(87,237)
\path(237,87)(312,12)
\end{picture}
}

%% file: draws/natUN.tex
%
\begingroup\makeatletter\ifx\SetFigFont\undefined%
\gdef\SetFigFont#1#2#3#4#5{%
  \reset@font\fontsize{#1}{#2pt}%
  \fontfamily{#3}\fontseries{#4}\fontshape{#5}%
  \selectfont}%
\fi\endgroup%
{\renewcommand{\dashlinestretch}{30}
\begin{picture}(324,339)(0,-10)
\path(12,312)(312,12)
\path(312,312)(237,237)
\path(87,87)(12,12)
\end{picture}
}

%% file: page1.tex
\title[A Diagrammatic Reshetikhin-Turaev Invariant]
  {Towards a Diagrammatic Analogue of the Reshetikhin-Turaev Link
   Invariants}

\author{Ami~Haviv}
\address{Institute of Mathematics\\
        The Hebrew University\\
        Giv'at-Ram, Jerusalem 91904\\
        Israel}
\email{haviv@math.huji.ac.il}

\date{September 1, 2002.}

\begin{abstract}
By considering spaces of directed Jacobi diagrams, we construct a diagrammatic 
version of the Etingof-Kazhdan quantization of complex semisimple Lie algebras.
This diagrammatic quantization is used to provide a construction of a
directed version of the Kontsevich integral, denoted $\Z_{\ek}$, in a 
way which is analogous to the construction of the Reshetikhin-Turaev 
invariants from the $R$-matrices of the Drinfel'd-Jimbo quantum groups.
Based on this analogy, we conjecture (and prove in a restricted sense) 
a formula for the value of the invariant $\Z_{\ek}$ on the unknot. This 
formula is simpler than the Wheels formula of 
~\cite{Bar-NatanGaroufalidisRozanskyThurston:WheelsWheeling}, 
but the precise relationship between the two is yet unknown.
\end{abstract}

\maketitle

\tableofcontents

%% file: intro.tex
\section{Introduction}

\subsection{Quick summary for the expert}

The algebra $\A$ of Jacobi diagrams modulo IHX-like relations is the
target space for the Kontsevich integral (of framed knots), which is a 
universal Vassiliev (finite-type) framed knot invariant. This algebra is 
fruitfully thought of as a diagrammatic model for the `invariant part' of the 
theory of complex semisimple Lie algebras. 
We propose a richer diagrammatic framework, involving {\em directed} 
Jacobi diagrams modulo directed versions of the IHX relations, which 
captures some features of the Cartan-Weyl theory. 
 
\smallskip
Within this framework we can provide a diagrammatic interpretation of
the Etingof-Kazhdan quantization of a complex semisimple Lie algebra $\g$ 
(more precisely, the quantization of the Manin triple corresponding to the
standard Lie bialgebra structure of $\g$). We are then able to
construct a framed link invariant $\Z_{\ek}$, taking values in suitable spaces 
of directed Jacobi diagrams, in a way which is analogous, to some extent, 
to the construction of the Reshetikhin-Turaev framed link invariants from the 
$R$-matrices of the Drinfel'd-Jimbo quantum groups. The invariant 
$\Z_{\ek}$ is closely related (and possibly equivalent) to the
Kontsevich integral, by means of the map $\iota$ sending an 
undirected Jacobi diagram to the sum over all possibilities of directing 
its arcs. The analogy with the Reshetikhin-Turaev invariants leads us to 
conjecture a formula, which holds `at the level of complex semisimple 
Lie algebras', for the value of $\Z_{\ek}$ on the unknot, 
in terms of the exponential of a certain simple directed Jacobi diagram:
\[
	\Z_{\ek}(\text{unknot})  
	= \tr(\exp(\frac{_1}{^2}\TP)) \, , 
\]
where the exponential is computed with respect to the juxtaposition 
product and $\tr$ is the operation of closing the line into a circle.

\smallskip
While our motivation for introducing the spaces of directed Jacobi 
diagrams is purely Lie algebraic, similar (but not identical) spaces
of `arrow diagrams' were considered by Polyak (\cite{Polyak:ArrowDgrmAlg}) in 
conjunction with the theory of virtual knots. We discuss the similarities and 
differences between our and Polyak's diagrammatic spaces.     

\subsection{Paper overview}

We provide here a more leisurely description of this paper, without 
strictly following the way it is organized. More information about
the organization of the paper can be obtained from the table
of contents, as well as from the short summaries at the beginning of 
each section.  

\subsubsection{\bf The universal Vassiliev invariant, $\protect\A$-spaces,
and invariant Lie tensors.}
Vassiliev invariants of  framed knots (and other  knot-like objects) are an 
important and wide family of knot invariants, including, for example, the 
coefficients of the  Conway, Jones, HOMFLY, and Kauffman polynomials. 

The fundamental theorem of the theory of Vassiliev invariants
asserts that all the (numerical) Vassiliev invariants of framed knots are 
encapsulated  in a single `universal' invariant $\ZK$, which takes its 
values in  $\A$, the {\em algebra of chord diagrams}. 
Every numerical Vassiliev invariant  is obtained by composing $\ZK$ with a 
{\em weight system} (a linear functional in the graded dual of $\A$). Similar 
results hold for other types of knotted objects.

There are several known  methods for constructing a universal Vassiliev
invariant (see ~\cite{Bar-NatanStoimenow:Fundamental} for a survey). While 
the defining universal property of $\ZK$ does not determine it uniquely,
all the constructions yield, at least conjecturally,  
the same result (up to minor normalization issues), which is usually referred 
to, after its first construction, as the {\em (framed) Kontsevich integral}.  

\smallskip
Recall that $\A$ is the graded space of formal linear combinations of
certain uni-trivalent diagrams, called
{\em Jacobi diagrams}, modulo certain linear relations --- the AS, IHX, and
STU relations. The diagrams in $\A$ are based on a {\em skeleton} consisting of
a single oriented interval, meaning that the univalent vertices of a diagram 
in $\A$ are ordered along the skeleton interval. The trivalent vertices of a 
diagram are {\em oriented}, meaning that one of the two possible orderings of 
the arcs meeting in a trivalent vertex is specified.
The space $\A$ is graded by half the total number of vertices in a diagram, 
and is often interchanged with its completion with respect to that grading.
Examples of diagrams in $\A$ and the local parts of the relations are shown
in Figure ~\ref{fig:Asketch}. 

By varying the skeleton, allowing any ordered set of oriented intervals, 
oriented circles, and `colors' as the set of skeleton components, we obtain a 
variety of spaces, collectively called {\em $\A$-spaces}. In particular, the 
space $\B$ is the $\A$-space corresponding
to a single color skeleton, namely the space spanned by Jacobi diagrams with 
unordered  univalent vertices  modulo AS and IHX relations. The $\A$-space 
based on an $n$-tuple of oriented intervals is denoted by $\A^{\bt n}$, and 
the one based on an $n$-tuple of colors is denoted $\B^{\bt n}$. Some
examples of diagrams in these spaces are also shown in Figure 
~\ref{fig:Asketch}.	  

\begin{figure}[htpb]
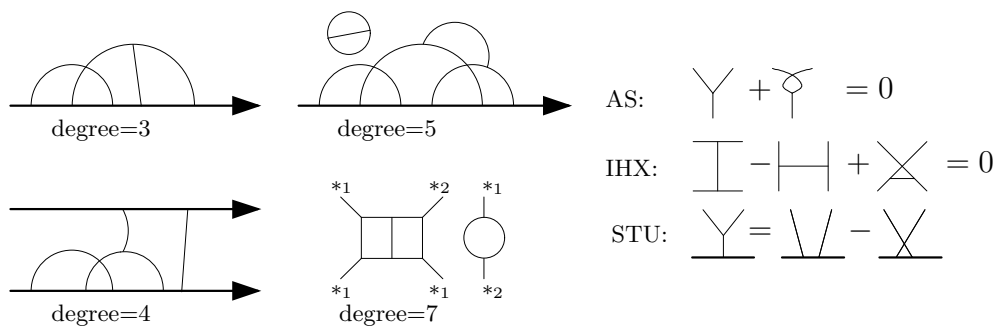

\[ 
  \printname{Asketch}
  \setlength{\unitlength}{0.85\standardunitlength}
  \begin{array}{c}  \hspace{-1.7mm}
    \raisebox{-8pt}{\input draws/Asketch.tex }
      \hspace{-1.9mm}
  \end{array}
 \]
\caption{Examples of diagrams in $\A$, $\A^{\bt 2}$, and $\B^{\bt 2}$, and 
the local parts of the AS, IHX, and STU relations. All the trivalent vertices 
are oriented counterclockwise. $\ast_1$ and $\ast_2$ are the colors of
univalent vertices of diagrams in $\B^{\bt 2}$.} \lbl{fig:Asketch}
\end{figure}

\smallskip
In many respects, $\A$ is a diagrammatic analogue of the 
$\ad$-invariant part of the universal enveloping algebras $U(\g)$ of 
{\em metrized} Lie algebras $\g$ (notably complex semisimple Lie algebras). 
This analogy is mediated by the `diagrams to invariant Lie tensors' maps, 
collectively denoted $\T_{\g}$, which assign to diagrams in $\A$-spaces 
elements in appropriate $\g$-spaces. 
For example, $\T_{\g}$ maps $\A^{\bt n}$  into $(U(\g)^{\ot n})^{\g}$ and maps
$\B$ into $S(\g)^{\g}$, the spaces of $\g$-invariants (with respect to the 
adjoint 
action) in the $n$-fold tensor product of $U(\g)$ and in the symmetric 
algebra $S(\g)$ of $\g$, respectively.  

$\A$ is equipped with the structure of a {\em diagrammatic Hopf algebra}, whose
components are maps between the  {\em diagrammatic tensor powers} 
$\A^{\bt n}$ of $\A$. All the Hopf structure components are intertwined by 
suitable $\T_{\g}$ maps with the parallel
components of the standard Hopf algebra structure of $U(\g)$. The 
product of $\A$ is done by juxtapositioning the skeleton intervals.
Similarly, $\B$ is endowed with the 
'two colors made one' product, which is
analogous to the product of $S(\g)$; when restricted to $\B\ot\B$  
(which naturally injects into $\B^{\bt 2}$), it equips $\B$ with the disjoint 
union product (which is a product operation in the usual algebraic sense).

The Poincar{\`e}-Birkhoff-Witt (PBW) isomorphism $\upchi: S(\g) \ra U(\g)$ 
restricts to an isomorphism $\upchi: S(\g)^{\g} \ra U(\g)^{\g}$. It has a
diagrammatic counterpart, the {\em diagrammatic PBW isomorphism} 
$\upchi: \B \ra \A$, which acts by averaging all the possible orderings
of the legs of a diagram in $\B$ along the oriented interval skeleton of $\A$.

Another instance of the $\A$ to $U(\g)^{\g}$ analogy is that the natural
`closing the interval' projection $\A = \A(\ua) \twoheadra \A(\bigcirc)$
is actually an isomorphism, just like its algebraic counterpart, which 
is the projection of $U(\g)^{\g}$ onto the space of its 
{\em $\g$-coinvariants}; the latter two spaces are trivially isomorphic.

\smallskip
Drinfel'd showed that the Hopf algebra structure of $U(\g)$ can be deformed,
in the space $\Ugh$ of formal power series in $\hbar$ over $U(\g)$, into a 
quasitriangular quasi-Hopf
structure\footnote{The basic definitions from the theory of quasitriangular 
quasi-Hopf algebras are recalled in Section ~\ref{_sect:qtrqH}}, 
sharing the product and coproduct of $U(\g)$. The $R$-matrix of the
deformed structure is the 
exponential $R_{\g} = e^{\hbar \Omega_{\g}/2} \in (U(\g)^{\ot 2})^{\g}\hb$ of 
half the metric tensor of $\g$, and its associator 
$\Phi_{\g} \in (U(\g)^{\ot 3})^{\g}\hb$ is the famous 
{\em Drinfel'd associator}, originating from the Knizhnik-Zamolodchikov 
(KZ) system of differential equations. We call this structure the 
{\em KZ structure} of $\Ugh$. Moreover, these KZ structures are
'universal' in $\g$: there are
elements $R\in \A^{\bt 2}$ and $\Phi\in \A^{\bt 3}$ such that 
$R_{\g} = \Tgh(R)$ and $\Phi_{\g} = \Tgh(\Phi)$ for any semisimple Lie 
algebra $\g$. (The maps $\Tgh$ are simple variants of $\T_{\g}$, multiplying 
the value of $\T_{\g}$ on a diagram of degree $d$ by $\hbar^d$.) 
The element $R$ is simple to describe: 
$R = \exp(\frac{1}{2}
  \setlength{\unitlength}{0.2\standardunitlength}
  \begin{array}{c}  \hspace{-1.7mm}
    \raisebox{-2pt}{\input draws/BCup.tex }
    \hspace{-1.9mm}
  \end{array}
)$. On the contrary, the Drinfel'd 
associator $\Phi$ is a mysterious object. The elements $R$ and $\Phi$
equip $\A$ with the structure of a {\em diagrammatic quasitriangular 
quasi-Hopf algebra} (which we  call the KZ structure of $\A$). 

\smallskip
The axioms of quasitriangular quasi-Hopf algebra $H$, either in the usual 
algebraic sense or in our diagrammatic setting, lead to a procedure for 
associating to such an algebra a topological invariant $\Z_H$ of 
{\em parenthesized framed tangles}. This procedure is explained in depth in 
Section ~\ref{sect:LinkInv}. It has the important property that the restriction
of $\Z_H$ to framed links is unchanged if the quasitriangular quasi-Hopf
structure of $H$ is replaced by a {\em twist equivalent} structure.  
The values of $\Z_H$ on framed links with $n$ components belong to the space 
of {\em coinvariants} of $H^{\ot n}$. For $H= \A$, this 
space is $\A(\bigcirc_1, \ldots, \bigcirc_n)$, whereas for $H = U(\g)$ it is
equal to $(U(\g)_{\g})^{\ot n}$, the $n$-fold tensor power of the space
of $\g$-coinvariants of $U(\g)$.

\smallskip
Applying the $\Z_H$-procedure to the KZ structure of $\A$ and 
restricting the resulting invariant to framed knots, we obtain a universal 
Vassiliev invariant taking values in $\A(\bigcirc) \cong \A$
(which is known to coincide with the Kontsevich integral $\ZK$). Applying
the same procedure to the accompanying structures on the enveloping
algebras $\Ugh$ yields a series of invariants of framed knots, which we 
denote by $\KZg$. The invariants $\KZg$ take values in the spaces 
$U(\g)_{\g}\hb$, and are $\T_{g}$-related to the Kontsevich integral: 
$\KZg = \Tgh \circ \ZK$.

\smallskip
The complicated nature of $\Phi$ makes it difficult to carry out explicit 
non-perturbative computations of the values of $\ZK$ for particular
knots (the situation is no better for other constructions of the
Kontsevich integral). In fact, even the 
value of $\ZK$ for the {\em unknot} was hard to compute. The following 
{\em Wheels formula} for this value, initially conjectured in 
~\cite{Bar-NatanGaroufalidisRozanskyThurston:WheelsWheeling}, was  
recently proved (see ~\cite{Thurston:Duflo}, 
~\cite{Bar-NatanLeThurston:TwoApplications},  and the references therein):
\begin{equation}
\label{eqn:BGRTformula}
	\upchi^{-1}\ZK(\text{unknot}) 
	= \exp_{\udot} \sum_{n=1}^{\infty} b_{2n} \omega_{2n}
    = \exp_{\udot}(\frac{\text{\Large$\twowheel$}}{48} 
                  - \frac{\text{\Large$\fourwheel$}}{5760}
                  + \frac{\text{\Large$\sixwheel$}}{362880} \mp \cdots) \, ,
\end{equation}
where the `modified Bernoulli numbers' $b_{2n}$ are defined by the power
series expansion
\[
\sum_{n=0}^{\infty} b_{2n} x^{2n} = \frac{1}{2}\log \frac{\sinh x/2}{x/2}  \, ,
\]
$\omega_{2n}$ is the `$2n$-wheel' (the uni-trivalent diagram made of a $2n$-gon
with $2n$ legs), and the exponential is computed using the disjoint
union product of $\B$.

\subsubsection{\bf The quest for a diagrammatic Reshetikhin-Turaev invariant.}
\label{__sect:DgrmRTquest}

When $\g$ is a complex semisimple Lie algebra, there is an alternative 
construction of the invariants $\KZg$, namely their values on framed links 
coincide with those of the {\em Reshetikhin-Turaev invariants} $\RTg$. In our 
language, the latter invariants are obtained by applying the $\Z_H$-procedure
to another quasitriangular quasi-Hopf structure on $\Ugh$, which is an
isomorphic image of the {\em Drinfel'd-Jimbo quantum group} $\Uhg$, 
a quantum universal enveloping algebra of $\g$ defined by means of  
`quantum  Chevalley generators and quantum Serre relations'. (All this is
reviewed in detail in Section ~\ref{_sect:DJ}.) We call this 
structure the {\em DJ structure} of $\Ugh$. The fact that the
framed link invariants $\KZg$ and $\RTg$ coincide follows from
a deep result of Drinfel'd, asserting that the underlying 
quasitriangular quasi-Hopf KZ and DJ structures are twist equivalent. 

\smallskip
We comment that in our approach the Reshetikhin-Turaev invariant  $\RTg(K)$ of 
a framed knot $K$ belongs to the space $U(\g)_{\g}\hb$. In the literature, 
starting with ~\cite{ReshetikhinTuraev:Ribbon}, the value
$\RTg(K)$ is identified with the equivalent collection of formal series
in $\hbar$, denoted $\RT_{\g; V}$, which are obtained by taking the trace
of the action of $\RTg(K)$ on a finite-dimensional module $V$ of $\g$. 
(A similar remark applies to framed links with any number of components.)  

\smallskip 
While the KZ construction is universal and manifestly $\g$-invariant,
the $\RTg$ construction has several merits:
\begin{itemize}
\item[(RT1)] The underlying DJ structure is actually Hopf (that is,
its associator is trivial), so that there is no need for the somewhat
artificial machinery of parenthesization for computing the values of $\RTg$.

\item[(RT2)] The $R$-matrix of $\Uhg$ has an {\em explicit} (albeit 
complicated) expression in terms of the generators of $\Uhg$, with 
{\em rational} coefficients,   
which implies that the values of $\RT_{\g; V}$ are rational power series in 
$\hbar$. (A similar rationality result can be proved directly for $\KZg$, and 
even in the `universal' framework of the KZ structure of $\A$,
but only by indirect means; see Section ~\ref{_sect:KontAndRT} for a review 
of such a proof.)

\item[(RT3)] There is a simple formula for the value of $\RTg$ on the
unknot. Fixing an arbitrary Cartan subalgebra of $\g$ and then choosing a 
set of positive roots for the associated root system, we have
\begin{equation}
\label{eqn:RTunknot}
	\RTg(\text{unknot}) = \tr(e^{\hbar \varrho_{\g}}) \, ,
\end{equation}
where $\varrho_{\g}$ is the element of the Cartan algebra dual to the half-sum
of the positive roots, and $\tr$ is the natural projection of $U(\g)\hb$
onto $U(\g)_{\g}\hb$.
\end{itemize}

\smallskip
These advantages make it appealing to look for a diagrammatic version of the
Reshetikhin-Turaev construction. However, there is an inherent difficulty
with trying to make such a construction in the diagrammatic framework 
offered by $\A$-spaces, in that the components of the DJ structure underlying 
$\RTg$ are {\em not $\g$-invariant}, and thus are inaccessible by the 
$\T_{\g}$-maps. Even if this initial difficulty is tackled somehow, the
generators and relations definition of the Drinfel'd-Jimbo quantum
groups is highly `non-universal' and does not seem to be directly
susceptible to a diagrammatic reformulation. 

\subsubsection{\bf A diagrammatic Etingof-Kazhdan quantization.}
\label{__sect:DgrmEK}
   
Our hint for (partially) solving both problems came from the 
Etingof-Kazhdan theory of quantization of Lie bialgebras,
which we quickly describe. (Section ~\ref{_sect:LBA} includes the needed 
definitions and facts from the theory of Lie bialgebras). 

\smallskip
A finite-dimensional Lie bialgebra $\g_+$ is part of a {\em Manin triple}
$(\tg, \g_+, \g_-)$. The ambient Lie bialgebra $\tg$ is the 
{\em Drinfel'd double} of the Lie bialgebra $\g_+$. It carries a natural 
metric, hence the enveloping
algebra $\Utgh$ has a KZ structure. By using {\em Verma modules} of $\tg$
and a fixed Drinfel'd associator, Etingof and Kazhdan 
(\cite{EtingofKazhdan:QuantizationI}) construct an element 
$J\in U(\tg)^{\ot 2}\hb$, such that twisting the KZ structure of $\Utgh$ 
by $J^{-1}$ leads to a {\em coassociative} structure. We shall call this 
quasitriangular (quasi-)Hopf structure of $\Utgh$ the {\em EK structure}. 

\smallskip
A complex semisimple Lie algebra $\g$ with a choice of Cartan subalgebra,
followed by a choice of the set of positive roots, has a {\em standard} Lie 
bialgebra structure. $\g$ equipped with this structure is not the ambient Lie 
bialgebra of a Manin triple, but is closely related to one, denoted $\tg$, 
which is obtained by `doubling' the Cartan subalgebra of $\g$, and is the 
Drinfel'd double of the positive Borel subbialgebra of $\g$ (see Section 
~\ref{__sect:CpxSimpleLA})\footnote{It is because of this difference between 
$\g$ and $\tg$ that we have chosen the notation $\tg$, rather than $\g$, for 
the `generic' ambient Lie bialgebra of a Manin triple, while reserving $\g$ 
for denoting most of the time a complex semisimple Lie algebra.}.

\smallskip
The EK structure of $U(\tg)$ is sent by the projection $p:\tg \twoheadra \g$ 
to a  quasitriangular Hopf structure on $\Ugh$, carrying the same name. By its
construction, the EK structure is twist equivalent to the KZ structure of 
$\Ugh$.
   
\smallskip
Since Etingof-Kazhdan theory is `universal' in the Manin triple 
$(\tg, \g_+, \g_-)$, it seemed to be a good candidate for a diagrammatic 
interpretation. Our results in doing so are the central theme of this paper.

\smallskip 
We begin by introducing a diagrammatic model of Manin triples 
(or, equivalently, 
double Lie bialgebras), by considering {\em directed} Jacobi diagrams, with 
suitable directed versions of the IHX-like relations. Most interestingly, 
there is a `5T relation'
\[ 

  \printname{ArrowI3}
  \setlength{\unitlength}{0.25\standardunitlength}
  \begin{array}{c}  \hspace{-1.7mm}
    \raisebox{-8pt}{\input draws/ArrowI3.tex }
      \hspace{-1.9mm}
  \end{array}
 = 
  \printname{ArrowH3a}
  \setlength{\unitlength}{0.25\standardunitlength}
  \begin{array}{c}  \hspace{-1.7mm}
    \raisebox{-8pt}{\input draws/ArrowH3a.tex }
      \hspace{-1.9mm}
  \end{array}
 - 
  \printname{ArrowX3a}
  \setlength{\unitlength}{0.25\standardunitlength}
  \begin{array}{c}  \hspace{-1.7mm}
    \raisebox{-8pt}{\input draws/ArrowX3a.tex }
      \hspace{-1.9mm}
  \end{array}
 
- 
  \printname{ArrowX3b}
  \setlength{\unitlength}{0.25\standardunitlength}
  \begin{array}{c}  \hspace{-1.7mm}
    \raisebox{-8pt}{\input draws/ArrowX3b.tex }
      \hspace{-1.9mm}
  \end{array}
 + 
  \printname{ArrowH3b}
  \setlength{\unitlength}{0.25\standardunitlength}
  \begin{array}{c}  \hspace{-1.7mm}
    \raisebox{-8pt}{\input draws/ArrowH3b.tex }
      \hspace{-1.9mm}
  \end{array}
 \, , 
\]
which is the diagrammatic analogue of the {\em cocycle identity} satisfied by 
the bracket and cobracket of $\g_+$.

The resulting diagrammatic spaces are called $\Aar$-spaces. In particular, 
$\Aar = \Aar(\ua)$ is called the {\em algebra of directed Jacobi diagrams}.

\smallskip
The vague notion of  `diagrammatic model' is made more precise
by constructing the `directed Jacobi diagrams to Lie tensors' maps, 
collectively denoted $\Tartg$, from $\Aar$-spaces to the paralleling 
tensor products of $U(\tg)$, $U(\tg)_{\tg}$, and $S(\tg)$. This
construction is carried out in Section ~\ref{__sect:DrctdJacobiToLie}.   
When $\g$ is a complex semisimple Lie algebra, we denote by $\Targ$ the 
composition
\[
\begin{CD}
\Aar\text{-space} @>{\Tartg}>> \tg\text{-space} @>{p}>> \g\text{-space} \, .
\end{CD}
\]
The $\hbar$-variants $\Tartgh$ and $\Targh$ of these maps are defined as in 
the undirected case.

\smallskip
The tensors in the images of the $\Targ$ maps are not necessarily 
$\g$-invariant. An important example is the element 
\begin{equation}
\label{eqn:DgrmVarrho}
	\varrho = \frac{_1}{^2}(\LHC - \RHC) \in \Aar \, , 
\end{equation}
which is mapped by $\Targ$ ($\g$ a complex semisimple Lie algebra) to 
$\varrho_{\g}$.

\smallskip
In the Lie algebraic setting, the $\g$-invariant tensors trivially form 
subspaces of their ambient tensor spaces. The diagrammatic situation 
is more involved, in that an $\A$-space is not directly a subspace of the 
$\Aar$-space based on the same skeleton. We do have a collection of natural 
maps between isoskeletal $\A$ and $\Aar$-spaces, called the 
$\iota$-maps, which act on a diagram in an $\A$-space via `directization': 
sending that diagram to the sum of all possible ways of directing its arcs.
We show that these maps are well-defined and `cover' the injections
$U(\tg)^{\tg} \ra U(\tg)$, for any double Lie bialgebra $\tg$. This motivates 
us to conjecture (see Conjecture ~\ref{conj:IotaIsInjective} in Section 
~\ref{__sect:ImgKerIota}) that all the $\iota$-maps are {\em injective}.

\smallskip
This sets the stage for a diagrammatic Etingof-Kazhdan construction. Using the
$\iota$-maps, the KZ structure is exported from $\A$ to $\Aar$. The 
{\em diagrammatic Verma modules} $\VM_+$ and $\VM_-$ are defined as the 
quotients of $\Aar$ by the RI (``right in'') and RO (``right out'') relations,
respectively, with natural projections $p_{\pm}:\Aar\ra \VM_{\pm}$. The RI and
RO relations are illustrated in the following figure. 
\[
	\text{ RI:\, } 
  \printname{RightLegIn}
  \setlength{\unitlength}{0.3\standardunitlength}
  \begin{array}{c}  \hspace{-1.7mm}
    \raisebox{-8pt}{\input draws/RightLegIn.tex }
      \hspace{-1.9mm}
  \end{array}
 = 0 \qquad; \qquad
	\text{ RO:\, } 
  \printname{RightLegOut}
  \setlength{\unitlength}{0.3\standardunitlength}
  \begin{array}{c}  \hspace{-1.7mm}
    \raisebox{-8pt}{\input draws/RightLegOut.tex }
      \hspace{-1.9mm}
  \end{array}
 = 0 \, .
\]

Diagrammatic versions of algebraic PBW arguments show that the (diagrammatic)
coproduct of $\Aar$ induces an isomorphism $\phi:\Aar \ra \VM_+ \bt \VM_-$.
The diagrammatic analogue of the Etingof-Kazhdan twist can then be defined
as
\[
	J = (\phi^{-1}\bt\phi^{-1})
		(p_+\bt p_-\bt p_+\bt p_-)
	\ZK\left(
  \setlength{\unitlength}{0.25\standardunitlength}
  \begin{array}{c}  \hspace{-1.7mm}
    \raisebox{-2pt}{\input draws/Jbraid.tex }
    \hspace{-1.9mm}
  \end{array}
\right) \, 
	\in \Aar^{\ot 2} \, ,
\]
and it is not too hard to prove that twisting by $J^{-1}$ establishes a twist 
equivalence between the KZ structure of $\Aar$ and a {\em coassociative} 
structure, called the {\em EK structure} of $\Aar$. The topological invariant 
associated to the latter structure, denoted $\Z_{\ek}$, takes its values in 
$\Aar$-spaces, and it rEadily follows that, when restricted to framed links, 
it is $\iota$-related to the Kontsevich integral (hence the two invariants are 
equivalent if Conjecture ~\ref{conj:IotaIsInjective} is true). 

\smallskip
In accordance with ~\eqref{eqn:RTunknot} and ~\eqref{eqn:DgrmVarrho}, we 
conjecture\footnote{We actually make a more refined conjecture about the
value of the {\em Drinfel'd element} of the EK structure of $\Aar$; see
Conjecture ~\ref{conj:Uek} in Section ~\ref{_sect:Zek}.} that
\begin{equation}
\label{eqn:EKunknot}
	\Z_{\ek}(\text{unknot})  
	= \tr(\exp(\frac{_1}{^2}\TP)) \, , 
\end{equation}
where $\tr$ is the `closing the circle' map $\Aar \twoheadra \Aar(\bigcirc)$,
and 
\[
	\TP = \TPD + \TPU \, .
\] 
We show that this conjecture holds at the level of complex 
semisimple Lie algebras: applying $\Targh$ to both sides, with $\g$ 
complex semisimple, yields an equality. 

\smallskip
This putative formula is simpler than the Wheels
formula \eqref{eqn:BGRTformula} in several respects. First, it 
gives the value on the unknot in the original target space of 
the invariant $\Z_{\ek}$, without having to resort to the PBW
map $\upchi$. Second, the diagrams and rational numbers appearing in it
are far less complicated than their counterparts in the Wheels formula.  

\subsection{Related work}
Spaces of directed Jacobi diagrams, similar to the $\Aar$-spaces, were 
independently studied by M.~Polyak in ~\cite{Polyak:ArrowDgrmAlg}, in
relation to the theory of virtual knots 
(see ~\cite{GoussarovPolyakViro:VirtualKnots}). We review Polyak's
work in Section ~\ref{_sect:Arrow}, and discuss the similarities and 
differences between his spaces and $\Aar$-spaces.

\subsection{Acknowledgments}
This paper is my Ph.D. thesis, which was  prepared in the Hebrew University 
of Jerusalem. I would like to thank my advisor,  Dror Bar-Natan, for 
suggesting the topic of this work and for numerous instructive talks. 

I also thank P.~Etingof, D.~Kazhdan, J.~Lieberum, T.~Ohtsuki, 
M.~Polyak, E. de-Shalit, and D.~P.~Thurston for some useful discussions
and helpful comments.

%% file: draws/BCup.tex
%
\begingroup\makeatletter\ifx\SetFigFont\undefined%
\gdef\SetFigFont#1#2#3#4#5{%
  \reset@font\fontsize{#1}{#2pt}%
  \fontfamily{#3}\fontseries{#4}\fontshape{#5}%
  \selectfont}%
\fi\endgroup%
{\renewcommand{\dashlinestretch}{30}
\begin{picture}(1086,1431)(0,-10)
\path(93,633)(993,633)
\thicklines
\path(93,33)(93,1383)
\blacken\path(153.000,1203.000)(93.000,1383.000)(33.000,1203.000)(153.000,1203.000)
\path(993,33)(993,1383)
\blacken\path(1053.000,1203.000)(993.000,1383.000)(933.000,1203.000)(1053.000,1203.000)
\end{picture}
}

%% file: draws/ArrowX3a.tex
%
\begingroup\makeatletter\ifx\SetFigFont\undefined%
\gdef\SetFigFont#1#2#3#4#5{%
  \reset@font\fontsize{#1}{#2pt}%
  \fontfamily{#3}\fontseries{#4}\fontshape{#5}%
  \selectfont}%
\fi\endgroup%
{\renewcommand{\dashlinestretch}{30}
\begin{picture}(2126,2115)(0,-10)
\path(695,678)(1444,675)
\path(2114,17)(1214,902)
\path(1661,592)(1864,262)(1531,462)
\path(12,12)(912,897)
\path(465,587)(262,257)(595,457)
\path(912,898)(2113,2088)
\path(1801,1648)(1464,1448)(1669,1774)
\path(324,1648)(661,1448)(456,1774)
\path(904,601)(1213,678)(904,749)
\path(1222,898)(12,2088)
\end{picture}
}

%% file: draws/ArrowX3b.tex
%
\begingroup\makeatletter\ifx\SetFigFont\undefined%
\gdef\SetFigFont#1#2#3#4#5{%
  \reset@font\fontsize{#1}{#2pt}%
  \fontfamily{#3}\fontseries{#4}\fontshape{#5}%
  \selectfont}%
\fi\endgroup%
{\renewcommand{\dashlinestretch}{30}
\begin{picture}(2126,2115)(0,-10)
\path(695,678)(1444,675)
\path(912,898)(2113,2088)
\path(1801,1648)(1464,1448)(1669,1774)
\path(324,1648)(661,1448)(456,1774)
\path(1222,898)(12,2088)
\path(1229,601)(920,678)(1229,749)
\path(12,12)(912,897)
\path(465,587)(262,257)(595,457)
\path(2114,17)(1214,902)
\path(1661,592)(1864,262)(1531,462)
\end{picture}
}

%% file: pre.tex
\section{Preliminaries} 
\label{sect:pre}

{\bf Section summary:} We recall some of the definitions, constructions, and 
results of the basic theory of Vassiliev invariants and $\A$-spaces (mostly 
from ~\cite{Bar-Natan:Vassiliev}), in a way which aims to ease the treatment of
their directed counterparts in the next section. After stating the 
definitions of Vassiliev invariants, weight systems, and chord 
diagrams, we describe the transition to Jacobi diagrams (allowing for 
trivalent vertices) and the fundamental construction of the maps from 
$\A$-spaces to spaces of (invariant) Lie tensors. We give a detailed 
review of the proof of the diagrammatic PBW theorem, and close this section
with a study of the basic algebraic operations on $\A$-spaces,
emphasizing that $\A$ has the structure of a diagrammatic Hopf algebra. 

\subsection{Vassiliev invariants}

\subsubsection{\bf Definition of Vassiliev invariants.}
\label{__sect:VassDef}

Suppose $V$ is an invariant of oriented knots taking values in an
abelian group. It can be inductively 
extended to an invariant $V^{(m)}$ of immersed oriented circles with exactly 
$m$ transversal self intersection points, using the formulas
\begin{gather}
V^{(0)}=V \, , \notag \\
V^{(m)}\left(\text{\Huge \raisebox{-6pt}{$\doublepoint$}}\right) =
	V^{(m-1)}\left(\text{\Huge \raisebox{-6pt}{$\overcrossing$}}\right) - 
	V^{(m-1)}\left(\text{\Huge \raisebox{-6pt}{$\undercrossing$}}\right)
        \, . \label{eqn:SingExtnd} 
\end{gather}

\begin{definition}
\lbl{def:VassInv}
$V$ is called a {\em Vassiliev invariant of type $m$} if $V^{(m+1)}$ 
vanishes identically. 
\end{definition}

Similar definitions can be made for a variety of `knot-like' objects, such as
braids, tangles and framed links. 

\subsubsection{\bf Chord diagrams and weight systems.}
\label{__sect:CDandWgtSys}

\begin{definition}
\lbl{def:sklt}
A {\em ($1$-dimensional) skeleton} is a finite {\em ordered} collection of 
oriented intervals and oriented circles, which are called the {\em components} 
of the skeleton. 
\end{definition}

\begin{definition}
\lbl{def:CD}
A {\em chord diagram} based on a skeleton $\M$ is a finite graph made of the
components of $\M$ and finitely many chords marked on those components.
The collection of all $\M$-based chord diagrams is denoted
by $\D^c(\M)$. This collection is graded by the number of chords in a diagram. 
Its degree $m$ part is denoted by $\G_m(\D^c(\M))$.
\end{definition}

An example of a chord diagram is depicted in Figure ~\ref{fig:CDexample}.

\begin{figure}[htpb]
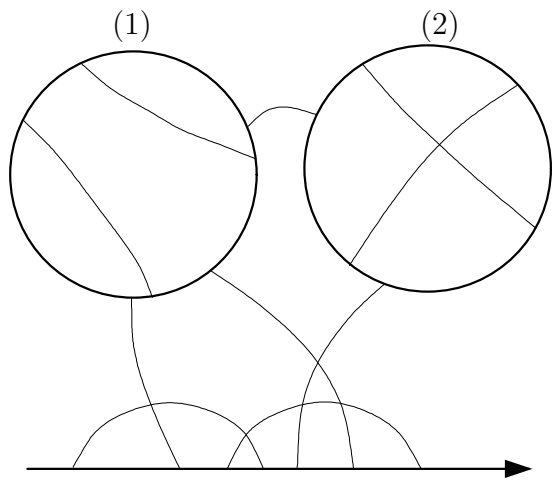

\[ 

  \printname{CDexample}
  \setlength{\unitlength}{0.8\standardunitlength}
  \begin{array}{c}  \hspace{-1.7mm}
    \raisebox{-8pt}{\input draws/CDexample.tex }
      \hspace{-1.9mm}
  \end{array}
 
\] 
\caption{Example of a chord  diagram based on a three component skeleton. By 
convention, skeleton circles are oriented counterclockwise, and the
order of components is such that the subset of oriented circles succeeds the 
subset of oriented intervals.} \lbl{fig:CDexample}
\end{figure}

Given an orientation preserving immersion $K: \M \ra \mathbb{R}^3$ 
whose only singularities are  $m$ transversal self intersections, we 
associate to $K$ an $\M$-based chord diagram $D$ of degree $m$, by the 
following procedure: For each double point $p$ of $K$, mark a chord in 
$D$ connecting the two points in $\M$ which are mapped to $p$ (see 
Figure ~\ref{fig:SingPtrn}). This chord diagram can be thought of as recording 
the pattern of singularities of the immersion.    

\begin{figure}[htpb]
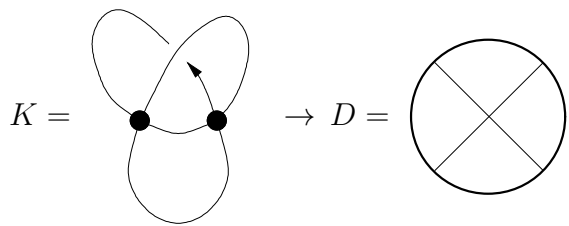

\[ K=
  \printname{TwoSingKnot}
  \setlength{\unitlength}{0.8\standardunitlength}
  \begin{array}{c}  \hspace{-1.7mm}
    \raisebox{-8pt}{\input draws/TwoSingKnot.tex }
      \hspace{-1.9mm}
  \end{array}
 \, \ra \, D=
  \printname{TwoSingDgrm}
  \setlength{\unitlength}{0.5\standardunitlength}
  \begin{array}{c}  \hspace{-1.7mm}
    \raisebox{-8pt}{\input draws/TwoSingDgrm.tex }
      \hspace{-1.9mm}
  \end{array}
 \]
\caption{Associating a chord diagram to a singular knot.} \lbl{fig:SingPtrn}
\end{figure}

Let $V$ be a Vassiliev knot invariant of type $m$ taking values in 
a field $\k$, and let $D$ be a 
(circle-based) chord diagram of degree $m$. $D$ can be {\em embedded} in
$\mathbb{R}^3$, that is, there is a knot 
$K_D:S^1 \ra \mathbb{R}^3$ (with $m$ double points) whose associated 
chord diagram is $D$.

Two different embeddings of $D$ are related by a sequence of moves of the form
$\text{\large \raisebox{-2pt}{$\overcrossing$}} \ra 
 \text{\large \raisebox{-2pt}{$\overcrossing$}}$. 
Such a move does 
not change the value of $V^{(m)}(K_D)$, because the change caused by the move
is equal to $V^{(m+1)}$ evaluated on a knot with $m+1$ double points, which
is zero by Definition ~\ref{def:VassInv}. Thus, we can unambiguously define the
{\em weight system of $V$}: 
\begin{equation} 
W_m(V):\G_m(\D^c(\bigcirc)) \ra \k \, , \quad 
W_m(V)(D)=V(K_D) \, . 
\end{equation}
$W_m(V)$ can be extended to $\D^c(\bigcirc)$ by setting it equal to zero on 
the other homogeneous components, and then can be extended by linearity to 
the $\k$-linear space formally spanned by chord diagrams, 
$\Span_{\k}(\D^c(\bigcirc))$. 

\smallskip
There is no difficulty in extending the definition of weight system to
Vassiliev invariants of tangles or framed tangles, yielding linear
functionals $W=W_m(V)$ on $\Span(\D^c(\M))$, for the appropriate
skeleton $\M$. In particular, it applies to Vassiliev invariants of framed 
links, where chord diagrams are based on a collection of circles. 

\smallskip
A simple topological argument (\cite[Section 2.2]{Bar-Natan:Vassiliev}) 
shows that the weight system $W$ of 
a Vassiliev invariant satisfies the {\em 4T relations}: $W$ vanishes 
on linear combinations of four chord diagrams of the form
\[

  \printname{4TW}
  \setlength{\unitlength}{0.5\standardunitlength}
  \begin{array}{c}  \hspace{-1.7mm}
    \raisebox{-8pt}{\input draws/4TW.tex }
      \hspace{-1.9mm}
  \end{array}
 - 
  \printname{4TE}
  \setlength{\unitlength}{0.5\standardunitlength}
  \begin{array}{c}  \hspace{-1.7mm}
    \raisebox{-8pt}{\input draws/4TE.tex }
      \hspace{-1.9mm}
  \end{array}

- 
  \printname{4TN}
  \setlength{\unitlength}{0.5\standardunitlength}
  \begin{array}{c}  \hspace{-1.7mm}
    \raisebox{-8pt}{\input draws/4TN.tex }
      \hspace{-1.9mm}
  \end{array}
 + 
  \printname{4TS}
  \setlength{\unitlength}{0.5\standardunitlength}
  \begin{array}{c}  \hspace{-1.7mm}
    \raisebox{-8pt}{\input draws/4TS.tex }
      \hspace{-1.9mm}
  \end{array}
 \, .
\]
(The oriented segments may belong to the same or different skeleton
components. The parts of the diagrams not shown in the figure are the
same for all four diagrams.) 

Thus, we may regard weight systems as linear functionals on the 
quotient space.
\[
\A^c(\M)=\Span_{\k}(\D^c(\M)) \left/
           \left\{\text{4T relations}\right\} \right. \, ,
\]
in which the above linear combinations vanish. $\A^c(\M)$ is a graded vector
space (as 4T relations are homogeneous). Its graded completion is also
denoted $\A^c(\M)$.

\begin{remark}
When we consider {\em unframed} knots, it is easy to derive another set of
relations satisfied by a weight system $W$ of a Vassiliev invariant --- the
{\em framing independence (FI) relations}: $W(D)=0$ if 
$D\in \G_m(\D^c(\bigcirc)$ has an isolated chord. Similar relations hold 
for unframed links or tangles (where the isolated chord has both endpoints 
on the {\em same} skeleton component).
\end{remark}

\begin{remark}
\lbl{rem:A-horizontal}
In the case of {\em pure braids} with $n$ strands, the suitable patterns of 
singularities are chord diagrams based on an $n$-interval skeleton with all 
the chords being {\em horizontal}. Denote this space by 
$\D^h(\ua_1\cdots \ua_n)$.
Weight systems of Vassiliev invariants of pure braids again satisfy the 4T 
relations. The corresponding quotient space is denoted by 
$\A^h(\ua_1\cdots \ua_n)$. 
\end{remark}

\subsection{The Algebra $\protect\A$ and invariant Lie tensors}
\subsubsection{\bf Jacobi diagrams.}

We extend the definition of a $1$-dimensional skeleton as follows.
\begin{definition}
\lbl{def:SkltWithColors}
A {\em skeleton} is a finite ordered collection of components, which are
either oriented intervals, oriented circles, or  
{\em colors}. (Colors are usually denoted $\ast_1, \ast_2, \ldots$)
\end{definition}

It is also useful to work with a family of diagrams extending the family of  
chord diagrams.
\begin{definition}
\lbl{def:diagram}
A {\em Jacobi diagram based on a skeleton $\M$} is a finite graph made of the
$1$-dimensional components of $\M$ 
and additional undirected arcs, with two types of vertices allowed:
\begin{enumerate}
\item {\em External} vertices in which an arc ends in a
component of $\M$. (Ending in $\ast_1$ means that the arc ends in a 
univalent vertex colored by $\ast_1$.) Arcs ending in external vertices 
are called {\em legs}.
\item {\em Internal} vertices in which three arcs meet. These vertices are 
oriented --- one of the two possible cyclic orderings of the arcs 
meeting in such a vertex is specified.
\end{enumerate}

The {\em body} of such a diagram is the graph obtained by removing all the 
components of $\M$ from it. (Removing the colors means making all the
univalent vertices colorless.) 
\end{definition}

\begin{figure}[htpb]
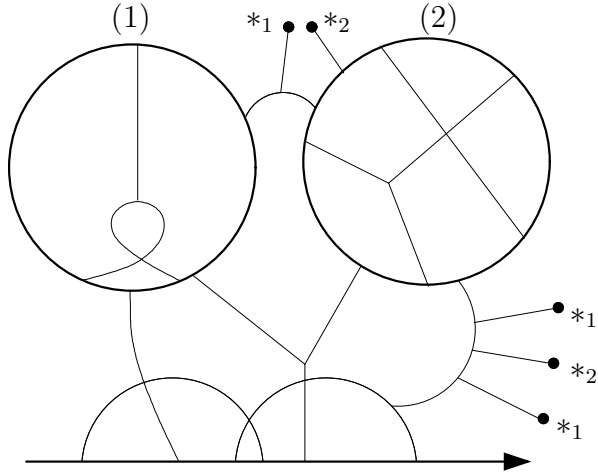

\[ 
  \printname{JDexample}
  \setlength{\unitlength}{0.8\standardunitlength}
  \begin{array}{c}  \hspace{-1.7mm}
    \raisebox{-8pt}{\input draws/JDexample.tex }
      \hspace{-1.9mm}
  \end{array}
 \]
\caption{Example of a Jacobi diagram based on a five component skeleton. By 
convention, internal trivalent vertices are oriented counterclockwise, and the
ordering of components is obtained by concatenating the separate 
orderings of isotypical component subsets, with the oriented intervals 
preceding the oriented circles, which are succeeded by the colors.} 
\end{figure}

The collection of all Jacobi diagrams based on a skeleton $\M$ will be denoted
by ${\D}(\M)$. Its subcollection consisting of diagrams having at least one 
external vertex in each connected component (called 
{\em boundary-connected diagrams}) will be denoted by $\D'(\M)$. Both $\D(\M)$ 
and $\D'(\M)$ are graded by half the number of vertices in a diagram.

\subsubsection{\bf Lie tensors.}

Let $\g$ be a finite-dimensional Lie algebra over a field $\k$ of 
characteristic zero. Let $\m$ be a component of a skeleton $\M$. We assign to 
$\m$ a $\g$-space, denoted\footnote{The Lie algebra $\g$ does not appear 
explicitly in the notation for the various tensor spaces, unless there is a 
risk of confusion.} $T(\m)$, which depends on the type of $\m$:
\begin{enumerate}
\item $\m$ is an oriented interval: $T(\m)$ is $T(\g)$, the tensor algebra 
over $\g$. $\g$ acts on $T(\m)$ via the adjoint representation.
\item $\m$ is an oriented circle: $T(\m)=T(\g)_{\g}$, where 
\[ T(\g)_{\g}=T(\g) \left/ \Span\{xv-vx \mid x \in \g, v \in T(\g)\} \right. \]
is the space of coinvariants for the action of $\g$ on $T(\g)$.
The action of $\g$ on $T(\g)_{\g}$ is the zero action (inherited from the 
action of $\g$ on $T(\g)$).
\item $\m$ is a color: set $T(\m)=S(\g)$, the symmetric algebra over $\g$,
with the corresponding adjoint action of $\g$. 
\end{enumerate}

\begin{remark}
A short reminder about the symmetric algebra may be helpful here.
Since we work in characteristic zero, there are two alternative definitions,
or `models', for the symmetric algebra: either as the {\em quotient} of $T(\g)$
by the two-sided ideal generated by $x\ot y - y\ot x \, (x, y \in \g)$, 
or as the {\em subspace} of symmetric tensors in $T(\g)$, which is the direct 
sum of the images of the symmetrization operators
\[
	\sym^k : T^k(\g) \lra T^k(\g)\, , \quad 
	sym^k(x_1\cdots x_k) = \frac{1}{k!} 
	                      \sum_{\pi \in S_k} x_{\pi(1)}\cdots x_{\pi(k)} \, ,
\]
acting on the homogeneous components of $T(\g)$. The adjoint action of $\g$ 
on $T(\g)$ leaves invariant the subspace model and descends to the quotient 
model. The linear isomorphism between
these two models is the restriction to the subspace model of the natural 
projection from $T(\g)$ onto the quotient model. This isomorphism is a 
$\g$-map (it intertwines the adjoint actions), and can be used to export
the product structure from the quotient model to the subspace model.
Another (equivalent) way to define the product in the subspace 
model is to multiply symmetric tensors as elements of $T(\g)$ and symmetrize
the result. In the sequel we use the subspace model of $S(\g)$.   
\end{remark}

\begin{definition}
\lbl{def:TM}
The tensor space $T(\M)$ associated to the skeleton $\M$ is the
tensor product (over $\k$) 
\[ 
T(\M)=\bigotimes_{\m \in \M} T(\m) \, , 
\]
where the order of factors is the order of components of $\M$.  
Being the tensor product of $\g$-spaces, $T(\M)$ carries an action
of $\g$.
\end{definition}

\subsubsection{\bf Mapping Jacobi diagrams into Lie tensors. } 
\label{__sect:DtoTensor}

Let $V$, $W$ be finite-dimensional linear spaces over $\k$. We 
make below a frequent use of the identification between
the spaces $\Hom (V , W)$ and $W \ot V^{\star}$, via the isomorphism
\begin{equation}
\label{eqn:HomEqTensor}
	\left(h:V \lra W\right) \ \  \mapsto \ \ \sum_i h(v_i) \ot v^i \, ,
\end{equation}
where $\{v_i\}$ is any basis of $V$, and $\{v^i\}$ is the dual basis of 
$V^{\star}$.

\smallskip
We now assume that the  Lie algebra $\g$ is equipped with a {\em metric} ---
a $\g$-invariant symmetric non-degenerate bilinear form $t$. (Prominent 
examples are semisimple Lie algebras with their Killing forms.) 
By \eqref{eqn:HomEqTensor} we may equivalently regard $t$ either as a map
$t:\g\ot \g \ra \k$, often denoted $(\cdot , \cdot)$, as a symmetric tensor 
$t \in \g^{\star} \ot \g^{\star}$, or as a map $t:\g \ra \g^{\star}$. The
non-degeneracy assumption implies that the last map is actually an isomorphism.
The inverse isomorphism,
$ t^{-1}:\g^{\star} \tilde{\ra} \g $,
may be identified with a tensor
$t^{-1} \in \g \ot \g^{\star \star}  \cong \g \ot \g$, which is again
symmetric. 

\smallskip
The Lie bracket of $\g$ is a map $\tilde{s}:\g \ot \g \ra \g$, that is, a 
tensor $\tilde{s} \in \g \ot \g^{\star} \ot \g^{\star}$. Applying the
map $1 \ot t^{-1} \ot t^{-1}$ to $\tilde{s}$ yields a tensor
$s \in \g \ot \g \ot \g$. It is simple to check that the $\g$-invariance
of the metric $t$ is equivalent to the tensor $s$ being invariant under
any {\em cyclic} permutation of its components.  

We note that $t$, $t^{-1}$, and $s$ are all $\g$-invariant tensors.

It is useful to have formulas for the above tensors in terms of
a basis $\{e_i\}$ of $\g$. Let $\{f^i\}$ be the $t$-dual basis of $\g$, so
that for all $x\in \g$, $x = (x, f^i) e_i$ (using the summation convention).
We have
\begin{align}
t^{-1} & = e_i \ot f^i = f^i \ot e_i \, , \label{eqn:CasInBasis} \\
s & = [e_i , e_j] \ot f^i \ot f^j \, . \label{eqn:BrcktInBasis}  
\end{align}

\begin{remark}
\label{rem:InvMetric}
The tensor $t^{-1}$ is often called the {\em Casimir element} of $\g$. We 
reserve this terminology for the closely related element of the universal
enveloping algebra of $\g$ (see Example ~\ref{exmpl:CasimirDgrm}), and call 
$t^{-1}$ the (inverse) metric tensor of $\g$. We sometimes denote it by  
$\Omega_{\g}$.
\end{remark}

We assign graphical representations to $t^{-1}$ and $s$. $t^{-1}$ is 
represented by the {\em atomic bivalent graph}
\[

  \printname{MetricGraph}
  \setlength{\unitlength}{0.5\standardunitlength}
  \begin{array}{c}  \hspace{-1.7mm}
    \raisebox{-8pt}{\input draws/MetricGraph.tex }
      \hspace{-1.9mm}
  \end{array}
 \, \raisebox{-5pt}{,}
\]
and $s$ is represented by the {\em atomic trivalent graph}
\[ 

  \printname{BracketGraph}
  \setlength{\unitlength}{0.5\standardunitlength}
  \begin{array}{c}  \hspace{-1.7mm}
    \raisebox{-8pt}{\input draws/BracketGraph.tex }
      \hspace{-1.9mm}
  \end{array}
 \, . 
\]
The graph with the opposite orientation of the trivalent vertex represents 
the tensor $(-s)$.

\smallskip
Let $D$ be a {\em Jacobi graph} --- a graph having either univalent vertices 
or oriented trivalent
vertices (basic example: the body of a Jacobi diagram based on some 
skeleton). We associate to $D$ the tensor space
\[ 
T(D)=\bigotimes_{\substack{\text{univalent vertices} \\ v \text{ of } D}} \g \,
.
\]
$T(D)$ is a $\g$-space, the $\g$-action being the tensor power of the adjoint 
action of $\g$ on itself.

\smallskip
Construct a tensor $\T_{\g}(D) \in T(D)$  by the following procedure:
\begin{itemize}
\item Separate the graph $D$ into a union of bivalent and trivalent atomic 
graphs, by putting separation points on some of the arcs. (It is possible to 
put more than one separation point on a single arc.)

\begin{figure}[htpb]
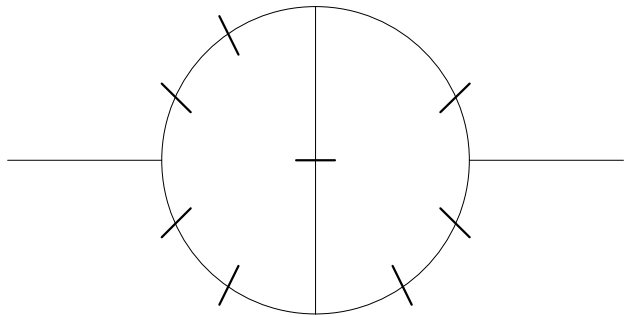

\[ 

  \printname{SeparatedGraph}
  \setlength{\unitlength}{0.8\standardunitlength}
  \begin{array}{c}  \hspace{-1.7mm}
    \raisebox{-8pt}{\input draws/SeparatedGraph.tex }
      \hspace{-1.9mm}
  \end{array}
 
\]
\caption{Separating a graph into atomic graphs.} \lbl{fig:SeparatedGraph}
\end{figure}

\item Attach to each atomic graph in the separation of $D$ the
tensor it represents, to get an intermediate tensor in 
\[ \bigotimes_{\substack{\text{atomic subgraphs } G \\
                          \text{making up } D}} T(G) . \]
Note that this step is unambiguous due to the symmetry of $t^{-1}$ and the 
cyclic invariance of $s$.
\item For each separation point, use the metric to contract the two components
of the intermediate  tensor belonging to the pair of spaces standing at the 
two sides of the point. The resulting tensor is in $T(D)$; this is 
$\T_{\g}(D)$.
\end{itemize}

\begin{proposition}
\lbl{prop:TgD}
Every Jacobi graph $D$ can be separated into atomic subgraphs, and 
$\T_{\g}(D)$ is independent of the separation used to compute it.
\end{proposition}
\begin{proof}
To get one possible separation, put one separation point on each arc. 
In order to prove that $\T_{\g}(D)$ is independent of the 
separation, it is enough to check that:
\begin{itemize}
\item
The tensor $s'$ computed according to the separation 
\[

  \printname{BracketGraphWithSepPt}
  \setlength{\unitlength}{0.5\standardunitlength}
  \begin{array}{c}  \hspace{-1.7mm}
    \raisebox{-8pt}{\input draws/BracketGraphWithSepPt.tex }
      \hspace{-1.9mm}
  \end{array}
 
\]
of the trivalent atomic graph is equal to $s$. 
\item
The tensor corresponding to the separation
\[

  \printname{MetricGraphWithSepPt}
  \setlength{\unitlength}{0.5\standardunitlength}
  \begin{array}{c}  \hspace{-1.7mm}
    \raisebox{-8pt}{\input draws/MetricGraphWithSepPt.tex }
      \hspace{-1.9mm}
  \end{array}
  
\]
is equal to $t^{-1}$.
\end{itemize}
Indeed, $s'$ is the contraction
of the second component of $[e_i, e_j] \ot f^i \ot f^j$ with the first
component of $e_k \ot f^k$, yielding
\[
(f^i, e_k) [e_i, e_j] \ot f^k \ot f^j = s \, ,
\]
since $(f^i, e_k) f^k = f^i$. The second assertion is easier.
\end{proof}

Since all the building blocks of $\T_{\g}(D)$ are $\g$-invariant and 
contractions are $\g$-maps, we have:
\begin{proposition}
\lbl{prop:TgImgInvar}
The tensors $\T_{\g}(D)$ are $\g$-invariant.
\end{proposition}

\begin{exercise}
Show that for the graph $D$ in Figure ~\ref{fig:SeparatedGraph},
\[
\T_{\g}(D) = [[e_i, e_j], e_k] \ot [f^i, [f^j, f^k]] \, .
\]
\end{exercise}

\smallskip
We apply the preceding construction to Jacobi diagrams. Let $\M$ be a skeleton 
and $D$ be a diagram in $\D(\M)$. We also denote its body by $D$. There is a 
natural $\g$-map $T(D) \ra T(\M)$, obtained by:
\begin{itemize}
\item Mapping the tensor product of $\g$'s along an oriented interval
into $T(\g)$. (The order of the factors in the tensor product is from 
the first vertex to the last vertex, according to the order dictated by
the orientation of the interval.)
\item For an oriented circle, fixing arbitrarily a starting vertex, mapping
into $T(\g)$ as for an oriented interval, and projecting to $T(\g)_{\g}$.
\item For a color, taking an arbitrary ordering of the vertices of this
color, mapping into $T(\g)$, and symmetrizing to get a tensor in $S(\g)$.
\end{itemize}
The resulting map is independent of the arbitrary choices 
that have been made. The composition
\[ D \stackrel{\T_{\g}}{\lra} T(D) \lra T(M) \]
yields a $\k$-linear map $\T_{\g}:\Span_{\k}(\D(\M)) \ra T(\M)$. 

\subsubsection{\bf Algebraic and diagrammatic Lie relations.}
\label{__sect:TgDef}

Let $\M$ be a skeleton. We form the $\g$-space $L_{\g}(\M)$ (usually 
denoted $L(\M)$), which is constructed like $T(\M)$ by dealing with each 
component $\m$ of $\M$ separately:
\begin{enumerate}
\item $\m$ is an oriented interval: $L(\m)$ is $\Ug$, the universal  
enveloping algebra of $\g$.
\item $\m$ is an oriented circle: $L(\m)=\Ug_{\g}$, where 
\[ 
\Ug_{\g}=\Ug \left/ \Span\{xv-vx \mid x \in \g, v \in \Ug\} \right.  
\]
is the space of coinvariants for the adjoint action of $\g$ on $U(\g)$.
\item $\m$ is a color: $L(\m)=S(\g)$ (no change).
\end{enumerate}
$\g$ acts on the each space $L(\m)$ by the appropriate adjoint actions, hence
also on
\[ 
L(\M) = \bigotimes_{\m \in \M} L(\m) \,. 
\]

\smallskip
The projection from $T(\M)$ to $L(\M)$ reflects moding out by the relations 
that make $\g$ into a Lie algebra, rather than just a vector space with a 
product. We next describe the analogous relations for Jacobi diagrams.

\smallskip
Consider the subspace of $\Span_{\k}(\D(\M))$ which is spanned by the following 
relations\footnote{As for the 4T relations, all the diagrams in a relation
share an identical part which is not shown in the figure.}:
\begin{enumerate}
\item Antisymmetry of internal vertices:
\[

  \printname{StraightY}
  \setlength{\unitlength}{0.25\standardunitlength}
  \begin{array}{c}  \hspace{-1.7mm}
    \raisebox{-8pt}{\input draws/StraightY.tex }
      \hspace{-1.9mm}
  \end{array}
 + 
  \printname{TwistedY}
  \setlength{\unitlength}{0.25\standardunitlength}
  \begin{array}{c}  \hspace{-1.7mm}
    \raisebox{-8pt}{\input draws/TwistedY.tex }
      \hspace{-1.9mm}
  \end{array}
 = 0 \, . 
\]
\item IHX relations:
\[

  \printname{Igraph}
  \setlength{\unitlength}{0.25\standardunitlength}
  \begin{array}{c}  \hspace{-1.7mm}
    \raisebox{-8pt}{\input draws/Igraph.tex }
      \hspace{-1.9mm}
  \end{array}
 - 
  \printname{Hgraph}
  \setlength{\unitlength}{0.25\standardunitlength}
  \begin{array}{c}  \hspace{-1.7mm}
    \raisebox{-8pt}{\input draws/Hgraph.tex }
      \hspace{-1.9mm}
  \end{array}
 + 
  \printname{Xgraph}
  \setlength{\unitlength}{0.25\standardunitlength}
  \begin{array}{c}  \hspace{-1.7mm}
    \raisebox{-8pt}{\input draws/Xgraph.tex }
      \hspace{-1.9mm}
  \end{array}
 = 0 \, . 
\]
\item STU relations:
\[

  \printname{Sgraph}
  \setlength{\unitlength}{0.25\standardunitlength}
  \begin{array}{c}  \hspace{-1.7mm}
    \raisebox{-8pt}{\input draws/Sgraph.tex }
      \hspace{-1.9mm}
  \end{array}
 = 
  \printname{Tgraph}
  \setlength{\unitlength}{0.25\standardunitlength}
  \begin{array}{c}  \hspace{-1.7mm}
    \raisebox{-8pt}{\input draws/Tgraph.tex }
      \hspace{-1.9mm}
  \end{array}
 - 
  \printname{Ugraph}
  \setlength{\unitlength}{0.25\standardunitlength}
  \begin{array}{c}  \hspace{-1.7mm}
    \raisebox{-8pt}{\input draws/Ugraph.tex }
      \hspace{-1.9mm}
  \end{array}
 \, . 
\]
(The bold line is part of a 1-dimensional component of $\M$.)
\end{enumerate}
 
\begin{definition}
$\A(\M)$ is the quotient space 
\[ \A(\M) = \Span_{\k}\left(\D(\M)\right) \left/
           \left\{\text{$AS$, $IHX$, and $STU$ relations}\right\} \right. . \]
$\A'(\M)$ is defined similarly as a quotient of $\Span_{\k}(\D'(\M))$. 
\end{definition}

Since all the relations are homogeneous (the diagrams appearing in a 
relation have the same number of vertices), $\A(\M)$ and $\A'(\M)$ are graded 
vector spaces. We usually work with their graded completions (with no change of
notation).

\begin{remark}
The spaces $\A'(\ua)$ ($\ua$ is a single oriented interval 
skeleton) and 
$\A'(\ast)$ ($\ast$ is a single color) are the $\A$  
and $\B$ spaces of ~\cite{Bar-Natan:Vassiliev}, 
respectively. We will reserve the symbols $\A$ and $\B$ for $\A(\ua)$ 
and $\A(\ast)$, respectively.
\end{remark}

The following proposition 
(compare ~\cite[Proposition 2.11]{Bar-Natan:Vassiliev}) 
gives a precise meaning to the intuitive analogy between the diagrammatic 
relations and the algebraic ones.
\begin{proposition}
\lbl{prop:TgWellDef}
The map  $\T_{\g}:\Span_{\k}(\D(\M)) \ra T(\M)$ 
descends to a map $\T_{\g}:\A(\M) \ra L(\M)$. The image of $\T_{\g}$ lies in 
$L(\M)^{\g}$, the space of $\g$-invariant tensors in $L(\M)$.
\end{proposition} 

\begin{proof}
We prove that each IHX relation is mapped to zero by $\T_{\g}$. 
Each diagram in the relation has an 'internal' part
(I, H, or X), and an 'external' part, which is the same
for all the participating diagrams. By construction, $\T_{\g}$ of such a 
diagram can be computed by computing first the tensor corresponding to the 
internal part (call it $T_I$, $T_H$, or $T_X$) and the one corresponding to the
external part (call it $T_E$), and then contracting these two tensors and 
mapping the resulting tensor to  $L(\M)$. Clearly, the operation of 
contraction against $T_E$ followed by mapping to $L(\M)$ is linear. Thus, we 
need only to show that the signed sum of tensors computed for the internal 
parts, $T_I - T_H + T_X$, is zero.

\smallskip
We draw the internal parts of an IHX relation in a slightly different way
than before, 
separate each graph into two trivalent atomic graphs, and number the 
univalent vertices and the separation point:
\[

  \printname{TgIgraph}
  \setlength{\unitlength}{0.5\standardunitlength}
  \begin{array}{c}  \hspace{-1.7mm}
    \raisebox{-8pt}{\input draws/TgIgraph.tex }
      \hspace{-1.9mm}
  \end{array}
 -  

  \printname{TgHgraph}
  \setlength{\unitlength}{0.5\standardunitlength}
  \begin{array}{c}  \hspace{-1.7mm}
    \raisebox{-8pt}{\input draws/TgHgraph.tex }
      \hspace{-1.9mm}
  \end{array}
 +

  \printname{TgXgraph}
  \setlength{\unitlength}{0.5\standardunitlength}
  \begin{array}{c}  \hspace{-1.7mm}
    \raisebox{-8pt}{\input draws/TgXgraph.tex }
      \hspace{-1.9mm}
  \end{array}
 = 0 \, .  
\]
We denote by $\g_i$ the copy of $\g$ standing at the $i$th position.
The I-graph is mapped to the contraction at $\g_5$ of the
tensors
\[
f^i \ot f^j \ot [e_i, e_j] \in \g_1 \ot \g_2 \ot \g_5 \, , \quad
f^k \ot f^l \ot [e_k, e_l] \in \g_5 \ot \g_3 \ot \g_4 \, ,
\]  
which is equal to 
\[
f^i \ot f^j \ot f^l \ot [[e_i, e_j], e_l] 
\in \g_1 \ot \g_2 \ot \g_3 \ot \g_4 \, .
\]
This is the tensor in $g^{\ot 4}$ representing the map 
$(x , y, z) \in \g^{\ot 3} \mapsto [[x, y], z]$ (more precisely, the tensor
represents the composition of this map with the map 
$(t^{-1})^{\ot 3}: (\g^{\star})^{\ot 3} \ra \g^{\ot 3}$.)
In the same fashion, the H-graph corresponds to the iterated bracket 
$[x, [y, z]]$ and  the X-graph corresponds to $[y, [x, z]]$. Therefore, 
the signed sum of the IHX relation is zero by the Jacobi identity in $\g$. 

Similarly, an STU relation is mapped to zero by $\T_{\g}$ because of the defining 
identity of the universal enveloping algebra, $[x , y] = xy - yx$, and an AS 
relation is mapped to zero because of the antisymmetry of the bracket.  

\smallskip
The second assertion of the proposition follows directly from 
Proposition ~\ref{prop:TgImgInvar}.
\end{proof}

\begin{example}
\lbl{exmpl:CasimirDgrm}
The diagram 
\[
	C = \raisebox{4pt}{$
  \printname{HalfCirc}
  \setlength{\unitlength}{0.25\standardunitlength}
  \begin{array}{c}  \hspace{-1.7mm}
    \raisebox{-8pt}{\input draws/HalfCirc.tex }
      \hspace{-1.9mm}
  \end{array}
$} \in \A(\ua)
\]
is mapped by $\T_{\g}$ to the {\em Casimir element}  
$C_{\g} = \sum_{i} e_i f^i \in \CUg$, while the diagram
\[
	\Omega = 
  \printname{BCup}
  \setlength{\unitlength}{0.25\standardunitlength}
  \begin{array}{c}  \hspace{-1.7mm}
    \raisebox{-8pt}{\input draws/BCup.tex }
      \hspace{-1.9mm}
  \end{array}
 \in \A(\ua_1, \ua_2)
\] 
is mapped to the (inverse) metric tensor 
$\Omega_{\g} \in (\g\ot\g)^{\g} \subset (U(\g)\ot U(\g))^{\g}$ 
(see Remark ~\ref{rem:InvMetric}). We call $C$ and $\Omega$ the diagrammatic
Casimir element and the diagrammatic metric tensor, respectively.
\end{example}

\begin{remark}
The construction of the maps $\T_{\g}$ can be extended with no change 
to the case when $\g$ is a {\em Lie superalgebra}.
\end{remark}

\subsection{Maps between $\A$-spaces} 
We discuss several maps between $\A$-spaces.

\subsubsection{\bf From chord diagrams to Jacobi diagrams.}
\label{__sect:ChordToJacobi} 

Let $\M$ be a $1$-dimensional skeleton. The obvious inclusion 
$\D^c(\M) \hookrightarrow \D(\M)$ descends to a map 
$\phi^c: \A^c(\M) \ra \A'(\M)$, because a 4T relation is seen to be equal to 
zero in $\A'(\M)$ by two applications of STU relations:
\[

  \printname{4TW}
  \setlength{\unitlength}{0.4\standardunitlength}
  \begin{array}{c}  \hspace{-1.7mm}
    \raisebox{-8pt}{\input draws/4TW.tex }
      \hspace{-1.9mm}
  \end{array}
 - 
  \printname{4TE}
  \setlength{\unitlength}{0.4\standardunitlength}
  \begin{array}{c}  \hspace{-1.7mm}
    \raisebox{-8pt}{\input draws/4TE.tex }
      \hspace{-1.9mm}
  \end{array}
 =

  \printname{Ygraph}
  \setlength{\unitlength}{0.4\standardunitlength}
  \begin{array}{c}  \hspace{-1.7mm}
    \raisebox{-8pt}{\input draws/Ygraph.tex }
      \hspace{-1.9mm}
  \end{array}
 =

  \printname{4TN}
  \setlength{\unitlength}{0.4\standardunitlength}
  \begin{array}{c}  \hspace{-1.7mm}
    \raisebox{-8pt}{\input draws/4TN.tex }
      \hspace{-1.9mm}
  \end{array}
 - 
  \printname{4TS}
  \setlength{\unitlength}{0.4\standardunitlength}
  \begin{array}{c}  \hspace{-1.7mm}
    \raisebox{-8pt}{\input draws/4TS.tex }
      \hspace{-1.9mm}
  \end{array}
 \, .
\]

\begin{proposition}
\lbl{prop:Achord-eq-Ajacobi}
The map $\phi^c: \A^c(\M) \ra \A'(\M)$ is an isomorphism. 
\end{proposition}
This is a slight generalization of ~\cite[Theorem 6]{Bar-Natan:Vassiliev}, 
where the case $\M = \bigcirc$ is proved. 

\subsubsection{\bf Circling an interval.} 
\label{__sect:IntrvlToCirc}

The map $\D(\ua) \ra \D(\bigcirc)$, defined by closing the interval to form 
a circle, clearly descends to a map $\phi^l: \A \ra \A(\bigcirc)$.
\begin{proposition}
\lbl{prop:Aline-eq-Acirc}
The map $\phi^l$ is an isomorphism.
\end{proposition}
We sketch the proof, following ~\cite[Lemma 3.1]{Bar-Natan:Vassiliev}. 
We need to show that the leftmost leg of a diagram in $\A(\ua)$ can be 
moved to the right without changing the diagram modulo the AS, IHX, and 
STU relations:
\begin{equation}
\label{eqn:legLtoR}

  \printname{LegLtoR}
  \setlength{\unitlength}{0.4\standardunitlength}
  \begin{array}{c}  \hspace{-1.7mm}
    \raisebox{-8pt}{\input draws/LegLtoR.tex }
      \hspace{-1.9mm}
  \end{array}
 
\end{equation}
Once this is known, it is immediate that the map $\D(\bigcirc)\ra \D(\ua)$ 
defined by `opening the circle at some point where no leg ends' descends 
to a well-defined map $\A(\bigcirc)\ra \A(\ua)$ (independent of the
opening point), which is an inverse of $\phi^l$. 

\smallskip
To prove \eqref{eqn:legLtoR}, attach three right turning hooks near each of 
the vertices except the leftmost vertex. Assign minus signs to hooks that 
are to the right of skeleton vertices. Disconnect the leftmost leg and 
consider the sum of diagrams obtained by connecting this leg to one of the 
hooks (taking the diagram with minus sign if the corresponding hook was 
assigned a minus sign). The terms in this sum can be grouped by arcs or by 
vertices. 
Grouping by arcs yields, using the sign choice and AS relations, the left 
hand side of \eqref{eqn:legLtoR}. Grouping by vertices, the sum in each 
group is seen to be zero (using IHX relations at internal vertices, and STU 
relations at skeleton vertices), thus proving our claim.  

\begin{remark}
\lbl{rem:PwrAcirc}
The same argument shows that if $\M = \ua \cup X$ and 
$\widetilde{\M} = \bigcirc \cup X$, where $X$ consists of 
circle components, then $\A(\M) \cong \A(\widetilde{\M})$. 
The situation when $\M$ has more than one interval component is more 
complicated; see ~\cite[Section 3.3]{Thurston:Duflo} for details.
\end{remark}

\subsubsection{\bf The PBW isomorphism.}

The PBW map
$\upchi:S(\g) \ra \Ug$ is the restriction to $S(\g)$ of the natural
projection of $T(\g)$ on $U(\g)$ (recall that we use the subspace model 
of $S(\g)$).
The PBW theorem asserts that $\upchi$ is a vector space isomorphism. 
Since $\upchi$ is a $\g$-map, it also induces
an isomorphism $\upchi: \CSg \tilde{\ra} \CUg$. 

\smallskip
We now construct the diagrammatic analogue of this map.
For a diagram $D \in \B$
having $k$ legs ending in vertices $v_1, \ldots, v_k$, we put
\begin{equation}
\label{eqn:DefChi} 
\upchi(D)=\frac{1}{k!}\sum_{\substack{\text{orderings }\pi \\
\text {of } v_1, \ldots, v_k}} \pi D \, , 
\end{equation}
where $\pi D \in \A$ is the diagram obtained from $D$ by arranging 
the vertices $v_1, \ldots, v_k$ along the skeleton interval according to $\pi$.
For example,
\[
\upchi\left(
  \printname{ThreeLegs}
  \setlength{\unitlength}{0.3\standardunitlength}
  \begin{array}{c}  \hspace{-1.7mm}
    \raisebox{-8pt}{\input draws/ThreeLegs.tex }
      \hspace{-1.9mm}
  \end{array}
\right) = 
\frac{1}{6}\left( 
\begin{array}{cccccc}
   & 
  \printname{ThreeLegs1}
  \setlength{\unitlength}{0.3\standardunitlength}
  \begin{array}{c}  \hspace{-1.7mm}
    \raisebox{-8pt}{\input draws/ThreeLegs1.tex }
      \hspace{-1.9mm}
  \end{array}
 & + &  
  \printname{ThreeLegs2}
  \setlength{\unitlength}{0.3\standardunitlength}
  \begin{array}{c}  \hspace{-1.7mm}
    \raisebox{-8pt}{\input draws/ThreeLegs2.tex }
      \hspace{-1.9mm}
  \end{array}
 & + & 
 
  \printname{ThreeLegs3}
  \setlength{\unitlength}{0.3\standardunitlength}
  \begin{array}{c}  \hspace{-1.7mm}
    \raisebox{-8pt}{\input draws/ThreeLegs3.tex }
      \hspace{-1.9mm}
  \end{array}
 \\
\\  
 + & 
  \printname{ThreeLegs4}
  \setlength{\unitlength}{0.3\standardunitlength}
  \begin{array}{c}  \hspace{-1.7mm}
    \raisebox{-8pt}{\input draws/ThreeLegs4.tex }
      \hspace{-1.9mm}
  \end{array}
 & + &  
  \printname{ThreeLegs5}
  \setlength{\unitlength}{0.3\standardunitlength}
  \begin{array}{c}  \hspace{-1.7mm}
    \raisebox{-8pt}{\input draws/ThreeLegs5.tex }
      \hspace{-1.9mm}
  \end{array}
 & + & 
 
  \printname{ThreeLegs6}
  \setlength{\unitlength}{0.3\standardunitlength}
  \begin{array}{c}  \hspace{-1.7mm}
    \raisebox{-8pt}{\input draws/ThreeLegs6.tex }
      \hspace{-1.9mm}
  \end{array}
 
\end{array}
\right) 
\, .
\]
It is easy to check that $\upchi:\B \ra \A$ is well-defined. The analogy 
between the Lie algebraic and diagrammatic maps is made precise by 
fitting them into a commutative square
\[
\begin{CD}
\B @>{\T_{\g}}>> S(\g)^{\g} \\
@V{\upchi}VV  @VV{\upchi}V \\
\A @>{\T_{\g}}>> \CUg \, . 
\end{CD}
\]
The analogy carries over to the isomorphism assertion.
\begin{proposition}
\lbl{prop:PBW-A}
The map $\upchi: \B \ra \A$ is an isomorphism. 
\end{proposition}
For later use, we review in detail the proof of this theorem, following  
~\cite[Section 5]{Bar-Natan:Vassiliev} with some small changes. We begin by 
introducing filtrations
\[
  \D_k(\M) = \Span_{k}\{D \in \D(\M) \text{ : $D$ has at most $k$ legs}\} 
\]
of the spaces $\Span_{\k}\{\D(\M)\}$, for any skeleton $\M$ (we need 
$M=\ua$ and $\M=\ast$ in this proof). Let $\I_k$ be
the subspace of $\D_k(\ua)$ which is spanned by AS, IHX, and STU relations
in which all the participating diagrams have at most $k$ legs. The quotient 
spaces $\A_k = \D_k(\ua)\left/\I_k\right.$ are naturally injected into $\A$ and
form a filtration of $\A$.

\smallskip
The plan of the proof is to construct a sequence of compatible 
maps $\sigma_k : \D_k(\ua) \ra \B$ satisfying:
\begin{enumerate}[($\Sigma$1)]
\item
$\sigma_k(\upchi(C)) = C$ for $C\in \D_k(\ast)$. \label{itm:sigma1}
\item
$\sigma_k(\I_k) = 0$. \label{itm:sigma2} 
\item
$ \upchi(\sigma_k(D)) - D \in \I_k$ for $D\in \D_k(\ua)$. \label{itm:sigma3} 
\end{enumerate}
These properties imply that the maps $\sigma_k$ descend to compatible maps 
$\sigma_k : \A_k \ra \B$, yielding a map 
$\sigma :\A \ra \B$ which is inverse to $\upchi : \B \ra \A$.

\smallskip
For $k=0, 1$ we define
\[
\sigma_0\left(
  \printname{D0LineDgrm}
  \setlength{\unitlength}{0.2\standardunitlength}
  \begin{array}{c}  \hspace{-1.7mm}
    \raisebox{-8pt}{\input draws/D0LineDgrm.tex }
      \hspace{-1.9mm}
  \end{array}
\right) 
= \raisebox{4pt}{$
  \printname{D0AstDgrm}
  \setlength{\unitlength}{0.2\standardunitlength}
  \begin{array}{c}  \hspace{-1.7mm}
    \raisebox{-8pt}{\input draws/D0AstDgrm.tex }
      \hspace{-1.9mm}
  \end{array}
$}
\, , \quad \text { and } \quad
\sigma_1\left(
  \printname{D1LineDgrm}
  \setlength{\unitlength}{0.2\standardunitlength}
  \begin{array}{c}  \hspace{-1.7mm}
    \raisebox{-8pt}{\input draws/D1LineDgrm.tex }
      \hspace{-1.9mm}
  \end{array}
\right) 
= 
  \printname{D1AstDgrm}
  \setlength{\unitlength}{0.2\standardunitlength}
  \begin{array}{c}  \hspace{-1.7mm}
    \raisebox{-8pt}{\input draws/D1AstDgrm.tex }
      \hspace{-1.9mm}
  \end{array}
 \, .
\]

\smallskip
Assume next that $k > 1$ and that $\sigma_{k-1}$ satisfying 
$\Sigma$\ref{itm:sigma1} -- $\Sigma$\ref{itm:sigma3} was already 
constructed.  

\smallskip
For $1 \leq i < k$ we denote by $u_i$ three related objects:
\begin{itemize}
\item
The transposition $u_i = (i, i+1) \in S_k$ swapping $i$ with $i+1$.
\item
The diagram
\begin{equation}
\label{eqn:uDgrm}
	u_i = 
  \printname{TranspDgrm}
  \setlength{\unitlength}{0.8\standardunitlength}
  \begin{array}{c}  \hspace{-1.7mm}
    \raisebox{-8pt}{\input draws/TranspDgrm.tex }
      \hspace{-1.9mm}
  \end{array}
 \in \D(\ua, \ua) \, .
\end{equation}
\item
The operator $u_i : \D_k(\ua) \ra \D_k(\ua)$ defined by `gluing $u_i$':
\begin{equation}
\label{eqn:uGlueOper}
u_i D = \left\{
\begin{array}{cl}

  \printname{GlueTranspDgrm}
  \setlength{\unitlength}{0.3\standardunitlength}
  \begin{array}{c}  \hspace{-1.7mm}
    \raisebox{-8pt}{\input draws/GlueTranspDgrm.tex }
      \hspace{-1.9mm}
  \end{array}
 & \text{if $D$ has $k$ legs}, \\
\\
0 & \text{otherwise}.
\end{array}
\right.
\end{equation}
\end{itemize}
In fact, such an operator can be associated to any $\pi \in S_k$,
by gluing the graphical representation of $\pi$. The result of applying
this operator to $D\in \D_k(\ua)$ is denoted by $\pi D$. These operators
define an action of $S_k$ on $\D_k(\ua)$: $(\pi \rho) D = \pi \rho D$
for $\pi, \rho\in S_k$..

In particular, the symbol $e$ stands for the identity element of $S_k$, 
its graphical representation
\begin{equation}
\label{eqn:eDgrm}
	e = 
  \printname{IdntDgrm}
  \setlength{\unitlength}{0.8\standardunitlength}
  \begin{array}{c}  \hspace{-1.7mm}
    \raisebox{-8pt}{\input draws/IdntDgrm.tex }
      \hspace{-1.9mm}
  \end{array}
 \in \D(\ua, \ua) \, ,
\end{equation}
and the `gluing $e$' operator on $\D_k(\ua)$ (which leaves unchanged diagrams
having precisely $k$ legs and maps other diagrams to zero).

In the same fashion, let $s_i$ ($1 \leq i < k$) denote the diagram
\begin{equation}
\label{eqn:sDgrm}
	s_i = 
  \printname{JoinDgrm}
  \setlength{\unitlength}{0.8\standardunitlength}
  \begin{array}{c}  \hspace{-1.7mm}
    \raisebox{-8pt}{\input draws/JoinDgrm.tex }
      \hspace{-1.9mm}
  \end{array}
 \in \D(\ua, \ua) \, ,	
\end{equation}
and its associated gluing operator, which maps $\D_k(\ua)$ into 
$\D_{k-1}(\ua)$. By STU, we have for $D\in \D_k(\ua)$ 
\begin{equation}
\label{eqn:OperSTU}
s_i D = e D - u_i D \text{ in } \A_k \, ,
\, 1\leq i < k \, .
\end{equation}

Let $\mathcal{M}$ be the free monoid on the $k-1$ letters 
$\{u_1, \ldots, u_{k-1}\}$.
Recall that the group ring $\mathbb{Z}[S_k]$ is the quotient of the monoid
ring $\mathbb{Z}[\mathcal{M}]$ by the two-sided ideal $\mathcal{K}$, which is 
generated by
\begin{align*}
\{u_i^2 - e \, : \, 1\leq i < k\} \cup 
\{u_i u_j - u_j u_i \, : \, |i-j| > 1\} \cup
\{u_i u_{i-1} u_i - u_{i-1} u_i u_{i-1} \, : \, 1 < i < k\} \, . 
\end{align*}
Fix $D\in \D_k(\ua)$ and define a map 
$\Gamma_D : \mathcal{M} \ra \D_{k-1}(\ua)$ by
\begin{equation}
\label{eqn:DefGammaD}
\Gamma_D(u_{i_1}\cdots u_{i_q}) 
= \sum_{p=1}^q s_{i_p} u_{i_{p+1}} \cdots u_{i_q} D\, .
\end{equation}

Using \eqref{eqn:OperSTU} in \eqref{eqn:DefGammaD} we get a telescopic
series 
\begin{equation}
\label{eqn:TelescopicTranspGammaD}
\begin{aligned}
\Gamma_D(u_{i_1}\cdots u_{i_q}) & = 
\sum_{p=1}^q (u_{i_{p+1}} \cdots u_{i_q} D 
             - u_{i_p} u_{i_{p+1}} \cdots u_{i_q} D)  \text{ mod } \I_k \\
& = D - u_{i_1} \cdots u_{i_q} D \, .
\end{aligned}
\end{equation}

Moding out the range of $\Gamma_D$ by $\I_{k-1}$ and extending by linearity
from $\mathcal{M}$ to $\mathbb{Z}[\mathcal{M}]$, 
we henceforth view $\Gamma_D$ as a map 
\[
\Gamma_D: \mathbb{Z}[\mathcal{M}]\ra \A_{k-1} \, .
\]
We observe that $\Gamma_D$ vanishes on the generators of 
$\mathcal{K}$. For example,
\begin{equation}
\label{eqn:VanishGammaD}
\begin{aligned}
\Gamma_D(u_i u_j - u_j u_i) & = s_i u_j D + s_j D - s_j u_i D - s_i D \\
&= 
  \printname{GammaD1a}
  \setlength{\unitlength}{0.3\standardunitlength}
  \begin{array}{c}  \hspace{-1.7mm}
    \raisebox{-8pt}{\input draws/GammaD1a.tex }
      \hspace{-1.9mm}
  \end{array}
 + 
  \printname{GammaD2a}
  \setlength{\unitlength}{0.3\standardunitlength}
  \begin{array}{c}  \hspace{-1.7mm}
    \raisebox{-8pt}{\input draws/GammaD2a.tex }
      \hspace{-1.9mm}
  \end{array}
 
- 
  \printname{GammaD1b}
  \setlength{\unitlength}{0.3\standardunitlength}
  \begin{array}{c}  \hspace{-1.7mm}
    \raisebox{-8pt}{\input draws/GammaD1b.tex }
      \hspace{-1.9mm}
  \end{array}
 - 
  \printname{GammaD2b}
  \setlength{\unitlength}{0.3\standardunitlength}
  \begin{array}{c}  \hspace{-1.7mm}
    \raisebox{-8pt}{\input draws/GammaD2b.tex }
      \hspace{-1.9mm}
  \end{array}
 \, .
\end{aligned}
\end{equation}
Apply twice STU relations from $\I_{k-1}$, once on the first and fourth 
diagrams and once on the two middle diagrams, to get
\[
\Gamma_D(u_i u_j - u_j u_i) = 
  \printname{GammaD3}
  \setlength{\unitlength}{0.3\standardunitlength}
  \begin{array}{c}  \hspace{-1.7mm}
    \raisebox{-8pt}{\input draws/GammaD3.tex }
      \hspace{-1.9mm}
  \end{array}
 - 
  \printname{GammaD3}
  \setlength{\unitlength}{0.3\standardunitlength}
  \begin{array}{c}  \hspace{-1.7mm}
    \raisebox{-8pt}{\input draws/GammaD3.tex }
      \hspace{-1.9mm}
  \end{array}
 
= 0 \, .
\]
Similarly, $\Gamma_D(u_i^2 - e) = 0$ by applying an AS relation and
$\Gamma_D(u_i u_{i-1} u_i - u_{i-1} u_i u_{i-1}) = 0$ by a combination
of STU and IHX relations. Together with the `cocycle property' of $\Gamma_D$,
\begin{equation}
\label{eqn:CocycleGammaD}
\Gamma_D(w_1 w_2) = \Gamma_{w_2 D}(w_1) + \Gamma_D(w_2)\, , \quad
w_1, w_2 \in \mathcal{M} \, ,
\end{equation}
which is easily derived from \eqref{eqn:DefGammaD}, it follows that
$\Gamma_D$ vanishes on any element of $\mathcal{K}$, hence it descends to 
a well-defined map 
$\Gamma_D : \mathbb{Z}[S_k] \ra \A_{k-1}$. This
means that we can unambiguously compute $\Gamma_D(\pi)$ for $\pi\in S_k$
by choosing any presentation of $\pi$ as a product of
transpositions and using \eqref{eqn:DefGammaD}. We can now write 
\eqref{eqn:TelescopicTranspGammaD} in the form
\begin{equation}
\label{eqn:TelescopicPermGammaD}
i(\Gamma_D(\pi)) = (D - \pi D) \text{ mod } \I_k \, , 
\quad D\in \D_k(\ua) \, ,  
\end{equation}
where $i: \A_{k-1}\ra\A_k$ is the natural injection.

\smallskip
As a last preparation for the definition of $\sigma_k$, we put, for 
$D\in \D_k(\ua)$, 
\begin{equation}
\label{eqn:DefLambdaD}
 \Lambda_D = \sigma_{k-1} \circ \Gamma_D : S_k \lra \B \, .
\end{equation}
Note that $\Lambda_D$ is well-defined as $\sigma_{k-1}(\I_{k-1}) = 0$
by $\Sigma$\ref{itm:sigma2}. By \eqref{eqn:CocycleGammaD},
\begin{equation}
\label{eqn:CocycleLambdaD}
\Lambda_D(\rho \pi) = \Lambda_{\pi D}(\rho) + \Lambda_D(\pi)\, , \quad
\rho, \pi \in S_k \, .
\end{equation}

\smallskip
Given a diagram $D\in \D_k(\ua)$ having precisely $k$ legs, let 
$D^{\ast} \in \D(\ast)$ be the diagram obtained by removing the interval
skeleton  of $D$ and coloring the detached legs with $\ast$. Set
\[
\sigma_k(D) = D^{\ast} + \frac{1}{k!}\sum_{\pi \in S_k} \Lambda_D(\pi) \, .
\]
If $D\in \D_k(\ua)$ has less than $k$ legs, then we simply set
$\sigma_k(D) = \sigma_{k-1}(D)$. We next check that $\sigma_k$ has the 
desired properties.

\smallskip
{\em Proof of $\Sigma$\ref{itm:sigma1}:} For $C \in \D_k(\ast)$ with $k$ legs, 
write $C=D^{\ast}$ for
some $D\in \D_k(\ua)$ (by picking an arbitrary ordering of the legs along
an added skeleton interval). Clearly, 
$\upchi(C) = \frac{1}{k!}\sum_{\pi \in S_k} \pi D$, so that
\begin{align*}
\sigma_k(\upchi(C)) 
& = \frac{1}{k!}\sigma_k\left(\sum_{\pi \in S_k} \pi D\right)
 = \frac{1}{k!}\left(\sum_{\pi \in S_k} (\pi D)^{\ast} + 
         \frac{1}{k!}\sum_{\pi, \rho \in S_k} \Lambda_{\pi D}(\rho)\right) \\
& = D^{\ast} +  \frac{1}{(k!)^2} 
               \sum_{\pi, \rho \in S_k} (\Lambda_D(\rho \pi) - \Lambda_D(\pi))
\qquad \text{by \eqref{eqn:CocycleLambdaD}} \\
& = C + \frac{1}{k!}\left(\sum_{\lambda \in S_k} \Lambda_D(\lambda) 
      - \sum_{\pi \in S_k} \Lambda_D(\pi)\right) = C \, .
\end{align*}
The case when $C$ has less than $k$ legs follows from $\Sigma$\ref{itm:sigma1}
for $\sigma_{k-1}$.

\smallskip
{\em Proof of $\Sigma$\ref{itm:sigma2}:} That 
$\sigma_k(\text{AS}) = 0$ and $\sigma_k(\text{IHX}) = 0$ essentially follows 
from the fact that the definitions of $\Gamma_D$ and the $D\mapsto D^{\ast}$ 
operation do not touch the internal vertices of $D$. To prove that 
$\sigma_k(\text{STU}) = 0$, it is enough to show that 
$\sigma_k(s_i D) = \sigma_k(D - u_i D)$, for a diagram $D\in \D_k(\ua)$ having 
$k$ legs. We compute:
\begin{align*}
\sigma_k(D) - \sigma_k(u_i D) & 
= D^{\ast} - (u_i D)^{\ast} + \frac{1}{k!}\sum_{\pi\in S_k} 
                              (\Lambda_D(\pi) - \Lambda_{u_i D}(\pi)) \\
& = \frac{1}{k!}\sum_{\pi\in S_k} 
	(\Lambda_D(\pi) - \Lambda_D(\pi u_i) + \Lambda_D(u_i)) \\
& = \Lambda_D(u_i) = \sigma_{k-1}(\Gamma_D(u_i)) = \sigma_k(s_i D) \, , 
\end{align*}
where the second equality used \eqref{eqn:CocycleLambdaD} and the last
equality used \eqref{eqn:DefGammaD}.

\smallskip
{\em Proof of $\Sigma$\ref{itm:sigma3}:} 
Working modulo $\I_k$ we have, for a diagram $\D\in \D_k(\ua)$ with $k$ legs,
\begin{alignat*}{2}
\upchi(\sigma_k(D)) & = 
	\upchi\left(D^{\ast} + \frac{1}{k!}\sum_{\pi\in S_k} \Lambda_D(\pi)\right) 
	&& \\
	& = \frac{1}{k!}\sum_{\pi\in S_k} \pi D 
	+ \frac{1}{k!}\sum_{\pi\in S_k} \upchi(\sigma_{k-1}(\Gamma_D(\pi))) 
	&& \qquad \text{by \eqref{eqn:DefLambdaD}} \\
	& = \frac{1}{k!}\sum_{\pi\in S_k} (\pi D + \Gamma_D(\pi))
	&& \qquad \text{by $\Sigma$\ref{itm:sigma3} for $\sigma_{k-1}$} \\
	& = \frac{1}{k!}\sum_{\pi\in S_k} (\pi D + D - \pi D) = D 
	&& \qquad \text{by \eqref{eqn:TelescopicPermGammaD}} \, ,
\end{alignat*}
as required. This completes our review of the proof of Proposition 
~\ref{prop:PBW-A}.

\begin{remark}
\lbl{rem:AfterPBW-A}
There is no difficulty in extending the definition of $\upchi$ to get a map 
from the space $\A(\M)$ based on a skeleton $\M = \ast \cup X$, to the space 
$\A(\widetilde{\M})$, where $\widetilde{\M} = \ua \cup X$ (a color component is
replaced by an oriented interval). The same proof works to show that this map 
is an isomorphism.
\end{remark}

\subsection{Algebraic structures on $\A$}

\subsubsection{\bf The Hopf algebra structure of $\Ug$.}

Recall that $\Ug$ has the structure of a Hopf algebra. The ingredients of this
 structure are:
\begin{enumerate}
\item
An associative product $\mu: \Ug^{\ot 2} \ra \Ug$ with a unit. The unit is 
considered as a linear 
map $\eta: \k \ra \Ug$ sending $1\in \k$ to $1\in \Ug$.
\item
A coassociative coproduct $\Delta: \Ug \ra \Ug^{\ot 2}$, defined by declaring 
elements of $\g$ to 
be primitive ($\Delta(x) = 1\ot x + x\ot 1, x\in \g$), with a counit 
$\eps: \Ug \ra \k$, determined 
by $\eps(1) = 1, \eps(x) = 0 \text{ for } x\in \g$. The coproduct and 
counit are morphisms of the algebra structure.
\item
The antipode map $S:\Ug \ra \Ug$, defined by 
\[
S(x_{1} \cdots x_{n}) = (-1)^{n} x_{n}\cdots x_{1} \, , 
\qquad x_{1} \cdots x_{n} \in \Ug \, . 
\]
It satisfies the identities
\begin{equation}
\label{eqn:antipod}
	\mu (S\ot\id)\Delta(x) = \eta \eps(x) = \mu (\id\ot S)\Delta(x) 
\, , x\in \Ug \, ,
\end{equation}
where $\id$ is the identity map.
\end{enumerate}
All the Hopf algebra axioms satisfied by these maps can be described by 
commutative diagrams of maps 
between tensor powers of $\Ug$. For example, the fact that the coproduct is 
an algebra morphism is described by the commutative diagram
\begin{equation}
\label{eqn:Ug-coprod-is-mult}
\begin{CD}
\Ug\ot \Ug @>{\Delta\ot \Delta}>> (\Ug\ot \Ug)\ot (\Ug\ot \Ug) \\
@V{\mu}VV  @VV{(\mu\ot \mu)(\id\ot \tau\ot \id)}V \\
\Ug @>{\Delta}>> \Ug\ot \Ug 
\end{CD} \, ,
\end{equation}
where $\tau:\Ug^{\ot 2} \ra \Ug^{\ot 2}$ is the flip morphism, 
$\tau(x\ot y) = y\ot x$.

\subsubsection{\bf Diagrammatic operations on $\A$-spaces.}
\label{__sect:DgrmOper}

Our aim is to define the analogous structure on $\A$. This structure, it turns
out, does not make $\A$ into a Hopf algebra in the usual algebraic 
sense\footnote{Thus, it should not be confused with the structure defined and 
discussed in ~\cite[Section 3]{Bar-Natan:Vassiliev}.}. 
The first step is to define a suitable 
notion of tensor powers of $\A$. Then we define diagrammatic maps between 
these tensor powers, which satisfy the Hopf algebra axioms in the sense that 
there exist commutative diagrams like \eqref{eqn:Ug-coprod-is-mult}.
 
\begin{definition}
Let $\M_1$, $\M_2$ be two skeletons. The {\em diagrammatic tensor 
product} of $\A(\M_1)$ and $\A(\M_2)$ is defined as $\A(\M_1 \udot \M_2)$ 
(disjoint union, in which the ordering of components is the concatenation of 
the separate orderings). It is denoted by $\A(\M_1) \bt \A(\M_2)$.
\end{definition}
There are associativity isomorphisms 
\[ (\A(\M_1) \bt \A(\M_2)) \bt \A(\M_3) \cong
   \A(\M_1) \bt (\A(\M_2) \bt \A(\M_3)) , \]
inherited from the identifications $(\M_1 \udot \M_2) \udot \M_3 = 
\M_1 \udot (\M_2 \udot \M_3)$. Thus, we can define unambiguously the 
diagrammatic tensor products of several $\A$-spaces. In particular, 
the diagrammatic tensor powers $\A(\M)^{\bt n}$ are well-defined.
These diagrammatic tensor powers are mapped by $\T_{\g}$ into (the 
$\g$-invariant part of) the algebraic tensor powers $L(\M)^{\ot n}$.
The space $\A(\emptyset)$ plays for diagrammatic tensor products the role
of the base field $\k$ in algebraic tensor products: 
$\A(\M) \bt \A(\emptyset) \cong \A(\M)$ for any skeleton $\M$.

\smallskip
We now define several general diagrammatic operations. Instances of these 
operations provide the Hopf algebra maps for $\A$.
 
\begin{definition}
Assume $\M$ is a skeleton and $\m_1$, $\m_2$ are two components of $\M$,
which are either both oriented intervals or both colors. The {\em juxtaposition}
map
\[ 
    \mu_{\m_1, \m_2} : \A(\M) \lra \A(\M-\{\m_2\}) 
\]
is defined by gluing the two components, with $\m_1$ before $\m_2$ (see Figure
~\ref{fig:gluing}).
\begin{figure}[htpb] 
\[ 
    
  \printname{DiagrammaticProduct}
  \setlength{\unitlength}{0.4\standardunitlength}
  \begin{array}{c}  \hspace{-1.7mm}
    \raisebox{-8pt}{\input draws/DiagrammaticProduct.tex }
      \hspace{-1.9mm}
  \end{array}
 
\]
\caption{Gluing oriented intervals ($\mu_{\ua_1, \ua_2}$) and gluing
colors ($\mu_{\ast_1, \ast_2}$). } \lbl{fig:gluing}
\end{figure} 
\end{definition}

In particular, we get a map $\mu:\A^{\bt 2} \ra \A$: 
\[ \mu=\mu_{\ua_1 , \ua_2}:\A(\ua_1 , \ua_2) \ra \A(\ua_1) \, . \]
The unit morphism $\eta:\A(\emptyset) \ra \A$ is the linear map sending a 
trivalent diagram to itself. The maps $\nu$ and $\eta$ make $\A$ into an 
'associative algebra with a unit'. For example, associativity of $\A$ is 
expressed by the commutativity of the diagram
\begin{equation}
\label{eqn:A-is-asoc}
\begin{CD}
\A(\ua_1, \ua_2, \ua_3) 
@>{\mu_{\ua_2 , \ua_3}}>>
\A(\ua_1, \ua_2) \\
@V{\mu_{\ua_1 , \ua_2}}VV  @VV{\mu_{\ua_2, \ua_3}}V \\
\A(\ua_1, \ua_3)
@>{\mu_{\ua_1 , \ua_3}}>>
\A(\ua_1)=\A(\ua_2)
\end{CD} \, ,
\end{equation}  
which holds already at the level of diagrams without relations: glueing the 
skeleton components in any of the two orders results in the same Jacobi 
diagram.

It is a simple observation that the product $\mu$ of $\A$ is indeed mediated 
by the maps $\T_{\g}$ to the product of $\Ug$, so that the square
\[
\begin{CD}
\A \bt \A @>{\T_{\g}}>> \Ug \ot \Ug \\
@V{\mu}VV  @VV{\mu}V \\
\A @>{\T_{\g}}>> \Ug
\end{CD}
\]
is commutative, and similarly for the unit maps.

\begin{definition}
\lbl{def:cabling}
Assume $\M$ is a skeleton and $\m$ is any component of $\M$. Denote by 
$\widetilde{\M}$ the skeleton obtained from $\M$ by adding a skeleton 
component $\tilde{\m}$ of the same type. (The ordering of $\widetilde{\M}$
is such that $\tilde{\m}$ is inserted immediately after $\m$.)  

The {\em cabling along $\m$},
\[ 
\Delta_{\m}:\A(\M) \lra \A(\widetilde{\M}) \, ,
\]
is defined by summing all the possible ways of splitting the legs ending on
$\m$ between $\m$ and $\tilde{\m}$.
\end{definition}
A proof that the map $\Delta_{\m}$, which is initially defined on Jacobi 
diagrams, does descend to $\A(\M)$ is demonstrated in Figure 
~\ref{fig:cabling}, when $\m$ is $1$-dimensional.
\begin{figure}[htpb]
\[ 
  \printname{CablingSTU}
  \setlength{\unitlength}{0.8\standardunitlength}
  \begin{array}{c}  \hspace{-1.7mm}
    \raisebox{-8pt}{\input draws/CablingSTU.tex }
      \hspace{-1.9mm}
  \end{array}
 \]
\caption{Cabling respects the STU relations.} \lbl{fig:cabling}
\end{figure}
(There is no problem with the AS and IHX relations, which pertain to internal 
vertices. In particular, the case when $\m$ is a color is trivial.) 

\begin{definition}
For $\m$ a component of a skeleton $\M$, the {\em component removal} map
\[ \eps_{\m}:\A(\M) \lra \A(\M-\{\m\}) \]
sends to zero any diagram having at least one external vertex on the component
$\m$, and leaves unchanged diagrams not touching $\m$, except for the removal 
of $\m$.  
\end{definition}

We put $\Delta = \Delta_{\ua_{1}}: \A(\ua_{1}) \ra A(\ua_{1}, \ua_{2})$, and 
$\eps = \eps_{\ua}:\A(\ua) \ra \A(\emptyset)$. $\Delta$ and 
$\eps$ make $\A$ 
into a coassociative coalgebra with a counit. The coassociativity of $\Delta$,
meaning that the square
\begin{equation}
\label{eqn:A-is-coasoc}
\begin{CD}
\A(\ua_1) @>{\Delta_{\ua_1}}>> \A(\ua_1, \ua_2) \\
@V{\Delta_{\ua_1}}VV  @VV{\Delta_{\ua_2}}V \\
\A(\ua_1, \ua_2)
@>{\Delta_{\ua_1}}>>
\A(\ua_1, \ua_{\tilde{1}}, \ua_2) = 
\A(\ua_1, \ua_2, \ua_{\tilde{2}})
\end{CD} 
\end{equation}  
commutes, follows from a simple argument: either way from the upper left 
corner to the lower right corner, starting with a Jacobi diagram $D$, the 
result is the sum of Jacobi diagrams obtained from $D$ by all possible 
splittings of its legs among the three skeleton intervals. Also, there is no 
difficulty in proving that the maps $\T_{\g}$ are coalgebra morphisms: 
\[ 
\Delta(\T_{\g}(D)) = \T_{\g}(\Delta(D)) \, , \, 
\eps(\T_{\g}(D)) = \T_{\g}(\eps(D)) \, \text{ for all } D\in \A \, . 
\]

Moreover, the coproduct of $\A$ is an algebra morphism, namely the following 
diagram of maps, which is analogous to the one in 
\eqref{eqn:Ug-coprod-is-mult}, is commutative:
\begin{equation}
\label{eqn:A-coprod-is-mult}
\begin{CD}
\A(\ua_1, \ua_2) @>{\Delta_{\ua_1}\bt\Delta_{\ua_2}}>> 
\A(\ua_1, \ua_2, \ua_3, \ua_4) \\
@V{\mu}VV                  @VV{\mu_{\ua_1, \ua_3}\bt\mu_{\ua_2, \ua_4}}V \\
\A(\ua_1) @>{\Delta}>>  \A(\ua_1, \ua_2) \, .
\end{CD}
\end{equation}
Here the diagrammatic tensor products of maps are informally defined as 
'applying each map to the legs ending on the corresponding skeleton 
component(s)'. This and similar definitions make sense because all the 
participating maps are defined by operations on legs of diagrams.

\begin{definition}
\lbl{def:CompReverse}
Let $\M$ be a skeleton and $\m$ be one of its components.
The {\em component reversal} map 
\[ S_{\m}:\A(\M) \lra \A(\M) \]
acts by reversing the orientation of $\m$ and multiplying the diagram by 
the parity ($\pm 1$) of the number of legs ending on $\m$.
(The orientation reversal applies only when $\m$ is a $1$-dimensional 
component.)
\end{definition}
It is easy to check that the maps $S_{\m}$ are well-defined and that the
map $S = S_{\ua}: \A \ra \A$ is $\T_{\g}$-compliant with the antipode of $\Ug$. 
Let us prove that $S$ satisfies the identity  
\begin{equation}
\label{eqn:A-antipode-identity} 
\mu(S \bt \id) \Delta (D) = \eta \eps(D) , \qquad D \in \A \, . 
\end{equation}
It is enough to show that if a diagram $D \in \A$ has at least one leg, then 
the left hand side of \eqref{eqn:A-antipode-identity} is zero. Assume that 
$D$ has $k$ legs and number the legs in ascending order from left to right. 
Write
\[ 
\Delta(D) = \sum_{J \subset \{1, \ldots, k\}} 
  \printname{AntipodeProof1}
  \setlength{\unitlength}{0.3\standardunitlength}
  \begin{array}{c}  \hspace{-1.7mm}
    \raisebox{-8pt}{\input draws/AntipodeProof1.tex }
      \hspace{-1.9mm}
  \end{array}
 \, ,
\]
where for $J = \{i_1, \ldots, i_r\}$, the $J$-box sends the legs numbered
$i_1, \ldots, i_r$ to the first component and the remaining legs to the
second component. Divide the sum into sums over pairs 
$\left( \{1\} \cup I \, , \, I \right)$, 
where $I$ is a subset of $\{2, \ldots, k\}$, so that
\[
\Delta(D) = \sum_{I}
\left( 
  \setlength{\unitlength}{0.3\standardunitlength}
  \begin{array}{c}  \hspace{-1.7mm}
    \raisebox{-2pt}{\input draws/AntipodeProof2a.tex }
    \hspace{-1.9mm}
  \end{array}
 + 

  \setlength{\unitlength}{0.3\standardunitlength}
  \begin{array}{c}  \hspace{-1.7mm}
    \raisebox{-2pt}{\input draws/AntipodeProof2b.tex }
    \hspace{-1.9mm}
  \end{array}
 \right) \, .
\]
Now our claim follows, as
\[
\mu(S \bt \id) \left( 
  \setlength{\unitlength}{0.3\standardunitlength}
  \begin{array}{c}  \hspace{-1.7mm}
    \raisebox{-2pt}{\input draws/AntipodeProof2a.tex }
    \hspace{-1.9mm}
  \end{array}
 + 
	             
  \setlength{\unitlength}{0.3\standardunitlength}
  \begin{array}{c}  \hspace{-1.7mm}
    \raisebox{-2pt}{\input draws/AntipodeProof2b.tex }
    \hspace{-1.9mm}
  \end{array}
 \right) =
(-1)^{|I| + 1} 
  \setlength{\unitlength}{0.3\standardunitlength}
  \begin{array}{c}  \hspace{-1.7mm}
    \raisebox{-2pt}{\input draws/AntipodeProof3a.tex }
    \hspace{-1.9mm}
  \end{array}
 + 
(-1)^{|I|}
  \setlength{\unitlength}{0.3\standardunitlength}
  \begin{array}{c}  \hspace{-1.7mm}
    \raisebox{-2pt}{\input draws/AntipodeProof3b.tex }
    \hspace{-1.9mm}
  \end{array}
 = 0 \, .
\]

\smallskip
We have proved most of   
\begin{proposition}
\lbl{prop:AisHopf} 
$(\A, \mu, \eta, \Delta, \eps, S)$ is a (diagrammatic) Hopf algebra.
For any metrized Lie algebra $\g$, the map $\T_{\g}:\A \ra \Ug$
is a morphism of Hopf algebras. Moreover, $\A$ is commutative and 
cocommutative.
\end{proposition}
\begin{proof}
Let $\tau : \A^{\bt 2} \ra \A^{\bt 2}$ be the flip morphism, which acts by 
interchanging the two skeleton intervals. The commutativity assertion means 
that $\mu \tau = \mu$. This follows from the argument in the proof of 
Proposition ~\ref{prop:Aline-eq-Acirc}: given $D \in \A(\ua_1, \ua_2)$, the 
legs of $\mu(D)$ which originally ended on $\ua_1$ can be moved to the right  
of the legs from $\ua_2$, without changing $\mu(D)$ modulo the relations of 
$\A$. In this process the ordering of the legs from $\ua_1$ along the interval
is preserved, hence the resulting diagram is precisely $\mu \tau(D)$. 
The cocommutativity of $\A$ is obvious.
\end{proof}

\begin{remark}
\lbl{rem:AlgDgrmCat}
Our algebro-diagrammatic terminology may be given a precise formal
meaning using the notion of algebraic structures in a $\k$-linear
tensor category. This notion is fairly clear: to define some sort of a 
structure on an object of the category (associative algebra object, Lie 
algebra object, Hopf algebra object, etc.), express the usual axioms for this 
structure as commutative diagrams of linear maps in the category 
$\Vect_{\k}$ of vector spaces over $\k$. These commutative diagrams make 
perfect sense in any  $k$-linear tensor category. 
(See ~\cite{EtingofKazhdan:QuantizationII} for more categorical details.)

\smallskip   
An appropriate category for the diagrammatic structures on $\A$ has as its 
objects the spaces
$\left\{\A(\M) : \text{$\M$ is a skeleton}\right\}$. 
A morphism from $\A(\M_1)$ to $\A(\M_{2})$ is a family 
$\{F_{k_1, k_2, \ldots }\}$ of elements of 
$\Span_{\k}\{\D(\M_1 \udot \M_{2})\}$ indexed by 
tuples of non-negative integers, one for each component of $\M_1$. 
$F_{k_1, k_2, \ldots}$ is required to be a sum of diagrams 
having $k_i$ legs on the $i$th component of $\M_1$. The morphism maps a 
Jacobi diagram $\D \in \D(\M_{1})$ with $k_i$ legs ending on the $i$th 
component of $\M_1$ to the image in $\A(\M_{2})$ of the sum of diagrams 
obtained by glueing  $\{F_{k_1, k_2, \ldots }\}$ to $D$ along $\M_1$. 

\smallskip
Glueing is different for each type of skeleton component of $\M_1$. For
an oriented interval component it is defined as in \eqref{eqn:uGlueOper}.
for the case of an oriented circle, see Figure ~\ref{fig:GlueAlongCirc}.
Glueing along a color component $\ast$ is done by summing over all ways of 
connecting the $\ast$-colored legs of $D$ to the $\ast$-colored legs of the 
glued diagrams (see 
~\cite[Definition 1.1]{Bar-NatanGaroufalidisRozanskyThurston:WheelsWheeling}).
\begin{figure}[htpb]
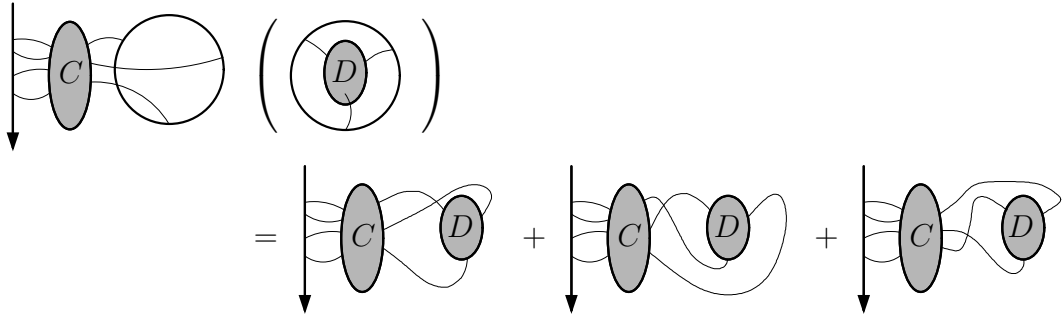

\begin{align*}

  \printname{GlueC}
  \setlength{\unitlength}{0.5\standardunitlength}
  \begin{array}{c}  \hspace{-1.7mm}
    \raisebox{-8pt}{\input draws/GlueC.tex }
      \hspace{-1.9mm}
  \end{array}
 & \left(
  \printname{GlueD}
  \setlength{\unitlength}{0.5\standardunitlength}
  \begin{array}{c}  \hspace{-1.7mm}
    \raisebox{-8pt}{\input draws/GlueD.tex }
      \hspace{-1.9mm}
  \end{array}
\right) \\
& = 
  \printname{CglueD1}
  \setlength{\unitlength}{0.5\standardunitlength}
  \begin{array}{c}  \hspace{-1.7mm}
    \raisebox{-8pt}{\input draws/CglueD1.tex }
      \hspace{-1.9mm}
  \end{array}
 + 
  \printname{CglueD2}
  \setlength{\unitlength}{0.5\standardunitlength}
  \begin{array}{c}  \hspace{-1.7mm}
    \raisebox{-8pt}{\input draws/CglueD2.tex }
      \hspace{-1.9mm}
  \end{array}
 + 
  \printname{CglueD3}
  \setlength{\unitlength}{0.5\standardunitlength}
  \begin{array}{c}  \hspace{-1.7mm}
    \raisebox{-8pt}{\input draws/CglueD3.tex }
      \hspace{-1.9mm}
  \end{array}
 
\end{align*}
\caption{Glueing diagrams along an oriented circle.} \lbl{fig:GlueAlongCirc}
\end{figure}

The resulting map from $\D(\M_{1})$ to $\A(\M_{2})$ 
is required to descend to $\A(\M_1)$, which amounts to respecting the STU 
relations for $1$-dimensional components. A typical example of such a 
morphism is the family $\{\sigma_k\}$ constructed in the proof of
Proposition ~\ref{prop:PBW-A}. 

Given two such morphisms, say
$F : \A(\M_1) \ra \A(\M_{2})$ and $G : \A(\M'_1) \ra \A(\M'_{2})$, it is clear
how to define their diagrammatic tensor product
$F\bt G: \A(\M_1)\bt\A(\M'_1) \ra \A(\M_{2})\bt \A(\M'_{2})$, 
by applying $F$ along $\M_1$ and $G$ along $\M'_1$.

\smallskip
We will not pursue this categorical approach. Nevertheless, we will 
freely use algebraic terminology in  connection with the diagrammatic objects 
and maps defined above and in the sequel.  
\end{remark}

\begin{exercise}
Find a family of elements of $\D(\ua_1, \ua_2)$ which gives rise to the
antipode $S:\A \ra \A$. 
\end{exercise}

\begin{remark}
\lbl{rem:CUgHopf}
The first part of the previous remark applies in a purely algebraic context: 
$\CUg$ is not a Hopf algebra with respect to the usual tensor product (because
the coproduct does not map $\CUg$ into $\CUg \ot \CUg$). However, it is 
a Hopf algebra object in the category whose objects are $(\Ug^{\ot n})^{\g}$. 
With this observation in mind, we see that $\T_{\g}: \A \ra \CUg$ is a morphism 
of Hopf algebra objects. 
\end{remark}

\smallskip
Having stressed that $\A$ is a diagrammatic Hopf algebra, it is nevertheless
useful to note that the diagrammatic product also makes $\A$ into an algebra in
the usual sense (that is, an algebra object in $\Vect_{\k}$)\footnote{This is 
not true for the diagrammatic coproduct; 
see  ~\cite[Definition 3.7]{Bar-Natan:Vassiliev} for the definition of a 
$\Vect_{\k}$-coproduct $\A\ra\A\ot\A$.}  
as follows. Starting
with the map $\D(\ua_1) \times \D(\ua_2) \ra \D(\ua_1, \ua_2)$ sending
a pair of diagrams to a `disconnected diagram' as in Figure 
~\ref{fig:DisconnectedDiagram}, we obtain
a map $\A\ot\A \ra \A\bt\A$. This map is injective, essentially 
because all the name-giving parts of diagrams participating in $\A$-relations
are connected. 
\begin{figure}[htpb]
\[ 

  \printname{DisconnectedDiagram}
  \setlength{\unitlength}{0.5\standardunitlength}
  \begin{array}{c}  \hspace{-1.7mm}
    \raisebox{-8pt}{\input draws/DisconnectedDiagram.tex }
      \hspace{-1.9mm}
  \end{array}
 
\]
\caption{An element of $\A\bt\A$ coming from $\A\ot\A$.}
\lbl{fig:DisconnectedDiagram}
\end{figure}
\\
The composition
\begin{equation}
\label{eqn:DisconnectedProd}
\begin{CD}
\A\ot\A \lra \A\bt\A @>{\mu}>> \A 
\end{CD}
\end{equation}
endows $\A$ with an algebra structure. 

\smallskip
Similarly, $A^{\bt n}$ is a $\Vect_{\k}$-algebra whose product is 
`juxtaposition along the $n$ skeleton components'; more precisely, it is the 
composition of maps
\[
\begin{CD}
\A(\ua_1, \ldots, \ua_n) \ot \A(\ua_{n+1}, \ldots, \ua_{2n})
@>>> \A(\ua_1, \ldots,  \ua_{2n}) 
@>{\mu_{\ua_n, \ua_{2n}}\circ\cdots\circ\mu_{\ua_1, \ua_{n+1}}}>> 
\A^{\bt n} \, . 
\end{CD}
\]
We denote by $D_1 D_2$ the product of $D_1, D_2 \in \A^{\bt n}$, and by
$[D_1, D_2]$ their commutator. 


\smallskip
Define the iterated coproduct map $\Delta^{(n)} : \A \ra \A^{\bt n}$ ($n\geq 1$)
as the composition
\begin{equation}
\label{eqn:AnVectProd}
\begin{CD}
\A(\ua_1) @>{\Delta_{\ua_1}}>> \A(\ua_1, \ua_2) 
@>{\Delta_{\ua_2}}>> \ \cdots \ @>{\Delta_{\ua_{n-1}}}>> 
\A(\ua_1, \ldots, \ua_n) \, .
\end{CD}
\end{equation}
We have the following commutativity result.
\begin{proposition}
\lbl{prop:AnInv}
For all $D\in \A$ and $E\in \A^{\bt n}$,
\[
	[\Delta^{(n)}(D), E] = 0 \, . 
\]
\end{proposition}
\begin{proof}
The proof is similar to that of Proposition ~\ref{prop:Aline-eq-Acirc} 
(which is closely related to the case $n=1$, where $\Delta^{(1)} = \id$). 
It is sufficient to give the proof for $n=2$. Assume $D$ is represented
by a diagram having $k$ legs. Write 
\[
\Delta(D) E = \sum_{I \subset \{1, \ldots, k-1\}} 
\left(
  \printname{DeltaDtimesEupL}
  \setlength{\unitlength}{0.3\standardunitlength}
  \begin{array}{c}  \hspace{-1.7mm}
    \raisebox{-8pt}{\input draws/DeltaDtimesEupL.tex }
      \hspace{-1.9mm}
  \end{array}
 + 
  \printname{DeltaDtimesEdownL}
  \setlength{\unitlength}{0.3\standardunitlength}
  \begin{array}{c}  \hspace{-1.7mm}
    \raisebox{-8pt}{\input draws/DeltaDtimesEdownL.tex }
      \hspace{-1.9mm}
  \end{array}

\right) \, .
\]

We attach three right-turning hooks around each vertex of $E$ (with minus 
signs on hooks to the right of skeleton vertices), and 
form the sum of all the diagrams obtained by connecting the free leg of
\[

  \printname{DeltaDtimesEfree}
  \setlength{\unitlength}{0.3\standardunitlength}
  \begin{array}{c}  \hspace{-1.7mm}
    \raisebox{-8pt}{\input draws/DeltaDtimesEfree.tex }
      \hspace{-1.9mm}
  \end{array}

\]
to one of the hooks. Grouping the summands by vertices shows that the sum is 
zero (in $\A^{\bt 2}$), and then regrouping by arcs shows that $\Delta(D) E$ 
is equal to
\[
\sum_{I \subset \{1, \ldots, k-1\}} 
\left(
  \printname{DeltaDtimesEupR}
  \setlength{\unitlength}{0.3\standardunitlength}
  \begin{array}{c}  \hspace{-1.7mm}
    \raisebox{-8pt}{\input draws/DeltaDtimesEupR.tex }
      \hspace{-1.9mm}
  \end{array}
 + 
  \printname{DeltaDtimesEdownR}
  \setlength{\unitlength}{0.3\standardunitlength}
  \begin{array}{c}  \hspace{-1.7mm}
    \raisebox{-8pt}{\input draws/DeltaDtimesEdownR.tex }
      \hspace{-1.9mm}
  \end{array}
\right) \, .
\]
We keep moving legs of $D$ across $E$ in this fashion, until
we get $E \Delta(D)$.  
\end{proof}  

For Lie algebras, the analogue of Proposition ~\ref{prop:AnInv} is  
just the easy fact that a $\g$-invariant tensor $X\in (U(\g)^{\ot n})^{\g}$
commutes with any tensor in the image of the iterated coproduct
$\Delta^{(n)}:U(\g) \ra U(\g)^{\ot n}$.

%% file: draws/MetricGraph.tex
%
\begingroup\makeatletter\ifx\SetFigFont\undefined%
\gdef\SetFigFont#1#2#3#4#5{%
  \reset@font\fontsize{#1}{#2pt}%
  \fontfamily{#3}\fontseries{#4}\fontshape{#5}%
  \selectfont}%
\fi\endgroup%
{\renewcommand{\dashlinestretch}{30}
\begin{picture}(1567,195)(0,-10)
\path(300,75)(1200,75)
\put(1275,0){\makebox(0,0)[lb]{\smash{{\mathmode{f^i  }}}}}
\put(0,0){\makebox(0,0)[lb]{\smash{{\mathmode{e_i}}}}}
\end{picture}
}

%% file: draws/BracketGraphWithSepPt.tex
%
\begingroup\makeatletter\ifx\SetFigFont\undefined%
\gdef\SetFigFont#1#2#3#4#5{%
  \reset@font\fontsize{#1}{#2pt}%
  \fontfamily{#3}\fontseries{#4}\fontshape{#5}%
  \selectfont}%
\fi\endgroup%
{\renewcommand{\dashlinestretch}{30}
\begin{picture}(1224,1239)(0,-10)
\path(12,1212)(612,612)
\path(1212,1212)(612,612)
\path(612,612)(612,12)
\thicklines
\path(237,837)(387,987)
\end{picture}
}

%% file: draws/MetricGraphWithSepPt.tex
%
\begingroup\makeatletter\ifx\SetFigFont\undefined%
\gdef\SetFigFont#1#2#3#4#5{%
  \reset@font\fontsize{#1}{#2pt}%
  \fontfamily{#3}\fontseries{#4}\fontshape{#5}%
  \selectfont}%
\fi\endgroup%
{\renewcommand{\dashlinestretch}{30}
\begin{picture}(924,381)(0,-10)
\path(12,183)(912,183)
\thicklines
\path(462,333)(462,33)
\end{picture}
}

%% file: draws/StraightY.tex
%
\begingroup\makeatletter\ifx\SetFigFont\undefined%
\gdef\SetFigFont#1#2#3#4#5{%
  \reset@font\fontsize{#1}{#2pt}%
  \fontfamily{#3}\fontseries{#4}\fontshape{#5}%
  \selectfont}%
\fi\endgroup%
{\renewcommand{\dashlinestretch}{30}
\begin{picture}(2424,2739)(0,-10)
\path(1212,762)(1212,12)
\path(1812,2112)(1212,1512)
\path(612,2112)(1212,1512)
\path(12,2712)(612,2112)
\path(2412,2712)(1812,2112)
\path(1212,762)(1212,1512)
\end{picture}
}

%% file: draws/Igraph.tex
%
\begingroup\makeatletter\ifx\SetFigFont\undefined%
\gdef\SetFigFont#1#2#3#4#5{%
  \reset@font\fontsize{#1}{#2pt}%
  \fontfamily{#3}\fontseries{#4}\fontshape{#5}%
  \selectfont}%
\fi\endgroup%
{\renewcommand{\dashlinestretch}{30}
\begin{picture}(2109,2435)(0,-10)
\path(1057,1809)(19,2408)
\path(1060,1809)(2097,2408)
\path(1055,1809)(1055,611)
\path(2095,12)(1057,611)
\path(1050,611)(12,12)
\end{picture}
}

%% file: draws/Hgraph.tex
%
\begingroup\makeatletter\ifx\SetFigFont\undefined%
\gdef\SetFigFont#1#2#3#4#5{%
  \reset@font\fontsize{#1}{#2pt}%
  \fontfamily{#3}\fontseries{#4}\fontshape{#5}%
  \selectfont}%
\fi\endgroup%
{\renewcommand{\dashlinestretch}{30}
\begin{picture}(2425,2117)(0,-10)
\path(611,1051)(13,12)
\path(611,1053)(12,2090)
\path(612,1049)(1809,1048)
\path(1814,1051)(2413,2088)
\path(2408,12)(1809,1049)
\end{picture}
}

%% file: draws/Xgraph.tex
%
\begingroup\makeatletter\ifx\SetFigFont\undefined%
\gdef\SetFigFont#1#2#3#4#5{%
  \reset@font\fontsize{#1}{#2pt}%
  \fontfamily{#3}\fontseries{#4}\fontshape{#5}%
  \selectfont}%
\fi\endgroup%
{\renewcommand{\dashlinestretch}{30}
\begin{picture}(2126,2117)(0,-10)
\path(695,680)(1444,677)
\path(912,900)(2113,2090)
\path(1222,900)(12,2090)
\path(916,900)(13,12)
\path(2114,19)(1214,904)
\end{picture}
}

%% file: draws/Sgraph.tex
%
\begingroup\makeatletter\ifx\SetFigFont\undefined%
\gdef\SetFigFont#1#2#3#4#5{%
  \reset@font\fontsize{#1}{#2pt}%
  \fontfamily{#3}\fontseries{#4}\fontshape{#5}%
  \selectfont}%
\fi\endgroup%
{\renewcommand{\dashlinestretch}{30}
\begin{picture}(2766,2820)(0,-10)
\path(1383,843)(1383,93)
\path(1983,2193)(1383,1593)
\path(783,2193)(1383,1593)
\thicklines
\path(33,93)(2733,93)
\blacken\path(2553.000,33.000)(2733.000,93.000)(2553.000,153.000)(2553.000,33.000)
\thinlines
\path(1383,843)(1383,1593)
\path(2583,2793)(1983,2193)
\path(183,2793)(783,2193)
\end{picture}
}

%% file: draws/Tgraph.tex
%
\begingroup\makeatletter\ifx\SetFigFont\undefined%
\gdef\SetFigFont#1#2#3#4#5{%
  \reset@font\fontsize{#1}{#2pt}%
  \fontfamily{#3}\fontseries{#4}\fontshape{#5}%
  \selectfont}%
\fi\endgroup%
{\renewcommand{\dashlinestretch}{30}
\begin{picture}(2766,2820)(0,-10)
\path(633,1443)(783,93)
\path(1833,1443)(1683,93)
\thicklines
\path(33,93)(2733,93)
\blacken\path(2553.000,33.000)(2733.000,93.000)(2553.000,153.000)(2553.000,33.000)
\thinlines
\path(483,2793)(633,1443)
\path(1983,2793)(1833,1443)
\end{picture}
}

%% file: draws/HalfCirc.tex
%
\begingroup\makeatletter\ifx\SetFigFont\undefined%
\gdef\SetFigFont#1#2#3#4#5{%
  \reset@font\fontsize{#1}{#2pt}%
  \fontfamily{#3}\fontseries{#4}\fontshape{#5}%
  \selectfont}%
\fi\endgroup%
{\renewcommand{\dashlinestretch}{30}
\begin{picture}(2235,820)(0,-10)
\put(1032.500,94.500){\arc{1401.000}{3.1409}{6.2839}}
\thicklines
\path(33,93)(2202,93)
\blacken\path(2022.000,33.000)(2202.000,93.000)(2022.000,153.000)(2022.000,33.000)
\end{picture}
}

%% file: Aarrow.tex
\section{The algebra $\protect\Aar$ of directed Jacobi diagrams}
\label{sect:Aarrow}

{\bf Section summary:} After recalling the definitions of Lie bialgebras,
the classical Yang-Baxter equation, the Drinfel'd double, and Manin triples,
we describe the standard Lie bialgebra structure of a complex semisimple 
Lie algebra. We then go to define $\Aar$-spaces as quotients of spaces spanned
by directed Jacobi diagrams modulo directed versions of the AS, IHX, and STU
relations, which are diagrammatic forms of the axioms of Lie bialgebra. We
give the construction of the maps from $\Aar$-spaces to spaces of (possibly
non-invariant) Lie tensors, which closely parallels the construction in the 
undirected case.

\smallskip
$\A$-spaces are mapped into $\Aar$-spaces by the $\iota$-maps, sending
an undirected diagram to the sum over all possible ways of directing 
its arcs. We verify that the $\iota$-maps are well-defined, show that 
elements in their images are invariant (in a suitable diagrammatic sense),
and conjecture that these maps are injective. Finally, we extend to the 
directed case the proof that the diagrammatic PBW map is an isomorphism.  

\subsection{Lie bialgebras} 
\label{_sect:LBA}

Let $\a$ be a finite-dimensional Lie algebra over $\k$. We recall the
basic definitions related to Lie bialgebras. We refer to 
~\cite[Sections 1.3--1.4, 2.1]{ChariPressley:Guide}
or ~\cite[Lectures 2--4]{EtingofSchiffmann:LectQuantumGrps} for full details.  

\subsubsection{\bf Definition of Lie bialgebras.} 

\begin{definition}
\lbl{def:lba}
A map $\delta: \a \ra \a \ot \a$ endows $\a$ with a Lie bialgebra
structure if it has the following properties:
\begin{enumerate}
\item $\delta$ satisfies the co-Jacobi identity
\[ 
(\id +\tau_3 +\tau_3^2)((\delta\ot\id)\delta(a)) =0 
\, , \quad a\in\a \, , 
\]
where $\tau_3:\a^{\ot 3} \ra \a^{\ot 3}$ is a
cyclic permutation of the tensor power.
\item $\delta$ is co-antisymmetric:
\[ 
\delta^{\op}(a) = -\delta(a) \, , \quad a\in \a \, .
\]
Here $\delta^{\op}$ is the composition of $\delta$ with the flip morphism 
$\tau: \a^{\ot 2} \ra \a^{\ot 2}$.
\item $\delta$ is a $1$-cocycle of $\a$ with coefficients in $\a \ot \a$;
that is, the following identity (the {\em cocycle identity}) holds\footnote{
The right hand side is computed in $U(\a)^{\ot 2}$, but the end result
belongs to $\a^{\ot 2}$.}:
\[ 
\delta([ab])=[1 \ot a + a \ot 1, \delta(b)] +
                [\delta(a), 1 \ot b + b \ot 1] \, , \quad a,b\in\a \, .
\]
\end{enumerate}
$\delta$ is called the {\em cobracket} of the Lie bialgebra $\a$. 
\end{definition}
The first two properties imply that the dual map 
$\delta^{\star}:\a^{\star} \ot \a^{\star} \ra \a^{\star}$
makes $\a^{\star}$ into a Lie algebra. The third property is a compatibility
criterion which the two Lie algebra structures, on $\a$ and $\a^{\star}$,
need to satisfy. 

\begin{example}
\label{exmpl:2D_lba}
We describe a non-trivial Lie bialgebra of dimension two, with
basis vectors $\{X, Y\}$. The bracket and cobracket are given by:
\[
[X, Y] = X \, , \quad \delta(X) = \frac{1}{2} X\land Y  \, , 
\delta(Y) = 0 \, ,
\]
where $X\land Y = X\ot Y - Y\ot X$.
\end{example}

\subsubsection{\bf Quasitriangular Lie bialgebras.}

For $r \in \a \ot \a$, the  
{\em classical Yang-Baxter equation (CYBE)} is the equation
\begin{equation}
\label{eqn:CYBE} 
[r^{12}, r^{13}] + [r^{12}, r^{23}] + [r^{13}, r^{23}] = 0 
\end{equation}
in $U(\a)^{\ot 3}$. ($r^{13}$, for example, is the image of $r$
under the injection 
$\a \ot \a \ra U(\a)^{\ot 3}, \, a \ot b \mapsto a \ot 1 \ot b$.) 

\begin{definition}
A Lie bialgebra $\a$ is called {\em quasitriangular} if there exists an 
element $r \in \a \ot \a$ (called the {\em $r$-matrix} of $\a$) such that
\begin{enumerate}
\item $r$ satisfies the CYBE.
\item $r+r^{21}$ is $\a$-invariant ($r^{21}$ denotes $\tau(r)$).  
\item $\delta$ is the coboundary of $r$, that is,
\[ 
 \delta(a)=[1 \ot a + a \ot 1, r] \text{ for all } a \in \a \, . 
\]
\end{enumerate}
\end{definition}
It can be checked that if $\a$ is a Lie algebra and $r \in \a \ot \a$ meets 
the first two conditions, 
then the last formula defines a cobracket on $\a$. 

\subsubsection{\bf The double of a Lie bialgebra.}
\label{__sect:DblLBA}

Let $\a$ be a Lie bialgebra. Recall that $\a$ acts on $\a^{\star}$ via the
coadjoint action: 
\[
(\ad^{\star}a)b = -b \circ \ad a \, , \quad
a \in \a \, , b \in \a^{\star} \, .
\]
The vector space $D\a= \a \oplus \a^{\star}$ has a natural structure of a 
quasitriangular Lie bialgebra, which is defined as follows:
\begin{itemize}
\item $\a$, $\a^{\star}$ are Lie subalgebras of $D\a$, and the bracket
of $a \in \a$, $b \in \a^{\star}$ is
\begin{equation} 
\label{eqn:MixedBracket}
[a , b] = (\ad^{\star}a)b - (1 \ot b)(\delta(a)) \, . 
\end{equation}
\item Let $r \in \a \ot \a^{\star} \subset D\a \ot D\a$ 
be the element corresponding to the identity map of $\a$. Then $r$ satisfies
the CYBE and  $r + r^{21}$ is $D\a$-invariant. The cobracket of 
$D\a$ is taken to be the coboundary of $r$. 
\end{itemize}
$D\a$ is  called the {\em (Drinfel'd) double} of the Lie bialgebra $\a$. It is
endowed with a metric\footnote{Recall from Section ~\ref{__sect:DtoTensor} 
that a metric is an invariant symmetric non-degenerate bilinear form.} given by
\[ 
(a+a', b+b') = a'(b) + b'(a) \, , \quad
a, b \in \a, \, a', b' \in \a^{\star} \, .
\] 

\subsubsection{\bf Manin triples.} 
\label{__sect:manin}

\begin{definition}
A triple $(\tg, \g_{+}, \g_{-})$ of Lie algebras is called a {\em Manin triple}
if $\tg$ is a metrized Lie algebra with a metric $(\cdot, \cdot)$, such that 
$\g_{\pm}$ are $(\cdot, \cdot)$-isotropic Lie subalgebras of $\tg$,
and $\tg = \g_{+} \oplus \g_{-}$. 
\end{definition}

The construction of the double of a Lie bialgebra $\a$ shows that
$(D\a, \a, \a^{\star})$ is a Manin triple. Conversely, if 
$(\tg, \g_{+}, \g_{-})$ is a (finite-dimensional) Manin triple, 
then $\g_{+}$ can be endowed with a Lie bialgebra structure such that 
\begin{equation}
\label{eqn:manin}
(\tg, \g_{+}, \g_{-}) \cong (D\g_{+}, \g_{+}, \g_{+}^{\star}) \, .  
\end{equation}

This is done in the following way: the non-degeneracy of the metric of $\tg$
and the isotropy of $\g_{\pm}$ imply that the paring 
$(\cdot, \cdot):\g_{+}\ot \g_{-} \ra \k$ is non-degenerate, hence 
it induces a vector space isomorphism $\phi:\g_{+}^{\star} \ra \g_{-}$. The Lie
bracket of $\g_{-}$ becomes, via this isomorphism, a putative cobracket
map $\delta: \g_{+} \ra \g_{+} \ot \g_{+}$, satisfying the first two 
properties of Definition ~\ref{def:lba}. The third property (the cocycle 
identity) can be verified by a direct computation. Furthermore, it can be 
checked that the vector space isomorphism 
\[
D\g_{+} \lra \tg : (x, y) \mapsto (x, \phi(y)) \, , \quad 
x\in \g_{+}\, , y\in \g_{+}^{\star} \, ,
\]
is a metric preserving Lie algebra map.

The isomorphism \eqref{eqn:manin} endows $\tg$, the ambient Lie
algebra of the Manin triple, with the structure of a quasitriangular Lie 
bialgebra. Its $r$-matrix (which is the image of the
$r$-matrix of $D\g_{+}$) is of the form $r = \sum_{i} e_i \ot f^i $, where
$\{e_i\}$ is a basis of $\g_{+}$ and $\{f^i\}$ is the $(\cdot, \cdot)$-dual
basis of $\g_{-}$.

\subsubsection{\bf The standard Lie bialgebra structure of a complex 
semisimple Lie algebra.}
\label{__sect:CpxSimpleLA}

Let $\g$ be a complex simple Lie algebra. (All the constructions and results
that are discussed here and in other sections below for simple Lie algebras,
can be trivially extended to semisimple Lie algebras by decomposing them into 
direct sums of simple Lie algebras.) We recall the rudiments of the 
structure theory of $\g$. Everything we need can be found in
~\cite[Sections 8--10]{Humphreys:LieAlgs}, but some care is needed
to adjust for our choice of a different normalization of the metric.

\smallskip
We fix a Cartan subalgebra $\h$ of $\g$.
For $\alpha \in \h^{\star} - \{0\}$, we denote by $\g_{\alpha}$ the root 
subspace
\begin{equation}
\label{eqn:RootSpaceDef}
\g_{\alpha} 
= \{x\in \g \, \mid \, [h, x] = \alpha(h) x  \text{ for all } h\in \h \} \, .
\end{equation}
The root subspaces are either zero or one-dimensional, and 
$[\g_{\alpha}, \g_{\beta}] = \g_{\alpha + \beta}$. 
The set of roots of $(\g, \h)$ is 
\[
R = \{\alpha \in \h^{\star} - \{0\} \, \mid \, \g_{\alpha} \neq 0 \} \, .
\]
The Cartan decomposition of $\g$ is
\[
\g = \h \oplus \bigoplus_{\alpha\in R} \g_{\alpha} \,.
\]

The set of roots $R$ can be partitioned (in several ways) into a union,
$R = R_{+} \cup R_{-}$, with $R_{-} = - \, R_{+}$ and 
$\alpha, \beta \in R_{+}, 
\alpha + \beta \in R \Rightarrow \alpha + \beta \in R_{+} $. We choose one
such a partition. $R_{+}$ ($R_{-}$)
is called the set of positive (negative) roots. There is a unique subset 
$\Lambda = \{\alpha_1, \ldots , \alpha_{n}\}$ of $R_{+}$ ($n = \dim \h)$, 
called the set of simple roots, such that every element of $R_{+}$
can be uniquely written as a linear combination with nonnegative
integral coefficients of the roots in $\Lambda$. Associated with the choice 
of positive roots are the positive and negative maximal nilpotent subalgebras
$\n_{\pm} = \oplus_{\alpha\in R_{\pm}}  \g_{\alpha}$, and the positive and
negative Borel (maximal solvable) subalgebras $\b_{\pm} = \h \oplus \n_{\pm}$.

\smallskip
There is a metric $(\cdot, \cdot)_{\g}$ on $\g$, which is unique up to a 
nonzero scalar. 
We have  $(\g_{\alpha} , \g_{\beta})_{\g} = 0$ if $\alpha + \beta \neq 0$, and 
the restriction of $(\cdot , \cdot)_{\g}$ to $\h$, denoted 
$(\cdot , \cdot)_{\h}$, is a metric on $\h$. We use the metrics to identify 
$\g$ and $\h$ with their dual spaces, and use the same notation for the 
metrics induced
on $\g^{\star}$ and $\h^{\star}$ by these identifications. Now
the metric can be normalized so that $(\alpha, \alpha)_{\g} > 0 $
for $\alpha \in R \subset \h^{\star}$, and then these squared root lengths 
take at most two values on $R$. We completely fix the normalization of 
$(\cdot , \cdot)_{\g}$ by scaling so that the {\em smaller} value (called the 
squared length of the short roots) is equal to $2$. 

For $\alpha\in R_{+}$ we can choose $e_{\alpha} \in \g_{\alpha}$ and
$f_{\alpha} \in \g_{-\alpha}$ such that $(e_{\alpha} , f_{\alpha})_{\g} = 1$.
Put $h_{\alpha} = [e_{\alpha}, f_{\alpha}] \in \h$. One checks that 
$\alpha \mapsto h_{\alpha}$ under the identification  $\h \cong \h^{\star}$ 
induced by the metric.
It follows that the following commutation relations 
hold\footnote{Suffixes of the metrics are sometimes omitted.}:
\begin{equation}
\label{eqn:CartanRootBracket}
[h_{\beta}, e_{\alpha}] = (\beta, \alpha) e_{\alpha} \, , \quad
[h_{\beta}, f_{\alpha}] = -(\beta, \alpha) f_{\alpha} \, , \quad
\alpha, \beta \in R_{+} \, . 
\end{equation}

\smallskip
We now turn to the description of the Lie bialgebra structure of $\g$, 
following the treatment of  
~\cite[Section 4.4]{EtingofSchiffmann:LectQuantumGrps}. 

We consider the Lie algebra obtained from $\g$ by `doubling' the
Cartan subalgebra: 
$\tg = \n_{+} \oplus \h^{(1)} \oplus \h^{(2)} \oplus \n_{-}$, where
$\h^{(1)} = \h^{(2)} = \h$, with the commutation relations
\begin{align*}
& [\h^{(1)}, \h^{(2)}] = 0 \, , \\
& [h^{(i)}, e_{\alpha}] = \alpha(h) e_{\alpha} \, , \
[h^{(i)}, f_{\alpha}] = -\alpha(h) f_{\alpha} \, \text{ for } h\in \h \, , \,
i = 1, 2 \, , \\
& [e_{\alpha},  f_{\alpha}]  
= \frac{1}{2} (h_{\alpha}^{(1)} + h_{\alpha}^{(2)}) \, .
\end{align*}  
The projection $p: \tg \twoheadra \g$ defined by 
\[
p(e_{\alpha}) = e_{\alpha} \, , \ p(f_{\alpha}) = f_{\alpha} \, , \
p(h^{(i)}) = h \, ,
\]
is a Lie algebra map. We next show that $\tg$ is a 
quasitriangular Lie bialgebra.

\smallskip
We define a bilinear form on $\tg$ by the formula
\[
(x + h^{(1)} + k^{(2)}, x' + h^{'(1)} + k^{'(2)})_{\tg} =
2((h, k') + (k, h')) + (x, x') \, , 
\]
for $h, k, h', k' \in \h, \, x, x' \in \n_{+} \oplus \n_{-}$. It is 
immediate that $(\cdot, \cdot)_{\tg}$ is a metric on $\tg$.  

Set $\tg_{+} = \n_{+} \oplus \h^{(1)}, \tg_{-} = \n_{-} \oplus \h^{(2)}$.  
Since $\n_{\pm}$ are $(\cdot, \cdot)_{\tg}$-isotropic, it follows that
$(\tg, \tg_{+}, \tg_{-})$ is a Manin triple. By the discussion in 
Section ~\ref{__sect:manin}, $\tg$ is a quasitriangular Lie bialgebra. To 
write its $r$-matrix, choose an orthonormal basis $\{k_i\}$ of $\h$ with 
respect to $(\cdot, \cdot)_{\h}$. Then 
$\{e_{\alpha}\}_{\alpha \in R_{+}}, \{k_i^{(1)}\}$ is a basis of
$\tg_{+}$, and its $(\cdot, \cdot)_{\tg}$-dual basis of $\tg_{-}$ is
$\{f_{\alpha}\}_{\alpha \in R_{+}}, \{\frac{1}{2} k_i^{(2)}\}$.
Thus, the $r$-matrix of $\tg$ is
\begin{equation}
\label{eqn:r_tg}
r_{\tg} = \sum_{\alpha\in R_{+}} e_{\alpha} \ot f_{\alpha} 
+\frac{1}{2} \sum_{i} k_i^{(1)} \ot k_i^{(2)} \, .
\end{equation}
The projection $p: \tg \twoheadra \g$ maps $r_{\tg}$ to the 
$r$-matrix  
\begin{equation}
\label{eqn:r_g}
r_{\g} = \sum_{\alpha\in R_{+}} e_{\alpha} \ot f_{\alpha} 
+\frac{1}{2} \sum_{i} k_i \ot k_i \, ,
\end{equation}
which endows $\g$ with a quasitriangular Lie bialgebra structure (called
the {\em standard} structure). We have $r_{\g} + r_{\g}^{21} = \Omega_{\g}$,
the inverse metric tensor of $\g$. With this structure, $\g$ is `almost' the
double of $\b_{+}$ (or $\b_{-}$), which is a Lie subbialgebra of $\g$.

\begin{remark}
There is another description of $\tg$ which avoids the slightly artificial
doubling of $\h$. It is the Lie algebra $g \oplus \h$ together with the
metric $(\cdot, \cdot)_{\g} - (\cdot, \cdot)_{\h}$. A metric 
preserving Lie algebra isomorphism $\tg \ra \g \oplus \h$ is given by 
\[
e_{\alpha} \mapsto (e_{\alpha}, 0), \ f_{\alpha} \mapsto (f_{\alpha}, 0), \
h^{(1)} \mapsto (h, h), \ h^{(2)} \mapsto (h, -h) \, .
\]
Note however that in the $\g \oplus \h$ picture, the description of the 
isotropic subalgebras is somewhat less symmetric.
\end{remark} 

It can be shown that the cobracket $\delta$ of $\g$ has the following
values when $\beta\in\Lambda$ is a simple root:
\[
\delta(e_{\beta}) = \frac{1}{2} e_{\beta} \land h_{\beta} \, \quad ; \quad
\delta(f_{\beta}) =  \frac{1}{2} f_{\beta} \land h_{\beta} \, .
\]
Since the Lie algebra $\n_{+}$ ($\n_{-}$) is generated by 
$\{e_{\beta}\}$ ($\{f_{\beta}\}$), $\beta \in \Lambda$, the cobracket of a
general element of $\g$ can be computed by combined applications of the above
formulas and the cocycle identity.  

\begin{example}
For $\g = \mathfrak{sl}(2)$, we take as usual 
\[
\h = \Span_{\mathbb C} \left\{ 
h = \begin{pmatrix} 1 & 0 \\ 0 & -1 \end{pmatrix}
\right\} \, .
\]
The roots are $\pm\alpha \in \h^{\star}$, where $\alpha(h) =2$, and we choose
$R_{+} = \{\alpha\}$. The corresponding root spaces are
\[
\n_{+} =  \Span_{\mathbb C} \left\{ 
e = \begin{pmatrix} 0 & 1 \\ 0 & 0 \end{pmatrix}
\right\} \, , \quad  
\n_{-} =  \Span_{\mathbb C} \left\{ 
f = \begin{pmatrix} 0 & 0 \\ 1 & 0 \end{pmatrix}
\right\} \, .
\]

The normalized metric of $\g$ is 
\[
(x , y) = \tr (xy) \, , \quad x, y \in \g \, .
\]
Since $(e,f) = 1$ and $[e, f] = h$, the $r$-matrix of $\g$ is
\[
r_{\g} = e \ot f + \frac{1}{4} h \ot h \, .
\]
We have $\delta(h) = 0$, $\delta(e) = \frac{1}{2} e\land h$, and
$\delta(f) = \frac{1}{2} f\land h$.
Note that the Lie subbialgebras $\b_{\pm}$ are isomorphic to the Lie bialgebra 
of Example ~\ref{exmpl:2D_lba}. 
\end{example}

\subsection{Directed Jacobi diagrams and directed Lie relations} 

\subsubsection{\bf The basic definitions.}
We want to construct diagrammatic objects which are related to Lie bialgebras 
in a way resembling the relationship between $\A$-spaces and Lie algebras.  

\begin{definition}
A {\em directed Jacobi diagram} based on a skeleton $\M$ is a Jacobi diagram 
in which each arc of the body has a prescribed direction.
\end{definition}

We denote the space of  directed Jacobi diagrams based on 
a skeleton $\M$ by $\Dar(\M)$. We consider the subspace of 
$\Span_{\k}(\Dar(\M))$ which is spanned by the following relations:

\begin{enumerate}
\item The no sink no source (NS) relations:
\[ 
  \printname{ArrowNS1}
  \setlength{\unitlength}{0.35\standardunitlength}
  \begin{array}{c}  \hspace{-1.7mm}
    \raisebox{-8pt}{\input draws/ArrowNS1.tex }
      \hspace{-1.9mm}
  \end{array}
 = 0 \qquad ; \qquad
   
  \printname{ArrowNS2}
  \setlength{\unitlength}{0.35\standardunitlength}
  \begin{array}{c}  \hspace{-1.7mm}
    \raisebox{-8pt}{\input draws/ArrowNS2.tex }
      \hspace{-1.9mm}
  \end{array}
 = 0 \, . \]  
\item The antisymmetry (AS) of internal vertices:
\[ 
  \printname{ArrowAS1a}
  \setlength{\unitlength}{0.35\standardunitlength}
  \begin{array}{c}  \hspace{-1.7mm}
    \raisebox{-8pt}{\input draws/ArrowAS1a.tex }
      \hspace{-1.9mm}
  \end{array}
 + 
  \printname{ArrowAS1b}
  \setlength{\unitlength}{0.35\standardunitlength}
  \begin{array}{c}  \hspace{-1.7mm}
    \raisebox{-8pt}{\input draws/ArrowAS1b.tex }
      \hspace{-1.9mm}
  \end{array}
 = 0 \qquad ; \qquad
   
  \printname{ArrowAS2a}
  \setlength{\unitlength}{0.35\standardunitlength}
  \begin{array}{c}  \hspace{-1.7mm}
    \raisebox{-8pt}{\input draws/ArrowAS2a.tex }
      \hspace{-1.9mm}
  \end{array}
 + 
  \printname{ArrowAS2b}
  \setlength{\unitlength}{0.35\standardunitlength}
  \begin{array}{c}  \hspace{-1.7mm}
    \raisebox{-8pt}{\input draws/ArrowAS2b.tex }
      \hspace{-1.9mm}
  \end{array}
 = 0 \, . \]  
\item IHX relations:
\[ 
  \printname{ArrowI1}
  \setlength{\unitlength}{0.4\standardunitlength}
  \begin{array}{c}  \hspace{-1.7mm}
    \raisebox{-8pt}{\input draws/ArrowI1.tex }
      \hspace{-1.9mm}
  \end{array}
 - 
  \printname{ArrowH1}
  \setlength{\unitlength}{0.4\standardunitlength}
  \begin{array}{c}  \hspace{-1.7mm}
    \raisebox{-8pt}{\input draws/ArrowH1.tex }
      \hspace{-1.9mm}
  \end{array}
 
+ 
  \printname{ArrowX1}
  \setlength{\unitlength}{0.4\standardunitlength}
  \begin{array}{c}  \hspace{-1.7mm}
    \raisebox{-8pt}{\input draws/ArrowX1.tex }
      \hspace{-1.9mm}
  \end{array}
 = 0 \, , \]
\[ 
  \printname{ArrowI2}
  \setlength{\unitlength}{0.4\standardunitlength}
  \begin{array}{c}  \hspace{-1.7mm}
    \raisebox{-8pt}{\input draws/ArrowI2.tex }
      \hspace{-1.9mm}
  \end{array}
 - 
  \printname{ArrowH2}
  \setlength{\unitlength}{0.4\standardunitlength}
  \begin{array}{c}  \hspace{-1.7mm}
    \raisebox{-8pt}{\input draws/ArrowH2.tex }
      \hspace{-1.9mm}
  \end{array}
 
+ 
  \printname{ArrowX2}
  \setlength{\unitlength}{0.4\standardunitlength}
  \begin{array}{c}  \hspace{-1.7mm}
    \raisebox{-8pt}{\input draws/ArrowX2.tex }
      \hspace{-1.9mm}
  \end{array}
 = 0 \, , \]
\[ 
  \printname{ArrowI3}
  \setlength{\unitlength}{0.4\standardunitlength}
  \begin{array}{c}  \hspace{-1.7mm}
    \raisebox{-8pt}{\input draws/ArrowI3.tex }
      \hspace{-1.9mm}
  \end{array}
 = 
  \printname{ArrowH3a}
  \setlength{\unitlength}{0.4\standardunitlength}
  \begin{array}{c}  \hspace{-1.7mm}
    \raisebox{-8pt}{\input draws/ArrowH3a.tex }
      \hspace{-1.9mm}
  \end{array}
 - 
  \printname{ArrowX3a}
  \setlength{\unitlength}{0.4\standardunitlength}
  \begin{array}{c}  \hspace{-1.7mm}
    \raisebox{-8pt}{\input draws/ArrowX3a.tex }
      \hspace{-1.9mm}
  \end{array}
 
- 
  \printname{ArrowX3b}
  \setlength{\unitlength}{0.4\standardunitlength}
  \begin{array}{c}  \hspace{-1.7mm}
    \raisebox{-8pt}{\input draws/ArrowX3b.tex }
      \hspace{-1.9mm}
  \end{array}
 + 
  \printname{ArrowH3b}
  \setlength{\unitlength}{0.4\standardunitlength}
  \begin{array}{c}  \hspace{-1.7mm}
    \raisebox{-8pt}{\input draws/ArrowH3b.tex }
      \hspace{-1.9mm}
  \end{array}
 \, . \]
\item STU relations:
\[ 
  \printname{ArrowS1}
  \setlength{\unitlength}{0.35\standardunitlength}
  \begin{array}{c}  \hspace{-1.7mm}
    \raisebox{-8pt}{\input draws/ArrowS1.tex }
      \hspace{-1.9mm}
  \end{array}
 = 
  \printname{ArrowT1}
  \setlength{\unitlength}{0.35\standardunitlength}
  \begin{array}{c}  \hspace{-1.7mm}
    \raisebox{-8pt}{\input draws/ArrowT1.tex }
      \hspace{-1.9mm}
  \end{array}
 
                         - 
  \printname{ArrowU1}
  \setlength{\unitlength}{0.35\standardunitlength}
  \begin{array}{c}  \hspace{-1.7mm}
    \raisebox{-8pt}{\input draws/ArrowU1.tex }
      \hspace{-1.9mm}
  \end{array}
 \, , \]
\[ 
  \printname{ArrowS2}
  \setlength{\unitlength}{0.35\standardunitlength}
  \begin{array}{c}  \hspace{-1.7mm}
    \raisebox{-8pt}{\input draws/ArrowS2.tex }
      \hspace{-1.9mm}
  \end{array}
 = 
  \printname{ArrowT2}
  \setlength{\unitlength}{0.35\standardunitlength}
  \begin{array}{c}  \hspace{-1.7mm}
    \raisebox{-8pt}{\input draws/ArrowT2.tex }
      \hspace{-1.9mm}
  \end{array}
 
                         - 
  \printname{ArrowU2}
  \setlength{\unitlength}{0.35\standardunitlength}
  \begin{array}{c}  \hspace{-1.7mm}
    \raisebox{-8pt}{\input draws/ArrowU2.tex }
      \hspace{-1.9mm}
  \end{array}
 \, , \]
\[ 
  \printname{ArrowS3a}
  \setlength{\unitlength}{0.35\standardunitlength}
  \begin{array}{c}  \hspace{-1.7mm}
    \raisebox{-8pt}{\input draws/ArrowS3a.tex }
      \hspace{-1.9mm}
  \end{array}
 + 
  \printname{ArrowS3b}
  \setlength{\unitlength}{0.35\standardunitlength}
  \begin{array}{c}  \hspace{-1.7mm}
    \raisebox{-8pt}{\input draws/ArrowS3b.tex }
      \hspace{-1.9mm}
  \end{array}
 
= 
  \printname{ArrowT3}
  \setlength{\unitlength}{0.35\standardunitlength}
  \begin{array}{c}  \hspace{-1.7mm}
    \raisebox{-8pt}{\input draws/ArrowT3.tex }
      \hspace{-1.9mm}
  \end{array}
 - 
  \printname{ArrowU3}
  \setlength{\unitlength}{0.35\standardunitlength}
  \begin{array}{c}  \hspace{-1.7mm}
    \raisebox{-8pt}{\input draws/ArrowU3.tex }
      \hspace{-1.9mm}
  \end{array}
 \, . \]
\[ 
  \printname{ArrowS4a}
  \setlength{\unitlength}{0.35\standardunitlength}
  \begin{array}{c}  \hspace{-1.7mm}
    \raisebox{-8pt}{\input draws/ArrowS4a.tex }
      \hspace{-1.9mm}
  \end{array}
 + 
  \printname{ArrowS4b}
  \setlength{\unitlength}{0.35\standardunitlength}
  \begin{array}{c}  \hspace{-1.7mm}
    \raisebox{-8pt}{\input draws/ArrowS4b.tex }
      \hspace{-1.9mm}
  \end{array}
 
= 
  \printname{ArrowT4}
  \setlength{\unitlength}{0.35\standardunitlength}
  \begin{array}{c}  \hspace{-1.7mm}
    \raisebox{-8pt}{\input draws/ArrowT4.tex }
      \hspace{-1.9mm}
  \end{array}
 - 
  \printname{ArrowU4}
  \setlength{\unitlength}{0.35\standardunitlength}
  \begin{array}{c}  \hspace{-1.7mm}
    \raisebox{-8pt}{\input draws/ArrowU4.tex }
      \hspace{-1.9mm}
  \end{array}
 \, . \]
\end{enumerate}

We refer to the first and second relation in each group as the incoming and
outgoing relation, respectively, to the third IHX relation as the 
{\em mixed IHX relation} or {\em 5T relation}, and to the last two relations 
as the {\em mixed STU relations}.
   
\begin{definition}
\lbl{def:Aar}
$\Aar(\M)$ is the quotient space 
\[ 
\Aar(\M) = \Span_{\k}\left(\Dar(\M)\right) \left/
	\left\{\text{$NS$, $AS$, $IHX$, and $STU$ relations}\right\} \right. \, . 
\]
\end{definition}
We put $\Aar = \Aar(\ua)$ and $\Bar = \Aar(\ast)$.

\begin{remark}
Let $\Dar_{NS}(\M)$ be the subset of $\Dar(\M)$ consisting of diagrams having
no sinks or sources. The quotient of $\Span_{\k}\left(\Dar(\M)\right)$ by the 
NS relations alone is clearly isomorphic to the subspace 
$\Span_{\k}\left(\Dar_{NS}(\M)\right)$. 
Thus, we may regard $\Aar(\M)$ as the quotient of the latter subspace by the
AS, IHX and STU relations.
\end{remark}

As in the undirected case, $\Dar(\M)$ is graded by half the number of 
vertices in a diagram. As all the directed relations are homogeneous, 
this grading is inherited by $\Aar(\M)$. We also denote by $\Aar(\M)$ its
graded completion.

\subsubsection{\bf Mapping directed Jacobi diagrams into Lie tensors.}
\label{__sect:DrctdJacobiToLie}

Given a finite-dimensional Lie bialgebra $\a$ and its double $\tg = D\a$ 
(or equivalently, a Manin triple $(\tg, \a, \a^{\star})$), we
are going to construct, for each skeleton $\M$, a map
$\Tartg:\Aar(\M) \ra L(\M)$, where $L(\M)$ is the $\tg$-space
attached to $\M$. (We stress that although $\a$ is suppressed from
the notation, the map $\Tartg$ depends on the Manin triple
$(\tg, \a, \a^{\star})$ and not just on the Lie bialgebra $\tg$.) 
This construction is the directed counterpart of the 
construction of the maps $\T_{\g}:\A(\M) \ra L(M)$ in Section 
~\ref{__sect:TgDef} for a metrized Lie algebra $\g$. 

\smallskip
We begin by describing the atomic graphs and the tensors corresponding to
them. There are four types of trivalent atomic graphs (up to the
choice of orientation of the vertex). The corresponding 
tensors all belong to $\tg^{\ot 3}$.
 
The first two atomic graphs are the sink and source graphs; to these we 
assign the zero tensor. The third has two incoming arcs;  the attached tensor
is $s_{+} \in \a \ot \a^{\star} \ot \a^{\star}  \subset \tg^{\ot 3}$ 
representing the bracket $[\cdot,\cdot]: \a \ot \a \ra \a$. The fourth graph, 
with two outgoing arcs, is assigned the tensor 
$s_{-} \in \a \ot \a \ot \a^{\star} \subset \tg^{\ot 3}$, which represents
the cobracket of $\a$. The 
same graphs with the reversed orientations are assigned the negated tensors.

Finally, there is a bivalent atomic graph with the corresponding
tensor being the $r$-matrix $r \in \a \ot \a^{\star} \subset \tg^{\ot 2}$.

Let $\{e_i\}$ be a basis of $\a$ and $\{f^i\}$ be the dual basis of 
$\a^{\star}$. In terms of the basis elements, we have the following
formulas: 
\begin{align*}
r & = e_i \ot f^i \, , \\
s_+ & = [e_i, e_j] \ot f^i \ot f^j \, , \\
s_- & = f^i \ot \delta(e_i)  \, .
\end{align*}

Figure ~\ref{fig:AtomicDirectedGraphs} shows the directed atomic graphs with 
their tensors.
\begin{figure}[htpb]
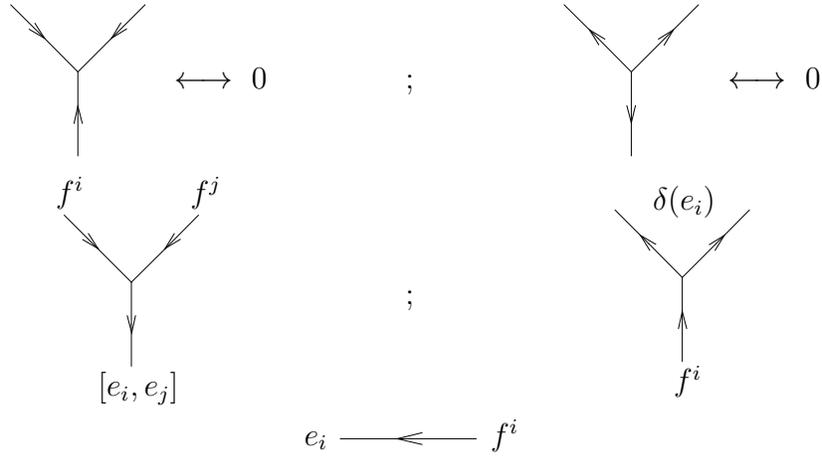

\[
\begin{array}{ccc}

  \printname{ArrowNS1}
  \setlength{\unitlength}{0.35\standardunitlength}
  \begin{array}{c}  \hspace{-1.7mm}
    \raisebox{-8pt}{\input draws/ArrowNS1.tex }
      \hspace{-1.9mm}
  \end{array}
 \longleftrightarrow \ 0 & ; & 

  \printname{ArrowNS2}
  \setlength{\unitlength}{0.35\standardunitlength}
  \begin{array}{c}  \hspace{-1.7mm}
    \raisebox{-8pt}{\input draws/ArrowNS2.tex }
      \hspace{-1.9mm}
  \end{array}
 \longleftrightarrow \ 0  \\
\\

  \printname{ArrAtomIn}
  \setlength{\unitlength}{0.35\standardunitlength}
  \begin{array}{c}  \hspace{-1.7mm}
    \raisebox{-8pt}{\input draws/ArrAtomIn.tex }
      \hspace{-1.9mm}
  \end{array}
  & ; & 

  \printname{ArrAtomOut}
  \setlength{\unitlength}{0.35\standardunitlength}
  \begin{array}{c}  \hspace{-1.7mm}
    \raisebox{-8pt}{\input draws/ArrAtomOut.tex }
      \hspace{-1.9mm}
  \end{array}
 \\
 & 
  \printname{ArrAtomBival}
  \setlength{\unitlength}{0.5\standardunitlength}
  \begin{array}{c}  \hspace{-1.7mm}
    \raisebox{-8pt}{\input draws/ArrAtomBival.tex }
      \hspace{-1.9mm}
  \end{array}
 &
\end{array}
\]
\caption{Attaching tensors to directed atomic graphs.} 
\lbl{fig:AtomicDirectedGraphs} 
\end{figure}  

As in the undirected case, we associate with a directed Jacobi graph $D$ the 
tensor power 
\[ 
 T(D)=\bigotimes_{\substack{\text{univalent vertices} \\
                              v \text{ of } D}} \tg \, .
\] 

\smallskip
We construct $\Tartg(D) \in T(D)$ following the same procedure as in
Section ~\ref{__sect:DtoTensor}: 
\begin{enumerate}
\item
$D$ is separated into atomic graphs. 
\item
For each atomic graph, the components of the corresponding tensor are
placed on its vertices, as illustrated in Figure 
~\ref{fig:AtomicDirectedGraphs}.
(That this process is unambiguous is even simpler than in the undirected case.)
\item
Tensor components from the two sides of each separation point are
contracted using the natural metric of $\tg$. 
\end{enumerate}

The proof that $\Tartg(D)$ is independent of the separation used is 
essentially identical to that of Proposition ~\ref{prop:TgD}. For example,
\[

  \printname{ArrAtomOutWithSepPt}
  \setlength{\unitlength}{0.4\standardunitlength}
  \begin{array}{c}  \hspace{-1.7mm}
    \raisebox{-8pt}{\input draws/ArrAtomOutWithSepPt.tex }
      \hspace{-1.9mm}
  \end{array}
 
 \mapsto \ (f^i, \delta(e_j)') \, f^j \ot e_i \ot \delta(e_j)''
 = f^j \ot \delta(e_j)' \ot \delta(e_j)'' = s_- \, ,
\]
where we used the notation (similar to the Sweedler notation) 
$\delta(x) = \delta(x)' \ot \delta(x)'' \, , x\in \a$, to indicate the 
decomposition of $\delta(x)$ into the sum of factorable tensors.

\smallskip
Let $D$ be a directed Jacobi diagram in $\Dar(\M)$.  Also denote by $D$ the
body of $D$ (which is a directed Jacobi graph). Recall the map 
$T(D) \ra T(\M)$ (defined at the end of Section ~\ref{__sect:DtoTensor}). 

\begin{proposition}
\lbl{prop:TarWellDef}
The series of maps
\[
 D \longmapsto \Tartg(D) \in T(D) \lra T(\M) \longmapsto L(\M)
\]
descends to a well-defined map $\Tartg: \Aar(\M) \ra L(\M)$.
\end{proposition}

\begin{proof}
We follow the lines of the proof of Proposition ~\ref{prop:TgWellDef}.
First, it is clear that the NS and AS relations are mapped to zero. 
The proofs for the incoming STU and IHX relations are
almost the same as in the undirected case, so we focus on the outgoing
and mixed relations.

\smallskip
We see from the figure below
\[ 

  \printname{ArrS2tnsr}
  \setlength{\unitlength}{0.35\standardunitlength}
  \begin{array}{c}  \hspace{-1.7mm}
    \raisebox{-8pt}{\input draws/ArrS2tnsr.tex }
      \hspace{-1.9mm}
  \end{array}
 = 
  \printname{ArrT2tnsr}
  \setlength{\unitlength}{0.35\standardunitlength}
  \begin{array}{c}  \hspace{-1.7mm}
    \raisebox{-8pt}{\input draws/ArrT2tnsr.tex }
      \hspace{-1.9mm}
  \end{array}
 
                        - 
  \printname{ArrU2tnsr}
  \setlength{\unitlength}{0.35\standardunitlength}
  \begin{array}{c}  \hspace{-1.7mm}
    \raisebox{-8pt}{\input draws/ArrU2tnsr.tex }
      \hspace{-1.9mm}
  \end{array}
 
 \] 
that the assertion for the outgoing STU relations
boils down to checking that 
\[
\delta(e_i) \ot f^i = e_i \ot e_j \ot [f^i, f^j] \, .
\]
This follows from the fact that the bracket (in $\tg$) of elements from 
$\a^{\star}$ is dual to the cobracket of $\a$; that is,
\[
([f^i, f^j], e_k) = (\delta(e_k), f^i \ot f^j) \, . 
\]
Hence
\[
 e_i \ot e_j \ot [f^i, f^j] =  (\delta(e_k), f^i \ot f^j)  e_i \ot e_j \ot f^k
= \delta(e_k) \ot f^k \, ,
\]
as required.

\smallskip
Similarly, the assertion for the mixed STU relations follows from 
\eqref{eqn:MixedBracket}, and that for the outgoing IHX relations from the
co-Jacobi identity.

\smallskip
We finally turn to the mixed IHX relations. Separating the I-graph
into two atomic graphs by placing one separation point on the inner arc, we 
find that it is mapped to the tensor
$f^i \ot f^j \ot \delta([e_i, e_j])$. This tensor corresponds to the map
$(x , y) \in \a^{\ot 2} \mapsto \delta([x, y])$, which is on the left hand
side of the cocycle identity. Doing the same for the other
four graphs, we find that each of them corresponds to 
one of the terms on the right hand side of the cocycle identity. Thus, the 
mixed IHX relations are mapped to zero as a consequence of the cocycle 
identity in $\a$.  
\end{proof}

\begin{example}
The diagram $r = 
  \printname{LeftArrow}
  \setlength{\unitlength}{0.25\standardunitlength}
  \begin{array}{c}  \hspace{-1.7mm}
    \raisebox{-8pt}{\input draws/LeftArrow.tex }
      \hspace{-1.9mm}
  \end{array}
$ is mapped by $\Tartg$
to the $r$-matrix of the Manin triple $(\tg, \a, \a^{\star})$, and
is accordingly called the {\em diagrammatic $r$-matrix}.
\end{example}

\begin{exercise}
\label{exer:DgrmCYBE}
Show that $r$ satisfies (a diagrammatic form of) the CYBE ~\eqref{eqn:CYBE}.
\end{exercise}

\smallskip
In the case of a complex semisimple Lie algebra $\g$ we introduce a variant of
the $\Tar$-maps. With the notation of Section ~\ref{__sect:CpxSimpleLA}, 
$\tg$ is the double of $\tg_{+}$, yielding a map  
$\Tartg : \Aar(\M) \ra L_{\tg}(\M)$.
The projection $p:\tg \twoheadra \g$ induces a map 
$p: L_{\tg}(\M) \ra L_{\g}(\M)$. We define the map
$\Targ: \Aar(\M) \ra L_{\g}(\M)$ as the composition 
$p \circ \Tartg$.

\begin{example}
\label{exmpl:DgrmHalfSum}
By the definition of $\Tartg$ and \eqref{eqn:r_tg},
\begin{align*}
\Tartg(\LHC) & = 
\sum_{\alpha\in R_{+}} e_{\alpha} f_{\alpha} 
+\frac{1}{2} \sum_{i} k_i^{(1)} k_i^{(2)} \, , \\ 
\Tartg(\RHC) & = 
\sum_{\alpha\in R_{+}} f_{\alpha} e_{\alpha} 
+\frac{1}{2} \sum_{i} k_i^{(2)} k_i^{(1)} \, . 
\end{align*}
Hence
\[
\Targ (\frac{_1}{^2}(\LHC - \RHC)) = 
\frac{1}{2} \sum_{\alpha\in R_{+}} [e_{\alpha}, f_{\alpha}] =
\frac{1}{2} \sum_{\alpha\in R_{+}} h_{\alpha}
\]
is the element of the Cartan algebra which is dual to the
half-sum of the positive roots; we denote this element by  $\varrho_{\g}$.
Thus, 
\[
  \varrho = \frac{_1}{^2}(\LHC - \RHC) \in \Aar
\]
is the diagrammatic analogue of $\varrho_{\g}$. By a mixed STU relation, 
we can also write $\varrho$ as
\[
  \varrho = \frac{_1}{^2}(\TPD + \TPU) \, . 
\]
\end{example}

\subsubsection{\bf Algebraic structures on $\protect\Aar$-spaces.}

$\Aar$-spaces can be endowed with the various algebraic structures
(products, coproducts, antipodes, \ldots), that were defined for 
$\A$-spaces; the definitions, proofs of maps being well-defined, and 
verifications of axioms are essentially identical. The maps $\Tartg$
interchange all these structures with the corresponding structures on the
Lie spaces $L_{\tg}(\M)$. We shall use the same notation as in the
undirected case for the various structure maps.

In particular, $\Aar$ is a diagrammatic Hopf algebra, and for any
Manin triple $(\tg, \a, a^{\star})$  the map $\Tartg:\Aar \ra U(\tg)$ is a 
morphism of Hopf algebras.

One basic difference from the undirected case is that $\Aar$, unlike $\A$,
is not commutative. For example,
\[
\mu \left(
  \printname{LeftArrow}
  \setlength{\unitlength}{0.25\standardunitlength}
  \begin{array}{c}  \hspace{-1.7mm}
    \raisebox{-8pt}{\input draws/LeftArrow.tex }
      \hspace{-1.9mm}
  \end{array}
\right) -
\mu\circ\tau\left(
  \printname{LeftArrow}
  \setlength{\unitlength}{0.25\standardunitlength}
  \begin{array}{c}  \hspace{-1.7mm}
    \raisebox{-8pt}{\input draws/LeftArrow.tex }
      \hspace{-1.9mm}
  \end{array}
\right) = 2 \varrho \neq 0 \, .
\]

We also endow the spaces $\lowupper{\Aar}{\bt n}$ with a $\Vect_{\k}$-algebra 
structure as in \eqref{eqn:AnVectProd}.

\subsection{Mapping $\protect\A$-spaces into $\protect\Aar$-spaces}

\subsubsection{\bf Definition of the $\iota$-maps.}

We define a map $\iota:\D(\M) \ra \Span_{\k}(\Dar(\M))$ by summing
over all the possible choices of directions for the arcs of a
diagram $D\in \D(\M)$.
\begin{proposition}
\lbl{prop:IotaWellDef}
$\iota$ descends to a well-defined map $\iota: \A(\M) \ra \Aar(\M)$.
\end{proposition}

\begin{proof}
As in the proof of Proposition ~\ref{prop:TgWellDef}, the local nature of the 
definition of $\iota$ makes it sufficient to check that the local
parts of the AS, STU, and IHX relations are mapped to zero modulo the
directed relations.

This is clear for the AS relations. The image of an STU relation is 
depicted below. (The two terms which are zero by NS relations are 
omitted from the image of the S-graph.)
\begin{align*}
+ 
  \printname{Sgraph}
  \setlength{\unitlength}{0.3\standardunitlength}
  \begin{array}{c}  \hspace{-1.7mm}
    \raisebox{-8pt}{\input draws/Sgraph.tex }
      \hspace{-1.9mm}
  \end{array}
   \mapsto\ & + 
  \printname{ArrowS1}
  \setlength{\unitlength}{0.3\standardunitlength}
  \begin{array}{c}  \hspace{-1.7mm}
    \raisebox{-8pt}{\input draws/ArrowS1.tex }
      \hspace{-1.9mm}
  \end{array}
  
                              + 
  \printname{ArrowS2}
  \setlength{\unitlength}{0.3\standardunitlength}
  \begin{array}{c}  \hspace{-1.7mm}
    \raisebox{-8pt}{\input draws/ArrowS2.tex }
      \hspace{-1.9mm}
  \end{array}

                              + \left(
  \printname{ArrowS3a}
  \setlength{\unitlength}{0.3\standardunitlength}
  \begin{array}{c}  \hspace{-1.7mm}
    \raisebox{-8pt}{\input draws/ArrowS3a.tex }
      \hspace{-1.9mm}
  \end{array}
 +  
                                
  \printname{ArrowS3b}
  \setlength{\unitlength}{0.3\standardunitlength}
  \begin{array}{c}  \hspace{-1.7mm}
    \raisebox{-8pt}{\input draws/ArrowS3b.tex }
      \hspace{-1.9mm}
  \end{array}
\right) \\
               & \hspace{7.55cm} + \left(
  \printname{ArrowS4a}
  \setlength{\unitlength}{0.3\standardunitlength}
  \begin{array}{c}  \hspace{-1.7mm}
    \raisebox{-8pt}{\input draws/ArrowS4a.tex }
      \hspace{-1.9mm}
  \end{array}
 + 
                                
  \printname{ArrowS4b}
  \setlength{\unitlength}{0.3\standardunitlength}
  \begin{array}{c}  \hspace{-1.7mm}
    \raisebox{-8pt}{\input draws/ArrowS4b.tex }
      \hspace{-1.9mm}
  \end{array}
\right) \\
- 
  \printname{Tgraph}
  \setlength{\unitlength}{0.3\standardunitlength}
  \begin{array}{c}  \hspace{-1.7mm}
    \raisebox{-8pt}{\input draws/Tgraph.tex }
      \hspace{-1.9mm}
  \end{array}
  \mapsto\ & - 
  \printname{ArrowT1}
  \setlength{\unitlength}{0.3\standardunitlength}
  \begin{array}{c}  \hspace{-1.7mm}
    \raisebox{-8pt}{\input draws/ArrowT1.tex }
      \hspace{-1.9mm}
  \end{array}
  
                                 - 
  \printname{ArrowT2}
  \setlength{\unitlength}{0.3\standardunitlength}
  \begin{array}{c}  \hspace{-1.7mm}
    \raisebox{-8pt}{\input draws/ArrowT2.tex }
      \hspace{-1.9mm}
  \end{array}
 
                                 - 
  \printname{ArrowT3}
  \setlength{\unitlength}{0.3\standardunitlength}
  \begin{array}{c}  \hspace{-1.7mm}
    \raisebox{-8pt}{\input draws/ArrowT3.tex }
      \hspace{-1.9mm}
  \end{array}
   
                                 - 
  \printname{ArrowT4}
  \setlength{\unitlength}{0.3\standardunitlength}
  \begin{array}{c}  \hspace{-1.7mm}
    \raisebox{-8pt}{\input draws/ArrowT4.tex }
      \hspace{-1.9mm}
  \end{array}
 \\
+ 
  \printname{Ugraph}
  \setlength{\unitlength}{0.3\standardunitlength}
  \begin{array}{c}  \hspace{-1.7mm}
    \raisebox{-8pt}{\input draws/Ugraph.tex }
      \hspace{-1.9mm}
  \end{array}
  \mapsto\ & +  
  \printname{ArrowU1}
  \setlength{\unitlength}{0.3\standardunitlength}
  \begin{array}{c}  \hspace{-1.7mm}
    \raisebox{-8pt}{\input draws/ArrowU1.tex }
      \hspace{-1.9mm}
  \end{array}
  
                                 +  
  \printname{ArrowU2}
  \setlength{\unitlength}{0.3\standardunitlength}
  \begin{array}{c}  \hspace{-1.7mm}
    \raisebox{-8pt}{\input draws/ArrowU2.tex }
      \hspace{-1.9mm}
  \end{array}
 
                                 +  
  \printname{ArrowU3}
  \setlength{\unitlength}{0.3\standardunitlength}
  \begin{array}{c}  \hspace{-1.7mm}
    \raisebox{-8pt}{\input draws/ArrowU3.tex }
      \hspace{-1.9mm}
  \end{array}
  
                                 +  
  \printname{ArrowU4}
  \setlength{\unitlength}{0.3\standardunitlength}
  \begin{array}{c}  \hspace{-1.7mm}
    \raisebox{-8pt}{\input draws/ArrowU4.tex }
      \hspace{-1.9mm}
  \end{array}
 
\end{align*}
The first three terms in the images of the undirected graphs sum to zero by 
an incoming STU relation, while the second terms do so by an 
outgoing STU relation. The sum of the third and fourth terms in the image of
the S-graph together with the third terms in the images of the T-graph and
U-graph is zero by a mixed STU relation, and so is the sum of the remaining 
four terms.

\smallskip
We turn to deal with the IHX relations.
Each of the graphs in an IHX relation is mapped to the sum of $32$ graphs. We
group the resulting $96$ signed graphs into $16$ groups according to the
directions of the four legs. We show that the signed sum of the six
terms in each group is equal to zero modulo the directed relations. 

There are two groups in which all the legs are incoming or outgoing. Each of
the graphs appearing in these groups has either a sink or a source, hence 
is equal to zero by an NS relation.

Consider one of the four groups having exactly one incoming leg. Three of the
graphs in this group have a source and are zero by NS relations; the 
signed sum of the remaining three graphs
(in which the direction of the middle arc 'follows' the direction of
the incoming leg) is zero by an outgoing IHX relation,
possibly in combination with AS relations. For example, below is the 
resulting signed sum when the incoming leg is arriving from the upper left 
corner:
\[

  \printname{OneIn1}
  \setlength{\unitlength}{0.35\standardunitlength}
  \begin{array}{c}  \hspace{-1.7mm}
    \raisebox{-8pt}{\input draws/OneIn1.tex }
      \hspace{-1.9mm}
  \end{array}
 - 
  \printname{OneIn2}
  \setlength{\unitlength}{0.35\standardunitlength}
  \begin{array}{c}  \hspace{-1.7mm}
    \raisebox{-8pt}{\input draws/OneIn2.tex }
      \hspace{-1.9mm}
  \end{array}
 + 
  \printname{OneIn3}
  \setlength{\unitlength}{0.35\standardunitlength}
  \begin{array}{c}  \hspace{-1.7mm}
    \raisebox{-8pt}{\input draws/OneIn3.tex }
      \hspace{-1.9mm}
  \end{array}
 \, .
\]
By making a half-turn clockwise and then applying AS at the two ends of the 
inner arc of the third graph, the sum becomes the outgoing IHX relation
in its standard form.
 
The four groups having exactly one outgoing leg are handled in a similar way.

We are left with six groups having two incoming and two outgoing legs. In
such a group, one term is zero by an NS relation, and the other five terms
make up (after using AS relations) a mixed IHX relation (see Exercise
~\ref{exer:2In2Out} below for an example). This completes the proof.  
\end{proof}

\begin{exercise}
\lbl{exer:2In2Out}
Show that the following signed sum is zero:
\[

  \printname{TwoInTwoOut1a}
  \setlength{\unitlength}{0.35\standardunitlength}
  \begin{array}{c}  \hspace{-1.7mm}
    \raisebox{-8pt}{\input draws/TwoInTwoOut1a.tex }
      \hspace{-1.9mm}
  \end{array}
 + 
  \printname{TwoInTwoOut1b}
  \setlength{\unitlength}{0.35\standardunitlength}
  \begin{array}{c}  \hspace{-1.7mm}
    \raisebox{-8pt}{\input draws/TwoInTwoOut1b.tex }
      \hspace{-1.9mm}
  \end{array}
 
- 
  \printname{TwoInTwoOut2a}
  \setlength{\unitlength}{0.35\standardunitlength}
  \begin{array}{c}  \hspace{-1.7mm}
    \raisebox{-8pt}{\input draws/TwoInTwoOut2a.tex }
      \hspace{-1.9mm}
  \end{array}
 - 
  \printname{TwoInTwoOut2b}
  \setlength{\unitlength}{0.35\standardunitlength}
  \begin{array}{c}  \hspace{-1.7mm}
    \raisebox{-8pt}{\input draws/TwoInTwoOut2b.tex }
      \hspace{-1.9mm}
  \end{array}
 
+ 
  \printname{TwoInTwoOut3a}
  \setlength{\unitlength}{0.35\standardunitlength}
  \begin{array}{c}  \hspace{-1.7mm}
    \raisebox{-8pt}{\input draws/TwoInTwoOut3a.tex }
      \hspace{-1.9mm}
  \end{array}
 + 
  \printname{TwoInTwoOut3b}
  \setlength{\unitlength}{0.35\standardunitlength}
  \begin{array}{c}  \hspace{-1.7mm}
    \raisebox{-8pt}{\input draws/TwoInTwoOut3b.tex }
      \hspace{-1.9mm}
  \end{array}
 \, .
\]
\end{exercise}

We note that the maps $\iota: \A(\M) \ra \Aar(\M)$ interchange 
all the diagrammatic operations defined in Section ~\ref{__sect:DgrmOper}. For 
example, for any skeleton $\M$ and one of its components $m$,
\[
\iota(\Delta_{\m}(D)) = \Delta_{\m}(\iota(D)) \, , \quad D\in \A(\M) \, .
\]
These properties hold already at the level of linear combinations of diagrams,
and are verified by simple inspections. In particular, the map
$\iota:\A \ra \Aar$ is a morphism of diagrammatic Hopf algebras.

\smallskip
The Lie algebraic counterpart of $\iota$ is described in the following
proposition.
\begin{proposition}
\lbl{prop:IotaIsLikeInject}
For any Manin triple $(\tg, \a, \a^{\star})$, the diagram
\begin{equation}
\label{cd:IotaIsLikeInject}
\begin{CD}
\A(\M) @>{\iota}>> \Aar(\M) \\
@V{\T_{\tg}}VV  @VV{\Tartg}V \\
L(\M)^{\tg} @>{i}>> L(\M) 
\end{CD}
\end{equation}
is commutative. (The metric of $\tg$ underlying $\T_{\tg}$ is its natural
metric as the double of $\a$.) 
\end{proposition}

\begin{proof}
It follows from the definitions of $\iota$, $\T_{\tg}$, and $\Tartg$ that
we only need to check that $\T_{\tg}(G) = \Tartg \circ \iota(G)$, when 
$G$ is one of the two undirected atomic graphs.

\smallskip
For the bivalent atomic graph, this amounts to the equality
\begin{equation}
\label{eqn:CasEqYBplusYBop}
t^{-1} = r + r^{21} \, .
\end{equation}
Using the basis $\{e_i\}$ of $\a$ and the dual basis $\{f^i\}$ of $\a^{\star}$
(in terms of which $r = e_i \ot f^i$), we note  that $\{e_i, f^i\}$ is a basis 
of $\tg$, and its
dual with respect to the natural metric of the double is simply
$\{f^i, e_i\}$. From \eqref{eqn:CasInBasis} we get
\[
t^{-1} = e_i \ot f^i + f^i \ot e_i \, ,
\]
proving \eqref{eqn:CasEqYBplusYBop}.

\smallskip
The image of the trivalent atomic graph under $\iota$ is the 
sum of six directed graphs, which are mapped by $\Tartg$ to the sum of the 
three cyclic permutations of $s_+$ and $s_-$:
\begin{align*}
s' 
= \Tartg \circ \iota\left(
  \printname{AtomicTrivalent}
  \setlength{\unitlength}{0.3\standardunitlength}
  \begin{array}{c}  \hspace{-1.7mm}
    \raisebox{-8pt}{\input draws/AtomicTrivalent.tex }
      \hspace{-1.9mm}
  \end{array}
\right) & =
[e_i, e_j] \ot f^i \ot f^j + f^j \ot [e_i, e_j] \ot f^i 
+ f^i \ot f^j \ot [e_i, e_j] \\
&   + \delta(e_i)'' \ot f^i \delta(e_i)' + \delta(e_i) \ot f^i
+ f^i \ot \delta(e_i)\, .
\end{align*}
$\T_{\tg}$ maps the trivalent atomic graph to the bracket tensor of the double,
which is
\[
s = [e_i, e_j] \ot f^i \ot f^j + [f^i, e_j] \ot e_i \ot f^j
  + [e_i, f^j] \ot f^i \ot e_j + [f^i, f^j] \ot e_i \ot e_j \, .
\]
We have to show that $s = s'$.
The first terms in $s$ and $s'$ are equal; so are their last terms
(see the proof of Proposition ~\ref{prop:TarWellDef}). We use
\eqref{eqn:MixedBracket} to compute the second term of $s$:
\begin{align*}
[f^i, e_j] \ot e_i \ot f^j & = 
((1 \ot f^i)\delta(e_j) - (\ad^{\star}e_j)f^i) \ot e_i \ot f^j \\
 & = (f^i, \delta(e_j)'')\, \delta(e_j)' + (f^i, [e_j, e_k])\, f^k) 
	\ot e_i \ot f^j \\
& = \delta(e_j)' \ot \delta(e_j)'' \ot f^j 
    	+ f^k \ot [e_j, e_k] \ot f^j \, , 
\end{align*}
which cancels out the second and fifth terms of $s'$. Similarly,
the third term of $s$ is equal to the sum of the remaining two terms
of $s'$, thus proving our claim.
\end{proof}

We elaborate on \eqref{cd:IotaIsLikeInject} in the example of a complex 
semisimple Lie algebra $\g$ (Section ~\ref{__sect:CpxSimpleLA}), where  
$\a = \tg_{+}$ and its double is $\tg$ ($\g$ with `doubled' Cartan 
subalgebra). We apply to the two $L$-spaces in \eqref{cd:IotaIsLikeInject} 
the maps induced by the projection $p:\tg \ra \g$. Recall that we defined 
$\Targ$ as the composition
\[
\begin{CD}
\Aar(\M) @>{\Tartg}>> L_{\tg}(\M) @>{p}>> L_{\g}(\M)  \, .
\end{CD}
\]
On the other hand, it is clear from the definition of $\tg$ that 
$L_{\tg}(\M)^{\tg}$ is mapped into $L_{\g}(\M)^{\g}$. We claim that the
composition
\[
\begin{CD}
\A(\M) @>{\T_{\tg}}>> L_{\tg}(\M)^{\tg} @>{p}>> L_{\g}(\M)^{\g}  
\end{CD}
\]
is equal to the map $\T_{\g}$ (defined with the metric of $\g$ being the one 
used in the construction of the Manin triple structure on $\tg$). As in the 
proof of the last proposition, it is enough to verify the equality for the
two atomic graphs. This is immediate for the bivalent atomic graph,
because by \eqref{eqn:r_tg} and \eqref{eqn:r_g}
\[
\Omega_{\g} = r_{\g} + r_{\g}^{21} = p(r_{\tg} + r_{\tg}^{21}
= p(\Omega_{\tg}) \, .
\]
The assertion for the trivalent atomic graph follows from a 
similar straightforward computation. We have thus proved the commutativity
of the following diagram, for any complex semisimple Lie algebra $\g$:
\begin{equation}
\label{cd:IotaVsCpxSimpleLA}
\begin{CD}
\A(\M) @>{\iota}>> \Aar(\M) \\
@V{\T_{\g}}VV  @VV{\Targ}V \\
L(\M)^{\g} @>{i}>> L(\M) \, .
\end{CD}
\end{equation}

\subsubsection{\bf The image and kernel of $\iota$.}
\label{__sect:ImgKerIota}

We continue to assume that $(\tg, \a, \a^{\star})$ is a Manin triple.
By the very definition, an element $u\in U(\tg)^{\tg}$ is (mapped
by $i:U(\tg)^{\tg} = L(\ua)^{\tg} \ra L(\ua) = U(\tg)$ to)  
a {\em $\tg$-invariant} element in $U(\tg)$, namely $xu-ux=0$ for all 
$x\in \tg$.  
It is easy to define the diagrammatic notion of invariance.
\begin{definition}
\lbl{def:InvarDgrm}
An element $D\in \Aar$ is called {\em invariant} if 
\begin{equation}
\label{eqn:DgrmInv}

  \printname{DirLegL}
  \setlength{\unitlength}{0.4\standardunitlength}
  \begin{array}{c}  \hspace{-1.7mm}
    \raisebox{-8pt}{\input draws/DirLegL.tex }
      \hspace{-1.9mm}
  \end{array}
 - 
  \printname{DirLegR}
  \setlength{\unitlength}{0.4\standardunitlength}
  \begin{array}{c}  \hspace{-1.7mm}
    \raisebox{-8pt}{\input draws/DirLegR.tex }
      \hspace{-1.9mm}
  \end{array}
 = 0 \,
\text{ in } \Aar(\ua, \ast) \, ,  
\end{equation}
and similarly with an outgoing leg (directed from the interval to the 
$\ast$-colored vertex).
The subspace of invariant elements in $\Aar(\ua)$ is denoted by 
$\lowupper{\Aar}{\inv}$.
\end{definition}

We can make a similar definition for $\A$, using the undirected version of 
\eqref{eqn:DgrmInv}; of course, any element of $\A$
is invariant, as follows from the argument in the proof of Proposition
~\ref{prop:Aline-eq-Acirc}.

\smallskip
An easy consequence of the invariance of $D\in \lowupper{\Aar}{\inv}$
is that   
\begin{equation}
\label{eqn:InvIsCntrl}
  \printname{DtimesE}
  \setlength{\unitlength}{0.4\standardunitlength}
  \begin{array}{c}  \hspace{-1.7mm}
    \raisebox{-8pt}{\input draws/DtimesE.tex }
      \hspace{-1.9mm}
  \end{array}
 =
  \printname{EtimesD}
  \setlength{\unitlength}{0.4\standardunitlength}
  \begin{array}{c}  \hspace{-1.7mm}
    \raisebox{-8pt}{\input draws/EtimesD.tex }
      \hspace{-1.9mm}
  \end{array}

\quad \text{ for all } E\in \Aar \, .
\end{equation}
This is proved, as in the proof of commutativity in 
Proposition ~\ref{prop:AisHopf}, by moving leg after leg of $E$ across $D$.

\smallskip
In view of Proposition \eqref{cd:IotaIsLikeInject}, the following result is no
surprise.
\begin{proposition}
\lbl{prop:ImgIotaInvar}
The image of $\iota : \A \ra \Aar$ lies in $\lowupper{\Aar}{\inv}$.
\end{proposition}
\begin{proof}
Let $D$ be a diagram in $\A$. We have to prove that \eqref{eqn:DgrmInv} 
holds for $\iota(D)$.
We consider the case of an incoming leg; the case of an outgoing leg is
completely analogous.

\smallskip
Led $C$ be a diagram in $\Aar$ (here we have in mind for $C$ one of the
summands making up $\iota(D)$). We examine first the effect of moving the
incoming leg across $C$. 
\begin{lemma}
\lbl{lemma:DirThreeHook}
For a diagram $C\in\Aar$,
\begin{equation}
\label{eqn:LegAcrossDirDgrm}

  \printname{DirLegLeftToC}
  \setlength{\unitlength}{0.4\standardunitlength}
  \begin{array}{c}  \hspace{-1.7mm}
    \raisebox{-8pt}{\input draws/DirLegLeftToC.tex }
      \hspace{-1.9mm}
  \end{array}
 - 
  \printname{DirLegRightToC}
  \setlength{\unitlength}{0.4\standardunitlength}
  \begin{array}{c}  \hspace{-1.7mm}
    \raisebox{-8pt}{\input draws/DirLegRightToC.tex }
      \hspace{-1.9mm}
  \end{array}
 
= \sum_{\bf a} (-1)^{\bf a} C_{\bf a} \, ,
\end{equation}
where:
\begin{itemize}
\item
The summation is over the internal arcs of $C$.
\item
The diagram $C_{\bf a}$ is obtained by connecting the incoming
leg to a new vertex in the middle of ${\bf a}$, directing the two newly
formed arcs as outgoing from the new vertex.
\item
The sign $(-1)^{\bf a}$ is $\pm 1$ according to the positioning of the 
internal arc ${\bf a}$ with respect to the incoming leg:
\[

  \printname{IntrnlArcPlus}
  \setlength{\unitlength}{0.4\standardunitlength}
  \begin{array}{c}  \hspace{-1.7mm}
    \raisebox{-8pt}{\input draws/IntrnlArcPlus.tex }
      \hspace{-1.9mm}
  \end{array}
 \ \longmapsto \ +1 \, , \qquad

  \printname{IntrnlArcMinus}
  \setlength{\unitlength}{0.4\standardunitlength}
  \begin{array}{c}  \hspace{-1.7mm}
    \raisebox{-8pt}{\input draws/IntrnlArcMinus.tex }
      \hspace{-1.9mm}
  \end{array}
 \ \longmapsto \ -1 \, .
\]
(Formally, the sign is the orientation type of the planar basis formed by the 
direction of the arc and the direction of the incoming leg.)
\end{itemize}
\end{lemma}
\begin{proof}
As in the proof of Proposition
~\ref{prop:Aline-eq-Acirc}, we attach three right-turning hooks near every
vertex of $C$, with minus signs under the hooks that are connected to a
skeleton arc leaving a vertex. An example of the resulting decorated
diagram is shown in the following figure:
\begin{figure}[htpb]
\[

  \printname{HookedDirDgrm}
  \setlength{\unitlength}{0.4\standardunitlength}
  \begin{array}{c}  \hspace{-1.7mm}
    \raisebox{-8pt}{\input draws/HookedDirDgrm.tex }
      \hspace{-1.9mm}
  \end{array}

\]
\end{figure}

We form the signed sum of the diagrams obtained by all possible ways of
connecting the incoming leg to one of the hooks and, if the incoming leg is
connected to a non-skeleton hook, specifying one out of the two possible 
directions of the small arc between the hook and its vertex. We group the 
terms of this signed sum in two ways.
\begin{itemize}
\item
By vertices:

Consider first the left-most skeleton vertex. There are four summands
in its group. They are identical outside a small neighborhood of the
vertex, and their sum locally looks like
\[

  \setlength{\unitlength}{0.4\standardunitlength}
  \begin{array}{c}  \hspace{-1.7mm}
    \raisebox{-2pt}{\input draws/SkltVrtxHook1.tex }
    \hspace{-1.9mm}
  \end{array}
 + 
  \setlength{\unitlength}{0.4\standardunitlength}
  \begin{array}{c}  \hspace{-1.7mm}
    \raisebox{-2pt}{\input draws/SkltVrtxHook2.tex }
    \hspace{-1.9mm}
  \end{array}
 +

  \setlength{\unitlength}{0.4\standardunitlength}
  \begin{array}{c}  \hspace{-1.7mm}
    \raisebox{-2pt}{\input draws/SkltVrtxHook3.tex }
    \hspace{-1.9mm}
  \end{array}
 - 
  \setlength{\unitlength}{0.4\standardunitlength}
  \begin{array}{c}  \hspace{-1.7mm}
    \raisebox{-2pt}{\input draws/SkltVrtxHook4.tex }
    \hspace{-1.9mm}
  \end{array}
 \, ,
\]
hence is zero by a mixed STU relation. The same argument applies for all
other skeleton vertices, possibly using NS and incoming STU relations
instead of the mixed STU relation.

For the internal vertex (in our example), the sum of the six terms in
its group is locally
\[

  \setlength{\unitlength}{0.35\standardunitlength}
  \begin{array}{c}  \hspace{-1.7mm}
    \raisebox{-2pt}{\input draws/IntrnlVrtxHook1.tex }
    \hspace{-1.9mm}
  \end{array}
 + 
  \setlength{\unitlength}{0.35\standardunitlength}
  \begin{array}{c}  \hspace{-1.7mm}
    \raisebox{-2pt}{\input draws/IntrnlVrtxHook2.tex }
    \hspace{-1.9mm}
  \end{array}
 +

  \setlength{\unitlength}{0.35\standardunitlength}
  \begin{array}{c}  \hspace{-1.7mm}
    \raisebox{-2pt}{\input draws/IntrnlVrtxHook3.tex }
    \hspace{-1.9mm}
  \end{array}
 + 
  \setlength{\unitlength}{0.35\standardunitlength}
  \begin{array}{c}  \hspace{-1.7mm}
    \raisebox{-2pt}{\input draws/IntrnlVrtxHook4.tex }
    \hspace{-1.9mm}
  \end{array}
 +

  \setlength{\unitlength}{0.35\standardunitlength}
  \begin{array}{c}  \hspace{-1.7mm}
    \raisebox{-2pt}{\input draws/IntrnlVrtxHook5.tex }
    \hspace{-1.9mm}
  \end{array}
 + 
  \setlength{\unitlength}{0.35\standardunitlength}
  \begin{array}{c}  \hspace{-1.7mm}
    \raisebox{-2pt}{\input draws/IntrnlVrtxHook6.tex }
    \hspace{-1.9mm}
  \end{array}
 \, ,
\]
which is zero by NS (applied to the sixth term), AS, and mixed IHX.
Similarly, if an internal vertex has two arcs emanating from it,
then three terms are zero by NS, and the remaining three sum
to zero by AS and incoming IHX.

\item
By arcs:

The two terms obtained by connecting the incoming leg to the left and right
`semi-infinite' skeleton arcs yield the difference created by moving
the incoming leg from left to right.

The group of an internal skeleton arc contains two terms which cancel
each other.

For an internal arc, we get (using AS and NS relations) either
\begin{equation}
\label{eqn:IncLegIntrnlArcUp}

  \setlength{\unitlength}{0.4\standardunitlength}
  \begin{array}{c}  \hspace{-1.7mm}
    \raisebox{-2pt}{\input draws/IntrnlArc1.tex }
    \hspace{-1.9mm}
  \end{array}
 + 
  \setlength{\unitlength}{0.4\standardunitlength}
  \begin{array}{c}  \hspace{-1.7mm}
    \raisebox{-2pt}{\input draws/IntrnlArc2.tex }
    \hspace{-1.9mm}
  \end{array}
 +

  \setlength{\unitlength}{0.4\standardunitlength}
  \begin{array}{c}  \hspace{-1.7mm}
    \raisebox{-2pt}{\input draws/IntrnlArc3.tex }
    \hspace{-1.9mm}
  \end{array}
 + 
  \setlength{\unitlength}{0.4\standardunitlength}
  \begin{array}{c}  \hspace{-1.7mm}
    \raisebox{-2pt}{\input draws/IntrnlArc4.tex }
    \hspace{-1.9mm}
  \end{array}
 =
- 
  \setlength{\unitlength}{0.4\standardunitlength}
  \begin{array}{c}  \hspace{-1.7mm}
    \raisebox{-2pt}{\input draws/IntrnlArc8.tex }
    \hspace{-1.9mm}
  \end{array}

\end{equation}
or
\begin{equation}
\label{eqn:IncLegIntrnlArcDown}

  \setlength{\unitlength}{0.4\standardunitlength}
  \begin{array}{c}  \hspace{-1.7mm}
    \raisebox{-2pt}{\input draws/IntrnlArc5.tex }
    \hspace{-1.9mm}
  \end{array}
 + 
  \setlength{\unitlength}{0.4\standardunitlength}
  \begin{array}{c}  \hspace{-1.7mm}
    \raisebox{-2pt}{\input draws/IntrnlArc6.tex }
    \hspace{-1.9mm}
  \end{array}
 +

  \setlength{\unitlength}{0.4\standardunitlength}
  \begin{array}{c}  \hspace{-1.7mm}
    \raisebox{-2pt}{\input draws/IntrnlArc7.tex }
    \hspace{-1.9mm}
  \end{array}
 + 
  \setlength{\unitlength}{0.4\standardunitlength}
  \begin{array}{c}  \hspace{-1.7mm}
    \raisebox{-2pt}{\input draws/IntrnlArc8.tex }
    \hspace{-1.9mm}
  \end{array}
 =

  \setlength{\unitlength}{0.4\standardunitlength}
  \begin{array}{c}  \hspace{-1.7mm}
    \raisebox{-2pt}{\input draws/IntrnlArc8.tex }
    \hspace{-1.9mm}
  \end{array}

\end{equation}

\end{itemize}

Now \eqref{eqn:LegAcrossDirDgrm} follows by comparing the two ways of 
grouping the terms. Note that the signs are dictated by
\eqref{eqn:IncLegIntrnlArcUp} and \eqref{eqn:IncLegIntrnlArcDown}. 
\end{proof}

We are now ready to prove the proposition. Apply the lemma with $C$ 
ranging over the summands of $\iota(D)$, and observe that the sum of the right
hand sides of \eqref{eqn:LegAcrossDirDgrm} is zero: 
the term $(-1)^{\bf a} C_{\bf a}$ is canceled against 
$(-1)^{\bf b} \tilde{C}_{\bf b}$, where 
$\tilde{C}$ is the summand of $\iota(D)$ which differs from $C$ only by the 
replacement of the internal arc ${\bf a}$ by an arc ${\bf b}$ with
the opposite direction, so that 
$\tilde{C}_{\bf b} = C_{\bf a}$ and
$(-1)^{\bf b} = - (-1)^{\bf a}$.   
\end{proof}

In view of the last proposition, it is natural to make the following
question.
\begin{question}
Are all the invariant elements in $\Aar$ coming from $\A$ ?
\end{question}

\smallskip
Combining the arguments of the proof of Proposition ~\ref{prop:ImgIotaInvar}  
with the method of proof of Proposition ~\ref{prop:AnInv}, we have the
following result.
\begin{proposition}
\lbl{prop:AarPwrInvar}
For all $D \in \Aar$ and $E \in \A^{\bt n}$,
\[
	[\Delta^{(n)}(D), \iota(E)] = 0 \, .
\]
\end{proposition}

\smallskip
{\bf The kernel of $\iota$:}
Another obvious fact about the map $i:U(\tg)^{\tg}\ra U(\tg)$ (and more 
generally, $i:L(\M)^{\tg} \ra L(\M)$ for any skeleton $\M$) is simply that it
is injective. Thus, the commutative diagram \eqref{cd:IotaIsLikeInject} 
motivates
\begin{conjecture}
\lbl{conj:IotaIsInjective}
For any skeleton $\M$, the map $\iota:\A(\M) \ra \Aar(\M)$ is injective. 
\end{conjecture}

For $M= \ua$ and $M=\bigcirc$, we can use \eqref{cd:IotaVsCpxSimpleLA} to 
prove Conjecture ~\ref{conj:IotaIsInjective} for low degrees, as
it is known from computer computations (see ~\cite{Kneissler:NumPrimVassiliev})
that the maps $\T_{\g}$, with $\g$ ranging over complex semisimple Lie 
algebras, separate the points of $\A(\ua) \ (= \A(\bigcirc))$ up to degree 
12.  This type of argument can not be extended to all degrees, because 
$\A(\ua)$ is known, by a result of P. Vogel (see ~\cite{Vogel:AlgStructs}), 
to contain non-zero elements that are  mapped to zero by all the maps 
$\T_{\g}$, where $\g$ is even allowed to be any {\em metrized} Lie 
(super-) algebra \footnote{In ~\cite{Vogel:AlgStructs} this result is 
explicitly proved for semisimple Lie (super-) algebras, but according 
to J.~Lieberum who quotes Vogel (private communication), it extends to all 
metrized Lie (super-) algebras.}. 

\smallskip
In any case, we will use all symbols and pictures of elements 
in $\A$-spaces to refer also to their images under $\iota$ in the 
corresponding $\Aar$-spaces.

\subsubsection{\bf The directed PBW isomorphism}
The PBW map $\upchi: \Bar \ra \Aar$ is defined exactly as in \eqref{eqn:DefChi}.
Using the $\iota$-maps, we can now prove without difficulty that this map is
an isomorphism (just like its algebraic counterpart 
$\upchi: S(\tg) \ra U(\tg)$).
\begin{proposition} 
\lbl{prop:PBW-Aar}
The map $\upchi: \Bar \ra \Aar$ is an isomorphism.
\end{proposition}
\begin{proof}
The proof is a simple adaptation of the proof of Proposition ~\ref{prop:PBW-A}.
We need to:
\begin{itemize}
\item
Replace all spaces ($\D_k(\M)$, $\I_k$, $\A_k$, etc.) with the 
corresponding directed spaces.
\item
Replace the graphical representations of $u_i$, $e$, and $s_i$ 
(\eqref{eqn:uDgrm}, \eqref{eqn:eDgrm}, and \eqref{eqn:sDgrm}), with their 
images under $\iota$ in $\Span_{\k}\{\Dar(\ua, \ua)\}$.
Also replace in this way the graphical representation of any $\pi\in S_k$. 
\item
For a directed diagram $C\in \Dar(\ua, \ua)$, define the `glueing $C$'
operator on $\Span_{\k}\{\Dar(\ua)\}$ by  
\begin{equation}
\label{eqn:GlueCtoD}
C(D) = \left\{
\begin{array}{cl}

  \printname{GlueCtoD}
  \setlength{\unitlength}{0.3\standardunitlength}
  \begin{array}{c}  \hspace{-1.7mm}
    \raisebox{-8pt}{\input draws/GlueCtoD.tex }
      \hspace{-1.9mm}
  \end{array}
 & \text{ $C$ and $D$ are compatible } \\
\\
0 & \text{ otherwise }
\end{array} 
\right. \, ,
\end{equation} 
where the compatibility criteria are clear from the figure: the number of legs
of $C$ on the upper skeleton is equal to the number of legs of $D$, and each 
incoming (outgoing) leg of $D$ is glued to an outgoing (incoming) leg of $C$.
\item
Extend the definition of glueing operators to elements of 
$\Span_{\k}\{\Dar(\ua, \ua)\}$ by linearity. We denote the operator
associated to $\iota(u_i)$ by $u_i$, and similarly for any $\pi \in S_k$ and 
$s_i$. For example, 
\[
s_2\left(
  \printname{MixedDgrm}
  \setlength{\unitlength}{0.3\standardunitlength}
  \begin{array}{c}  \hspace{-1.7mm}
    \raisebox{-8pt}{\input draws/MixedDgrm.tex }
      \hspace{-1.9mm}
  \end{array}
\right) =
	
  \printname{MixedDgrmGlueS2a}
  \setlength{\unitlength}{0.3\standardunitlength}
  \begin{array}{c}  \hspace{-1.7mm}
    \raisebox{-8pt}{\input draws/MixedDgrmGlueS2a.tex }
      \hspace{-1.9mm}
  \end{array}
 + 
  \printname{MixedDgrmGlueS2b}
  \setlength{\unitlength}{0.3\standardunitlength}
  \begin{array}{c}  \hspace{-1.7mm}
    \raisebox{-8pt}{\input draws/MixedDgrmGlueS2b.tex }
      \hspace{-1.9mm}
  \end{array}
 \, . 
\]
\item
Notice that \eqref{eqn:OperSTU} and the vanishing of $\Gamma_D$
on the generators of $\mathcal{K}$, as in \eqref{eqn:VanishGammaD},
are all consequences of relations only between the elements of 
$\Span_{\k}\{\Dar(\ua, \ua)\}$ representing the permutations $\pi \in S_k$ 
($u_i$ and $e$ in particular) 
and the elements representing $s_i$. The fact that 
$\iota$ is well-defined (Proposition ~\ref{prop:IotaWellDef}) guarantees that 
the needed relations continue to hold in the directed spaces. A slightly 
closer look reveals that the filtrations used are also preserved.  
\end{itemize}
The rest of the proof consists now of formal computations identical to 
the ones carried out in the proof of the undirected case.
\end{proof}

%% file: draws/ArrowNS1.tex
%
\begingroup\makeatletter\ifx\SetFigFont\undefined%
\gdef\SetFigFont#1#2#3#4#5{%
  \reset@font\fontsize{#1}{#2pt}%
  \fontfamily{#3}\fontseries{#4}\fontshape{#5}%
  \selectfont}%
\fi\endgroup%
{\renewcommand{\dashlinestretch}{30}
\begin{picture}(2424,2739)(0,-10)
\path(340,2276)(612,2113)(445,2382)
\path(1977,2385)(1814,2113)(2083,2280)
\path(1290,593)(1213,901)(1141,593)
\path(1812,2112)(1212,1512)
\path(612,2112)(1212,1512)
\path(12,2712)(612,2112)
\path(2412,2712)(1812,2112)
\path(1212,762)(1212,1512)
\path(1212,762)(1212,12)
\end{picture}
}

%% file: draws/ArrowNS2.tex
%
\begingroup\makeatletter\ifx\SetFigFont\undefined%
\gdef\SetFigFont#1#2#3#4#5{%
  \reset@font\fontsize{#1}{#2pt}%
  \fontfamily{#3}\fontseries{#4}\fontshape{#5}%
  \selectfont}%
\fi\endgroup%
{\renewcommand{\dashlinestretch}{30}
\begin{picture}(2424,2739)(0,-10)
\path(740,2099)(468,2262)(635,1993)
\path(1800,1990)(1963,2262)(1694,2095)
\path(1136,911)(1213,603)(1285,911)
\path(1812,2112)(1212,1512)
\path(612,2112)(1212,1512)
\path(12,2712)(612,2112)
\path(2412,2712)(1812,2112)
\path(1212,762)(1212,1512)
\path(1212,762)(1212,12)
\end{picture}
}

%% file: draws/ArrowAS1a.tex
%
\begingroup\makeatletter\ifx\SetFigFont\undefined%
\gdef\SetFigFont#1#2#3#4#5{%
  \reset@font\fontsize{#1}{#2pt}%
  \fontfamily{#3}\fontseries{#4}\fontshape{#5}%
  \selectfont}%
\fi\endgroup%
{\renewcommand{\dashlinestretch}{30}
\begin{picture}(2424,2739)(0,-10)
\path(1136,916)(1213,608)(1285,916)
\path(340,2276)(612,2113)(445,2382)
\path(1977,2385)(1814,2113)(2083,2280)
\path(1812,2112)(1212,1512)
\path(612,2112)(1212,1512)
\path(12,2712)(612,2112)
\path(2412,2712)(1812,2112)
\path(1212,762)(1212,1512)
\path(1212,762)(1212,12)
\end{picture}
}

%% file: draws/ArrowAS2a.tex
%
\begingroup\makeatletter\ifx\SetFigFont\undefined%
\gdef\SetFigFont#1#2#3#4#5{%
  \reset@font\fontsize{#1}{#2pt}%
  \fontfamily{#3}\fontseries{#4}\fontshape{#5}%
  \selectfont}%
\fi\endgroup%
{\renewcommand{\dashlinestretch}{30}
\begin{picture}(2424,2739)(0,-10)
\path(740,2099)(468,2262)(635,1993)
\path(1800,1990)(1963,2262)(1694,2095)
\path(1290,598)(1213,906)(1141,598)
\path(1812,2112)(1212,1512)
\path(612,2112)(1212,1512)
\path(12,2712)(612,2112)
\path(2412,2712)(1812,2112)
\path(1212,762)(1212,1512)
\path(1212,762)(1212,12)
\end{picture}
}

%% file: draws/ArrowX1.tex
%
\begingroup\makeatletter\ifx\SetFigFont\undefined%
\gdef\SetFigFont#1#2#3#4#5{%
  \reset@font\fontsize{#1}{#2pt}%
  \fontfamily{#3}\fontseries{#4}\fontshape{#5}%
  \selectfont}%
\fi\endgroup%
{\renewcommand{\dashlinestretch}{30}
\begin{picture}(2126,2117)(0,-10)
\path(695,680)(1444,677)
\path(912,900)(2113,2090)
\path(1801,1650)(1464,1450)(1669,1776)
\path(1222,900)(12,2090)
\path(324,1650)(661,1450)(456,1776)
\path(2114,19)(1214,904)
\path(1661,594)(1864,264)(1531,464)
\path(916,900)(13,12)
\path(179,307)(516,507)(310,181)
\path(904,603)(1213,680)(904,751)
\end{picture}
}

%% file: draws/ArrowX2.tex
%
\begingroup\makeatletter\ifx\SetFigFont\undefined%
\gdef\SetFigFont#1#2#3#4#5{%
  \reset@font\fontsize{#1}{#2pt}%
  \fontfamily{#3}\fontseries{#4}\fontshape{#5}%
  \selectfont}%
\fi\endgroup%
{\renewcommand{\dashlinestretch}{30}
\begin{picture}(2176,2133)(0,-10)
\path(727,661)(1476,658)
\path(1253,584)(944,661)(1253,732)
\path(12,2106)(1222,915)
\path(2164,2105)(963,915)
\path(1220,915)(2107,12)
\path(1811,178)(1611,516)(1938,310)
\path(769,1492)(434,1691)(638,1367)
\path(1416,1492)(1751,1691)(1547,1367)
\path(653,463)(321,261)(523,595)
\path(78,13)(965,913)
\end{picture}
}

%% file: draws/ArrowT1.tex
%
\begingroup\makeatletter\ifx\SetFigFont\undefined%
\gdef\SetFigFont#1#2#3#4#5{%
  \reset@font\fontsize{#1}{#2pt}%
  \fontfamily{#3}\fontseries{#4}\fontshape{#5}%
  \selectfont}%
\fi\endgroup%
{\renewcommand{\dashlinestretch}{30}
\begin{picture}(2766,2820)(0,-10)
\path(1791,1765)(1834,1450)(1943,1746)
\path(675,1765)(632,1450)(523,1746)
\path(1983,2793)(1833,1443)
\path(1833,1443)(1683,93)
\path(633,1443)(783,93)
\path(483,2793)(633,1443)
\thicklines
\path(33,93)(2733,93)
\blacken\path(2553.000,33.000)(2733.000,93.000)(2553.000,153.000)(2553.000,33.000)
\end{picture}
}

%% file: draws/ArrowT2.tex
%
\begingroup\makeatletter\ifx\SetFigFont\undefined%
\gdef\SetFigFont#1#2#3#4#5{%
  \reset@font\fontsize{#1}{#2pt}%
  \fontfamily{#3}\fontseries{#4}\fontshape{#5}%
  \selectfont}%
\fi\endgroup%
{\renewcommand{\dashlinestretch}{30}
\begin{picture}(2766,2820)(0,-10)
\path(555,1428)(598,1743)(707,1447)
\path(1907,1428)(1864,1743)(1755,1447)
\path(483,2793)(633,1443)
\path(633,1443)(783,93)
\path(1829,1443)(1679,93)
\path(1979,2793)(1829,1443)
\thicklines
\path(33,93)(2733,93)
\blacken\path(2553.000,33.000)(2733.000,93.000)(2553.000,153.000)(2553.000,33.000)
\end{picture}
}

%% file: draws/ArrowT3.tex
%
\begingroup\makeatletter\ifx\SetFigFont\undefined%
\gdef\SetFigFont#1#2#3#4#5{%
  \reset@font\fontsize{#1}{#2pt}%
  \fontfamily{#3}\fontseries{#4}\fontshape{#5}%
  \selectfont}%
\fi\endgroup%
{\renewcommand{\dashlinestretch}{30}
\begin{picture}(2766,2820)(0,-10)
\path(673,1765)(630,1450)(521,1746)
\path(1905,1437)(1862,1752)(1753,1456)
\path(1829,1443)(1679,93)
\path(631,1443)(781,93)
\path(481,2793)(631,1443)
\path(1977,2793)(1827,1443)
\thicklines
\path(33,93)(2733,93)
\blacken\path(2553.000,33.000)(2733.000,93.000)(2553.000,153.000)(2553.000,33.000)
\end{picture}
}

%% file: draws/ArrowT4.tex
%
\begingroup\makeatletter\ifx\SetFigFont\undefined%
\gdef\SetFigFont#1#2#3#4#5{%
  \reset@font\fontsize{#1}{#2pt}%
  \fontfamily{#3}\fontseries{#4}\fontshape{#5}%
  \selectfont}%
\fi\endgroup%
{\renewcommand{\dashlinestretch}{30}
\begin{picture}(2766,2820)(0,-10)
\path(555,1428)(598,1743)(707,1447)
\path(1787,1765)(1830,1450)(1939,1746)
\path(483,2793)(633,1443)
\path(633,1443)(783,93)
\path(1829,1443)(1679,93)
\path(1979,2793)(1829,1443)
\thicklines
\path(33,93)(2733,93)
\blacken\path(2553.000,33.000)(2733.000,93.000)(2553.000,153.000)(2553.000,33.000)
\end{picture}
}

%% file: draws/LeftArrow.tex
%
\begingroup\makeatletter\ifx\SetFigFont\undefined%
\gdef\SetFigFont#1#2#3#4#5{%
  \reset@font\fontsize{#1}{#2pt}%
  \fontfamily{#3}\fontseries{#4}\fontshape{#5}%
  \selectfont}%
\fi\endgroup%
{\renewcommand{\dashlinestretch}{30}
\begin{picture}(1086,1431)(0,-10)
\path(93,633)(993,633)
\thicklines
\path(93,33)(93,1383)
\blacken\path(153.000,1203.000)(93.000,1383.000)(33.000,1203.000)(153.000,1203.000)
\path(993,33)(993,1383)
\blacken\path(1053.000,1203.000)(993.000,1383.000)(933.000,1203.000)(1053.000,1203.000)
\thinlines
\path(693,560)(384,637)(693,708)
\end{picture}
}

%% file: draws/OneIn3.tex
%
\begingroup\makeatletter\ifx\SetFigFont\undefined%
\gdef\SetFigFont#1#2#3#4#5{%
  \reset@font\fontsize{#1}{#2pt}%
  \fontfamily{#3}\fontseries{#4}\fontshape{#5}%
  \selectfont}%
\fi\endgroup%
{\renewcommand{\dashlinestretch}{30}
\begin{picture}(2144,2119)(0,-10)
\path(695,678)(1444,675)
\path(324,1648)(661,1448)(456,1774)
\path(1222,898)(12,2088)
\path(1229,601)(920,678)(1229,749)
\path(465,587)(262,257)(595,457)
\path(2114,17)(1214,902)
\path(1661,592)(1864,262)(1531,462)
\path(12,12)(912,897)
\path(2132,2092)(931,902)
\path(1384,1479)(1719,1678)(1515,1354)
\end{picture}
}

%% file: draws/TwoInTwoOut3a.tex
%
\begingroup\makeatletter\ifx\SetFigFont\undefined%
\gdef\SetFigFont#1#2#3#4#5{%
  \reset@font\fontsize{#1}{#2pt}%
  \fontfamily{#3}\fontseries{#4}\fontshape{#5}%
  \selectfont}%
\fi\endgroup%
{\renewcommand{\dashlinestretch}{30}
\begin{picture}(2126,2131)(0,-10)
\path(695,694)(1444,691)
\path(12,28)(912,913)
\path(465,603)(262,273)(595,473)
\path(324,1664)(661,1464)(456,1790)
\path(1222,914)(12,2104)
\path(886,617)(1195,694)(886,765)
\path(1227,915)(2114,12)
\path(1818,178)(1618,516)(1945,310)
\path(2112,2097)(911,907)
\path(1364,1484)(1699,1683)(1495,1359)
\end{picture}
}

%% file: draws/TwoInTwoOut3b.tex
%
\begingroup\makeatletter\ifx\SetFigFont\undefined%
\gdef\SetFigFont#1#2#3#4#5{%
  \reset@font\fontsize{#1}{#2pt}%
  \fontfamily{#3}\fontseries{#4}\fontshape{#5}%
  \selectfont}%
\fi\endgroup%
{\renewcommand{\dashlinestretch}{30}
\begin{picture}(2126,2131)(0,-10)
\path(695,694)(1444,691)
\path(12,28)(912,913)
\path(465,603)(262,273)(595,473)
\path(324,1664)(661,1464)(456,1790)
\path(1222,914)(12,2104)
\path(1227,915)(2114,12)
\path(1818,178)(1618,516)(1945,310)
\path(2112,2097)(911,907)
\path(1364,1484)(1699,1683)(1495,1359)
\path(1203,619)(894,696)(1203,767)
\end{picture}
}

%% file: draws/AtomicTrivalent.tex
%
\begingroup\makeatletter\ifx\SetFigFont\undefined%
\gdef\SetFigFont#1#2#3#4#5{%
  \reset@font\fontsize{#1}{#2pt}%
  \fontfamily{#3}\fontseries{#4}\fontshape{#5}%
  \selectfont}%
\fi\endgroup%
{\renewcommand{\dashlinestretch}{30}
\begin{picture}(1224,1239)(0,-10)
\put(612.000,799.500){\arc{375.000}{3.7851}{5.6397}}
\path(495.158,964.207)(462.000,912.000)(515.825,942.461)
\path(12,1212)(612,612)
\path(1212,1212)(612,612)
\path(612,612)(612,12)
\end{picture}
}

%% file: draws/SkltVrtxHook1.tex
%
\begingroup\makeatletter\ifx\SetFigFont\undefined%
\gdef\SetFigFont#1#2#3#4#5{%
  \reset@font\fontsize{#1}{#2pt}%
  \fontfamily{#3}\fontseries{#4}\fontshape{#5}%
  \selectfont}%
\fi\endgroup%
{\renewcommand{\dashlinestretch}{30}
\begin{picture}(1566,1106)(0,-10)
\path(333,763)(473,670)
\path(473,670)(463,838)
\path(955,536)(1114,589)
\path(1114,589)(1061,430)
\put(-679.500,183.000){\arc{2493.116}{5.4765}{6.4038}}
\put(1514.250,64.250){\arc{1313.987}{3.0940}{4.5113}}
\thicklines
\path(33,33)(1533,33)
\end{picture}
}

%% file: qtrqH.tex
\section{Link invariants from quasitriangular quasi-Hopf algebras}
\label{sect:LinkInv}

{\bf Section summary:} After reviewing the basic theory of ribbon
quasitriangular quasi-Hopf algebras and their twist transformations, 
we go on to describe the examples of interest to us. These are the
Knizhnik-Zamolodchikov (KZ)  structures, having both Lie algebraic and 
diagrammatic versions, and the Drinfel'd-Jimbo (DJ) quantum groups.

\smallskip
We then work out in complete generality the construction of an invariant
$\Z_H$ of framed links associated to a ribbon quasitriangular quasi-Hopf 
algebra $H$ (possibly a diagrammatic one), showing that the natural target 
spaces of $\Z_H$ are spaces of coinvariants of $H$. (Close relatives of this
construction were considered by many authors; see Section 
~\ref{_sect:HinvOfLinks} for an explanation of our reasons for reworking it 
here.) We prove that $\Z_H$ remains unchanged when twisting the quasitriangular
quasi-Hopf structure of $H$.

\smallskip
In the last part of this section, we fit the Kontsevich integral and the 
Reshetikhin-Turaev invariants into the framework of $\Z_H$-invariants
corresponding to the KZ and DJ structures, respectively. This 
motivates the search for a diagrammatic structure on $\Aar$ having
some of the features of the DJ structures, which is the subject matter
of Section ~\ref{sect:EK}.  

\subsection{Quasitriangular quasi-Hopf algebras}
\label{_sect:qtrqH}

Let $H$ be an associative algebra over $\k$ with a product 
$\mu: H\ot H \ra H$ and a unit $\eta:\k \ra H$. For brevity we refer to $H$
simply as an algebra, and accordingly morphisms of algebras are 
always unital. Even though we have much in mind diagrammatic algebras, we do  
use in this section the element-wise notation $\mu(a, b) = ab$ and write all 
identities as equalities between $H$-tensors, rather than as commutative 
diagrams of maps. The translations to the diagrammatic
scenarios are usually straightforward, and are occasionally indicated.

\smallskip
We use throughout the following notational convention. Let 
$X = \sum_i x_i^{(1)} \ot\cdots\ot x_i^{(p)} \in H^{\ot p}$ be a tensor. For 
any $p$-tuple $(k_1, \ldots, k_p)$ of distinct indices from $\{1,\ldots, n\}$, 
$n\geq p$ (not necessarily in ascending order), we set
\[
X^{k_1\ldots k_p} = \sum_i y_i^{(1)} \ot\cdots\ot y_i^{(n)} \in H^{\ot n} 
\, \quad  \text{where} \quad \ y_i^{(k)} = 
	\begin{cases} x_i^{(j)} & k = k_j \text{ for some } j\leq p \, , \\
                      1 & \text{otherwise.}
	\end{cases}
\]
As $n$ does not appear explicitly in the notation, there is in principle
some danger of ambiguity (if $X\in H^{\ot 3}$, for example, then $X^{312}$
may stand for a tensor in any of the spaces $H^{\ot n}$, $n\geq 3$). In 
practice, however, it is always possible to deduce $n$ from the 
context.

Another convention we often employ is to drop the summation signs; summations 
over all indices of tensor components are implicitly implied.

As an example of both conventions, if $R = s_i \ot t_i \in H^{\ot 2}$, then 
$R^{31} \in H^{\ot 3}$ is given by
\[
	R^{31} = t_i\ot 1 \ot s_i \, .  
\]

\smallskip
If $\H$ is a diagrammatic algebra ($\A$ or $\Aar$), then the parallel operation
on a diagram $D\in \H^{\bt p}$ is adding $n-p$ oriented intervals to
the skeleton of $D$ (with no leg attached to them), and then relabeling the 
$n$ skeleton intervals as dictated by the $p$-tuple $(k_1, \ldots, k_p)$:
the $j$th skeleton component of $D$ ($j\leq p$) is labeled by $k_j$,
and the remaining $n-p$ labels are arbitrarily assigned to the added
skeleton components. Below is an example, with $D\in \A^{\bt 2}$.
\[
D = 
  \printname{A2Dgrm}
  \setlength{\unitlength}{0.4\standardunitlength}
  \begin{array}{c}  \hspace{-1.7mm}
    \raisebox{-8pt}{\input draws/A2Dgrm.tex }
      \hspace{-1.9mm}
  \end{array}
 \quad ; \quad 
D^{31} = 
  \printname{A2Dgrm31toA3}
  \setlength{\unitlength}{0.4\standardunitlength}
  \begin{array}{c}  \hspace{-1.7mm}
    \raisebox{-8pt}{\input draws/A2Dgrm31toA3.tex }
      \hspace{-1.9mm}
  \end{array}
 \, .
\]

Note that the algebraic and diagrammatic notations are compatible:
$\T_g(D^{31}) = (\T_g(D))^{31}$ for any metrized Lie algebra $\g$, and
similarly with $\Aar$-diagrams and $\Tartg$-maps, for any Manin
triple  $(\tg, \g_+, \g_-)$. 

\smallskip
In the algebraic framework we may write expressions like $1\ot X$ for 
the image of $X$ under the mapping $H^p \ra H^{p+1}$ determined
by the tuple $(2, \ldots, p+1)$. Accordingly, we shall write
$1\bt D$ for the image of $D\in H^{\bt p}$ under the corresponding
diagrammatic map.  
 
\subsubsection{\bf Definition of quasitriangular quasi-Hopf algebras.}
\label{__sect:QuasiHopf}

We recall the basics of the theory of quasitriangular quasi-Hopf 
algebras. Missing proofs for which explicit references are not given
may be found in the fundamental paper ~\cite{Drinfeld:QuasiHopf}. 

\smallskip
Assume an algebra $H$ is equipped with linear maps $\Delta: H \ra H\ot H$,
$\eps: H \ra \k$, and possibly $S:H\ra H$. We let 
\[
\Delta_i = \id \ot \cdots \ot \Delta \ot \cdots \ot \id 
  : H^{\ot n} \lra H^{\ot (n+1)} \qquad \text{($\Delta$ in the $i$th slot})
\]
be the map obtained by applying $\Delta$ to the $i$th component of
a tensor in $H^{\ot n}$, and similarly define 
$\eps_i : H^{\ot n} \ra H^{\ot (n-1)}$, $S_i : H^{\ot n} \ra H^{\ot n}$.

\begin{definition}
\lbl{def:QuasiBiAlg}
A {\em quasi-bialgebra} is an algebra $H$ equipped with 
an algebra morphism $\Delta: H \ra H\ot H$ (the {\em coproduct}), 
an algebra morphism $\eps: H \ra \k$ (the {\em counit}), and an invertible
element $\Phi\in H^{\ot 3}$ (the {\em associator}) such that
\begin{align}
(\Delta_2 \Delta)(x) &= \Phi (\Delta_1 \Delta)(x) \Phi^{-1}\, ,
\quad x\in H\, , \label{eqn:QuasiCoAssoc} \\
\Delta_3(\Phi) \Delta_1(\Phi) &= 
\Phi^{234} \Delta_2(\Phi) \Phi^{123} \, , \label{eqn:Pentagon} \\
\eps_1 \Delta &= \id = \eps_2 \Delta \, , 
\label{eqn:EpsForDelta} \\
 \eps_2(\Phi) &= 1\ot 1 \, . \label{eqn:EpsForPhi}
\end{align}    
\end{definition}
The associator $\Phi$ measures the failure of the coproduct to be 
coassociative. If $\Phi = 1\ot 1\ot 1$, then $\Delta$ is coassociative 
and $H$ is a bialgebra. \eqref{eqn:Pentagon} is called the
{\em pentagon relation} for $\Phi$.

\smallskip
It can be checked that \eqref{eqn:Pentagon}--\eqref{eqn:EpsForPhi}
imply
\begin{equation}
\label{eqn:PhiND}
\eps_1(\Phi) = 1\ot 1 = \eps_3(\Phi) \, .   
\end{equation}
In general, a tensor $X\in H^{\ot n}$ is called {\em non-degenerate}
if $\eps_i(X)= 1^{\ot (n-1)}$ for $i=1, \ldots ,n$. Thus, the 
associator $\Phi$ is non-degenerate.

\smallskip
We shall use a variant of Sweedler notation and write  
$\Delta(x) = x'\ot x''$ for the coproduct of $x\in H$, which implicitly
indicates a sum of factorable tensors. 
If
\begin{equation}
\label{eqn:PhiComps}
\Phi = X_i \ot Y_i \ot Z_i \quad \text{and} \quad 
\Phi^{-1} = \bar{X}_i \ot \bar{Y}_i \ot \bar{Z}_i \, ,
\end{equation}
then we may write \eqref{eqn:QuasiCoAssoc}, for example, 
in the form
\[
x'\ot (x'')' \ot (x'')'' =
X_i (x')' \bar{X}_j \ot Y_i (x')'' \bar{Y}_j \ot Z_i x'' \bar{Z}_j \, .
\] 

\smallskip
A morphism of quasi-bialgebras
\[
f: (H, \Delta, \eps, \Phi) \lra (H', \Delta', \eps', \Phi')
\]
is an algebra morphism satisfying
\[
	(f\ot f) \circ \Delta = \Delta' \circ f \, ,\quad
	\eps = \eps' \circ f \, , \quad
	(f\ot f\ot f)(\Phi) = \Phi' \, .	 
\] 
Definitions of morphisms for the more refined structures below are similar and
will not be given explicitly.

\smallskip
Let $H$ be a quasi-bialgebra.
The opposite coproduct of $H$ is 
$\Delta^{\op}(x) = \Delta(x)^{21} = x''\ot x'$. 
\begin{definition}
$H$ is {\em quasitriangular} if there exists an invertible element 
$R\in H^{\ot 2}$, called the {\em universal $R$-matrix} of $H$, such that
\begin{equation}
\Delta^{\op}(x) = R \Delta(x) R^{-1} \, , \quad x\in H \, ,
\label{eqn:QuasiCoComm} 
\end{equation}
\begin{align}
(\id\ot\Delta)(R) &= 
(\Phi^{231})^{-1} R^{13} \Phi^{213} R^{12} (\Phi^{123})^{-1} \, , 
\label{eqn:Hexagon1} \\  
(\Delta\ot\id)(R) &= 
\Phi^{312} R^{13} (\Phi^{132})^{-1} R^{23} \Phi^{123} \, .
\label{eqn:Hexagon2}   
\end{align}
\end{definition}
The universal $R$-matrix measures how far is the coproduct from
being cocommutative. \eqref{eqn:Hexagon1} -- \eqref{eqn:Hexagon2} are 
called the {\em hexagon relations}. It is not hard to see that $R$
is non-degenerate:
\begin{equation}
\label{eqn:EpsRmat}
\eps_1(R) = 1 = \eps_2(R) \, .
\end{equation}
We write the components of $R$ and its inverse as
\begin{equation}
\label{eqn:RmatComps}
R = s_i \ot t_i \quad \text{and} \quad 
R^{-1} = \bar{s}_i \ot \bar{t}_i \, .
\end{equation} 

\smallskip
The universal $R$-matrix satisfies the {\em quantum quasi-Yang-Baxter 
equation} 
\begin{equation}
\label{eqn:QQYBE}
R^{12} \Phi^{312} R^{13} (\Phi^{132})^{-1} R^{23} \Phi^{123} =
\Phi^{321} R^{23} (\Phi^{231})^{-1} R^{13} \Phi^{213} R^{12}\, ;
\end{equation}
see ~\cite[Section 3]{Drinfeld:QuasiHopf} or 
~\cite[Corollary XV.2.3]{Kassel:QuantumGroups} for a proof.   

\begin{definition}
\lbl{def:quasiHopf}
A {\em quasi-Hopf algebra} is a quasi-bialgebra $H$
together with an anti-automorphism $S:H\ra H$ of the underlying algebra 
($S(xy) = S(y)S(x)$), called the {\em antipode} of $H$, and elements 
$\alpha, \beta \in H$, such that 
\begin{align}
S(x') \alpha x'' &= \eps(x) \alpha \, , 
\label{eqn:AntipodeId1} \\
x' \beta S(x'') &= \eps(x) \beta \, , 
\label{eqn:AntipodeId2} 
\end{align}
for all $x\in H$, and, with the notation of \eqref{eqn:PhiComps},
\begin{align}
X_i \beta S(Y_i) \alpha Z_i &= 1 \, , 
\label{eqn:AntipodePhi} \\  
S(\bar{X}_i) \alpha \bar{Y}_i \beta S(\bar{Z}_i) &= 1 \, . 
\label{eqn:AntipodePhiInv}   
\end{align}
A quasi-Hopf algebra is quasitriangular if the underlying quasi-bialgebra
is.
\end{definition}
We will also assume that the element $\alpha$ is invertible, as is the case 
in all of our examples.

\eqref{eqn:AntipodeId1} -- \eqref{eqn:AntipodeId2} are the
quasi-versions of the antipodal identities \eqref{eqn:antipod}. Two
consequences of the conditions satisfied by $S, \alpha, \beta$ are
\[
	\eps(\alpha)\eps(\beta) = 1 \, , \quad \eps \circ S = \eps \, .
\]

\smallskip
If $(S, \alpha, \beta)$ is an antipodal triple for a quasi-Hopf algebra $H$,
then another antipodal triple can be obtained by an
{\em antipodal rescaling} transformation of the form 
\begin{equation}
\label{eqn:AntipodalRescale}
S \mapsto \hat{S} \ , \hat{S}(x)=cS(x)c^{-1} \, ; \quad 
\alpha \mapsto c \alpha\, ; \quad \beta \mapsto \beta c^{-1} \, ,
\end{equation}
where $c\in H$ is an invertible element.
In fact, the antipodal triple is unique up to antipodal rescaling 
(\cite[Proposition 1.1]{Drinfeld:QuasiHopf}). 

\smallskip
If $H$ is a {\em coassociative} quasi-Hopf algebra, then  
\eqref{eqn:AntipodePhi} -- \eqref{eqn:AntipodePhiInv} show that 
$\beta = \alpha^{-1}$, hence rescaling with $c=\beta$ gives an 
antipodal triple $(\hat{S}, 1, 1)$. This means that $\hat{S}$ satisfies
the axioms of a {\em Hopf algebra antipode}, so that $H$ is a Hopf algebra.
We always use this normalization of the antipodal triple
for Hopf algebras. 

\subsubsection{\bf The square of the antipode.}

Let $H$ be a quasitriangular quasi-Hopf algebra, with $\Phi$ and $R$ as 
in \eqref{eqn:PhiComps} and 
\eqref{eqn:RmatComps}. 
\begin{definition}
The {\em (Drinfel'd) $u$-element} of $H$ is given by
\begin{equation}
\label{eqn:uElement}
	u = S(\bar{Y}_i \beta S(\bar{Z}_i))S(t_j) \alpha s_j \bar{X}_i \, .
\end{equation}
\end{definition}

The following proposition is due to Altsch{\"u}ler and Coste 
(\cite[Section 3]{AltschulerCoste:QuasiQuantum}),
generalizing a result of Drinfel'd for quasitriangular Hopf algebras.
\begin{proposition}
The $u$-element is invertible, and the inner automorphism of $H$ it determines
is the square of the antipode: for any $x\in H$,
\begin{equation}
\label{eqn:SqrAntipode}
	S^2(x) = uxu^{-1} \, .
\end{equation}
We also have
\begin{equation}
\label{eqn:Salpha_u}
	S(\alpha)u = S(t_i) \alpha s_i \, .
\end{equation}
\end{proposition}
An easy corollary of \eqref{eqn:SqrAntipode} is that $uS(u) = S(u)u$.

\begin{definition}
\lbl{def:vElement}
$H$ is a {\em ribbon} quasitriangular quasi-Hopf if there exists
a central element $v \in H$ (called the {\em ribbon element} of $H$) 
such that
\begin{equation}
\label{eqn:vProps}
v^2 = uS(u) = S(u)u \, , \quad
S(v) = v \, , \quad
\eps(v) = 1 \, ,
\end{equation} 
and
\begin{equation}
\label{eqn:Delta_v}
\Delta(v) = (v\ot v) (R^{21} R)^{-1} \, .
\end{equation}
\end{definition}

\begin{remark}
We have used the assumption that $\alpha$ is invertible in order to be
able to write the condition \eqref{eqn:Delta_v} in the current form. A more
complicated condition is needed in the general case. (See 
~\cite[Section 4.1]{AltschulerCoste:QuasiQuantum} for the general form of the  
condition and the end of ~\cite[Section 4.4]{AltschulerCoste:QuasiQuantum} for
a proof of its equivalence to \eqref{eqn:Delta_v} when $\alpha$ is invertible.)
\end{remark}

\begin{remark}
\lbl{rem:Central-vs-Invar}
For diagrammatic examples, we replace the requirement of
centrality in Definition ~\ref{def:vElement} by that of {\em invariance} 
(see Definition ~\ref{def:InvarDgrm}).
 
\smallskip
We might have attempted to define centrality with respect to the 
$\Vect_{\k}$-products of $\A$ and $\Aar$, and use this notion of 
centrality in Definition ~\ref{def:vElement}. The problem
is that in several $\Vect_{\k}$ computations, we use the fact that if $x\in H$ 
is central, then so is, for example, $1\ot x \in H^{\ot 2}$. (See the end
of the proof of Proposition ~\ref{prop:RibbonTwist} below for an example). The
hypothetical $\Vect_{\k}$-notion of diagrammatic centrality may not be enough 
to guarantee that if $D$ is central in $\Aar$, then $1\bt D$ commutes with 
every element of $\lowupper{\Aar}{\bt 2}$; invariance does guarantee that, and 
allows us to use the $\Vect_{\k}$ results in diagrammatic computations. 

\smallskip
The only diagrammatic ribbon element we shall actually encounter in $\Aar$  
will belong to the image of $\iota:\A \ra \Aar$, and is thus
invariant by Proposition ~\ref{prop:ImgIotaInvar}.
\end{remark}

\smallskip
The following proposition describes the effect of antipodal rescaling
on the $u$-element and the ribbon element.
\begin{proposition}
\lbl{prop:RibbonRescale}
Let $H$ be a ribbon quasitriangular quasi-Hopf algebra. If the antipodal
triple $(S, \alpha, \beta)$ is rescaled by an invertible element $c$
into an antipodal triple $(\hat{S}, \hat{\alpha}, \hat{\beta})$, then
\[
	\hat{u} = cS(c^{-1})u\, , \quad \hat{v} = v \, .
\] 
\end{proposition}
\begin{proof}
The formula for $u$ follows by substituting 
\eqref{eqn:AntipodalRescale} into \eqref{eqn:uElement}. To show that 
$\hat{v} = v$, we need to check that $\hat{S}(\hat{u})\hat{u} = S(u)u$.
Using \eqref{eqn:SqrAntipode} and the centrality (or invariance) of $v^2$,
\[
   \hat{S}(\hat{u})\hat{u} = cS(u)S^2(c^{-1})u = cS(u)uc^{-1}u^{-1}u
	= c v^2 c^{-1} = v^2 = S(u)u \, ,
\]
as required.
\end{proof}

\subsubsection{\bf Twisting a quasitriangular quasi-Hopf algebra.}
\label{__sect:twist}

\begin{definition}
Let $H = (H, \Delta, \eps, \Phi, R, S)$ be a quasitriangular 
quasi-algebra. An element $F\in H\ot H$ is a {\em twist} (alternatively
called {\em gauge transformation}) if $F$ is invertible and non-degenerate. 
\end{definition}

Write
\begin{equation}
\label{eqn:TwistComps}
F = f_i\ot g_i \quad \text{and} \quad F^{-1} = \bar{f}_i \ot \bar{g}_i \, ,
\end{equation}
and define a new set of structure ingredients as follows:
\begin{align}
\Delta_F(x) &= F\Delta(x) F^{-1} \, , \quad x\in H 
	\, , \label{eqn:TwistDelta}\\
\eps_F &= \eps \, , \\
\Phi_F &= F^{23}\Delta_2(F) \Phi \Delta_1(F^{-1}) (F^{12})^{-1}
	\, , \label{eqn:TwistPhi} \\
\R_F &= F^{21} R F^{-1} \, , \label{eqn:TwistR}
\end{align}

\begin{proposition}
$H_F = (H, \Delta_F, \eps, \Phi_F, R_F)$ is 
a quasitriangular quasi-bialgebra. 
If $F_1, F_2$ are both twists, then so is their
product $F_1 F_2$, and $H_{F_1 F_2} = (H_{F_2})_{F_1}$.
\end{proposition}
We say that $H_F$ is obtained from $H$ by {\em twisting via $F$} (or that 
$F$ twists $H$ into $H_F$), and that
$H$, $H_F$ are {\em twist equivalent}.   

\smallskip
If $H$ is a quasitriangular quasi-Hopf algebra, with antipodal triple 
$(S, \alpha, \beta)$,
then $H_F$ is also quasi-Hopf; an antipodal triple for $H_F$ is given by
\begin{equation}
\label{eqn:TwistAntipode}
\begin{aligned} 
S_F & = S \, , \\
\alpha_F &= S(\bar{f}_i) \alpha \bar{g}_i \, , \\
\beta_F &= f_i \beta S(g_i) \, . 
\end{aligned}
\end{equation}
To adhere to our assumption that the $\alpha$-element of a quasi-Hopf
algebra is invertible, we impose on $F$ the additional condition that 
$\alpha_F$ is invertible; this will be trivially satisfied in our examples.

\smallskip
If $H, H'$ are quasitriangular quasi-Hopf algebras, we say that they are
twist equivalent if the underlying quasitriangular quasi-bialgebras are. Note
that if  $H'$ is obtained from $H$ by twisting via $F$, then the antipodal 
triple $(S', \alpha', \beta')$ of $H'$ need not coincide with the antipodal
triple given by \eqref{eqn:TwistAntipode}, but, as discussed at 
the end of Section ~\ref{__sect:QuasiHopf}, the two antipodal triples are
related by an antipodal rescaling transformation.

\begin{proposition}
\lbl{prop:RibbonTwist}
Let $H, H'$ be quasitriangular quasi-Hopf algebras.
If $H$ has a ribbon element $v$ and $H'$ is twist equivalent to $H$ 
(via $F \in H^{\ot 2}$), then $v$ is also a ribbon element for $H'$.
\end{proposition}
\begin{proof}
We may assume, by Proposition ~\ref{prop:RibbonRescale}, that the antipodal
triple of $H'$ is given by \eqref{eqn:TwistAntipode}. We denote by $u_F$ 
the $u$-element of $H'$.  

By applying $S$ to \eqref{eqn:Salpha_u}, and then using 
\eqref{eqn:SqrAntipode} and the centrality of $v$, we obtain
\begin{equation}
\label{eqn:Su_u_alpha}
	S(u)u\alpha = S(t_i)S(\alpha)u t_i = S(s_i) S(t_j) \alpha s_j t_i \, .
\end{equation}
We apply this formula to $H'$, to compute $S_F(u_F)u_F \alpha_F$. 
Substituting $R_F = g_i s_j \bar{f}_k \ot f_i t_j \bar{g}_k$, 
$\alpha_F = S(\bar{f}_l) \alpha \bar{g}_l$, and using $S_F=S$, 
a straightforward computation (utilizing once more each of \eqref{eqn:Salpha_u}
and \eqref{eqn:Su_u_alpha}) leads to
\[
	S_F(u_F)u_F \alpha_F = S(\bar{f}_l) S(u)u\alpha \bar{g}_l 
		= v^{2} \alpha_F \, .		
\]
Since $\alpha_F$ is invertible, we have $v^2 = S_F(u_F)u_F$.

It remains to check that
\[
	\Delta_F(v) = (v\ot v) (R_F^{21} R_F)^{-1} \, .
\]
As $(R_F^{21} R_F)^{-1} = F (R^{21} R)^{-1} F^{-1}$ and 
$\Delta_F(v) = F\Delta(v) F^{-1}$, this amounts to
\[
	\Delta(v) = F^{-1}(v\ot\ v) F (R^{21} R)^{-1} \, ,
\]
which follows from the centrality (or invariance, for diagrammatic
algebras) of $v\ot v$ and \eqref{eqn:Delta_v}.
\end{proof}
See ~\cite[Section 5]{LinksGouldZhang:TwistInvarLinkPoly} for a proof
which does not rely on the assumption that $\alpha$ is invertible.

\subsection{The KZ structures}
\label{_sect:KZ}

We turn to describe a set of examples of quasitriangular quasi-Hopf
structures originating from the monodromy of the Knizhnik-Zamolodchikov (KZ) 
system of differential equations.

\subsubsection{\bf The algebra of horizontal chords and KZ associator.}
We recall the spaces $\A^h_n = \A^h(\ua_1, \ldots, \ua_n)$ 
that were defined in Remark ~\ref{rem:A-horizontal}.  $A^h_n$ is the 
quotient by 4T relations of the formal span of the space $\D^h_n$ of chord 
diagrams based on the skeleton $\M= \{\ua_1, \ldots, \ua_n\}$, all of their 
chords are horizontal.
It has the usual grading by half the number of vertices (= number of
chords). As with $\A$ and $\Aar$, we almost always work with the graded 
completion of $\A^h_n$, without introducing a separate notation.

It is customary to refer to the oriented intervals making up the skeleton
of $\A^h_n$ as {\em strands}. We denote by $t^{ij} \in \A^h_n$ 
($1 \leq i \neq j \leq n$) the element
represented by a diagram with a single chord connecting the $i$th strand
to the $j$th strand:
\[
t^{ij} = 
  \printname{TijDgrm}
  \setlength{\unitlength}{0.4\standardunitlength}
  \begin{array}{c}  \hspace{-1.7mm}
    \raisebox{-8pt}{\input draws/TijDgrm.tex }
      \hspace{-1.9mm}
  \end{array}
 \, .
\]

As for $\A^{\bt n}$, the juxtaposition product $D_1 D_2$ of diagrams 
$D_1, D_2 \in \A^h_n$ is defined
by glueing the head of the $i$th strand of $D_1$ to the tail of
the $i$th strand of $D_2$, for $i=1, \ldots, n$; 
it endows $\A^h_n$ with the 
structure of an algebra (in $\Vect_{\k}$). 
The algebra $\A^h_n$ is generated by the elements $t^{ij}$, which satisfy the
relations 
\begin{equation}
\label{eqn:AhorRels} 
t^{ij} = t^{ji}\, , \qquad [t^{ij}, t^{kl}] = 0 \quad \text{if} \
\left| \{i, j, k, l\} \right| = 4\, , 
\end{equation}
and
\begin{equation}
\label{eqn:AhorRels4T}  
[t^{jk}, t^{ij}+t^{ik}] = 0 \quad \text{if} \
\left| \{i, j, k\} \right| = 3\, . 
\end{equation}
The last relations follow from the 4T relations; for example,
\begin{align*}
[t^{34}, t^{13} &+ t^{14}] =  \\
& 
  \printname{Ahor4T}
  \setlength{\unitlength}{0.3\standardunitlength}
  \begin{array}{c}  \hspace{-1.7mm}
    \raisebox{-8pt}{\input draws/Ahor4T.tex }
      \hspace{-1.9mm}
  \end{array}

= 0 \quad \text{in} \ \A^h_4 \, .
\end{align*}
(The strands are labeled $1$--$4$ from left to right.)
It is not hard to show that $\A^h_n$ is isomorphic to the graded
completion of the algebra generated over $\k$ by symbols $t^{ij}$ 
(of degree $1$) subject to the relations \eqref{eqn:AhorRels} and 
\eqref{eqn:AhorRels4T}.  
 
\smallskip
We can define the cabling and strand removal operations
\[
	\Delta_i = \Delta_{\ua_i} : \A^h_n \ra \A^h_{n+1} \quad ; \quad
	\eps_i = \eps_{\ua_i} : \A^h_n \ra \A^h_{n-1} \, ,
\]
just as for $\A^{\bt n}$. The injection $A^h_p \ra \A^h_n$ $(p\leq n)$,
corresponding to a tuple $(k_1, \ldots, k_p)$, is defined as in
Section ~\ref{_sect:qtrqH}, by adding and relabeling strands. 
Thus, the notion of non-degeneracy and the pentagon and hexagon relations
make sense in the sequence of spaces\footnote{These spaces may be thought of
as the diagrammatic tensor powers of $\A^h_1$, even though the latter
space is spanned by the single diagram having no chords.} $\{\A^h_n\}$, and
we can ask whether there exist a pair of non-degenerate, invertible  elements
$(R\in A^h_2$, $\Phi\in \A^h_3)$,  solving \eqref{eqn:Pentagon}, 
\eqref{eqn:Hexagon1}, and \eqref{eqn:Hexagon2}. We are particularly
interested in solutions $(R, \Phi)$ in which the $R$-matrix is given by
\begin{equation}
\label{eqn:Rkz}
R = e^{t^{12}/2} = \sum_{m=0}^{\infty} \frac{1}{2^m m!}
	{
  \printname{T12Pwr}
  \setlength{\unitlength}{0.3\standardunitlength}
  \begin{array}{c}  \hspace{-1.7mm}
    \raisebox{-8pt}{\input draws/T12Pwr.tex }
      \hspace{-1.9mm}
  \end{array}
}_{\ \text{($m$ chords)}}	\, ,
\end{equation}
where the power series, having finitely many terms in each degree, is 
convergent in (the graded completion of) $\A^h_2$. Such an element 
$\Phi$ is then called a {\em Drinfel'd associator}. 

\smallskip
For $k = \mathbb{C}$, Drinfel'd (\cite[Section 2]{Drinfeld:GalQQ}) 
proved the existence of such an associator, denoted $\Phi_{\kz}$, 
by studying the monodromy of the KZ system
\[
\frac{\partial W}{\partial z_i} = \frac{1}{2\pi i}\sum_{j\neq i}
\frac{t^{ij}}{z_i - z_j} W\, , \quad 
W:\{(z_1, \ldots, z_n) \in \mathbb{C}^n \, \mid \, 
	z_i \neq z_j \text{ for } i\neq j\} \ra \A^h_n \, .
\]
The element $\Phi_{\kz}\in A^h_3$, called the {\em KZ associator}, is  
highly complicated; it can be written explicitly in terms of
{\em multiple zeta values} (see ~\cite[Section 4]{LeMurakami:Universal}).

Drinfel'd further established the following rationality result.
\begin{theorem}{~\cite[Theorem A']{Drinfeld:GalQQ}}
\lbl{thm:RationalPhi}
There exists a Drinfel'd associator defined over $\mathbb{Q}$. 
\end{theorem}

This theorem implies that there exists a Drinfel'd associator
over any field $\k$ (of characteristic zero). We fix one such an
associator $\Phi$ from now on. 
We note that the non-degeneracy of $\Phi$ 
implies that  
\begin{equation}
\label{eqn:HorPhi}
	\Phi = 
  \printname{A3hUnit}
  \setlength{\unitlength}{0.2\standardunitlength}
  \begin{array}{c}  \hspace{-1.7mm}
    \raisebox{-8pt}{\input draws/A3hUnit.tex }
      \hspace{-1.9mm}
  \end{array}
 + \sum_i c_i D_i \, , \quad
	c_i\in \k\, , D_i \in \D^h_3\, , 
\end{equation}
where every strand of a diagram $D_i$ touches at least one chord; in 
particular,
$\deg D_i \geq 2$, because all the elements of degree $1$ in $A^h_3$ 
(linear combinations of $t_{12}, t_{13}, t_{23}$) fail to meet this condition. 

\begin{remark}
The existence of the associator $\Phi_{\kz}$ (over $\mathbb{C}$) plays an
essential role in the proof of Theorem ~\ref{thm:RationalPhi}. Accordingly,
we refer to the structures involving $\Phi$ as {\em KZ structures}, even when 
$\k \neq \mathbb{C}$ or if another Drinfel'd associator is used over 
$\mathbb{C}$. 
\end{remark}

\subsubsection{\bf KZ structures on $\protect\A$ and $\protect\Aar$.}

The inclusion maps $\D^h_n \ra \D(\ua_1, \ldots, \ua_n)$ descend to maps
$\phi^h : \A^h_n \ra \A^{\bt n}$; see Section ~\ref{__sect:ChordToJacobi} for
the proof that 4T relations map to zero. These maps evidently
interchange all the algebraic structures on $\A^h$-spaces with the
parallel structures on the spaces $\A^{\bt n}$.

The maps $\phi^h$ are actually injective 
(see ~\cite{Bar-Natan:Braids}, in particular Corollary 4.4). We will not need 
this fact in the sequel, but we do retain the notation $\Phi$ for
$\phi^h(\Phi)$.
Notice that $t^{12}$ is mapped by $\phi^h$ to the diagrammatic metric tensor 
$\Omega$ (see Example ~\ref{exmpl:CasimirDgrm}), and hence 
$R = e^{t^{12}/2}$ is mapped to $R = e^{\Omega/2}$.  

The elements $R \in \A^{\bt 2}$ and 
$\Phi \in \A^{\bt 3}$ are the major components of a quasitriangular
quasi-Hopf structure on $\A$. 
\begin{proposition}
\lbl{prop:Akz}
$(\A, \Delta, \eps, \Phi, R)$ is a diagrammatic quasitriangular 
quasi-bialgebra. 
\end{proposition}

We stress that this structure shares the same coproduct (and counit) with 
the usual bialgebra structure $(A, \Delta, \eps, 1^{\bt 3}, 1^{\bt 2})$.
\begin{proof}
The pentagon and hexagon relations follow from the fact that they hold in 
$\A^h$-spaces, and so does \eqref{eqn:EpsForPhi}. We know that $\Delta$ and 
$\eps$ are algebra morphisms and that \eqref{eqn:EpsForDelta} holds,
because these properties are enjoyed by
the usual Hopf structure. Thus, it remains to check 
\eqref{eqn:QuasiCoAssoc} and \eqref{eqn:QuasiCoComm}. As $\Delta$ is
coassociative and cocommutative, they take the form
\[
	[(\Delta_2\Delta)(D), \Phi] = 0 \quad \text{and} \quad
	[\Delta(D), R] = 0 \quad 
	\text{for all} \ D\in \A \, ,
\]
which follows at once from Proposition ~\ref{prop:AnInv}.
\end{proof}

\smallskip
We denote $\A$ equipped with the above described structure by $\A_{\kz}$.
It remains to describe an antipodal triple of $\A_{\kz}$. 
For that we need a short preparation.
\begin{definition}
An element $D\in A^{\bt n}$ is called a {\em perturbation of the identity} if
$D = 1^{\bt n} + D_{> 0}$, where $1^{\bt n}$ is the no-chord diagram on
$n$ strands, and $D_{>0}$ is $D$ with its degree $0$ piece removed.
\end{definition}
Clearly, the $R$-matrix and associator of $\A_{\kz}$ are perturbations of 
the identity.

A perturbation of the identity in $\A^{\bt n}$ is 
invertible\footnote{The notion
of an element being invertible is always with respect to the 
$\Vect_{\k}$-product.},
because a two-sided inverse can be defined using the power series
of $(1+x)^{-1}$. We note for a later purpose that, similarly,  
we can use the power series of $\sqrt{1+x}$ to define the square root 
$D^{\frac{1}{2}}$ of a  perturbation of the identity $D$. 

\smallskip
We return to describing the antipodal triple of $\A_{\kz}$. 
Its antipode is the same as the antipode of the usual Hopf structure of 
$\A$ (see Definition ~\ref{def:CompReverse}). Put 
\begin{equation}
\label{eqn:NuDef}
\nu= (X_i S(Y_i) Z_i)^{-1} \in \A \, .
\end{equation}
Since all elements of $\A$ are invariant, 
\eqref{eqn:AntipodeId1} and \eqref{eqn:AntipodeId2} hold for any choice of 
$\alpha, \beta$, and we can satisfy \eqref{eqn:AntipodePhi} by choosing
\begin{equation}
\label{eqn:AkzAlpha}
\alpha = \nu ; \quad \beta = 1\, . 
\end{equation}
In diagrammatic notation, $\nu$ is the inverse of
$(\mu_{\ua_1, \ua_3}\mu_{\ua_1, \ua_2} S_{\ua_2})(\Phi)$; the latter 
element is indeed invertible, because it is a perturbation of the identity 
in $\A$. If we depict (a summand of) $\Phi$ as
\[

  \printname{PhiPict}
  \setlength{\unitlength}{0.4\standardunitlength}
  \begin{array}{c}  \hspace{-1.7mm}
    \raisebox{-8pt}{\input draws/PhiPict.tex }
      \hspace{-1.9mm}
  \end{array}
\, ,
\]
then $\nu^{-1}$ is (omitting the parity sign in the definition of $S$)
\[

  \printname{NuInvPict}
  \setlength{\unitlength}{0.4\standardunitlength}
  \begin{array}{c}  \hspace{-1.7mm}
    \raisebox{-8pt}{\input draws/NuInvPict.tex }
      \hspace{-1.9mm}
  \end{array}
\, ,
\]
and the antipodal identity \eqref{eqn:AntipodePhi}, with 
$\beta = 1 = \text{empty diagram}$, looks like
\[

  \printname{AntipodeIdPict}
  \setlength{\unitlength}{0.4\standardunitlength}
  \begin{array}{c}  \hspace{-1.7mm}
    \raisebox{-8pt}{\input draws/AntipodeIdPict.tex }
      \hspace{-1.9mm}
  \end{array}
\, .
\]

We will use below more algebraic expressions as 
in \eqref{eqn:NuDef} for elements of $\A$-spaces;
a little thought shows in each case how to translate those expressions
and algebraic manipulations with them into diagrammatic arguments.  

\smallskip
We claim that the triple $(S, \alpha = \nu, \beta = 1)$ also satisfies 
\eqref{eqn:AntipodePhiInv}.
Here is a sketch of the proof. We rewrite the pentagon relation 
\eqref{eqn:Pentagon} as
\[
(\id\bt\id\bt\Delta)(\Phi) \cdot (\Delta\bt\id\bt\id)(\Phi) 
\cdot (\Phi^{-1}\bt 1) = (1\bt\Phi) \cdot (\id\bt\Delta\bt\id)(\Phi) \, .
\]
We apply to both sides the map $\A^{\bt 4} \ra \A$ taking 
$D\in \A^{\bt 4}$ to 
\[
\mu^{(7)} (S_1\bt S_5)(D\bt\alpha\bt\beta\bt\alpha)^{1357246} \, , 
\]
where $S_i = S_{\ua_i}$, and $\mu^{(7)}: \A^{\bt 7} \ra \A$ is the iterated 
product map, concatenating the seven skeleton components, with $\ua_1$ before 
$\ua_2$ before \ldots  $\ua_7$. (Just for orientation, the corresponding 
map in a $\Vect_k$ notation is: 
$a \ot b \ot c \ot d \ra S(a) \alpha b \beta S(c) \alpha d$.) The 
resulting expressions on both sides can be simplified in several steps, by 
using \eqref{eqn:AntipodeId1}, \eqref{eqn:AntipodeId2}, and the 
non-degeneracy of $\Phi$, until we get
\[
(S(\bar{X}_i \alpha \bar{Y}_i \beta S(\bar{Z}_i)) \alpha =
\alpha (X_j \beta S(Y_j) \alpha Z_j) = \alpha \, ,
\]
where the last equality follows from \eqref{eqn:AntipodePhi}. Since 
$\alpha$ is invertible, we see that \eqref{eqn:AntipodePhiInv} also holds 
true. Thus, $(S, \nu, 1)$ is an antipodal triple for $\A_{\kz}$.

\begin{remark}
The argument above is an adaptation of an excerpt from the proof of 
~\cite[Proposition 1.3]{Drinfeld:QuasiHopf}, which basicly states that for
a quasitriangular quasi-bialgebra in $\Vect_{\k}$, conditions
\eqref{eqn:AntipodeId1} -- \eqref{eqn:AntipodePhi} imply
condition \eqref{eqn:AntipodePhiInv}. (The proof there is somewhat more
complicated, as it avoids the need to assume that $\alpha$ is invertible.)
\end{remark}

\smallskip
We now want to compute the $u$-element of $\A_{\kz}$. Before doing that,
we state a useful property of the associator and comment on its proof.
\begin{proposition}
The associator of $\A_{\kz}$ satisfies
\begin{equation}
\label{eqn:NormalizedPhi}
\Phi^{-1} = \Phi^{321} \, .
\end{equation}
\end{proposition}
Indeed, Drinfel'd (\cite[Proposition 3.5]{Drinfeld:QuasiHopf}) proves that
for a certain class of graded quasitriangular quasi-Hopf algebras, called 
QUE algebras (see Definition ~\ref{def:QUE} below),
the symmetry of the $R$-matrix ($R^{21} = R$) 
implies \eqref{eqn:NormalizedPhi}. An inspection of the proof shows that 
it carries over to $\A_{\kz}$; except for algebraic manipulations valid in any
quasitriangular quasi-Hopf ($\Vect_{\k}$ or diagrammatic) algebra, the proof 
uses only the ability to extract square roots of perturbations of the 
identity. 

Another version of the proof, already in a diagrammatic framework, 
can be found in ~\cite[Proposition 3.7]{Bar-Natan:NAT}, which utilizes the 
restriction to braids of the invariants to be described in Section 
~\ref{_sect:HinvOfLinks}.

\smallskip
Plugging \eqref{eqn:NormalizedPhi} into \eqref{eqn:NuDef}, we find that
\[
\nu = (\bar{Z}_i S(\bar{Y}_i) \bar{X}_i)^{-1}\, .
\]
We compute $u$ according to \eqref{eqn:uElement}. Using the 
commutativity of $\A$, that $S^2 = \id$, and the last formula for $\nu$, we 
have
\begin{align*}
u & = S(\bar{Y}_i \beta S(\bar{Z}_i))S(t_j) \alpha s_j \bar{X}_i  \\
& = (\bar{Z}_i S(\bar{Y}_i)  \bar{X}_i) \nu S(t_j) s_j  \\
& = S(t_j) s_j = \mu(S\bt \id)(R^{21}) \, .
\end{align*}
From \eqref{eqn:Rkz} we get
\[
u = \sum_{m=1}^{\infty} \frac{(-1)^m}{2^m m!} 
	{
  \printname{RainBow}
  \setlength{\unitlength}{0.3\standardunitlength}
  \begin{array}{c}  \hspace{-1.7mm}
    \raisebox{-8pt}{\input draws/RainBow.tex }
      \hspace{-1.9mm}
  \end{array}
}_{\ \text{($m$ chords)}} \, .	
\]
Repeatedly moving the rightmost leg across the half circles below it until
it is immediately to the right of its origin (which we may do by the  
invariance of the diagram made of the lower circles), we obtain
\[
	
  \printname{RainBow}
  \setlength{\unitlength}{0.3\standardunitlength}
  \begin{array}{c}  \hspace{-1.7mm}
    \raisebox{-8pt}{\input draws/RainBow.tex }
      \hspace{-1.9mm}
  \end{array}
 = 
  \printname{HalfCircPwr}
  \setlength{\unitlength}{0.3\standardunitlength}
  \begin{array}{c}  \hspace{-1.7mm}
    \raisebox{-8pt}{\input draws/HalfCircPwr.tex }
      \hspace{-1.9mm}
  \end{array}
 = 
	\left( 
  \printname{HalfCirc}
  \setlength{\unitlength}{0.3\standardunitlength}
  \begin{array}{c}  \hspace{-1.7mm}
    \raisebox{-8pt}{\input draws/HalfCirc.tex }
      \hspace{-1.9mm}
  \end{array}
 \right)^m \, ,
\]
and we conclude that
\begin{equation}
\label{eqn:uElmntKZ}
u = e^{-C/2} \, ,
\end{equation}
where $C$ is the diagrammatic Casimir element (Example 
~\ref{exmpl:CasimirDgrm}).

The coproduct of $C$ is
\begin{equation}
\label{eqn:DeltaCas}
\Delta(C) = 
  \printname{1btC}
  \setlength{\unitlength}{0.25\standardunitlength}
  \begin{array}{c}  \hspace{-1.7mm}
    \raisebox{-8pt}{\input draws/1btC.tex }
      \hspace{-1.9mm}
  \end{array}
 + 
  \printname{Cbt1}
  \setlength{\unitlength}{0.25\standardunitlength}
  \begin{array}{c}  \hspace{-1.7mm}
    \raisebox{-8pt}{\input draws/Cbt1.tex }
      \hspace{-1.9mm}
  \end{array}
 + 2 
  \printname{BCup}
  \setlength{\unitlength}{0.25\standardunitlength}
  \begin{array}{c}  \hspace{-1.7mm}
    \raisebox{-8pt}{\input draws/BCup.tex }
      \hspace{-1.9mm}
  \end{array}
 =
	1\bt C + C\bt 1 + 2\Omega\, .
\end{equation}
The elements on the right hand side commute with each other,
thanks to the invariance of $C$ (compare Remark ~\ref{rem:Central-vs-Invar}).
Therefore,
\[
\Delta(u) = e^{-\Delta(C)/2} =  
  (1\bt e^{-C/2}) (e^{-C/2}\bt 1) e^{-\Omega} = (u\bt u) (R^{21} R)^{-1} \, .
\] 
As $S(u) = u$, it follows that $v=u$ meets all the conditions of 
Definition ~\ref{def:vElement}, and is a ribbon element for $\A_{\kz}$.

\smallskip
We record another consequence of \eqref{eqn:DeltaCas}. As $\Delta(C)$
commutes with all $D\in \A^{\bt 2}$ (by Proposition ~\ref{prop:AnInv}),
and so do $1\bt C$ and $C\bt 1$, we have:
\begin{proposition}
\lbl{prop:OmegaA2Central}
$[\Omega, D] = 0$ for all $D\in \A^{\bt 2}$.
\end{proposition}

\smallskip   
The quasitriangular quasi-Hopf structure of $\A_{\kz}$ depends on the
choice of the Drinfel'd associator $\Phi$. Suppose $\tilde{\Phi}\in \A^h_3$ 
is another Drinfel'd associator, and denote by $\tilde{\A}_{\kz}$ the
structure corresponding to $\tilde{\Phi}$.

The following result is due to Le and Murakami 
(\cite[Theorem 8]{LeMurakami:Universal}). It is the diagrammatic 
analogue of a result of Drinfel'd for universal enveloping algebras
(\cite[Theorem 3.15]{Drinfeld:QuasiHopf}).
\begin{theorem}
There exists a symmetric twist $F\in \A^{\bt 2}$ such that 
$\tilde{\Phi} = \Phi_F$, where $\Phi_F$ is given by the formula of Equation 
~\ref{eqn:TwistPhi}.
\end{theorem}

We note that twisting by $F$ leaves the $R$-matrix unchanged, because 
$F$ commutes with $R=e^{\Omega/2}$ (by Proposition ~\ref{prop:OmegaA2Central}),
and therefore
\[
	R_{F} = F^{21} R F^{-1} = F^{21} F^{-1} R = R \, ,
\]
as $F$ is symmetric.
  
By Proposition ~\ref{prop:AnInv}, $\Delta_F = \Delta$. We have thus
proved: 
\begin{proposition}
\lbl{prop:AkzIndepPhi}
If $\Phi, \tilde{\Phi}$ are Drinfel'd associators, then the
quasitriangular quasi-Hopf algebras $\A_{\kz}$ and
$\tilde{\A}_{\kz}$ are twist equivalent.
\end{proposition}

\smallskip
The $\A_{\kz}$ structure can be exported to $\Aar$ with the help
of the $\iota$-maps. We recall that $\iota$-images in $\Aar$-spaces 
of $\A$-elements are denoted by the same symbols as their preimages.
\begin{proposition}
\lbl{prop:AarKZ}
$\Aar_{\kz} = (\Aar, \Delta, \eps, \Phi, R = e^{\Omega/2},                    
       S, \alpha=\nu, \beta=1)$
is a ribbon quasitriangular quasi-Hopf algebra. Its $u$-element and ribbon
element are equal to $e^{-C/2}$.
\end{proposition}
\begin{proof}
As in the proof of Proposition ~\ref{prop:Akz}, we need only verify axioms
pertaining to a general element $D\in \Aar$, namely the 
quasi-coassociativity \eqref{eqn:QuasiCoAssoc}, quasi-cocommutativity
\eqref{eqn:QuasiCoComm}, and the quasi-antipodal
identities \eqref{eqn:AntipodeId1} -- \eqref{eqn:AntipodeId2}. The first two 
follow from Proposition \ref{prop:AarPwrInvar},
and the other two follow immediately from the fact that elements from $\A$
are (mapped by $\iota$ to) invariant elements in $\Aar$.
\end{proof}

\subsubsection{\bf KZ structures on universal enveloping algebras.}
\label{__sect:UghKZ}

We want to use the maps $\T_{\g}$ to transfer the KZ structure of $\A$ into
$U(\g)$, the universal enveloping algebra of a metrized Lie algebra 
$\g$. Here we face a problem: The maps $\T_{\g}$, initially defined for single
diagrams, can only be extended to {\em finite} linear combinations of such 
diagrams, whereas $R = e^{\Omega/2}$, for example, is a linear combination of
infinitely many diagrams.

The remedy to this problem is replacing $U(\g)$ by $\Ugh$, the space
of formal power series in $\hbar$ with coefficients from $U(\g)$, which is
a (left) module over the ring $\kh$. More generally, for any skeleton $\M$, 
we replace the space $L(\M)$ with $L(\M)\hb$. The action of $\g$
on $L(\M)$ extends to an action on $L(\M)\hb$, and
$(L(\M)\hb)^{\g} = L(\M)^{\g}\hb$. For a homogeneous element $D\in \A(\M)$
of degree $d$, we define
\[
	\Tgh(D) = \T_{\g}(D)\hbar^d \, ,
\]
and extend by linearity to (the graded completion of) $\A(\M)$, obtaining 
a map $\Tgh:\A(\M) \ra  L(\M)^{\g}\hb$. In the same 
way we define $\Tartgh: \Aar(\M) \ra L(\M)\hb$, when $\tg$ is the ambient 
Lie bialgebra of a Manin triple, and $\Targh: \Aar(\M) \ra L(\M)\hb$ 
for a complex semisimple Lie algebra $\g$.

\smallskip
A $\kh$-module $M$ is endowed with the $\hbar$-adic topology, which is 
defined by declaring $\{\hbar^n M\}_n$ as a basis of open neighborhoods of $0$.
For $x, y \in M$, we write 
\[
	x\equiv y \mod \hbar^n\quad\text{ if } x-y\in \hbar^n M \, . 
\]

$M$ is called a {\em topologically free} module if there exist a $\k$-linear
space $V$ such that $M\cong V\hb$ as $\kh$-modules.

The {\em topological tensor product} of two topologically free $\kh$-modules of
the form $M=V\hb$, $N=W\hb$, is
\[
	M \hat{\ot} N = (V\ot W)\hb \, .
\]
More intrinsically, it can be defined as the 
{\em $\hbar$-adic completion} of the algebraic tensor product $M\ot N$ (see
~\cite[XVI.2--3]{Kassel:QuantumGroups} for details). The latter definition
is applicable to any pair of $\kh$-modules.
  
With the notion of tensor product at our disposal, we can carry over 
the definitions of Section ~\ref{_sect:qtrqH} to the category of 
$\kh$-modules; all the results mentioned there hold in this category as well, 
with identical proofs. 
We refer to objects in this category as `{\em topological}', partly to 
follow common terminology, and partly as a reminder that topological tensor
products are used throughout.

\smallskip
In particular, $\Ugh$ is a topological Hopf algebra, with its coproduct, for
example, given by
\[
	\Delta\left(\sum_{n=0}^{\infty} u_n \hbar^n \right) 
	= \Delta(u_n) \hbar^n \, .
\]
The maps $\Tgh$ and $\Tartgh$ intertwine all the structure maps of
$\Ugh$-spaces with their parallels for $\A$-spaces and $\Aar$-spaces,
respectively, essentially because all the diagrammatic structure maps
are degree preserving.

\smallskip
We now describe the KZ structure on $\Ugh$. Set 
$\Phi_{\g} = \Tgh(\Phi) \in (\Ugh^{\hat{\ot} 3})^{\g}$, 
$R_{\g} = \Tgh(R) = e^{\hbar\Omega_{\g}/2} \in (\Ugh^{\hat{\ot} 2})^{\g}$, 
and
$\nu_{\g} = \Tgh(\nu) \in (\Ugh)^{\g}$. 
\begin{proposition}
\lbl{prop:UghKZ}
$\Ugh_{\kz} = (\Ugh, \Delta, \eps, \Phi_{\g}, R_{\g}, S, \nu_{\g}, 1)$ is
a topological ribbon quasitriangular quasi-Hopf algebra. Its $u$-element
is $u_{\g} = e^{-\hbar C_{\g}/2}$, which is also the ribbon element. 
\end{proposition}
\begin{proof}
Similarly to the proof of Proposition ~\ref{prop:AarKZ}, we only need to 
comment that the properties of quasi-coassociativity, quasi-cocommutativity, 
and the quasi-antipodal identities hold for all $x\in \Ugh$, due to the 
$\g$-invariance of $\Phi_{\g}$, $R_{\g}$, and $\nu_{\g}$.
\end{proof}   

We note that by \eqref{eqn:HorPhi},
\[
	\Phi_{\g} = 1\hat{\ot} 1 \hat{\ot} 1  \mod{\hbar^2}.
\]

The algebras $\Ugh_{\kz}$ are examples of {\em quantum enveloping algebras}:
\begin{definition}
\lbl{def:QUE}
A quasitriangular quasi-Hopf quantum enveloping (QUE) algebra for a Lie 
algebra $\g$ is a topologically free  quasitriangular quasi-Hopf algebra
$H = (H, \bar{\Delta}, \bar{\eps}, \bar{\Phi}, \bar{R})$, 
such that `$H$ is isomorphic to $U(\g)$ modulo $\hbar$', namely
$H\left/ \hbar H \right. \cong (U(\g), \Delta, \eps, 1^{\ot 3}, 1^{\ot 2})$. 
In addition, $\bar{\Phi}$ is required to satisfy
\[
	\frac{1}{3}(\bar{\Phi}^{123} + \bar{\Phi}^{312} + \bar{\Phi}^{231})
	\equiv 1^{\hat{\ot} 3} \mod{\hbar^2} \, . 
\]
\end{definition}
More explicitly, the topological freeness condition implies that there exists 
an invertible $\k\hb$-linear map $f:H \ra V\hb$, for a certain $k$-linear space
$V$; necessarily then,
$V\cong U(\g)$. Hence there is an invertible 
$\k\hb$-linear map $f:H \ra \Ugh$. That $f$ is a `morphism modulo $\hbar$'
means, for example, that 
\[
	f(xy) \equiv f(x) f(y) \mod{\hbar} \, ,
\]
for all $x, y \in H$, and
\[
	(f\hat{\ot} f) (\tilde{R}) \equiv 1\hat{\ot} 1 \mod{\hbar} \, .	
\]

\begin{remark}
We have omitted the antipodal triple from the displayed components of the 
structure of $H$ because the relevant conditions, such as 
\[
   f(\tilde{S}(x)) \equiv S(f(x)) \mod{\hbar} \, ,	
\]
are automatically satisfied; this can be deduced from 
~\cite[Theorem 1.6]{Drinfeld:QuasiHopf}.
\end{remark}

\subsection{Drinfel'd-Jimbo quantum groups}
\label{_sect:DJ}

Let $\g$ be a complex (semi-)simple Lie algebra. Our aim is to describe an
important construction of a {\em coassociative} QUE algebra for $\g$.

\smallskip
We begin by stating a technical definition.  
\begin{definition}
Let $X$ be a set and $R$ be a subset of $\k\{X\}\hb$, where $\k\{X\}$ is the
$\k$-linear space formally spanned by $X$.
A topological algebra $H$ is said to be {\em (topologically) generated} by the
set $X$ of generators and the set $R$ of relations if $H$ is isomorphic to 
the quotient of $\k\{X\}\hb$ by the $\hbar$-adic closure of the two-sided
ideal generated by $R$.  
\end{definition}

We recall the definitions of several $q$-symbols. For 
an invertible $q\in \mathbb{C}\hb$, the $q$-integers are
\[
	[n]_q = \frac{q^n - q^{-n}}{q - q^{-1}} 
	= q^{n-1} + q^{n-3} + \cdots + q^{-n+3} + q^{-n+1} \, .
\]
For a non-negative integer $k$, the $q$-factorials and $q$-binomial 
coefficients are defined as
\[
	[k]_q ! = \prod_{i=1}^k [i]_q \, , \qquad 
	\bbinom{n}{k}_q = \frac{[n]_q !}{[n-k]_q ! [k]_q !} \, ,
\]
where these expressions are considered as elements of $\mathbb{Z}[q, q^{-1}]$;
that $\bbinom{n}{k}_q \in \mathbb{Z}[q, q^{-1}]$ follows from  
the recursion formula
\[
	\bbinom{n+1}{k}_q 
	= q^{-k} \bbinom{n}{k}_q + q^{n-k+1} \bbinom{n}{k-1}_q \, .
\]

\smallskip
We use the notation and results that were reviewed in Section
~\ref{__sect:CpxSimpleLA}. We further recall that the entries of the Cartan
matrix of $\g$ are 
\[
	a_{ij} = 2 \frac{(\alpha_i, \alpha_j)}{(\alpha_i, \alpha_i)}\, , \quad 
		i, j = 1,\ldots, n\, ,
\]
and set
\[
	d_i = \frac{(\alpha_i, \alpha_i)}{2}\, , \quad i = 1,\ldots, n\, .
\]
The numbers $a_{ij}$ are non-positive integers for $i\neq j$, 
while  $d_i \in \{1, 2, 3\}$ (due to our choice of normalization of
the metric).

For a simple root $\alpha_i\in\Lambda$, let
\[
	e_i = e_{\alpha_i}\, ,  \quad f_i = d_i^{-1} f_{\alpha_i} \, \quad 
	h_i = [e_i, f_i] \, (= d_i^{-1} h_{\alpha_i}) \, .
\]
These elements satisfy {\em Serre's relations}:
\[
	[h_i, h_j] = 0\, , \quad 
	[h_i, e_j] = a_{ij} e_j \, , \quad
	[h_i, f_j] = - a_{ij} f_j \, ,  \quad
	[e_i, f_i] = h_i \, ,
\]
and for $i\neq j$, 
\[
	(\ad e_i)^{1 - a_{ij}} (e_j) = 0 \, ,  \quad 
	(\ad f_i)^{1 - a_{ij}} (f_j) = 0 \, .
\]
By a theorem of Serre (see ~\cite[Section 18]{Humphreys:LieAlgs}), 
$U(\g)$ is isomorphic to the algebra 
generated by a set of $3n$ generators $\{e_i, f_i, h_i\}_{1\leq i\leq n}$ 
subject to the above relations.

In terms of the generators $\{e_i, f_i, h_i\}_{1\leq i\leq n}$ , the standard 
cobracket of $\g$ is determined by
\[
\delta(h_i) = 0\, , \quad \delta(e_i) = \frac{d_i}{2} e_i \land h_i \, , \quad
\delta(f_i) = \frac{d_i}{2} f_i \land h_i \, .
\]

\smallskip
We are now ready to define the Drinfel'd-Jimbo algebra $\Uhg$. Put
$q = e^{\hbar / 2}$. 
\begin{definition}
$\Uhg$ is the topological algebra generated by the set of generators \\
$\{E_i, F_i, H_i\}_{1\leq i\leq n}$ subject to the relations
\[
	[H_i, H_j] = 0\, , \quad 
	[H_i, E_j] = a_{ij} E_j \, , \quad
	[H_i, F_j] = - a_{ij} F_j \, , 
\]
\[
	[E_i, F_i] = \delta_{ij}\frac{q^{d_i H_i} - q^{-d_i H_i}}
                         {q^{d_i} - q^{-d_i}} \, ,
\]
and for $i\neq j$, 
\[
	\sum_{k=0}^{1-a_{ij}} (-1)^k \bbinom{1-a_{ij}}{k}_{q^{d_i}} 
			E_i^k E_j E_i^{1-a_{ij}-k} = 0 \, ,
\]
\[
	\sum_{k=0}^{1-a_{ij}} (-1)^k \bbinom{1-a_{ij}}{k}_{q^{d_i}} 
			F_i^k F_j F_i^{1-a_{ij}-k} = 0 \, .
\]
\end{definition}

It is not hard to see that as $\hbar \ra 0$, the above relations
become Serre's relations for $U(\g)$. Thus, 
$\Uhg\left/ \hbar\Uhg \right. \cong U(\g)$. 

\smallskip
We endow $\Uhg$ with putative coproduct $\Delta_{\hbar}$,
counit $\eps_{\hbar}$, and antipode $S_{\hbar}$, which are determined by
\begin{align*}
	\Delta_{\hbar}(H_i) & = H_i \hat{\ot} 1 + 1 \hat{\ot} H_i \, , \\
	\Delta_{\hbar}(E_i) & = E_i \hat{\ot} q^{d_i H_i /2} 
                         + q^{-d_i H_i /2} \hat{\ot} E_i \, , \\
	\Delta_{\hbar}(F_i) & = F_i \hat{\ot} q^{d_i H_i /2} 
                         + q^{-d_i H_i /2} \hat{\ot} F_i \, ;	
\end{align*}
\[
	\eps_{\hbar}(H_i) = \eps_{\hbar}(E_i) = \eps_{\hbar}(F_i) = 0 \, ;
\]
and
\[
S_{\hbar}(H_i) = -H_i\, , \quad 
S_{\hbar}(E_i) = - q^{d_i} E_i\, , \quad 
S_{\hbar}(F_i) = - q^{-d_i} F_i\, .
\]

The following theorem is proved in ~\cite{Drinfeld:QuantumGroups} (see
also ~\cite[Section 6.5]{ChariPressley:Guide}).
\begin{theorem}
There is a universal $R$-matrix\footnote{The Lie algebra $\g$ is 
suppressed from the notation for elements of $\Uhg$-spaces.} 
$R_{\hbar} \in \Uhg^{\hat{\ot} 2}$ such that  
\[
\Uhg = (\Uhg, \Delta_{\hbar}, \eps_{\hbar}, R_{\hbar}, S_{\hbar}) 
\]
is a quasitriangular Hopf QUE algebra for $\g$.
\end{theorem}
Further properties and explicit expressions for $R_{\hbar}$ can be found in 
~\cite[Section 8.3]{ChariPressley:Guide}.

\smallskip
The theorem implicitly asserts that $\Uhg$ is topologically free, so that
$\Uhg$ and $\Ugh$ are isomorphic as 
$\mathbb{C}\hb$-modules (it is a difficult part of the proof). In fact, much 
more is proved in ~\cite{Drinfeld:AlmostCocomm} and 
~\cite[Proposition 6.5.7]{ChariPressley:Guide}. 
\begin{theorem}
\lbl{thm:Uhg-Iso-Ugh}
There is an isomorphism of topological {\em algebras} $\varphi:\Uhg \ra \Ugh$ 
satisfying 
\[
\varphi(E_i) = e_i \mod{\hbar}\, , \ 
\varphi(F_i) = f_i \mod{\hbar}\, , \
\varphi(H_i) = h_i \mod{\hbar}\, , \ \  i=1,\ldots, n \, .
\]
The restriction of $\varphi$ to $\Zhg$, the center of $\Uhg$, is 
uniquely determined by these conditions and induces a canonical isomorphism 
of topological algebras $\varphi: \Zhg \ra \CUgh$. 
Moreover, $\varphi$ can be chosen in such a way that 
$\varphi(H_i) = h_i$ for $i=1,\ldots, n$ (strict equality).
\end{theorem}

In particular, we can define the {\em Quantum Casimir element}
$C_{\hbar} \in \Zhg$ as the inverse image of the (classical) 
Casimir element $C_{\g} \in \CUg \subset \CUgh$. Write $\varrho_{\g}$, the 
(element of $\h$ dual to the) half-sum of positive roots, as
$\varrho_{\g} = \sum_{i=1}^n c_i h_i$, and set 
\[
	\varrho_{\hbar} = \sum_{i-1}^n c_i H_i \in \Uhg \, .
\]
(It can be checked that the coefficients are given by 
$c_i = \sum_{j=1}^n a^{ji} d_j$, where $(a^{ij})$ is the inverse of the Cartan
matrix.)
\begin{proposition}
Let $v_{\hbar} = e^{-\hbar C_{\hbar}/2}$. Then the $u$-element of $\Uhg$ is
$u_{\hbar} = e^{\hbar \varrho_{\hbar}} v_{\hbar}$, and $v_{\hbar}$ is a ribbon 
element for $\Uhg$. 
\end{proposition}
See ~\cite[Section XVII.3]{Kassel:QuantumGroups} for a proof, based on the 
results of ~\cite[Section 5]{Drinfeld:AlmostCocomm}.

\smallskip
Using an isomorphism $\varphi$ as in Theorem ~\ref{thm:Uhg-Iso-Ugh}, we obtain
a ribbon quasitriangular Hopf structure on $\Ugh$ by using $\varphi$ to export
the quasitriangular Hopf structure from $\Uhg$ to $\Ugh$. We call such a
structure a {\em DJ structure} of $\Ugh$ and denote it by $\Ugh_{\qg}$. We
write
\[
\Ugh_{\qg} = 
  (\Ugh, \Delta_{\hbar}, \eps_{\hbar}, R_{\g, \hbar}, S_{\hbar}) \, ,
\]
where the coproduct, the counit, and the antipode are denoted by the same 
symbols as their $\varphi$-related maps for $\Uhg$, and the $R$-matrix is
\[
	R_{\g, \hbar} = \varphi \hat{\ot} \varphi (R_{\hbar}) \, .
\]
It can be shown (see ~\cite[Lemma XIX.5.1]{Kassel:QuantumGroups})
that $\eps_\hbar = \Ugh \ra \mathbb{C}\hb$ is equal to the usual 
counit of $\Ugh$.

\smallskip
It is immediate from the properties of $\varphi$ stated in Theorem 
~\ref{thm:Uhg-Iso-Ugh}, that the $u$-element and ribbon 
element of $\Ugh_{\qg}$ are
\begin{equation}
\label{eqn:uDJ}
	u_{\g,\hbar} = e^{\hbar \varrho_{\g}} e^{-\hbar C_{\g}/2} \, , \quad
	v_{\g,\hbar} = e^{-\hbar C_{\g}/2} \, .
\end{equation}

\smallskip
As $\Ugh_{\qg}$ has a trivial associator, Equation ~\ref{eqn:QQYBE} for 
$R_{\g,\hbar}$ becomes the {\em quantum Yang-Baxter equation (QYBE)},
\[
	R_{\g,\hbar}^{12} R_{\g,\hbar}^{13} R_{\g,\hbar}^{23} = 
	R_{\g,\hbar}^{23} R_{\g,\hbar}^{13} R_{\g,\hbar}^{12} 
\, ,
\] 
which is a `quantum form' of the CYBE \eqref{eqn:CYBE}.  
It can be proved (see the remark following 
~\cite[Theorem 8.3.9]{ChariPressley:Guide}) that
\begin{equation}
\label{eqn:UhgRmatrix}
	R_{\g,\hbar} = 1\hat{\ot} 1 + \hbar r_{\g} \ \mod{\hbar^2}\, , 
\end{equation}
where $r_{\g}$ is the classical $r$-matrix of $\g$, given by \eqref{eqn:r_g}.
Thus, $R_{\g,\hbar}$ is a {\em quantization} of $r_{\g}$.

\smallskip
The next theorem is an important ingredient in Drinfel'd's celebrated proof 
of the Drinfel'd-Kohno theorem (see ~\cite[XIX.4--5]{Kassel:QuantumGroups}
for a streamlined exposition). 
\begin{theorem}
\lbl{thm:QGtwistKZ}
$\Ugh_{\qg}$ is twist equivalent to  $\Ugh_{\kz}$.
\end{theorem}
We stress that this result is independent of the particular
isomorphism $\varphi$ which is used to define a DJ structure on $\Ugh$.

\subsection{$H$-invariants of framed links}
\label{_sect:HinvOfLinks}

Coupled with some important examples of algebraic and diagrammatic 
ribbon quasitriangular quasi-Hopf algebras, we are going to describe 
how to construct from such algebras topological invariants of 
{\em framed parenthesized tangles} (restricting, in particular, to invariants
of framed links). Similar constructions appear in many texts, including
~\cite{AltschulerCoste:QuasiQuantum}, ~\cite{Bar-Natan:NAT},
~\cite{Cartier:Construction}, ~\cite[Chapter XX]{Kassel:QuantumGroups},
~\cite{LinksGouldZhang:TwistInvarLinkPoly}, 
~\cite{LeMurakami:Universal}, ~\cite{Piunikhin:Weights}, and
~\cite{ReshetikhinTuraev:Ribbon}. Unfortunately, none of those
treatments is sufficient for our purposes. This is because constructions in the
usual algebraic framework (most notably ~\cite{AltschulerCoste:QuasiQuantum}
and ~\cite{LinksGouldZhang:TwistInvarLinkPoly}) work with finite-dimensional
modules, which are unavailable in diagrammatic cases, while diagrammatic
constructions (such as  ~\cite{Bar-Natan:NAT} and ~\cite{LeMurakami:Universal})
utilize specialized properties of $\A$ and its KZ structure, which we can not
use when working with the EK structure of $\Aar$, to be described in Section 
~\ref{sect:EK}.

\subsubsection{\bf The space of coinvariants.}
\label{__sect:CoInvSpace}

\begin{definition}
\lbl{def:CoInv}
The {\em space of coinvariants} of an algebra $H$ is the quotient space 
\[
H^{\coinv} = H / \Span\{xy - yx \, \mid \, x, y \in H\} \, .  
\]
We often use the alternative notation $\bar{H}$. We denote by 
$\tr:H \ra \bar{H}$ the natural projection.
\end{definition}

For $H=U(\g)$, with $\g$ a metrized Lie algebra, an elementary argument 
shows that $\bar{H} = U(\g)_{\g}$, the
space of coinvariants for the adjoint action of $\g$ on $U(\g)$ 
(see Section ~\ref{__sect:TgDef}). If $V$ is a finite-dimensional 
representation of $\g$, then the associated algebra morphism 
$U(\g) \ra \End(V)$ descends upon composition with the trace map 
$\tr:\End(V) \ra \k$ to a map $\tr_V:U(\g)_{\g} \ra \k$. Therefore, 
$U(\g)_{\g}$ is in a sense a 'trace object'. If $\g$ is a complex semisimple 
Lie algebra, then it can be proved\footnote{The crucial ingredient in the 
proof is a well-known theorem of Chevalley 
(\cite[Section 23.1]{Humphreys:LieAlgs}), asserting
that the trace maps $\{\tr_V\}$ span the space of
linear functionals on $U(\g)^{\g}$.} that the trace maps $\{\tr_V\}$, with 
$V$ ranging over finite-dimensional representations of $\g$, separate
the points of $U(\g)_{\g}$, so that $U(\g)_{\g}$ is a `universal
trace object'.

\smallskip
For the diagrammatic algebra $\Aar$, the subspace corresponding to 
$\Span\{xy - yx\}$ is the image of the map 
\[ 
\gamma: \lowupper{\Aar}{\bt 2}\ra \Aar\, , \quad  D \mapsto 
	\mu(D) - \mu(D^{21}) \, .
\]
Let $\phi^{l}:\Aar \ra \Aar(\bigcirc)$ be the `circling the interval' map, 
defined as in Section ~\ref{__sect:IntrvlToCirc}. It is clear that 
$\img(\gamma) \subset \ker(\phi^l)$, so $\phi^{l}$ descends to a map
$\phi^l:\Aar \left/ \img(\gamma) \right. \ra \Aar(\bigcirc)$. By an argument
similar to the one in the beginning of the proof of Proposition
~\ref{prop:Aline-eq-Acirc}, the latter map has an `opening the circle' 
inverse. We conclude that
\[
	\Aar \left/ \img(\gamma) \right. \cong \Aar(\bigcirc) \, , 
\]
so it is reasonable to define  $\lowupper{\Aar}{\coinv} = \Aar(\bigcirc)$.
Accordingly, we henceforth denote the natural projection from $\Aar$ to 
$\Aar(\bigcirc)$ by $\tr$ instead of $\phi^l$. The spaces of coinvariants and 
trace maps for $\Aar$ and $U(\g)$ ($\g$ a complex semisimple Lie algebra) fit 
into the commutative diagram
\begin{equation}
\label{eqn:AarVsUgCoinv}
\begin{CD}
\Aar @>{\tr}>> \Aar(\bigcirc) \\
@V{\Targ}VV @VV{\Targ}V \\
U(\g) @>{\tr}>> U(\g)_{\g} \,    
\end{CD} \, .
\end{equation}

\smallskip
For a $\Vect_{\k}$-algebra $H$, the space of coinvariants of the tensor
power $H^{\ot n}$ (with its usual $\Vect_k$-algebra structure) is clearly
equal to $\bar{H}^{\ot n}$. To extend \eqref{eqn:AarVsUgCoinv}, we
define
\[
(\lowupper{\Aar}{\bt n})^{\coinv} = \Aar(\bigcirc_1, \ldots, \bigcirc_n) \, ,
\]
and make an analogous definition for $\A$:
\[
(\A^{\bt n})^{\coinv} = \A(\bigcirc_1, \ldots, \bigcirc_n) \, ,
\]
Note that $\A^{\coinv} = \A(\bigcirc) = \A$,
but for $n\geq 2$, $(\A^{\bt n})^{\coinv} \neq \A^{\bt n}$ (see Remark 
~\ref{rem:PwrAcirc}). That these are the correct diagrammatic definitions 
will become clear below when we discuss in Section ~\ref{__sect:Z_H} the 
$\HS$ categories. 
  
\subsubsection{\bf The category of parenthesized framed tangles.}
\label{__sect:CatPFT}

\begin{definition}
A {\em parenthesized string of length $l$} over the alphabet $\{\ua, \da\}$ 
is either the empty string $\ew$ ($l=0$), one of the two single-letter strings
$\ua$ and $\da$ ($l=1$), or a string of the form $(W_1 W_2)$, 
where $W_1, W_2$ are parenthesized strings of lengths $0 < l_1, l_2 < l$ such 
that $l=l_1 + l_2$. The length of a parenthesized string $W$ is denoted by 
$|W|$. 
\end{definition}

We denote by $\W$ the set of parenthesized strings over $\{\ua, \da\}$.
Some elements of $\W$ are 
$\da, (\ua\da)$, $(\ua(\da\ua))$,  $((\ua\da)\ua)$, $((\ua\da)((\ua\ua)\da))$.
Such elements are often represented by a collection of vertical arrows along
a horizontal line, whose distances from each other reflect the 
parenthesization; for example, 
\[
	((\ua\da)((\ua\ua)\da)) \qquad \longLRa \qquad
	\ua\da \ \ \ \ \ua\ua \ \ \da \, .
\]

We endow $\W$ with the concatenation product $\ew \cdot W = W\cdot\ew = W$, 
$W_1 \cdot W_2 = (W_1 W_2)$ if $|W_1|, |W_2| \geq 1$. Henceforth, we drop
the dots  and reserve the expression $W_1 W_2$ for the  product of 
$W_1, W_2 \in \W$, rather than for their concatenation.
The product of $\W$ is not associative; we follow the usual habit of 
indicating the order of products by parentheses, which should not
(and are unlikely to) be confused with the left and right parenthesis letters
appearing in elements of $\W$.    

\smallskip
We similarly define the set $\W^{\star}$ as the set of parenthesized strings 
over the alphabet $\{\ua, \da, \star\}$, in which the letter $\star$ appears 
exactly once.
\begin{definition}
Let $W\in \W^{\star}$ and $A\in \W$. The {\em substitution} of $A$ into 
$W$ is the operation of replacing the $\star$ letter in $W$ by $A$. The 
resulting element in $\W$  is denoted by $W/\{\star \ra A\}$.
\end{definition}
A little care is needed when substituting the empty word $A=\ew$, in order to 
avoid redundant parentheses. The rule is that $\star$ is erased, and 
(if $|W| > 1$) a matching pair of parentheses, one of which was adjacent to
$\star$ in $W$, is also erased. For example,
\[
((\ua\da)(\star((\ua\ua)\da)))\, /\{\star\ra \ew\} = 
((\ua\da)((\ua\ua)\da)) \, .	
\]

\smallskip
We briefly recall 
some definitions and facts about tangles and framed tangles 
(consult, for example, 
~\cite[Section X.5]{Kassel:QuantumGroups} for more details). 
A {\em tangle} $T$ is a smooth oriented embedding  
of a $1$-dimensional skeleton $\M$ into $\mathbb{R}^2 \times [0, 1]$, 
considered up to isotopy relative to the boundary, such that $T$ intersects 
the two boundary planes of $\mathbb{R}^2 \times [0, 1]$ transversally, and the 
intersection points (necessarily the images of ends of oriented interval 
components of $\M$) are along the lines $\mathbb{R} \times \{0\} \times \{0\}$ 
and $\mathbb{R} \times \{0\} \times \{1\}$. The {\em domain} of $T$ is the 
(non-parenthesized) string over $\{\ua, \da\}$ indicating the orientation types
of the intersection points along the line 
$\mathbb{R} \times \{0\} \times \{0\}$, while the {\em target} of $T$ is 
defined in the same way with respect to the line 
$\mathbb{R} \times \{0\} \times \{1\}$.

A tangle can be represented by a {\em tangle diagram}, which is a 
projections of the tangle (more precisely, of a certain embedding in its
isotopy class) on the plane $\mathbb{R} \times \{0\} \times \mathbb{R}$,
whose only singularities are transversal double points. The topological
notion of isotopic tangles can be expressed in terms of a suitable set 
of {\em Reidemeister moves} on tangle diagrams (see 
~\cite[Theorem X.5.9]{Kassel:QuantumGroups}).
 
Figure ~\ref{fig:Tangle}  shows an example of a tangle diagram. The domain 
of this tangle is $\ua\ua\da$ and its target is $\ua\da\ua$.

\begin{figure}[htpb]
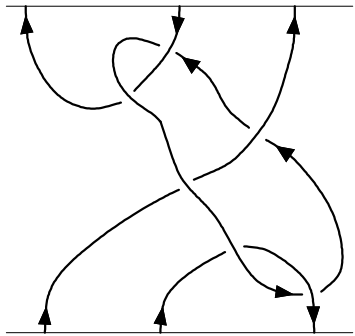

\[ 

  \printname{Tangle}
  \setlength{\unitlength}{0.8\standardunitlength}
  \begin{array}{c}  \hspace{-1.7mm}
    \raisebox{-8pt}{\input draws/Tangle.tex }
      \hspace{-1.9mm}
  \end{array}
 
\]
\caption{Example of a tangle diagram.} \lbl{fig:Tangle}
\end{figure}

We define a {\em normal vector field} on a tangle $T$ as a smooth vector
field on $T$ that is nowhere tangent to $T$, and that is given by the vector
$(0, 1, 0)$ at the points of $T$ along its domain and target lines. A 
{\em framing} of $T$ is a normal vector field on $T$, considered up to 
homotopy (within the class of normal vector fields). 

A {\em framed tangle} is a tangle equipped with a framing, considered up
to framing preserving isotopy relative to the boundary. A framed tangle may
be thought of as a {\em tangled ribbon}, whose rims are the tangle itself and
another tangle obtained by a small translation along a normal vector field
in the framing homotopy class. 

Framed tangles can be represented by tangle diagrams equipped with the
{\em blackboard framing} --- the framing determined by the constant 
vector field $(0, 1, 0)$, which is perpendicular to the diagram plane and
points away from the reader. Figure ~\ref{fig:Ribbon} shows an example
of a vertical framed tangle whose normal vector field makes a counterclockwise
turn by an angle of $2\pi$, and a tangle diagram representing it. It can be 
proved (see ~\cite[Section X.8]{Kassel:QuantumGroups}) that isotopy classes
of framed tangles are the same as equivalence classes of tangle diagrams
with respect to certain Reidemeister moves.

\begin{figure}[htpb]
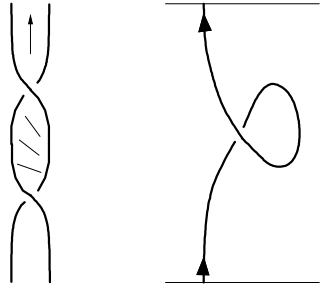

\[ 

  \printname{Ribbon}
  \setlength{\unitlength}{0.8\standardunitlength}
  \begin{array}{c}  \hspace{-1.7mm}
    \raisebox{-8pt}{\input draws/Ribbon.tex }
      \hspace{-1.9mm}
  \end{array}
 
\]
\caption{A framed tangle and its tangle diagram.} \lbl{fig:Ribbon}
\end{figure}

\begin{definition}
The {\em category $\catFT$ of framed tangles} is the category whose objects 
are (non-parenthesized) strings over $\{\ua, \da\}$, and whose
morphisms are framed tangles (a framed tangle $T$ belongs to the set
of morphisms between its domain and target strings). 

The {\em category $\catPFT$ of parenthesized framed tangles} is defined in the 
same way, only with elements of $\W$ as its objects. Note that the 
endomorphisms of the object $\ew$ in both categories are simply 
{\em framed links}. 
\end{definition}

We proceed to provide a description of the category $\catPFT$ in terms of 
`generators and relations', following the treatment in
~\cite[Section 2]{Bar-Natan:NAT}. We begin by defining an auxiliary category
which is `freely generated' by a certain collection of morphisms, and
whose `quotient' by suitable relations is $\catPFT$. We refer to
~\cite[Chapter XII]{Kassel:QuantumGroups} for precise definitions 
of all the categorical terms, and for detailed discussions  
of presentations of several categories of (non-parenthesized) knot-like 
objects by generators and relations.  

\begin{definition} 
\lbl{def:NAT}
Let $\catNAT$ be the category whose set of objects is $\W$, and whose
morphisms are freely generated by the following morphisms:

\par\noindent
\begin{itemize}
\item[(G1)] {\em Associativity morphisms: } A pair of morphisms
\[
\begin{aligned}
 	W^{\NATra}_{A,B,C}\in\mor(W/\{\star\ra (AB)C\}\, ,\, 
                              W/\{\star\ra A(BC)\})\, ;
	\\
	W^{\NATla}_{A,B,C}\in\mor(W/\{\star\ra A(BC)\}\, ,\, 
                              W/\{\star\ra (AB)C)\}) \, ,
\end{aligned}
\]
for each $W\in \W^{\star}$, and each triple $(A,B,C)$ of $\W$-strings.
Graphically, we represent these morphisms as follows:
\[ 
	\star^{\NATra}_{(\ua\da),(\da\ua),\ua}
	\quad \longLRa \quad
	
  \printname{natRightAssoc}
  \setlength{\unitlength}{0.75\standardunitlength}
  \begin{array}{c}  \hspace{-1.7mm}
    \raisebox{-8pt}{\input draws/natRightAssoc.tex }
      \hspace{-1.9mm}
  \end{array}
; \qquad
     \star^{\NATla}_{(\ua\da),(\ua\da),\ua}
	\quad \longLRa \quad
	
  \printname{natLeftAssoc}
  \setlength{\unitlength}{0.75\standardunitlength}
  \begin{array}{c}  \hspace{-1.7mm}
    \raisebox{-8pt}{\input draws/natLeftAssoc.tex }
      \hspace{-1.9mm}
  \end{array}
 \, .
\]

\item[(G2)] {\em Braiding morphisms: } A pair of morphisms 
\[
	W^{\NATov}_{A,B}, W^{\NATun}_{A,B} \in 
	\mor(W/\{\star\ra AB\}\, ,\, W/\{\star\ra BA\}) \, ,
\]
for each $W\in\W^{\star}$, and each pair $A$, $B$ of $\W$-strings.
The graphical representations of these morphisms are:
\[ 
((\ua\star)\da)^{\NATov}_{\da,(\ua\ua)} \qquad \longLRa \qquad
        
  \printname{natOver}
  \setlength{\unitlength}{0.5\standardunitlength}
  \begin{array}{c}  \hspace{-1.7mm}
    \raisebox{-8pt}{\input draws/natOver.tex }
      \hspace{-1.9mm}
  \end{array}
 \, ;
\]
\[ ((\da\star)\ua)^{\NATun}_{\da,(\ua\ua)} \qquad \longLRa \qquad
        
  \printname{natUnder}
  \setlength{\unitlength}{0.5\standardunitlength}
  \begin{array}{c}  \hspace{-1.7mm}
    \raisebox{-8pt}{\input draws/natUnder.tex }
      \hspace{-1.9mm}
  \end{array}
 \, .
\]

\item[(G3)] {\em Pair creations and annihilations: } A quadruple of morphisms
\[
\begin{aligned}
 W^{\NATcp} \in \mor(W/\{\star\ra \ew\}\, ,\, W/\{\star\ra (\da\ua) \})\, ; \\
 W^{\NATcn} \in \mor(W/\{\star\ra \ew\}\, ,\, W/\{\star\ra (\ua\da) \})\, ; \\
 W^{\NATap} \in \mor(W/\{\star\ra (\da\ua)\}\, ,\, W/\{\star\ra \ew \})\, ; \\
 W^{\NATan} \in \mor(W/\{\star\ra (\ua\da)\}\, ,\, W/\{\star\ra \ew \})\, ,
\end{aligned}
\]
for each $W\in \W^{\star}$.
Graphically, we will represent a morphism $W^{\NATcp}$, for example, as 
follows:
\[ 
	((\ua\da)\star)^{\NATcp} \qquad \longLRa \qquad
	
  \printname{natCreatePos}
  \setlength{\unitlength}{0.5\standardunitlength}
  \begin{array}{c}  \hspace{-1.7mm}
    \raisebox{-8pt}{\input draws/natCreatePos.tex }
      \hspace{-1.9mm}
  \end{array}
 \, .
\]

\end{itemize} 
\end{definition}
\par\noindent
\parbox{3.4in}{
The graphical representation of the composition $W_1\cdot W_2$ of any
two composable morphisms $W_1$ and $W_2$ is obtained by stacking the graphical
representation of $W_2$ on top of that of $W_1$, as shown on the right.}
$\hspace{-1.5in}$
\parbox{3.5in}{
\[ {\def\arg{W_1}
  \printname{natBoxMorphism}
  \setlength{\unitlength}{0.5\standardunitlength}
  \begin{array}{c}  \hspace{-1.7mm}
    \raisebox{-8pt}{\input draws/natBoxMorphism.tex }
      \hspace{-1.9mm}
  \end{array}
}\cdot\,\,
   {\def\arg{W_2}
  \printname{natBoxMorphism}
  \setlength{\unitlength}{0.5\standardunitlength}
  \begin{array}{c}  \hspace{-1.7mm}
    \raisebox{-8pt}{\input draws/natBoxMorphism.tex }
      \hspace{-1.9mm}
  \end{array}
} =
   \begin{array}{c}
	{\def\arg{W_2}
  \printname{natBoxMorphism}
  \setlength{\unitlength}{0.5\standardunitlength}
  \begin{array}{c}  \hspace{-1.7mm}
    \raisebox{-8pt}{\input draws/natBoxMorphism.tex }
      \hspace{-1.9mm}
  \end{array}
} \\
	\vspace{-9mm} \\
	{\def\arg{W_1}
  \printname{natBoxMorphism}
  \setlength{\unitlength}{0.5\standardunitlength}
  \begin{array}{c}  \hspace{-1.7mm}
    \raisebox{-8pt}{\input draws/natBoxMorphism.tex }
      \hspace{-1.9mm}
  \end{array}
}
   \end{array}
\]
}

\smallskip
There is an obvious forgetful functor
$\catFor:\catNAT\ra\catPFT$ mapping each object $W\in \W$ of $\catNAT$ to 
itself, and mapping a morphism $T$ of $\catNAT$ to the framed tangle whose
tangle diagram is the graphical representation of $T$, as dictated by 
Definition ~\ref{def:NAT}. We list several groups 
of relations belonging to the `kernel' of $\catFor$. 

\begin{myitemize}

\item {\em Associativity relations: }
\begin{enumerate}
\item[(R1)]
$W^{\NATra}_{A,B,C}$ is equal to an identity morphism if any of $A$, $B$, or
$C$ is the empty string. Similarly for $W^{\NATla}_{A,B,C}$.
\item[(R2)]
$W^{\NATla}_{A,B,C}$ is a two-sided inverse of $W^{\NATra}_{A,B,C}$:
\[ 
	
  \printname{natRala}
  \setlength{\unitlength}{0.5\standardunitlength}
  \begin{array}{c}  \hspace{-1.7mm}
    \raisebox{-8pt}{\input draws/natRala.tex }
      \hspace{-1.9mm}
  \end{array}
 \quad ; \qquad 
  \printname{natLara}
  \setlength{\unitlength}{0.5\standardunitlength}
  \begin{array}{c}  \hspace{-1.7mm}
    \raisebox{-8pt}{\input draws/natLara.tex }
      \hspace{-1.9mm}
  \end{array}
 
\]
(Strands in the above figure are actually bundles of strands, determined by
$A$, $B$, or $C$).
Notice that only the `active' part of each relation is drawn, omitting the
`far away' strands of $W^{\NATra}$ and $W^{\NATla}$, which connect $\ua$ and 
$\da$ arrows in $W$ that are to the left or to the right of the $ABC$ part.

\par\noindent
\parbox{1.2in}{
\item[(R3)]
{\em The pentagon: }
}
$\hspace{-1.0in}$
\parbox{4in}{
\[ 
	
  \printname{natPentagon}
  \setlength{\unitlength}{0.5\standardunitlength}
  \begin{array}{c}  \hspace{-1.7mm}
    \raisebox{-8pt}{\input draws/natPentagon.tex }
      \hspace{-1.9mm}
  \end{array}
 
\]
}
\end{enumerate}

\item {\em Locality relations: }
\begin{enumerate}
\parbox{2.5in}{
\item[(R4)]
{\em Locality in space:} (for any pair $W_A$ and $W_B$ of morphisms)
}
$\hspace{-1.9in}$
\parbox{3.8in}{
\[ 
	
  \printname{natLocSpace}
  \setlength{\unitlength}{0.575\standardunitlength}
  \begin{array}{c}  \hspace{-1.7mm}
    \raisebox{-8pt}{\input draws/natLocSpace.tex }
      \hspace{-1.9mm}
  \end{array}
 
\]
}
\item[(R5)]
{\em Locality in scale: } for any morphisms $W_A$, $W_B$, and $W_C$,
\[	
	\begin{array}{c}
		
  \printname{natLocRight}
  \setlength{\unitlength}{0.575\standardunitlength}
  \begin{array}{c}  \hspace{-1.7mm}
    \raisebox{-8pt}{\input draws/natLocRight.tex }
      \hspace{-1.9mm}
  \end{array}
, \\
		\rule{0in}{0.5cm} \\
		
  \printname{natLocLeft}
  \setlength{\unitlength}{0.575\standardunitlength}
  \begin{array}{c}  \hspace{-1.7mm}
    \raisebox{-8pt}{\input draws/natLocLeft.tex }
      \hspace{-1.9mm}
  \end{array}
,
\end{array}
\]
\[ 
	
  \printname{natLocOver}
  \setlength{\unitlength}{0.575\standardunitlength}
  \begin{array}{c}  \hspace{-1.7mm}
    \raisebox{-8pt}{\input draws/natLocOver.tex }
      \hspace{-1.9mm}
  \end{array}
 \quad\text{ and }\qquad 
	
  \printname{natLocUnder}
  \setlength{\unitlength}{0.575\standardunitlength}
  \begin{array}{c}  \hspace{-1.7mm}
    \raisebox{-8pt}{\input draws/natLocUnder.tex }
      \hspace{-1.9mm}
  \end{array}
. 
\]
\end{enumerate}

\item {\em Braiding relations: }
\begin{enumerate}
\item[(R6)]
$W^{\NATov}_{A,B}$ is equal to an identity morphism if $A$ or $B$ is the
empty string. Similarly for $W^{\NATun}_{A,B}$.
\item[(R7)]
$W^{\NATun}_{B,A}$ is a two-sided inverse of $W^{\NATov}_{A,B}$:
$\quad 
  \printname{natOvun}
  \setlength{\unitlength}{0.5\standardunitlength}
  \begin{array}{c}  \hspace{-1.7mm}
    \raisebox{-8pt}{\input draws/natOvun.tex }
      \hspace{-1.9mm}
  \end{array}
 \ \ ; \ \ 
  \printname{natUnov}
  \setlength{\unitlength}{0.5\standardunitlength}
  \begin{array}{c}  \hspace{-1.7mm}
    \raisebox{-8pt}{\input draws/natUnov.tex }
      \hspace{-1.9mm}
  \end{array}
 $
\item[(R8)]
{\em The hexagons: } 
\[ 
	
  \printname{natPHex}
  \setlength{\unitlength}{0.5\standardunitlength}
  \begin{array}{c}  \hspace{-1.7mm}
    \raisebox{-8pt}{\input draws/natPHex.tex }
      \hspace{-1.9mm}
  \end{array}
 \quad;\quad 
  \printname{natNHex}
  \setlength{\unitlength}{0.5\standardunitlength}
  \begin{array}{c}  \hspace{-1.7mm}
    \raisebox{-8pt}{\input draws/natNHex.tex }
      \hspace{-1.9mm}
  \end{array}
 
\]
\end{enumerate}

\item {\em Creation/annihilation relations: }
\begin{enumerate}
\parbox{2.7in}{
\item[(R9)]
$ \quad
  \printname{natCA}
  \setlength{\unitlength}{0.5\standardunitlength}
  \begin{array}{c}  \hspace{-1.7mm}
    \raisebox{-8pt}{\input draws/natCA.tex }
      \hspace{-1.9mm}
  \end{array}
, $
}
\parbox{2.8in}{
\item[(R10)]
$\quad
  \printname{natDoubleLoop}
  \setlength{\unitlength}{0.5\standardunitlength}
  \begin{array}{c}  \hspace{-1.7mm}
    \raisebox{-8pt}{\input draws/natDoubleLoop.tex }
      \hspace{-1.9mm}
  \end{array}
\hspace{-12pt}, $ 
}
\\
where $a\in\{\ua,\da\}$ is an arrow and $\bar{a}$ is its opposite.

\end{enumerate}
\end{myitemize}

\smallskip
By Theorem 1 in ~\cite{Bar-Natan:NAT} (see also 
~\cite{AltschulerCoste:QuasiQuantum}), the kernel of 
$\catFor:\catNAT\ra\catPFT$ is generated by the above relations. 
Thus, $\catPFT$ is generated by the morphisms (G1)--(G3) subject to the 
relations (R1)--(R10). 

\subsubsection{\bf Construction of the $H$-invariant.}
\label{__sect:Z_H}

The description of $\catPFT$ in terms of generators and relations provides
the following method for constructing topological invariants of parenthesized 
framed tangles:
\begin{itemize}
\item 
Construct a functor $\Z:\catNAT\ra\mathbf{C}$, where $\mathbf{C}$ is some 
arbitrary `target category'. This can be done by simply stating the values 
of $\Z$ on the generators of $\catNAT$.
\item Check that $\Z$ respects relations (R1)--(R10).
\item If it does, it descends to a functor (denoted by the same symbol)
$\Z:\catPFT\ra \mathbf{C}$, which is an invariant of parenthesized framed 
tangles.
\end{itemize}


We carry out this procedure in the case where $\mathbf{C}=\HS$ is the category 
of $H$-decorated skeletons\footnote{All skeletons considered in this
section are $1$-dimensional.}, which we now turn to define. 

\smallskip
Let $\M$ be a skeleton. For $\m$, a component of $\M$, we put $H(\m) = H$ if 
$\m = \ua$ (an oriented interval), and $H(\m) = \bar{H}$ if $\m=\bigcirc$
(an oriented circle). If $H$ is a $\Vect_{\k}$-algebra, then we define the 
$H$-space associated to $\M$ as 
(compare Definition ~\ref{def:TM} and Section ~\ref{__sect:TgDef})
\[
	H(\M) = \bigotimes_{\m \in \M} H(\m) \, . 
\]
(Recall from Definition ~\ref{def:SkltWithColors} that the set of 
components of $\M$ is ordered, hence the tensor product is well-defined.)
For the diagrammatic algebras, $\A$ and $\Aar$, the spaces $H(\M)$ are taken to
be the familiar spaces $\A(\M)$ and $\Aar(\M)$, respectively.  
\begin{definition}
Let $\HS$ be the category whose objects are (non-parenthesized) strings of
$\ua$ and $\da$ arrows (same as the objects of $\catFT$), and whose 
morphisms are pairs $(\M, X)$, where $\M$ is a skeleton and $X\in H(\M)$.
The tails and heads of the oriented intervals in $\M$ are in a bijective
direction preserving correspondence with the $\ua$ and $\da$ arrows in the
domain and target objects of the morphism $(\M, X)$. We call $\M$ the
{\em morphism skeleton}, and say that $X$ {\em decorates} $\M$, or that
$X$ is the {\em decorating tensor} of the morphism. When $H$ is a 
$\Vect_{\k}$-algebra, we graphically represent a morphism $(\M, X)$ (with 
$X$ a factorable tensor) by drawing $\M$ with each component decorated by a 
pair consisting of an order tag and the corresponding factor of $X$. (For 
diagrammatic algebras, the graphical representations of morphisms are simply 
the usual pictures of elements in $\A(\M)$ or $\Aar(\M)$.)

\smallskip
The rule for composing two composable morphisms $(\M_1, X_1)$ and $(\M_2, X_2)$
is illustrated by stacking the graphical representation of $(\M_2, X_2)$ on 
top of that of $(\M_1, X_1)$: when several oriented interval components are 
glued, the elements of $H$ decorating them are multiplied left to right from 
tail to head; if an oriented circle is formed by glueing two oriented 
intervals at both ends, then the two elements decorating it are multiplied in 
any order and the result is mapped  by $\tr$ to $\bar{H}$. (These 
multiplication rules refer to a $\Vect_{\k}$-algebra; for $\A$ and $\Aar$, 
stacking the graphical representations one on top of the other automatically 
provides the skeleton of the composed morphism with its decoration.) 
The order of the components of the skeleton of the composed morphism is 
determined by first declaring the components of $\M_1$ as preceding those 
of $\M_2$, and then marking each component of the composition skeleton
with the smallest order tag it carries; see Figure ~\ref{fig:HSexample} 
for an example.
\end{definition}     

\begin{figure}[htpb]
\[ 
	
  \printname{HSexample1}
  \setlength{\unitlength}{0.8\standardunitlength}
  \begin{array}{c}  \hspace{-1.7mm}
    \raisebox{-8pt}{\input draws/HSexample1.tex }
      \hspace{-1.9mm}
  \end{array}
 \circ 
  \printname{HSexample2}
  \setlength{\unitlength}{0.8\standardunitlength}
  \begin{array}{c}  \hspace{-1.7mm}
    \raisebox{-8pt}{\input draws/HSexample2.tex }
      \hspace{-1.9mm}
  \end{array}
 =
	
  \printname{HSexampleRslt}
  \setlength{\unitlength}{0.8\standardunitlength}
  \begin{array}{c}  \hspace{-1.7mm}
    \raisebox{-8pt}{\input draws/HSexampleRslt.tex }
      \hspace{-1.9mm}
  \end{array}
 
\]
\caption{The composition in $\HS$ of a morphism in 
$\mor(\ua\ua\da, \ua\ua\ua\da\da)$ whose decorating tensor is 
$a_1\ot \a_2\ot a_3\ot a_4 \ot a_5 \in H^{\ot 4}\ot\bar{H}$, with a 
morphism in $\mor(\ua\ua\ua\da\da, \ua\da\ua)$ decorated by
$b_1\ot b_2\ot b_3\ot b_4 \in H^{\ot 4}$. The general composition
law of $\HS$ is obtained from the one indicated here by multilinear 
extension.}
\lbl{fig:HSexample}
\end{figure} 

\begin{remark}
When $H=\A$, we retrieve the category $\mathbf{AT}$ of 
~\cite[Definition 3.1]{Bar-Natan:NAT}. 
\end{remark}

In order to define the functor $\Z:\catNAT\ra \HS$, we need a
collection of auxiliary maps $\Delta_W:H \ra H^{\ot |W|}$, indexed by
$W\in \W$.
\begin{definition}
We set $\Delta_W = S_W \Delta_W^0$, where 
the maps $\Delta_W^0:H\ra  H^{\ot |W|}$ and $S_W:H^{\ot |W|} \ra H^{\ot |W|}$
are recursively defined as
\begin{alignat*}{2}
& \Delta_{\ew}^0 = \eps \, &; \quad  & S_{\ew} = \id_{\k} \, , \\
& \Delta_{\ua}^0 = \Delta_{\da}^0 = \id \, &; \quad
& S_{\ua} = \id \, , \ \ S_{\da} = S \, , \\
\intertext{and for $W=W_1 W_2$ with $|W_1|, |W_2|\geq 1$,}
& \Delta_W^0 = (\Delta_{W_1}^0 \ot \Delta_{W_2}^0)\Delta \, &; \quad 
& S_{W_1 W_2} = S_{W_1} \ot S_{W_2} \, .
\end{alignat*}
\end{definition}
For example, $\Delta_{(\ua(\da\ua))} = S_2\Delta_2\Delta$.
Note that $\Delta_W^0$ is independent of the directions of the arrows in $W$
and is an algebra morphism, 
while $S_W$ is independent of the parenthesization of $W$.

\begin{definition}
\lbl{def:ZcatNAT}
The action of $\Z:\catNAT\ra \HS$ on objects is just dropping the 
parentheses. The action of $Z$ on the generating morphisms of $\catNAT$
is divided into two stages:
\begin{itemize}
\item
First, a generating morphism is sent to the skeleton read from the graphical
representation of the morphism, with the skeleton components ordered from
left to right. For example, the
skeletons of $\Z((((\ua\ua)\star)\da)^{\NATcp})$ and 
$\Z(\star_{(\da\ua),\da}^{\NATov})$ are
\[
	
  \printname{SkltNATcp}
  \setlength{\unitlength}{0.5\standardunitlength}
  \begin{array}{c}  \hspace{-1.7mm}
    \raisebox{-8pt}{\input draws/SkltNATcp.tex }
      \hspace{-1.9mm}
  \end{array}
 \quad \text{and} \quad 
	
  \printname{SkltNATov}
  \setlength{\unitlength}{0.5\standardunitlength}
  \begin{array}{c}  \hspace{-1.7mm}
    \raisebox{-8pt}{\input draws/SkltNATov.tex }
      \hspace{-1.9mm}
  \end{array}
 \, ,
\]
respectively.
\item
Second, the skeleton is decorated by a tensor specified by the following
table. Here $W\in W^{\star}$, $\, A,B,C \in W$, and we assume that there
are $n_1$ arrows in $W$ to the left of the $\star$, and $n_2$ arrows
to its right. 
\begin{align*}
 	W^{\NATra}_{A,B,C} \, : \quad & 
		1^{\ot n_1} \ot 
		(\Delta_A\ot\Delta_B\ot\Delta_C)(\Phi^{-1}) \ot 1^{\ot n_2} \\
	W^{\NATla}_{A,B,C} \, : \quad & 
		1^{\ot n_1} \ot 
		(\Delta_A\ot\Delta_B\ot\Delta_C)(\Phi) \ot 1^{\ot n_2} \\
	W^{\NATov}_{A,B}  \, : \quad & 
		1^{\ot n_1} \ot 
		(\Delta_A\ot\Delta_B)(R^{21}) \ot 1^{\ot n_2} \\
	W^{\NATun}_{A,B}  \, : \quad &
		1^{\ot n_1} \ot  
		(\Delta_A\ot\Delta_B)(R^{-1}) \ot 1^{\ot n_2} \\
 	W^{\NATcp} \, : \quad & 
		1^{\ot n_1} \ot \alpha \ot 1^{\ot n_2} \\
 	W^{\NATcn} \, : \quad & 
		1^{\ot n_1} \ot S(\alpha)uv^{-1} \ot 1^{\ot n_2} \\
 	W^{\NATap} \, : \quad & 
		1^{\ot n_1} \ot u^{-1}vS(\beta) \ot 1^{\ot n_2} \\
 	W^{\NATan} \, : \quad & 
		1^{\ot n_1} \ot \beta \ot 1^{\ot n_2} 
\end{align*}
\end{itemize}
\end{definition}

We will graphically represent the values of $\Z$ on generating morphisms as
in the following examples, using the notation of \eqref{eqn:PhiComps}
and \eqref{eqn:RmatComps}. 
\begin{align*}
	\Z\left(
  \printname{natLeftAssoc}
  \setlength{\unitlength}{0.7\standardunitlength}
  \begin{array}{c}  \hspace{-1.7mm}
    \raisebox{-8pt}{\input draws/natLeftAssoc.tex }
      \hspace{-1.9mm}
  \end{array}
\right) & =
	
  \printname{natLeftAssocZval}
  \setlength{\unitlength}{0.7\standardunitlength}
  \begin{array}{c}  \hspace{-1.7mm}
    \raisebox{-8pt}{\input draws/natLeftAssocZval.tex }
      \hspace{-1.9mm}
  \end{array}
 \, , \\
	\Z\left(
  \printname{natOver1}
  \setlength{\unitlength}{0.53\standardunitlength}
  \begin{array}{c}  \hspace{-1.7mm}
    \raisebox{-8pt}{\input draws/natOver1.tex }
      \hspace{-1.9mm}
  \end{array}
\right) & =
	
  \printname{natOver1Zval}
  \setlength{\unitlength}{0.53\standardunitlength}
  \begin{array}{c}  \hspace{-1.7mm}
    \raisebox{-8pt}{\input draws/natOver1Zval.tex }
      \hspace{-1.9mm}
  \end{array}
 \, , \\
	\Z\left(
  \printname{natCreatePos1}
  \setlength{\unitlength}{0.7\standardunitlength}
  \begin{array}{c}  \hspace{-1.7mm}
    \raisebox{-8pt}{\input draws/natCreatePos1.tex }
      \hspace{-1.9mm}
  \end{array}
\right) & =
	
  \setlength{\unitlength}{0.7\standardunitlength}
  \begin{array}{c}  \hspace{-1.7mm}
    \raisebox{-2pt}{\input draws/natCreatePosZval.tex }
    \hspace{-1.9mm}
  \end{array}
 \hspace{-0.8cm} \, .
\end{align*}
(Order tags are omitted and the convention of implicit summation over 
indices of tensor components is employed.)

\smallskip
The interpretation of Definition ~\ref{def:ZcatNAT} for diagrammatic algebras 
is simple. We content with a single example for $H= \A_{\kz}$:
\[
	\Z\left(
  \printname{natOver2}
  \setlength{\unitlength}{0.7\standardunitlength}
  \begin{array}{c}  \hspace{-1.7mm}
    \raisebox{-8pt}{\input draws/natOver2.tex }
      \hspace{-1.9mm}
  \end{array}
\right) =
	
  \printname{A3hUnit}
  \setlength{\unitlength}{0.4\standardunitlength}
  \begin{array}{c}  \hspace{-1.7mm}
    \raisebox{-8pt}{\input draws/A3hUnit.tex }
      \hspace{-1.9mm}
  \end{array}
 + \frac{1}{2} 
  \printname{A3hT13}
  \setlength{\unitlength}{0.4\standardunitlength}
  \begin{array}{c}  \hspace{-1.7mm}
    \raisebox{-8pt}{\input draws/A3hT13.tex }
      \hspace{-1.9mm}
  \end{array}

	+ \frac{1}{2} 
  \printname{A3hT23}
  \setlength{\unitlength}{0.4\standardunitlength}
  \begin{array}{c}  \hspace{-1.7mm}
    \raisebox{-8pt}{\input draws/A3hT23.tex }
      \hspace{-1.9mm}
  \end{array}
 + \cdots \, .	
\]

\begin{theorem}
\lbl{thm:ZcatPFT}
The functor $\Z:\catNAT\ra \HS$ respects the relations (R1)--(R10), and 
descends to an invariant of parenthesized framed tangles 
$\Z:\catPFT\ra \HS$.
\end{theorem}
\begin{proof}
(R1) essentially follows from the non-degeneracy of $\Phi$. In more detail,
we prove that $\Z(W^{\NATra}_{A,B,C}) = \Z(\id)$ if $A=\ew$ as follows. Start
from the equality $(\eps\ot\id\ot\id)(\Phi^{-1}) = 1\ot 1$, which is an 
immediate consequence of \eqref{eqn:PhiND}, and apply to it the map 
$(\Delta_B\ot\Delta_C)$, followed by addition of $1^{\ot n_1}$ to the
left and $1^{\ot n_2}$ to the right ($n_1, n_2$ are as in Definition
~\ref{def:ZcatNAT}). The result on the left hand side is 
$\Z(W^{\NATra}_{\ew,B,C})$, while on the right hand side we obtain, using 
$\Delta_U(1)=1^{\ot |U|}$ for all $U\in \W$,  the value
assigned by $\Z$ to the identity morphism of $W/\{\star\ra (\ew B)C\}$.
Similar considerations apply to all the other cases of (R1), as well as to
(R6), using the non-degeneracy of the $R$-matrix \eqref{eqn:EpsRmat}.

A lesson we learn from the above argument and keep for the rest of the proof, 
is that it is sufficient to deal only with cases in which $W=\star$ 
(so that $n_1=n_2=0$).

\smallskip
We turn to (R7) and prove that 
$\Z(\star^{\NATov}_{A,B})\Z(\star^{\NATun}_{B,A}) = \Z(\id)$.
We note that it is enough to prove this
equality for $A=B=\ua$, for then we can apply the map $\Delta_A\ot\Delta_B$
and conclude the equality for any $A, B$, using that $\Delta$ and $\eps$ are 
algebra morphisms and $S$ is an algebra anti-morphism.

The proof for $A=B=\ua$ is simple:
\[
	\Z\left(
  \printname{R7pfOv}
  \setlength{\unitlength}{0.5\standardunitlength}
  \begin{array}{c}  \hspace{-1.7mm}
    \raisebox{-8pt}{\input draws/R7pfOv.tex }
      \hspace{-1.9mm}
  \end{array}
\right)
	\Z\left(
  \printname{R7pfUn}
  \setlength{\unitlength}{0.5\standardunitlength}
  \begin{array}{c}  \hspace{-1.7mm}
    \raisebox{-8pt}{\input draws/R7pfUn.tex }
      \hspace{-1.9mm}
  \end{array}
\right) =
	
  \printname{R7pfOvun}
  \setlength{\unitlength}{0.5\standardunitlength}
  \begin{array}{c}  \hspace{-1.7mm}
    \raisebox{-8pt}{\input draws/R7pfOvun.tex }
      \hspace{-1.9mm}
  \end{array}
 = 
  \printname{R7pfId}
  \setlength{\unitlength}{0.5\standardunitlength}
  \begin{array}{c}  \hspace{-1.7mm}
    \raisebox{-8pt}{\input draws/R7pfId.tex }
      \hspace{-1.9mm}
  \end{array}
 \, .
\]
(R2) is verified in a similar way.

\smallskip
Reasoning as for the previous relations, we find that (R3) and (R8) basicly
reduce to the pentagon relation and hexagon relations, respectively. For
example, it is clear that both sides of the first hexagon are mapped by $\Z$ 
to morphisms having the same skeleton; the decoration for the $\Z$-value of 
the right hand side (with $A=B=C=\ua$) is 
\[
(\Phi^{123})^{-1} R^{32} \Phi^{132} R^{31} (\Phi^{312})^{-1} \, ,
\]
which, by applying the map $X\mapsto X^{231}$ to \eqref{eqn:Hexagon1}, 
is seen to be equal to $(\Delta\ot \id)(R^{21})$, the decorating tensor of 
$\Z(\star^{\NATov}_{(\ua\ua), \ua})$.

\smallskip
Locality in space (R4) obviously follows from the fact that tensors living on
different components commute with each other. We next examine locality in 
scale (R5). It is clearly enough to verify the relations when
the participating morphisms are generators of $\catNAT$, because if a morphism
is the composition of several generating morphisms, then the generating  
morphisms can be `slided up one by one'. Similarly, we may also restrict to 
cases in which one of $W_A$ and $W_B$ (and $W_C$, when relevant) is one of the 
generators of $\catNAT$, while the other one (or two) are the identity 
morphisms. A final reduction is that, like in the proofs of previous
relations, it is enough to check (R4) when the identity morphism(s) are
$\id_{\ua}\in\mor(\ua, \ua)$, and the non-identity morphism is one
of the fundamental generators: 
$\star^{\NATra}_{\ua, \ua, \ua}$, $\star^{\NATla}_{\ua, \ua, \ua}$,
$\star^{\NATov}_{\ua, \ua}$, $\star^{\NATun}_{\ua, \ua}$,
$\star^{\NATcp}$, $\star^{\NATcn}$, $\star^{\NATap}$, $\star^{\NATan}$;
the required equalities with identity morphisms of arbitrary objects and 
all generating morphisms then follow by applying appropriate
tensor products of $\Delta_W$-maps.

\smallskip
We give several representative examples of the cases that remain to 
be checked. In all cases, the skeletons of the morphisms being compared 
are evidently equal, so we focus on the decorating tensors.

\smallskip
Take in the second group of (R4) relations
$W_A=W_C=\id_{\ua}$, $W_B=\star^{\NATla}_{\ua, \ua, \ua}$. 
The decorating tensors for the $\Z$-values of the two sides are
\[
	\Phi^{234}(\id\ot(\Delta_2\Delta_2)\ot\id)(\Phi) \quad \text{and} \quad
	(\id\ot(\Delta_3\Delta_2)\ot\id)(\Phi)\Phi^{234}\, ,
\]
which are equal by the quasi-coassociativity axiom  
\eqref{eqn:QuasiCoAssoc}, applied to the second component of $\Phi$.
Similarly, relations in the third or forth groups with such $W_A$ and 
$W_B$ follow from the quasi-coassociativity applied to components of the 
$R$-matrix or its inverse.  

\smallskip
A relation in the fourth group with $W_A=\star^{\NATun}_{\ua, \ua}$ and 
$W_B = \id_{\ua}$ is equivalent to the equality
\[
  (\Delta^{21}\ot\id)(R^{-1}) = R^{12}(\Delta\ot\id)(R^{-1})(R^{12})^{-1},
\]
which is the quasi-cocommutativity axiom \eqref{eqn:QuasiCoComm}
applied to the first component of $R^{-1}$. Plugging a fundamental braiding
morphism (together with two identity morphisms) into relations of the first 
or second group, we encounter equalities that follow by applying the 
quasi-cocommutativity axiom to components of the associator or its inverse.

\smallskip
As a last example of verifying locality in scale, we put   
$W_A=\star^{\NATap}$, $W_B=W_C=\id_{\ua}$ in the second group
of relations. Using that by (R1)
$\Z(\star^{\NATla}_{\ew, \ua. \ua}) = \Z(\id)$, we need to check that
\[
	X_i^{''} u^{-1} v S(\beta) X_i^{'}\ot Y_i\ot Z_i 
	= u^{-1} v S(\beta)\ot 1 \ot 1 \, .
\]
By applying $S$ to \eqref{eqn:AntipodeId2}, and then using 
\eqref{eqn:SqrAntipode} and the centrality of $v$, the first component
is seen to be equal to $\eps(X_i) u^{-1} v S(\beta)$; our claim now follows
from the non-degeneracy of $\Phi$.

\smallskip
The four (R9) relations follow from \eqref{eqn:AntipodePhi} and 
\eqref{eqn:AntipodePhiInv}, either directly or after applying the antipode to
these equations. We finish the proof of the theorem by verifying the
relation (R10) with $a=\da$. Mapping to $\HS$ each of the seven generating 
morphisms composing the morphism on the left hand side, and then multiplying   
the factors of the decorating tensors along the single skeleton component, 
we are left with proving that
\begin{equation}
\label{eqn:R10}
  	S(\bar{X}_j) S(s_k) S(\alpha) u v^{-1} t_k 
	\bar{Y}_j X_i \beta S(Y_i) S(\bar{Z}_j^{'}) 
	S(\bar{t}_l) S(\alpha) u v^{-1} \bar{s}_l \bar{Z}_j^{''} Z_i
	= \alpha \, .
\end{equation}
By \eqref{eqn:Su_u_alpha},
\[
	S(s_k)S(\alpha) u t_k = S(u)u \alpha = v^2 \alpha \, .
\]
We also have
\[
	S(\bar{t}_l) S(\alpha) u \bar{s}_l = \alpha \, ,
\]
which is established by substituting \eqref{eqn:Salpha_u} into it, and
comparing with the result of the application of the map 
$X \mapsto ((\id\ot S)(X)\ot\alpha)^{312}$ to both
sides of the equation $s_i \bar{s}_k \ot t_i \bar{t}_k = 1\ot 1$.
Utilizing the centrality of $v$, the last two equations simplify the left 
hand side of \eqref{eqn:R10} into
\[
	S(\bar{X}_j) \alpha \bar{Y}_j X _i \beta S(Y_i) 
	S(\bar{Z}_j^{'}) \alpha \bar{Z}_j^{''} Z_i \, ,		
\]
which, by \eqref{eqn:AntipodeId1} and then the non-degeneracy of 
$\Phi^{-1}$, is equal to 
\[
	\alpha X_i \beta S(Y_i) \alpha Z_i \, .
\]
The proof of \eqref{eqn:R10} is now completed by invoking 
\eqref{eqn:AntipodePhi}. 
\end{proof}

\smallskip
Let $T$ be a parenthesized framed tangle whose underlying tangle is the 
(isotopy class of an) embedding of a skeleton $\M_T$.
The invariant $\Z(T)$ is then of the form 
$(\M, X)$, where the skeleton $\M$ has the same components as $\M_T$, but with
a possibly different ordering. The two orderings of components determine
an isomorphism $H(\M)\tilde{\ra} H(\M_T)$, and the image of $X\in H(\M_T)$
under this isomorphism may be regarded as a tensor decorating the components 
of $T$. (In terms of the graphical representation of $\Z$, 
this simply means that the auxiliary order tags used in composing the
images in $\HS$ of the generating morphisms making up $T$, are omitted  
from the final result.) Henceforth, we sometimes consider $\Z(T)$ in this way 
as an element of $H(\M_T)$. 

For a framed link $L$, the invariant $\Z(L)$ takes its value 
in the space $\bar{H}^{\ot m}$, where $m$ is the number of components of $L$.
Suppose that the components of $L$ are {\em colored} by finite-dimensional
$H$-modules $V_1, \ldots, V_m$. We obtain a {\em numerical} invariant 
(taking values in $\k$), denoted $\Z_{V_1, \ldots, V_m}(L)$, by picking a
representative in $H^{\ot m}$ of $\Z(L)$, and computing the trace of
its action on $V_1\ot\cdots\ot V_m$; by the definition of $\bar{H}$, the
value of the trace is independent of the choice of a representative. The
invariants $\Z_{V_1, \ldots, V_m}$ coincide with the invariants constructed 
in ~\cite{AltschulerCoste:QuasiQuantum} (and previously in 
~\cite{ReshetikhinTuraev:Ribbon} for the case of quasitriangular Hopf 
algebras).
  
\begin{example}
{\bf The $H$-invariant of the unknot:} Presenting the unknot as a composition
of generating morphisms in two ways, 
\[
	
  \printname{unknotANCN}
  \setlength{\unitlength}{0.5\standardunitlength}
  \begin{array}{c}  \hspace{-1.7mm}
    \raisebox{-8pt}{\input draws/unknotANCN.tex }
      \hspace{-1.9mm}
  \end{array}
 \quad \text{and} \quad 
	
  \printname{unknotAPCP}
  \setlength{\unitlength}{0.5\standardunitlength}
  \begin{array}{c}  \hspace{-1.7mm}
    \raisebox{-8pt}{\input draws/unknotAPCP.tex }
      \hspace{-1.9mm}
  \end{array}
 \, ,	
\]
we obtain
\begin{equation}
\label{eqn:Zunknot}
\Z(\text{unknot}) = \tr(u^{-1}vS(\beta)\alpha) 
	= \tr(\beta S(\alpha)uv^{-1}) \in \bar{H} \, .
\end{equation}
\end{example}

\begin{exercise}
Show that the $H$-invariant of the framed tangle in Figure ~\ref{fig:Ribbon}
is equal to the ribbon element of $H$.
\end{exercise}

\subsubsection{\bf Behavior of the $H$-invariant under twisting.}
\label{__sect:TwistZ}
Let $F\in H^{\ot 2}$ be a twist, and 
$H_F = (H, \Delta_F, \eps, \Phi_F, R_F, S, \alpha_F, \beta_F)$ the
corresponding twist equivalent quasitriangular quasi-Hopf algebra, which for
the time being we endow with the 
antipodal triple specified by \eqref{eqn:TwistAntipode}. Recall 
(Proposition ~\ref{prop:RibbonTwist}) that $H_F$ has a ribbon element
(same as that of $H$),
hence we may construct the $H_F$-invariant $\Z_F:\catPFT \ra \HS$. We want to 
describe the relationship between $\Z$ and $\Z_F$.  

\smallskip
We recursively define a collection of elements 
$\{F_W^0\in H^{\ot |W|}\}_{W\in\W}$ by
\begin{align*}
	F_{\ew}^0 & = 1_{\k} \\ 
	F_{\ua}^0 & = F_{\da}^0 = 1 \, , \\
\intertext{and for $W=W_1 W_2$ with $|W_1|, |W_2|\geq 1$,}
	F_W^0 & 
	= (F_{W_1}^0 \ot F_{W_2}^0)(\Delta_{W_1}^0\ot \Delta_{W_2}^0)(F) \, .
\end{align*}
We then set 
\[
	F_W = S_W(F_W^0) \in H^{\ot |W|}\, , \quad W\in \W \, .
\]
Similarly, we define elements $G^0_W \in H^{\ot |W|}$ by
\begin{align*}	
  G_{\ew}^0 & = 1_{\k} \\
  G_{\ua}^0 & = G_{\da}^0 = 1 \, , \\
  G_{W_1 W_2}^0 & 
  = (\Delta_{W_1}^0\ot \Delta_{W_2}^0)(F^{-1}) (G_{W_1}^0 \ot G_{W_2}^0)\, ,
\end{align*}
and set $G_W = S_W(G_W^0)$. The elements $F_W^0$ and $G_W^0$ are inverse to 
each other (but this is not necessarily true for $F_W$ and $G_W$).
Two examples are 
\begin{equation}
\label{eqn:LeftAssocTwist}
F_{(\ua(\ua\ua))} = F^{23}\Delta_2(F) \quad \text{and} \quad
G_{((\ua\ua)\ua)} = \Delta_1(F^{-1}) (F^{12})^{-1} \, .      
\end{equation}

\smallskip
We sometimes denote the twisted coproduct $\Delta_F$ by $\Delta'$, for 
typographical convenience. 
It is not hard to prove by induction on $|W|$ the following generalization
of \eqref{eqn:TwistDelta}:
\begin{equation}
\label{eqn:GenrlTwistDelta}
  \Delta_W^{'0}(x) = F_W^0 \Delta_W^0(x) G_W^0 \, , \quad \text{for all} \
	W\in \W\, , x\in H\, .
\end{equation}

If $W\in \W^{\star}$ has $n_1$ arrows left to the $\star$ and $n_2$ arrows
to its right, then for all $A\in \W$ the following formulas hold: 
\begin{equation}
\label{eqn:FwSubst}
\begin{aligned}
F_{W/\{\star\ra A\}}^0 & = 
  (1^{\ot n_1} \ot F_A^0 \ot 1^{\ot n_2})\cdot
  (\id^{\ot n_1} \ot \Delta_A^0 \ot \id^{\ot n_2})(F_{W/\{\star\ra \ua\}}^0)
  \, , \\
G_{W/\{\star\ra A\}}^0 & = 
  (\id^{\ot n_1} \ot \Delta_A^0 \ot \id^{\ot n_2})(G_{W/\{\star\ra \ua\}}^0)
  \cdot (1^{\ot n_1} \ot G_A^0 \ot 1^{\ot n_2})\, .
\end{aligned}
\end{equation}
These are proved by induction on the structure of $W$.

\smallskip
For $W\in \W$, we also denote by $W$ the same string after erasing the
parentheses, which is an object in $\HS$. We let $\F_W$ be the morphism 
in $\mor_{\HS}(W, W)$ whose skeleton consists of $|W|$ vertical strands 
(ordered from left to right), and whose decorating tensor is $F_W$. For 
example, 
\[
	F_{(\ua(\da\ua))} = 
  \printname{FwExample}
  \setlength{\unitlength}{0.75\standardunitlength}
  \begin{array}{c}  \hspace{-1.7mm}
    \raisebox{-8pt}{\input draws/FwExample.tex }
      \hspace{-1.9mm}
  \end{array}
 \hspace{-1cm} \, .
\]

By applying the map $S_W$ to the equation
\[
	F_W^0 G_W^0 = 1^{\ot |W|} = G_W^0 F_W^0 \, ,
\]
and recalling the composition rule in $\HS$, we see that $\F_W$ has a 
two-sided inverse $\F_W^{-1}$, which is 
the morphism with the same skeleton as $\F_W$ and decorating tensor $G_W$. 

\begin{proposition}
\lbl{prop:LMtwist}
Let $D, U \in \W$, and let $T$ be a parenthesized framed tangle belonging to 
$\mor_{\catPFT}(D, U)$. Then
\begin{equation}
\label{eqn:TwistZ}
	\Z_F(T) = \F_U^{-1} \Z(T) \F_D \, .
\end{equation}
\end{proposition}
\begin{proof}
It is clearly enough to prove the proposition when $T$ is a generator of 
$\catPFT$. The morphisms on both sides of \eqref{eqn:TwistZ} have the same
skeleton, so it is left to establish the equality of their decorating 
tensors, which we denote in this proof by $\Z(T)$ and 
$\Z_F(T)$.

\smallskip
Before explaining the general proof, it may be instructive to check 
\eqref{eqn:TwistZ} for several fundamental generators. For 
$T = \star^{\NATla}_{\ua, \ua, \ua}$ we have $\Z(T)= \Phi$, $\Z_F(T)= \Phi_F$,
and the result is immediate from
\eqref{eqn:LeftAssocTwist} and \eqref{eqn:TwistPhi};
for $T= \star^{\NATun}_{\ua, \ua}$ the required equality is
\[
	R_F^{-1} = F R^{-1} (F^{21})^{-1}
\]
which is obtained by inverting \eqref{eqn:TwistR}, and for 
$T = \star^{\NATan}$, using $F_{(\ua\da)} = S_2(F)$, it is 
\[
	\beta_F = f_i \beta S(g_i) \, , 
\]
which is indeed the twisting rule for $\beta$.

\smallskip
As a typical example of the proof for general associativity and braiding
morphisms, we consider an $(A,B)$-braiding extended by a single $C$-bundle
of vertical strands on the right; that is, $T=W^{\NATun}_{A,B}$, where
$W=\star C$ and $A, B, C\in \W$. We may assume that there are only $\ua$ 
arrows in $A, B, C$; once \eqref{eqn:TwistZ} is established in those 
cases, it can be extended to any combination of arrow directions by applying
$S$ (which is also $S_F$) on both sides of  \eqref{eqn:TwistZ} to the tensor 
factors corresponding to $\da$ arrows. Note that with only $\ua$ arrows 
around, we have $\Delta_A = \Delta_A^0$, $F_B = F_B^0$, etc.; in particular, 
$\Delta_A$, $\Delta_B$, and $\Delta_C$ are algebra morphisms. 

\smallskip
We turn to the computation itself. By Definition ~\ref{def:ZcatNAT},
\[
  \Z_F(T) = (\Delta_A^{'} \ot \Delta_B^{'})(R_F^{-1}) \ot 1^{\ot |C|} \, ,
\]
which by \eqref{eqn:GenrlTwistDelta} is equal to
\begin{equation}
\label{eqn:TwistZ_braidLHS}
  (F_A \ot F_B)\cdot (\Delta_A \ot \Delta_B)(R_F^{-1})\cdot (G_A \ot G_B) 
		\ot 1^{\ot |C|} \, .
\end{equation}
The decorating tensor for $\F_{(BA)C}^{-1} \Z(T) \F_{(AB)C)}$ is
\[
	F_{(AB)C}\cdot (\Delta_A \ot \Delta_B)(R^{-1} \ot 1^{|C|})
	\cdot G_{(BA)C}^{\pi} \, ,
\]
where $\pi$ is the permutation of tensor factors `swapping $A$ and $B$';
more precisely, $\pi$ leaves unchanged the last $|C|$ factors, and cycles
the first $|A|+|B|$ factors by $|A|$ places to the right.
Repeatedly using \eqref{eqn:FwSubst}, the above expression becomes
\begin{equation}
\label{eqn:TwistZ_braidRHS}
\begin{aligned}
  & (F_A\ot F_B\ot F_C) \cdot 
    (\Delta_A\ot \Delta_B \ot \Delta_C)(F_{((\ua\ua)\ua)})\cdot
    (\Delta_A\ot \Delta_B\ot \Delta_C)(R^{-1} \ot 1) \\
  & \cdot (\Delta_A\ot \Delta_B \ot \Delta_C)(G_{((\ua\ua)\ua)}^{213})\cdot
    (G_A \ot G_B \ot G_C) \, .
\end{aligned}
\end{equation}
By the defining property \eqref{eqn:QuasiCoComm} of the $R$-matrix,
\[
\begin{aligned}
	F_{((\ua\ua)\ua)}\cdot (R^{-1} \ot 1)\cdot G_{((\ua\ua)\ua)}^{213} & =
	F^{12}\cdot (\Delta\ot \id)(F)\cdot (R^{-1} \ot 1)\cdot
	(\Delta^{21}\ot \id)(F^{-1})\cdot (F^{21})^{-1} \\
	& = R_F^{-1} \ot 1 \,	
\end{aligned}
\]
which simplifies \eqref{eqn:TwistZ_braidRHS} into
\[
	(F_A\ot F_B\ot F_C)\cdot
	(\Delta_A\ot \Delta_B\ot \Delta_C)(R_F^{-1} \ot 1)\cdot
	(G_A \ot G_B \ot G_C) \, ;
\]
this last expression is clearly equal to \eqref{eqn:TwistZ_braidLHS}.

\smallskip
We move to demonstrate the proof for creation and annihilation morphisms,
and consider the case $T=(\star A)^{\NATan}$. As in the previous case, we may 
assume that all arrows in $A$ are pointing up. By Definition 
~\ref{def:ZcatNAT}, $\Z_F(T) = \beta_F \ot 1^{\ot |A|}$. The decorating
tensor for $\F_{(\ua\da)A}$ is, using \eqref{eqn:FwSubst},
\[
  F_{(\ua\da)A} = S_2((1^{\ot 2} \ot F_A)\cdot 
	(\id^{\ot 2} \ot \Delta_A)(F^{12}\cdot (\Delta\ot \id)(F))) \, ,
\] 
and $\F_A^{-1}$ is decorated by $G_A$. The composition of morphisms 
on the right hand side of \eqref{eqn:TwistZ} is therefore
\[
	
  \printname{TwistZpfAN}
  \setlength{\unitlength}{0.5\standardunitlength}
  \begin{array}{c}  \hspace{-1.7mm}
    \raisebox{-8pt}{\input draws/TwistZpfAN.tex }
      \hspace{-1.9mm}
  \end{array}
 \, ,
\] 
so that its decorating tensor is
\[
    f_i f_j^{'} \beta S(f_j^{''}) S(g_i) \ot F_A \Delta_A(g_j) G_A =
    \eps(f_j) \beta_F \ot F_A \Delta_A(g_j) G_A = \beta_F \ot 1^{\ot |A|} \, ,
\]
where we have used \eqref{eqn:AntipodeId2}, \eqref{eqn:TwistAntipode}, and 
finally the non-degeneracy of $F$.
\end{proof}

\begin{remark}
Proposition ~\ref{prop:LMtwist} is a generalization of Theorem 7 of
~\cite{LeMurakami:Universal}, which was formulated for the case $H = \A_{\kz}$.
\end{remark}

\smallskip
The use of the particular antipodal triple $(S, \alpha_F, \beta_F)$ is too 
restrictive for our applications, 
so we need to understand the effect of antipodal rescaling on $\Z$. 
Suppose that the antipodal triple $(S, \alpha, \beta)$ of $H$ is 
rescaled by an invertible element $c$ into an antipodal triple 
$(\hat{S}, \hat{\alpha}, \hat{\beta})$ given by \eqref{eqn:AntipodalRescale},
and denote by $\hat{\Z}$ the corresponding invariant of $\catPFT$.

\smallskip
For $W\in \W$, write $W$ after erasing the parentheses as 
$W = a_1 a_2 \ldots a_{|W|}$, where $a_i\in \{\ua, \da\}$. Put
$c_{\ua} = 1$, $c_{\da} = c$, and let $\C_W \in \mor_{\HS}(W, W)$ be the 
morphism with the same skeleton as $\F_W$, decorated by the tensor 
$\otimes_{i=1}^{|W|} c_{a_i}$.
The inverse morphism $\C_W^{-1}$ is defined by replacing $c$ with $c^{-1}$.
\begin{proposition}
For $D, U \in \W$ and $T\in \mor_{\catPFT}(D, U)$, 
\begin{equation}
\label{eqn:RescaleZ}
	\hat{\Z}(T) = \C_U \Z(T) \C_D^{-1} \, .
\end{equation}
\end{proposition}
\begin{proof}
It is sufficient to verify \eqref{eqn:RescaleZ} for generating morphisms. The
only non-trivial cases are the morphisms of types $W^{\NATcn}$ and 
$W^{\NATap}$, which are settled with the help of Proposition 
~\ref{prop:RibbonRescale}.
\end{proof}

Combining \eqref{eqn:TwistZ} and \eqref{eqn:RescaleZ} when $T$ is a 
{\em framed link} $L$ (so that $D = U = \ew$), we conclude:
\begin{corollary}
\lbl{corr:ZlinkIndepTwist}
If  $H_1$ and $H_2$ are twist equivalent ribbon quasitriangular quasi-Hopf
algebras, then $\Z_{H_1}(L) = \Z_{H_2}(L)$ for any framed link L.
\end{corollary}

\begin{remark}
For another proof of Corollary ~\ref{corr:ZlinkIndepTwist} (when  
the invariants are evaluated on finite-dimensional modules coloring the 
components of $L$), see ~\cite{LinksGouldZhang:TwistInvarLinkPoly}. The proof 
there uses the difficult machinery of the presentation of framed links as 
closures of braids modulo Markov moves.
\end{remark}

\subsection{The Kontsevich integral and the Reshetikhin-Turaev invariants}
\label{_sect:KontAndRT}

We now discuss the invariants of framed links obtained by applying the 
$H$-invariant construction to the ribbon quasitriangular quasi-Hopf
algebras described in Sections ~\ref{_sect:KZ} and ~\ref{_sect:DJ}.

\smallskip
The invariant corresponding to $\A_{\kz}$ is denoted by $\ZK$. Note that by 
Proposition ~\ref{prop:AkzIndepPhi} and Corollary ~\ref{corr:ZlinkIndepTwist},
the invariant $\ZK$ is independent of the Drinfel'd
associator used for defining $\A_{\kz}$. 

$\ZK$, when restricted to framed knots\footnote{The discussion can be
easily extended to framed links.}, is a 
{\em universal Vassiliev invariant}: for a chord diagram $D$ of degree $m$, 
\[
	\ZK(K_D) = D + (\text{terms of degree } > m) \, , 
\]
where $K_D$ is an embedding of $D$ as in Section ~\ref{__sect:CDandWgtSys}, 
and $\ZK$ is extended to singular framed knots by \eqref{eqn:SingExtnd}. 
The basic observation in the proof is that by Definition ~\ref{def:ZcatNAT}, 
with $R = e^{\Omega/2}$,
\[
	\ZK\left(\text{\Huge \raisebox{-6pt}{$\overcrossing$}}\right) - 
	\ZK\left(\text{\Huge \raisebox{-6pt}{$\undercrossing$}}\right)
        = R^{21} - R^{-1} = 
	
  \printname{BCup}
  \setlength{\unitlength}{0.25\standardunitlength}
  \begin{array}{c}  \hspace{-1.7mm}
    \raisebox{-8pt}{\input draws/BCup.tex }
      \hspace{-1.9mm}
  \end{array}
 + \text{higher degree terms} \, .
\]
See ~\cite[Theorem 4]{Bar-Natan:NAT} for more details.

\smallskip
Le and Murakami (\cite{LeMurakami:KontIterIntgrl}) proved that when
$\Phi = \Phi_{\kz}$, the invariant $\Z_K$ coincides with 
the iterated integral formula of Kontsevich 
(\cite[Section 4]{Bar-Natan:Vassiliev}), which was the first (rigorous)
construction of a universal Vassiliev invariant. An interesting consequence
of the independence of $\ZK$ of the Drinfel'd associator and Theorem 
~\ref{thm:RationalPhi}, is that the values taken by Kontsevich's highly 
transcendental formula are linear combinations of chord diagrams with 
{\em rational} coefficients.

Since the $u$-element and ribbon element of $\A_{\kz}$ are equal and
$\beta = 1$, we have by \eqref{eqn:Zunknot} that
\[
	\Z_K(\text{unknot}) = \nu \in \A(\bigcirc)\, .
\]
(We have omitted the trace projection, as $\A(\bigcirc)$ may be 
identified with $\A$.) An explicit formula for this element is given by
\eqref{eqn:BGRTformula}.

\smallskip
It is clear that the invariants associated to the KZ structures $\Aar_{\kz}$
and $\Ugh_{\kz}$ ($\g$ a metrized Lie algebra) are simply the compositions 
$\iota \circ \ZK$ and $\Tgh \circ \ZK$,  respectively. We denote the latter
invariant by $\KZg$.

\smallskip
Let $\g$ be a complex semisimple Lie algebra.
The invariant associated to $\Ugh_{\qg}$ is denoted by $\RTg$; 
the numerical invariants of colored framed links derived from $\RTg$ are the
{\em Reshetikhin-Turaev invariants} $\RT_{\g; V_1, \ldots, V_m}$ 
(\cite{ReshetikhinTuraev:Ribbon}).
From Theorem ~\ref{thm:QGtwistKZ} and Corollary ~\ref{corr:ZlinkIndepTwist}
we deduce that following equality holds in $\CoUgh^{\ot m}$, for any framed 
link $L$ with $m$ components:
\begin{equation}
\label{eqn:RTgEqualKZg}
	\RTg(L) = \KZg(L) \, .
\end{equation}
See ~\cite[Section 5]{Piunikhin:Weights}, 
~\cite[Theorem 10]{LeMurakami:Universal}, and 
~\cite[Section XX.8]{Kassel:QuantumGroups}
for other proofs of this equality, at the (equivalent) level of the
associated numerical invariants of framed links colored by finite-dimensional
$\g$-modules.  
  
\smallskip
Since $\Ugh_{\qg}$ is a Hopf algebra, we get from \eqref{eqn:Zunknot}
and \eqref{eqn:uDJ} that (compare ~\cite[Section 7]{ReshetikhinTuraev:Ribbon})
\begin{equation}
\label{eqn:RTgUnknot}
	\RTg(\text{unknot}) = \tr(e^{\hbar \varrho_{\g}}) 
	= \tr(e^{-\hbar \varrho_{\g}}) \in \CoUgh \, .
\end{equation}

%% file: draws/1btC.tex
%
\begingroup\makeatletter\ifx\SetFigFont\undefined%
\gdef\SetFigFont#1#2#3#4#5{%
  \reset@font\fontsize{#1}{#2pt}%
  \fontfamily{#3}\fontseries{#4}\fontshape{#5}%
  \selectfont}%
\fi\endgroup%
{\renewcommand{\dashlinestretch}{30}
\begin{picture}(936,1431)(0,-10)
\thicklines
\path(93,33)(93,1383)
\blacken\path(153.000,1203.000)(93.000,1383.000)(33.000,1203.000)(153.000,1203.000)
\path(843,33)(843,1383)
\blacken\path(903.000,1203.000)(843.000,1383.000)(783.000,1203.000)(903.000,1203.000)
\thinlines
\put(843.000,633.000){\arc{600.000}{1.5708}{4.7124}}
\end{picture}
}

%% file: draws/Cbt1.tex
%
\begingroup\makeatletter\ifx\SetFigFont\undefined%
\gdef\SetFigFont#1#2#3#4#5{%
  \reset@font\fontsize{#1}{#2pt}%
  \fontfamily{#3}\fontseries{#4}\fontshape{#5}%
  \selectfont}%
\fi\endgroup%
{\renewcommand{\dashlinestretch}{30}
\begin{picture}(936,1431)(0,-10)
\thicklines
\path(93,33)(93,1383)
\blacken\path(153.000,1203.000)(93.000,1383.000)(33.000,1203.000)(153.000,1203.000)
\path(843,33)(843,1383)
\blacken\path(903.000,1203.000)(843.000,1383.000)(783.000,1203.000)(903.000,1203.000)
\thinlines
\put(93.000,633.000){\arc{600.000}{4.7124}{7.8540}}
\end{picture}
}

%% file: draws/SkltNATov.tex
%
\begingroup\makeatletter\ifx\SetFigFont\undefined%
\gdef\SetFigFont#1#2#3#4#5{%
  \reset@font\fontsize{#1}{#2pt}%
  \fontfamily{#3}\fontseries{#4}\fontshape{#5}%
  \selectfont}%
\fi\endgroup%
{\renewcommand{\dashlinestretch}{30}
\begin{picture}(1674,789)(0,-10)
\path(12,762)(1587,12)
\path(1465.759,36.506)(1587.000,12.000)(1491.555,90.678)
\path(837,12)(1662,762)
\path(1593.387,659.081)(1662.000,762.000)(1553.027,703.477)
\path(80.613,114.919)(12.000,12.000)(120.973,70.523)
\path(12,12)(837,762)
\end{picture}
}

%% file: draws/natOver2.tex
%
\begingroup\makeatletter\ifx\SetFigFont\undefined%
\gdef\SetFigFont#1#2#3#4#5{%
  \reset@font\fontsize{#1}{#2pt}%
  \fontfamily{#3}\fontseries{#4}\fontshape{#5}%
  \selectfont}%
\fi\endgroup%
{\renewcommand{\dashlinestretch}{30}
\begin{picture}(924,789)(0,-10)
\path(773.360,655.934)(837.000,762.000)(730.934,698.360)
\path(837,762)(87,12)
\path(462,462)(162,762)
\path(268.066,698.360)(162.000,762.000)(225.640,655.934)
\path(762,12)(462,312)
\path(537,387)(912,12)
\path(387,387)(12,762)
\path(118.066,698.360)(12.000,762.000)(75.640,655.934)
\end{picture}
}

%% file: draws/R7pfId.tex
%
\begingroup\makeatletter\ifx\SetFigFont\undefined%
\gdef\SetFigFont#1#2#3#4#5{%
  \reset@font\fontsize{#1}{#2pt}%
  \fontfamily{#3}\fontseries{#4}\fontshape{#5}%
  \selectfont}%
\fi\endgroup%
{\renewcommand{\dashlinestretch}{30}
\begin{picture}(1126,1539)(0,-10)
\path(42,12)(42,1512)
\path(72.000,1392.000)(42.000,1512.000)(12.000,1392.000)
\path(867,12)(867,1512)
\path(897.000,1392.000)(867.000,1512.000)(837.000,1392.000)
\put(117,687){\makebox(0,0)[lb]{\smash{{\mathmode{\scriptstyle 1}}}}}
\put(942,687){\makebox(0,0)[lb]{\smash{{\mathmode{\scriptstyle 1}}}}}
\end{picture}
}

%% file: draws/unknotANCN.tex
%
\begingroup\makeatletter\ifx\SetFigFont\undefined%
\gdef\SetFigFont#1#2#3#4#5{%
  \reset@font\fontsize{#1}{#2pt}%
  \fontfamily{#3}\fontseries{#4}\fontshape{#5}%
  \selectfont}%
\fi\endgroup%
{\renewcommand{\dashlinestretch}{30}
\begin{picture}(935,1008)(0,-10)
\put(468.000,459.000){\arc{900.000}{6.2832}{9.4248}}
\path(64.420,344.348)(18.000,459.000)(5.005,335.991)
\put(468.000,534.000){\arc{900.000}{3.1416}{6.2832}}
\path(871.580,648.652)(918.000,534.000)(930.995,657.009)
\end{picture}
}

%% file: draws/unknotAPCP.tex
%
\begingroup\makeatletter\ifx\SetFigFont\undefined%
\gdef\SetFigFont#1#2#3#4#5{%
  \reset@font\fontsize{#1}{#2pt}%
  \fontfamily{#3}\fontseries{#4}\fontshape{#5}%
  \selectfont}%
\fi\endgroup%
{\renewcommand{\dashlinestretch}{30}
\begin{picture}(935,1008)(0,-10)
\put(468.000,459.000){\arc{900.000}{6.2832}{9.4248}}
\path(930.995,335.991)(918.000,459.000)(871.580,344.348)
\put(468.000,534.000){\arc{900.000}{3.1416}{6.2832}}
\path(5.005,657.009)(18.000,534.000)(64.420,648.652)
\end{picture}
}

%% file: EK.tex
\section{The Etingof-Kazhdan quantization of $\protect\Aar$}
\label{sect:EK}

{\bf Section summary:} The first two subsections are devoted to our
diagrammatic reformulation of the Etingof-Kazhdan quantization of the
ambient Lie bialgebra of a finite-dimensional Manin triple. We then
consider the resulting framed link invariant $\Z_{\ek}$, and examine the 
extent to which it has the nice features of the Reshetikhin-Turaev invariants.
A particular attention is given to the value of $\Z_{\ek}$ on the unknot.

\smallskip
After giving a diagram of spaces and maps, which summarizes 
some of the algebraic and  diagrammatic quasitriangular quasi-Hopf structures
(together with the related framed link invariants) that we have considered 
in this paper, we provide a short review of the work of Polyak on directed 
versions of the chord  diagram models of $\A$-spaces, which we call 
$\Par$-spaces, and ask whether it is possible to carry out our construction 
within these spaces.

\subsection{Verma modules}
\subsubsection{\bf The algebras $\protect\Aar_{\pm}$.} 
\label{__sect:Aar-pm}

Let $(\tg, \g_+, \g_-)$ be a finite-dimensional Manin triple 
(see Section ~\ref{__sect:manin}). We saw that the map 
$\Tartg: \Aar \ra U(\tg)$ makes $\Aar$ into a `diagrammatic version'
of the Hopf algebra $U(\tg)$. The algebras $U(\g_+)$ and $U(\g_-)$ are
naturally Hopf subalgebras of $U(\tg)$. These subalgebras have diagrammatic 
companions inside $\Aar$.

\begin{definition}
\lbl{def:Aar-pm}
Let $\Dar_+(\ua)$ be the subset of $\Dar(\ua)$ consisting of diagrams
whose legs are all {\em incoming} into the skeleton interval. The quotient
of $\Span_{\k}\{\Dar_+(\ua)\}$ by NS, AS, STU, and IHX relations
(in which all the participating diagrams belong to $\Dar_+(\ua)$) is denoted
by $\Aar_+$. Symbolically,
\[
	\Aar_+ = \{
  \printname{InDgrm}
  \setlength{\unitlength}{0.3\standardunitlength}
  \begin{array}{c}  \hspace{-1.7mm}
    \raisebox{-8pt}{\input draws/InDgrm.tex }
      \hspace{-1.9mm}
  \end{array}
\} 
		\left/ \{\text{ NS, AS, STU, IHX }\}\right. \, . 
\]
Similarly, restricting to diagrams having all their legs {\em outgoing} from 
the skeleton interval, we define
\[
	\Aar_- = \{
  \printname{OutDgrm}
  \setlength{\unitlength}{0.3\standardunitlength}
  \begin{array}{c}  \hspace{-1.7mm}
    \raisebox{-8pt}{\input draws/OutDgrm.tex }
      \hspace{-1.9mm}
  \end{array}
\} 
		\left/ \{\text{ NS, AS, STU, IHX }\}\right. \, .
\]
\end{definition}

Analogous definitions can be made for other skeletons, which enables
us to define the diagrammatic tensor powers of $\Aar_{\pm}$. For example,
$\Aar_{++}(\ua, \ua) = \Aar_+ \bt \Aar_+$ is the quotient of the (formal
span of)
\[
\Dar_{++}(\ua_1, \ua_2) = \biggl\{
  \printname{InDgrmDblSklt}
  \setlength{\unitlength}{0.3\standardunitlength}
  \begin{array}{c}  \hspace{-1.7mm}
    \raisebox{-8pt}{\input draws/InDgrmDblSklt.tex }
      \hspace{-1.9mm}
  \end{array}
 \biggr\} \, .
\]

We can also consider mixed diagrammatic tensor products, such as
$\Aar \bt \Aar_+$ and $\Aar_+ \bt \Aar_-$, which arise from the diagram
collections
\[
\Dar_{\bullet +}(\ua_1, \ua_2) = \biggl\{
  \printname{MixedInDgrm}
  \setlength{\unitlength}{0.3\standardunitlength}
  \begin{array}{c}  \hspace{-1.7mm}
    \raisebox{-8pt}{\input draws/MixedInDgrm.tex }
      \hspace{-1.9mm}
  \end{array}
 \biggr\}
\text{ and }
\Dar_{+-}(\ua_1, \ua_2) = \biggl\{
  \printname{InOutDgrm}
  \setlength{\unitlength}{0.3\standardunitlength}
  \begin{array}{c}  \hspace{-1.7mm}
    \raisebox{-8pt}{\input draws/InOutDgrm.tex }
      \hspace{-1.9mm}
  \end{array}
 \biggr\}
\, \text{, respectively.}
\]

The next proposition lists some basic properties of $\Aar_{\pm}$.
\begin{proposition}
$\Aar_{\pm}$ are diagrammatic Hopf subalgebras of $\Aar$, and are mapped by 
$\Tartg: \Aar \ra U(\tg)$ into $U(\g_{\pm})$. 
\end{proposition}
\begin{proof}
The proof is a series of simple observations, which we demonstrate by two 
samples:
\begin{itemize}
\item
$\Aar_+$ is closed under the coproduct of $\Aar$, because the coproduct of
a diagram $D \in \Dar_+(\ua)$ is, by Definition ~\ref{def:cabling}, the sum
over diagrams obtained by splitting the legs of $D$ among the two skeleton 
intervals. Clearly, all the resulting diagrams belong to 
$\Dar_{++}(\ua_1, \ua_2)$.
\item
When mapping a diagram in $\Dar_-(\ua)$ to $U(\tg)$, all the elements of $\tg$
placed on the vertices along the skeleton interval belong to $\g_-$ 
(see Figure ~\ref{fig:AtomicDirectedGraphs}). Hence 
$\Tartg(\Aar_-) \subset U(\g_-)$.
\end{itemize}
\end{proof}

Taking a single color skeleton, we obtain the Hopf subalgebras $\Bar_{\pm}$
of $\Bar$, which are thought of as diagrammatic versions of $S(\g_{\pm})$ (and 
are mapped into them by $\Tartg$). In consistency with the algebraic 
situation, we have:
\begin{proposition}
\lbl{prop:PBW-Aar-pm}
The PBW isomorphism
$\upchi : \Bar \ra \Aar$ (Proposition ~\ref{prop:PBW-Aar}) restricts to
isomorphisms $\upchi_{\pm}: \Bar_{\pm} \ra \Aar_{\pm}$. 
\end{proposition}
\begin{proof}
It is clear that  $\upchi$ maps $\Bar_+$ into $\Aar_+$. 
The map $\sigma: \Bar \ra \Aar$ is constructed by glueing various
linear combinations of directed diagrams. By the definition of the
glueing operators in \eqref{eqn:GlueCtoD}, the result of glueing any directed 
diagram to a diagram in $\Dar_+(\ua)$ is again an element of 
$\Dar_+(\ua)$. Thus, $\sigma$ restricts to a map 
$\sigma_+: \Bar_+ \ra \Aar_+$, which is obviously an inverse of
 $\upchi_+: \Bar_+ \ra \Aar_+$.   
\end{proof}

Let $\upchi_+ \bt \upchi_- : \Bar_+\bt \Bar_- \ra \Aar_+ \bt \Aar_-$ 
be the diagrammatic tensor product
of the isomorphisms in Proposition ~\ref{prop:PBW-Aar-pm}. We recall that
it is informally defined as applying $\upchi$ separately on each skeleton 
component (see Remark ~\ref{rem:AlgDgrmCat} for a more formal discussion). 
Explicitly, if $D\in\Dar_{+-}(\ast_1, \ast_2)$ has $m_1$ vertices 
$v_1, \ldots, v_{m_1}$ colored by $\ast_1$ and $m_2$ vertices 
$w_1, \ldots, w_{m_2}$ colored by $\ast_2$, then
\[
\upchi\bt\upchi(D) = \frac{1}{m_1! m_2!} 
	\sum_{\pi_1\in S_{m_1}} \sum_{\pi_2\in S_{m_2}} (\pi_1\bt\pi_2)D\, ,
\]
where $(\pi_1\bt\pi_2)D$ is the diagram obtained from $D$ by placing
$v_{\pi_1(1)}, \ldots, v_{\pi_1(m_1)}$ (in this order) along a first 
skeleton interval, and placing $w_{\pi_2(1)}, \ldots, w_{\pi_2(m_2)}$ along
a second skeleton interval. 

Since the inverse isomorphisms 
$\sigma_{\pm}: \Aar_{\pm} \ra \Bar_{\pm}$ also act at the level of skeleton
components, we may form the map 
$\sigma_+ \bt \sigma_- : \Aar_+\bt \Aar_- \ra \Bar_+ \bt \Bar_-$.
A little thought shows that it is an inverse of $\upchi_+ \bt \upchi_-$. We
conclude:
\begin{proposition}
The map 
$\upchi_+ \bt \upchi_- : \Bar_+\bt \Bar_- \ra \Aar_+ \bt \Aar_-$ 
is an isomorphism.
\end{proposition}

The PBW theorem implies that the product map 
$\mu: U(\g_+) \ot U(\g_-) \ra U(\tg)$ is a linear isomorphism.
Accordingly, we have the following diagrammatic result.
\begin{proposition}
\lbl{prop:AarPlus-bt-AarMinusEqAar}
The map $\mu:\Aar_+ \bt \Aar_- \ra \Aar$ is an isomorphism of
linear spaces.
\end{proposition}
\begin{proof}
Consider the diagram of maps
\begin{equation}
\label{cd:Aar-pm-Bar-pm-prod}
\begin{CD}
\Aar_+ \bt \Aar_- @>{\mu}>> \Aar \\
@A{\upchi_+ \bt \upchi_-}AA  @AA{\upchi}A \\
\Bar_+ \bt \Bar_- @>{\mu}>> \Bar 
\end{CD}
\end{equation}
This diagram is {\em not} commutative. Nevertheless, we claim that it does
commute up to lower filtration levels. More precisely, as in the proofs
of Proposition ~\ref{prop:PBW-A} and Proposition ~\ref{prop:PBW-Aar}, 
we filter the diagram collections $\Dar_{+-}(\ast_1, \ast_2)$, $\Dar(\ast)$, 
$\Dar_{+-}(\ua_1, \ua_2)$, and $\Dar(\ua)$
according to the number of legs; this induces filtrations of the spaces
appearing in \eqref{cd:Aar-pm-Bar-pm-prod}, which are preserved
by the maps between them. Our claim is that if 
$D\in \Dar_{+-}(\ast_1, \ast_2)$ has $k$ legs, then
the element $\upchi\circ\mu(D) - \mu\circ(\upchi_+\bt\upchi_-)(D)$
is equivalent, modulo $\Aar$-relations, to a linear combination of diagrams
having no more than $k-1$ legs. 

\smallskip
Assume that $D$ has $m_1$ legs ending in $\ast_1$ and $m_2$ legs ending in 
$\ast_2$ (so that $m_1 + m_2 = k$). Write
\begin{align*}
(\upchi\circ\mu & - \mu\circ(\upchi_+\bt\upchi_-))
	 \left(
  \printname{InOutAstDgrmD}
  \setlength{\unitlength}{0.35\standardunitlength}
  \begin{array}{c}  \hspace{-1.7mm}
    \raisebox{-8pt}{\input draws/InOutAstDgrmD.tex }
      \hspace{-1.9mm}
  \end{array}
\right) = \\
\\
& \frac{1}{k!}\sum_{\pi\in S_k}
  \printname{AstDgrmPi}
  \setlength{\unitlength}{0.3\standardunitlength}
  \begin{array}{c}  \hspace{-1.7mm}
    \raisebox{-8pt}{\input draws/AstDgrmPi.tex }
      \hspace{-1.9mm}
  \end{array}

	- \frac{1}{m_1! m_2!}
	\sum_{\substack{\pi_1\in S_{m_1}, \\ \pi_2\in S_{m_2}}}
		
  \printname{AstDgrmPi1Pi2}
  \setlength{\unitlength}{0.3\standardunitlength}
  \begin{array}{c}  \hspace{-1.7mm}
    \raisebox{-8pt}{\input draws/AstDgrmPi1Pi2.tex }
      \hspace{-1.9mm}
  \end{array}
 \, ,
\end{align*}
where the permutation boxes determine the orderings of vertices along
the skeleton interval. In a $\pi$-summand, we can move all the outgoing legs
to the right by a repeated use of (mixed) STU relations. Thus, each
$\pi$-summand is equivalent, modulo $\Aar$-relations and diagrams
having less than $k$ legs, to a $(\pi_1, \pi_2)$-summand. Every
$(\pi_1, \pi_2)$-summand is obtained in this way from 
$\binom{k}{m_1} = \frac{k!}{m_1! m_2!}$ different $\pi$-summands, because $\pi$
is uniquely determined by the orderings of incoming legs ($\pi_1$) and
outgoing legs ($\pi_2$), and the choice of $m_1$ out of $k$ places for the
incoming legs. This proves our claim.

\smallskip
We next show that the map $\mu: \Bar_+ \bt \Bar_- \ra \Bar$ is an
isomorphism. At the level of
diagrams we have $\Dar_{+-}(\ast_1, \ast_2) \cong \Dar(\ast)$  via the 
pair of inverse maps
\[
\Dar_{+-}(\ast_1, \ast_2) \ni 
  \printname{InOutAstDgrm}
  \setlength{\unitlength}{0.35\standardunitlength}
  \begin{array}{c}  \hspace{-1.7mm}
    \raisebox{-8pt}{\input draws/InOutAstDgrm.tex }
      \hspace{-1.9mm}
  \end{array}
 
\quad \begin{array}{c}\stackrel{\mu}{\lra} \\ 
               \stackrel{\xi}{\lla} \end{array} \quad

  \printname{AstDgrm}
  \setlength{\unitlength}{0.35\standardunitlength}
  \begin{array}{c}  \hspace{-1.7mm}
    \raisebox{-8pt}{\input draws/AstDgrm.tex }
      \hspace{-1.9mm}
  \end{array}
 \in \Dar(\ast) \, ,
\]
where $\xi$ colors the univalent end vertices of incoming legs
with $\ast_1$ and those of outgoing legs with $\ast_2$. Both maps descend 
to $\Bar$-spaces (as there are no STU relations to worry about; this is what
makes the $\Bar$ situation easier), and our assertion is proved.

\smallskip
We can now use an elementary lemma from linear algebra.
\begin{lemma}
\lbl{lemma:IsoDeformToIso}
Let $V$ and $W$ be filtered linear spaces, with filtrations
$\{V_k\}$ and $\{W_k\}$, respectively. Let $\mu, \eta :V\ra W$
be filtration preserving linear maps. Assume that $\eta$ 
is an isomorphism and that
\[
\eta(x) - \mu(x) \in V_{k-1} \text{ for all } x\in V_k \quad (k\geq 1) \, .
\]
Then $\mu$ is also an isomorphism.
\end{lemma}

The proof of the proposition is completed by applying Lemma 
~\ref{lemma:IsoDeformToIso} with
$\mu: \Aar_+ \bt \Aar_- \ra \Aar$ and
$\eta = \upchi\circ\mu\circ(\upchi_+\bt\upchi_-)^{-1}$.
\end{proof}

\begin{exercise}
Prove Lemma ~\ref{lemma:IsoDeformToIso}. Observe that if 
$\mu, \eta, \eta^{-1}$ are all diagrammatic morphisms, then so is
$\mu^{-1}$.  
\end{exercise}

\subsubsection{\bf Definition of Verma modules.}
\label{__sect:VermaDef} 

Let $\mathbf{1}$ denote the trivial $1$-dimensional representation
of both $\g_+$ and $\g_-$. We consider Verma modules of $U(\tg)$ which are 
the induced representations $M_+ = \Ind_{\g_+}^{\tg}\mathbf{1}$, 
$M_- = \Ind_{\g_-}^{\tg}\mathbf{1}$. Recall that
$\Ind_{\g_{\pm}}^{\tg}\mathbf{1} = U(\tg) \ot_{U(\g_{\pm})} \mathbf{1}$
(see ~\cite[Section 20.3]{Humphreys:LieAlgs}). Therefore, $M_{\pm}$ are 
generated (as $U(\tg)$-modules) by vectors $1_{\pm}$ such that 
$\g_{\pm} 1_{\pm} = 0$. This leads us to the following diagrammatic
definition.
\begin{definition}
Let $\J_+$ be the subspace of $\Aar$ which is spanned by the images in 
$\Aar$ of diagrams from $\Dar(\ua)$ whose rightmost leg is incoming.
The {\em diagrammatic Verma module} $\VM_+$ is the quotient of $\Aar$ by 
$\J_+$.
Thus, $\VM_+$ is the quotient of $\Span_{\k}\{\Dar(\ua)\}$ by the usual
directed relations (Definition ~\ref{def:Aar}) and the additional RI relations
\[
	
  \printname{RightLegIn}
  \setlength{\unitlength}{0.3\standardunitlength}
  \begin{array}{c}  \hspace{-1.7mm}
    \raisebox{-8pt}{\input draws/RightLegIn.tex }
      \hspace{-1.9mm}
  \end{array}
 = 0 \, .
\]
We similarly define $\J_-$ to be the span of images of diagrams whose
rightmost leg is outgoing, and put $\VM_- = \Aar \left/ J_- \right.$. The
corresponding RO relations are
\[
	
  \printname{RightLegOut}
  \setlength{\unitlength}{0.3\standardunitlength}
  \begin{array}{c}  \hspace{-1.7mm}
    \raisebox{-8pt}{\input draws/RightLegOut.tex }
      \hspace{-1.9mm}
  \end{array}
 = 0 \, .
\]
We denote by $1_{\pm}$ the elements of $\VM_{\pm}$ corresponding to the
empty diagram.
\end{definition}

Since RI and RO relations are homogeneous, $\VM_{\pm}$ have the usual
grading by half the number of vertices. We denote by 
$p_{\pm}:\Aar \ra \VM_{\pm}$ the natural projections. When we view 
a diagram in $\Dar(\ua)$ as an element of $\VM_+$ we draw it in the form
\[

  \printname{MplusDgrm}
  \setlength{\unitlength}{0.3\standardunitlength}
  \begin{array}{c}  \hspace{-1.7mm}
    \raisebox{-8pt}{\input draws/MplusDgrm.tex }
      \hspace{-1.9mm}
  \end{array}
 \, ,
\]
and similarly for $\VM_-$. 

An RI diagram is mapped by $\Tartg$ to a linear combination of
elements of $U(\tg)$ of the form $u = yx, y\in U(\tg), x\in \g_+$; hence
$u 1_+ = 0$. Thus, the composition
\[
\begin{CD}
\Aar @>{\Tartg}>> U(\tg) @>{u \mapsto u 1_+}>> M_+
\end{CD}
\]
descends to a map $\Tartg:\VM_+ \ra M_+$. A twin argument establishes the
existence of a well-defined map $\Tartg:\VM_- \ra M_-$. 

The definitions of diagrammatic tensor products involving $\VM_{\pm}$,
such as $\Aar\bt\VM_+$ and $\VM_+\bt\VM_-$, should be clear. For example,
$\VM_+\bt\VM_-$ is the quotient of $\Aar\bt\Aar$ by 
$\J_+\bt\Aar + \Aar\bt\J_-$ , which is the span of (images of) 
diagrams having either an ingoing rightmost leg on the first skeleton 
interval, or an outgoing rightmost leg on the second one. The accompanying
projection map is denoted $p_+\bt p_- : \Aar\bt\Aar \ra \VM_+\bt\VM_-$.

\smallskip
We can now justify the module qualifications of $\VM_{\pm}$. The 
maps
\[
\begin{CD}
	\Aar\bt \Aar @>{\mu}>> \Aar @>{p_{\pm}}>> \VM_+ 
\end{CD}
\]
send $\Aar\bt\J_{\pm}$ to zero, because $\J_{\pm}$ are `diagrammatic
ideals' as illustrated (for $\J_+$) in the figure below:
\[
	\mu\left(
  \printname{AarJplusDgrm}
  \setlength{\unitlength}{0.3\standardunitlength}
  \begin{array}{c}  \hspace{-1.7mm}
    \raisebox{-8pt}{\input draws/AarJplusDgrm.tex }
      \hspace{-1.9mm}
  \end{array}
\right) 
	= 
  \printname{AarJplusProdDgrm}
  \setlength{\unitlength}{0.3\standardunitlength}
  \begin{array}{c}  \hspace{-1.7mm}
    \raisebox{-8pt}{\input draws/AarJplusProdDgrm.tex }
      \hspace{-1.9mm}
  \end{array}
 \, .
\]
Hence these maps  descend to $\Aar$-module structure
maps $\mu_{\pm}:\Aar\bt\VM_{\pm} \ra \VM_{\pm}$. The basic module
axiom, namely the commutativity of the diagram
\[
\begin{CD}
\Aar\bt\Aar\bt\VM_{\pm} @>{\id\bt\mu_{\pm}}>> \Aar\bt\VM_{\pm} \\
@V{\mu\bt 1}VV                              @VV{\mu_{\pm}}V \\
\Aar\bt\VM_{\pm}        @>{\mu_{\pm}}>>     \VM_{\pm} 
\end{CD} \, ,
\]
follows from the associativity of the diagrammatic product. The analogy of
the diagrammatic module structures to the algebraic ones is conveyed by
the commutative diagram
\[
\begin{CD}
\Aar\bt\VM_{\pm} @>{\mu_{\pm}}>> \VM_{\pm} \\
@V{\Tartg}VV               @VV{\Tartg}V \\
U(\tg)\ot M_{\pm}        @>>>    M_{\pm} 
\end{CD} \, .
\] 

\smallskip
It is clear that the product $\mu:\lowupper{\Aar}{\bt 2} \ra \Aar$ makes 
$\Aar$ itself into an $\Aar$-module (which is analogous to the left regular 
action of $U(\tg)$ on itself). The maps $p_{\pm}:\Aar\ra\VM_{\pm}$ are then 
morphisms of $\Aar$-modules (or {\em $\Aar$-maps}, for short), in the sense 
that the diagram
\[
\begin{CD}
\Aar\bt\Aar @>{\id\bt p_{\pm}}>> \Aar\bt\VM_{\pm} \\
@V{\mu}VV                      @VV{\mu_{\pm}}V \\
\Aar        @>{p_{\pm}}>>     \VM_{\pm}  
\end{CD}
\]
is commutative.

\subsubsection{\bf The PBW isomorphisms 
$p_{\pm}:\protect\Aar_{\mp}\ra \protect\VM_{\pm}$.}
\label{__sect:VermaPBW} 
The PBW isomorphism
$\mu: U(\g_+)\ot U(\g_-) \ra U(\tg)$ implies that the map
$U(\g_+) \ra M_- \, :\, u\mapsto u1_-$ is an isomorphism of linear spaces
(not of $U(\tg)$-modules). We state and prove the corresponding 
diagrammatic result.
\begin{proposition}
The maps $p_{\pm}:\Aar \ra \VM_{\pm}$ restrict to linear
isomorphisms $p_{\pm}:\Aar_{\mp}\ra \VM_{\pm}$.
\end{proposition}
\begin{proof}
We treat the case of $p_-:\Aar_+ \ra \VM_-$. We follow the pattern of proof of 
Proposition ~\ref{prop:PBW-A}, and freely use the objects employed there, or
their adaptations to the directed case, as indicated in the proof
of Proposition ~\ref{prop:PBW-Aar}. We filter $\Dar(\ua)$
by the number of legs, and denote the induced filtrations of $\Aar_+$ and  
$\VM_-$ by $(\Aar_+)_k$ and $(\VM_-)_k$.

\smallskip
We shall inductively construct a sequence of compatible maps 
$q_k:\Dar_k(\ua) \ra \Aar_+$ satisfying:
\begin{enumerate}[(Q1)]
\item
$q_k(E) = E$ for all $E\in (\Dar_+(\ua))_k$. \label{itm:qpEqId}
\item
$q_k(\Iar_k) = 0$, where $\Iar_k$ is the subspace of $\Dar_k(\ua)$ spanned by
NS, AS, IHX, and STU relations which belong to $\Dar_k(\ua)$ (that is, all the
diagrams participating in a relation have at most $k$ legs). \label{itm:qSTU} 
\item
$q_k((\J_-)_k) = 0$, where $(\J_-)_k$ is spanned by RO diagrams having
no more than $k$ legs. \label{itm:qRO}
\item
$p_-(q_k(D)) = p_-(D)$ for all $D\in\Dar_k(\ua)$. \label{itm:pqEqId}
\end{enumerate}
By Q\ref{itm:qSTU} and Q\ref{itm:qRO}, the maps $q_k$ descend to
compatible maps $q_k:(\VM_-)_k \ra \Aar_+$. The resulting map 
$q:\VM_- \ra \Aar_+$ is an inverse of $p_-:\Aar_+ \ra \VM_- $ by 
Q\ref{itm:qpEqId} and Q\ref{itm:pqEqId}.

\smallskip
For $k=0$ we set
\[
q_0\left(
  \printname{D0LineDgrm}
  \setlength{\unitlength}{0.2\standardunitlength}
  \begin{array}{c}  \hspace{-1.7mm}
    \raisebox{-8pt}{\input draws/D0LineDgrm.tex }
      \hspace{-1.9mm}
  \end{array}
\right) = 
  \printname{D0LineDgrm}
  \setlength{\unitlength}{0.2\standardunitlength}
  \begin{array}{c}  \hspace{-1.7mm}
    \raisebox{-8pt}{\input draws/D0LineDgrm.tex }
      \hspace{-1.9mm}
  \end{array}
 \, ,
\]
and for $k=1$
\[
q_1\left(
  \printname{D1LineDgrmIn}
  \setlength{\unitlength}{0.2\standardunitlength}
  \begin{array}{c}  \hspace{-1.7mm}
    \raisebox{-8pt}{\input draws/D1LineDgrmIn.tex }
      \hspace{-1.9mm}
  \end{array}
\right) = 
  \printname{D1LineDgrmIn}
  \setlength{\unitlength}{0.2\standardunitlength}
  \begin{array}{c}  \hspace{-1.7mm}
    \raisebox{-8pt}{\input draws/D1LineDgrmIn.tex }
      \hspace{-1.9mm}
  \end{array}

\quad ; \quad 
q_1\left(
  \printname{D1LineDgrmOut}
  \setlength{\unitlength}{0.2\standardunitlength}
  \begin{array}{c}  \hspace{-1.7mm}
    \raisebox{-8pt}{\input draws/D1LineDgrmOut.tex }
      \hspace{-1.9mm}
  \end{array}
\right) = 0 \, .
\]
Let $k>1$ and assume that a map $q_{k-1}:\Dar_{k-1}(\ua) \ra \Aar_+$ with the
desired properties was already constructed. Let $D\in \Dar_k(\ua)$ be 
a diagram with $k$ legs. It is convenient to call $D$ a {\em fully incoming} 
diagram if $D\in (\Dar_+(\ua))_k$. If $D$ is not fully incoming, then there
exists $\pi\in S_k$ such that $\pi D$ is an RO diagram ; we 
say that such a permutation is {\em admissible} for $D$. 

\smallskip
Recall that the element 
$\Gamma_D(\pi) \in \Aar_{k-1}$ is defined by
presenting $\pi$ as a product of transpositions, 
$\pi = u_{i_1}\cdots u_{i_m}$, and putting 
\begin{equation}
\label{eqn:RptDefGammaD}
\Gamma_D(\pi) 
= \sum_{l=1}^m s_{i_l} u_{i_{l+1}} \cdots u_{i_m} D\, ;
\end{equation} 
it was shown to be independent of the presentation of $\pi$ being used.
We now mod out the target space of $\Gamma_D$ by $(\J_-)_{k-1}$ and regard it
for the rest of this proof as a map $\Gamma_D : S_k  \ra (\VM_-)_{k-1}$.

\smallskip
The crux of the proof is showing that if $D$ is not fully incoming, then 
\begin{equation}
\label{eqn:AdmIndepGammaD}
\Gamma_D(\pi_1) = \Gamma_D(\pi_2) \
\text{ in $(\VM_-)_{k-1}$, when both $\pi_1$, $\pi_2$ are admissible for $D$} 
\, .
\end{equation}
Granting this fact for the moment, we denote by $\pi_D$ any permutation
admissible for $D$ (unless $D$ is fully incoming), and define
\begin{equation}
\label{eqn:DefQk}
q_k(D) = \left\{
\begin{array}{cl}
D & \text{if $D$ is fully incoming}, \\
\\
q_{k-1}(\Gamma_D(\pi_D)) & \text{otherwise}.
\end{array}
\right.
\end{equation}
(Q\ref{itm:qSTU} and Q\ref{itm:qRO} for $q_{k-1}$ make this a legitimate 
definition.)
For a diagram $D$ with less than $k$ legs, we define $q_k(D) = q_{k-1}(D)$.

We note that that if $D$ is not fully incoming, then by 
\eqref{eqn:TelescopicPermGammaD}
\[
q_k(D) - q_k(\pi_D D) = q_{k-1}(\Gamma_D(\pi_D)) \, . 
\]
Since $\pi_D D$ is an RO diagram, we are led to the definition in 
\eqref{eqn:DefQk}. 

\smallskip
{\em Proof of Q\ref{itm:qpEqId}:} Immediate for $E\in (\Dar_{+})_k$ with exactly
$k$ legs, and follows from the induction hypothesis about $q_{k-1}$
for $E$ with less than $k$ legs. We omit similar remarks about diagrams
having less than $k$ legs in the proofs of the remaining properties.

\smallskip
{\em Proof of Q\ref{itm:qSTU}:} As usual, the only difficulty is verifying that
STU relations are sent to zero, which amounts to
\[
q_k(D - u_i D) = q_k(s_i D) \quad \text{for } \, D\in \Dar_k(\ua)\, , 
1\leq i < k \, .
\]
For a fully incoming diagram $D$, this is just an incoming STU relation. 
Otherwise, let $\pi_D$ be an admissible permutation for $D$. Then
$\pi_D u_i$ is an admissible permutation for $u_i D$. Using
\eqref{eqn:AdmIndepGammaD} and the cocycle property \eqref{eqn:CocycleGammaD}, 
we obtain
\begin{align*}
q_k(D - u_i D) & = q_{k-1}(\Gamma_D(\pi_D)) 
    - q_{k-1}(\Gamma_{u_i D}(\pi_D u_i)) \\
    & =  q_{k-1}(\Gamma_D(\pi_D)) 
    - (q_{k-1}(\Gamma_{u_i u_i D}(\pi_D)) + q_{k-1}(\Gamma_{u_i D}(u_i))) \\
    & = - q_{k-1}(\Gamma_{u_i D}(u_i)) = - q_{k-1}(s_i u_i D) \\
    & = q_{k-1}(s_i D) = q_k(s_i D) \, ,
\end{align*}
where the last three equalities follow from \eqref{eqn:RptDefGammaD},
an AS relation, and the definition of $q_k$, respectively.
 
\smallskip
{\em Proof of Q\ref{itm:qRO}:} If $D$ is an RO diagram with $k$ legs, then we 
can take $\pi_D=e$, the identity permutation. Since $\Gamma_D(e)=0$ 
(presenting $e$ as an empty product of transpositions), we get
$q_k(D) = 0$. 

\smallskip
{\em Proof of Q\ref{itm:pqEqId}:} This is clear for a fully incoming  
diagram. If $D$ is not fully incoming, then by Q\ref{itm:pqEqId} for 
$q_{k-1}$ we have
\[
p_-(q_k(D)) = (p_- \circ q_{k-1})(\Gamma_D(\pi_D)) = p_-(\Gamma_D(\pi_D)) \, . 
\]
By \eqref{eqn:TelescopicPermGammaD}, 
\[
p_-(\Gamma_D(\pi_D)) = p_-(D - \pi_D D) = p_-(D) \, ,
\]
because $\pi_D D$ is an RO diagram. Thus $p_-(q_k(D)) = p_-(D)$, as required. 

\smallskip
It remains to prove \eqref{eqn:AdmIndepGammaD}. We embed $S_{k-1}$ in
$S_k$ as the subgroup of permutations fixing $k$. We observe that if $D$ is an
RO diagram and $\rho\in S_{k-1}$, then $\Gamma_D(\rho)$ is a linear
combination of RO diagrams, and is thus zero in $(\VM_-)_{k-1}$. To establish 
this fact, present $\rho$ as a product of transpositions from 
$\{u_1, \ldots, u_{k-2}\}$. Looking at \eqref{eqn:RptDefGammaD}, we see that 
none of the glueing operators in $\Gamma_D(\rho)$ alters the rightmost leg of 
$D$, hence all terms are RO diagrams. 

Let $\pi_1$, $\pi_2$
be two admissible permutations for a diagram $D\in \Dar_k(\ua)$ with $k$
legs (such that at least one of the legs is outgoing). There are two
cases:
\begin{enumerate}
\item
$\pi_1$ and $\pi_2$ move the same (outgoing) leg of $D$ to the right end. In
this case,
\[
	\pi_2 = \rho \pi_1 \ \text{ for some } \rho\in S_{k-1} \, .
\]
\item
$\pi_1$ and $\pi_2$ move different legs of $D$ to the right end, say
$\pi_1(i_1) = k$ and $\pi_2(i_2) = k $, where $i_1, i_2 < k$ and 
$i_1 \neq i_2$. Pick $\rho_1 \in S_{k-1}$ such that 
$\rho_1\pi_1(i_2) = k-1$. Then
$\pi_2 \pi_1^{-1} \rho_1^{-1} u_{k-1}(k) = k$, so that 
\[
  \pi_2 = \rho_2 u_{k-1} \rho_1 \pi_1 \quad (\rho_1, \rho_2 \in S_{k-1}) \, .
\] 
\end{enumerate}

In the first case we find, using \eqref{eqn:CocycleGammaD},
\begin{align*}
\Gamma_D(\pi_2) & = \Gamma_{\pi_1 D}(\rho) + \Gamma_D(\pi_1)  \\
  & = \Gamma_D(\pi_1) \quad \text{(since $\pi_1 D$ is an RO diagram)} \, . 
\end{align*}
In the second case,
\[
\Gamma_D(\pi_2) = \Gamma_{u_{k-1}\rho_1\pi_1 D}(\rho_2)
	+ \Gamma_D(u_{k-1}\rho_1\pi_1) \, .
\]
Since $u_{k-1}\rho_1\pi_1 D$ is an RO diagram, the first term on the right 
hand side is zero. Hence
\[
\Gamma_D(\pi_2) = \Gamma_{\rho_1\pi_1 D}(u_{k-1}) + \Gamma_D(\rho_1\pi_1) \, .
\]
By \eqref{eqn:RptDefGammaD}, 
$\Gamma_{\rho_1\pi_1 D}(u_{k-1}) = s_{k-1} (\rho_1\pi_1 D)$. But
$\rho_1$ was chosen in such a way that the two rightmost legs of 
$\rho_1\pi_1 D$ are outgoing. Therefore, applying
$s_{k-1}$ results in an RO diagram, which is equivalent to zero. Thus,
\[
\Gamma_D(\pi_2) = \Gamma_D(\rho_1 \pi_1) \, ,
\]  
which is equal to $\Gamma_D(\pi_1)$ by the argument of the first case.
\end{proof}

As with other PBW maps, the proof shows that 
the maps $p_{\pm}:\Aar_{\mp} \ra \VM_{\pm}$ can be inverted 
`skeleton by skeleton'. In particular, the map
$p_+\bt p_- : \Aar_- \bt \Aar_+ \ra \VM_+ \bt \VM_-$ is an isomorphism.

\subsubsection{\bf The map $\phi:\protect\Aar \ra 
			\protect\VM_+ \bt \protect\VM_-$.} 
\label{__sect:phi}

The tensor product of two $U(\tg)$-modules $V, W$ is given the structure
of a $U(\tg)$-module with the help of the coproduct of $U(\tg)$:
\[
  u(v\ot w) = \Delta(u)(v\ot w)\, , \quad u\in U(\g) \ , v\in V, w\in W \, .
\] 
We can imitate this construction without difficulty for the diagrammatic
tensor product of $\Aar$-modules. In particular, $\VM_+\bt \VM_-$ is
an $\Aar$-module; the module map is the composition
\[
\begin{CD}
\Aar\bt\VM_+\bt\VM_-  @>{\Delta_{\ua_1}}>> 
\Aar\bt\Aar\bt\VM_+\bt\VM_- @>{\mu_{\ua_1, \ua_3}\bt \mu_{\ua_2, \ua_4}}>>
\VM_+\bt\VM_- \, ,
\end{CD}
\]
which can be pictured as
\[

  \printname{AarMplusMminusDgrm}
  \setlength{\unitlength}{0.4\standardunitlength}
  \begin{array}{c}  \hspace{-1.7mm}
    \raisebox{-8pt}{\input draws/AarMplusMminusDgrm.tex }
      \hspace{-1.9mm}
  \end{array}
 \longmapsto

  \printname{CoprodMplusMminusDgrm}
  \setlength{\unitlength}{0.4\standardunitlength}
  \begin{array}{c}  \hspace{-1.7mm}
    \raisebox{-8pt}{\input draws/CoprodMplusMminusDgrm.tex }
      \hspace{-1.9mm}
  \end{array}
 \, .
\]

\begin{definition}
\lbl{def:phi}
The map $\phi:\Aar \ra \VM_+ \bt \VM_-$ is defined as the composition
\[
\begin{CD}
\Aar @>{\Delta}>> \Aar\bt\Aar @>{p_+\bt p_-}>> \VM_+ \bt \VM_- \, .
\end{CD}
\]
\end{definition}

\begin{proposition}
$\phi$ is an isomorphism of $\Aar$-modules.
\end{proposition}
\begin{proof}
That $\phi$ is an $\Aar$-map, namely
\[
  \phi(DE) = \Delta(D) \phi(E) \, \text{ for } D, E \in \Aar \, ,
\]
is immediate from the compatibility of the coproduct of $\Aar$ with its 
product (see 
\eqref{eqn:A-coprod-is-mult}). To prove
that $\phi$ is an isomorphism, we use the method of proof of Proposition
~\ref{prop:AarPlus-bt-AarMinusEqAar}: We show that $\phi$ coincides, up to 
lower filtration levels, with the composition of isomorphisms
\[
\begin{CD}
\Aar @>{\mu^{-1}}>> \Aar_-\bt\Aar_+ @>{p_+\bt p_-}>> \VM_+ \bt \VM_- \, ,
\end{CD}
\]
that is
\begin{equation}
\label{eqn:PhiDeformToIso}
\phi(D) - (p_+\bt p_-)\mu^{-1}(D) \in (\VM_+\bt\VM_-)_{k-1} 
\quad \text{ for } D\in \Aar_k \, .
\end{equation}
We can then apply Lemma ~\ref{lemma:IsoDeformToIso} to finish the proof.

\smallskip
To establish \eqref{eqn:PhiDeformToIso}, we first claim that
for a diagram $D\in\Aar_k$, $\mu^{-1}(D)$ is equal modulo 
$(\Aar_-\bt\Aar_+)_{k-1}$ to the diagram obtained from $D$ by placing
its outgoing legs on the first skeleton interval and placing
its incoming legs on the second.
For example,
\[
\mu^{-1}\left(
  \printname{MixedDgrm2}
  \setlength{\unitlength}{0.3\standardunitlength}
  \begin{array}{c}  \hspace{-1.7mm}
    \raisebox{-8pt}{\input draws/MixedDgrm2.tex }
      \hspace{-1.9mm}
  \end{array}
\right) -  

  \printname{MuInvMixedDgrm2}
  \setlength{\unitlength}{0.3\standardunitlength}
  \begin{array}{c}  \hspace{-1.7mm}
    \raisebox{-8pt}{\input draws/MuInvMixedDgrm2.tex }
      \hspace{-1.9mm}
  \end{array}
 \, \in (\Aar_-\bt\Aar_+)_4 \, .
\]
This is because when applying $\mu$ to the diagram on the right hand side we 
get a diagram of the form $\pi D$, for some $\pi \in S_k$, and we know that
$D - \pi D \in \Aar_{k-1}$ (by the directed version of 
\eqref{eqn:TelescopicPermGammaD}). Thus, modulo $(\VM_+\bt\VM_-)_{k-1}$ we
have
\begin{equation}
\label{eqn:PhiTopLevel}
(p_+\bt p_-)\mu^{-1}\left(
  \printname{MixedDgrm2}
  \setlength{\unitlength}{0.3\standardunitlength}
  \begin{array}{c}  \hspace{-1.7mm}
    \raisebox{-8pt}{\input draws/MixedDgrm2.tex }
      \hspace{-1.9mm}
  \end{array}
\right) = 

  \printname{PhiTopDgrm}
  \setlength{\unitlength}{0.3\standardunitlength}
  \begin{array}{c}  \hspace{-1.7mm}
    \raisebox{-8pt}{\input draws/PhiTopDgrm.tex }
      \hspace{-1.9mm}
  \end{array}
 \ .
\end{equation}

On the other hand, $\phi(D)$ is the sum over splittings of the legs of
$D$ among the $\VM_+$ and $\VM_-$ skeletons. All splittings in which at least 
one outgoing leg is sent to the $\VM_-$ skeleton belong to
$(\VM_+\bt\VM_-)_{k-1}$: We can bring such a splitting to become zero 
(in $\VM_+\bt\VM_-$)
by a sequence of STU moves pushing the outgoing leg to the right, and the 
price for these moves is a linear combination of diagrams with less than
$k$ legs. Arguing similarly about an incoming leg sent to the $\VM_+$ skeleton,
we find that the only top level contribution to $\phi(D)$ comes from
the unique splitting in which all outgoing legs are sent to
the  $\VM_+$ skeleton and all incoming legs are sent to the $\VM_-$
skeleton. This is precisely the diagram as in the right hand side of
\eqref{eqn:PhiTopLevel}.
\end{proof}

\subsubsection{\bf The maps $i_{\pm}:\protect \VM_{\pm} \ra 
			\protect\VM_{\pm} \bt \protect\VM_{\pm}$.} 
Two more $\Aar$-maps we will use below are the maps  
$i_{\pm}:\VM_{\pm} \ra \VM_{\pm} \bt \VM_{\pm}$. The map $i_+$ is defined as 
follows:
\[
   i_+(p_+(D)) = (p_+\bt p_+)(\Delta(D)) \, , \quad D\in \Aar \, . 
\]
It is immediate to check that $i_+$ is well-defined and is an $\Aar$-map,
namely
\[
   i_+(D E) = \Delta(D) i_+(E) \, , 
\]
for all $D\in \Aar, E\in \VM_+$. Similar definition and properties are 
shared by the map $i_-$.

\subsection{The diagrammatic Etingof-Kazhdan twist} 
\label{_sect:EKtwist}

Following ~\cite{EtingofKazhdan:QuantizationI}, we construct an element
$J\in \lowupper{\Aar}{\bt 2}$ as follows. Put
\[
\wtJ = 
   \ZK\Big(\raisebox{2pt}{$\hs{-3pt}
  \printname{Jbraid}
  \setlength{\unitlength}{0.3\standardunitlength}
  \begin{array}{c}  \hspace{-1.7mm}
    \raisebox{-8pt}{\input draws/Jbraid.tex }
      \hspace{-1.9mm}
  \end{array}
\hs{-6pt}$}\Big) 
   \in \A^{\bt 4} \, ,
\]
and use the same symbol to denote  $\iota(\wtJ) \in \lowupper{\Aar}{\bt 4}$. 
Set
\[
  J = (\phi^{-1}\bt\phi^{-1}) (p_+\bt p_-\bt p_+\bt p_-)(\wtJ)  
  \in \lowupper{\Aar}{\bt 2} \, .
\] 

By expressing the parenthesized braid in the formula for $\wtJ$ as the
composition of generating associativity and braiding morphisms, and
recalling \eqref{eqn:HorPhi}, we see that
\[
\wtJ = 
  \printname{A4hUnit}
  \setlength{\unitlength}{0.25\standardunitlength}
  \begin{array}{c}  \hspace{-1.7mm}
    \raisebox{-8pt}{\input draws/A4hUnit.tex }
      \hspace{-1.9mm}
  \end{array}
 + \frac{_1}{^2}
  \printname{A4hT23}
  \setlength{\unitlength}{0.25\standardunitlength}
  \begin{array}{c}  \hspace{-1.7mm}
    \raisebox{-8pt}{\input draws/A4hT23.tex }
      \hspace{-1.9mm}
  \end{array}
 + 
            \text{ terms of degree \ $>1$} \, ,
\] 
and it is then easy to check that
\begin{equation}
\label{eqn:JuptoDeg1}
J = 
  \printname{A2hUnit}
  \setlength{\unitlength}{0.25\standardunitlength}
  \begin{array}{c}  \hspace{-1.7mm}
    \raisebox{-8pt}{\input draws/A2hUnit.tex }
      \hspace{-1.9mm}
  \end{array}
 + \frac{_1}{^2}
  \printname{LeftArrow}
  \setlength{\unitlength}{0.25\standardunitlength}
  \begin{array}{c}  \hspace{-1.7mm}
    \raisebox{-8pt}{\input draws/LeftArrow.tex }
      \hspace{-1.9mm}
  \end{array}
 + 
            \text{ terms of degree \ $>1$} \, .
\end{equation}

\begin{theorem}
\lbl{thm:EKtwist}
Twisting $\Aar_{\kz}$ by $F = J^{-1}$ yields a quasitriangular {\em Hopf} 
structure on $\Aar$, denoted $\Aar_{\ek}$, whose $R$-matrix is a 
quantization of the diagrammatic $r$-matrix: 
\begin{equation}
\label{eqn:RekUptoDeg1}
	R_{\ek} = 
  \printname{A2hUnit}
  \setlength{\unitlength}{0.25\standardunitlength}
  \begin{array}{c}  \hspace{-1.7mm}
    \raisebox{-8pt}{\input draws/A2hUnit.tex }
      \hspace{-1.9mm}
  \end{array}
 + 
  \printname{LeftArrow}
  \setlength{\unitlength}{0.25\standardunitlength}
  \begin{array}{c}  \hspace{-1.7mm}
    \raisebox{-8pt}{\input draws/LeftArrow.tex }
      \hspace{-1.9mm}
  \end{array}
 + 
            \text{ terms of degree \ $>1$} \, .
\end{equation}
\end{theorem}
\begin{proof}
By the transformation formula \eqref{eqn:TwistPhi} of  
the associator under twisting, we have to show that 
\begin{equation}
\label{eqn:TwistToHopf}
     \Phi \cdot (\Delta\bt\id)(J) \cdot J^{12} =
	(\id\bt\Delta)(J) \cdot J^{23} \, ,
\end{equation}
in order to prove that $\Aar_{\ek}$ is Hopf (that is, its coproduct
is coassociative). 

Our starting point is the topological equivalence
of the following two parenthesized braids:
\[
   
  \printname{Jproof1}
  \setlength{\unitlength}{0.6\standardunitlength}
  \begin{array}{c}  \hspace{-1.7mm}
    \raisebox{-8pt}{\input draws/Jproof1.tex }
      \hspace{-1.9mm}
  \end{array}
 = 
  \printname{Jproof2}
  \setlength{\unitlength}{0.6\standardunitlength}
  \begin{array}{c}  \hspace{-1.7mm}
    \raisebox{-8pt}{\input draws/Jproof2.tex }
      \hspace{-1.9mm}
  \end{array}
 \, .
\]
Applying $\ZK$ to both sides, using Definition ~\ref{def:ZcatNAT} and
the definition of $\wtJ$, we find that 
\[
   \Delta^{\bt 3}(\Phi) \cdot (\wtJ\bt 1^{\bt 2}) \cdot
    ((\Delta^{\bt 2} \bt \id^{\bt 2})(\wtJ))^{(23)} = \\
    (1^{\bt 2}\bt \wtJ) \cdot ((\id^{\bt 2} \bt \Delta^{\bt 2})(\wtJ))^{(45)}  
    \cdot (\Phi^{135} \bt \Phi^{246})
\]
in $\A^{\bt 6}$. Here we have used the following notation: for an element 
$X$ in some (diagrammatic) tensor product, $X^{(ij)}$ is the element 
obtained from $X$ by transposing its $i$th and $j$th factors.

We map the last equality into $\lowupper{\Aar}{\bt 6}$ (using 
$\iota: \A^{\bt 6} \ra \lowupper{\Aar}{\bt 6}$) and then apply the map
$(p_+ \bt p_-)^{\bt 3} : \lowupper{\Aar}{\bt 6} \ra (\VM_+\bt \VM_-)^{\bt 3}$.
Since $p_{\pm}$ are $\Aar$-maps we get:
\begin{multline*}
   \Delta^{\bt 3}(\Phi) \cdot (p_+ \bt p_-)^{\bt 3} ((\wtJ\bt 1^{\bt 2}) 
    \cdot ((\Delta^{\bt 2} \bt \id^{\bt 2})(\wtJ))^{(23)}) = \\
    (1^{\bt 2}\bt \wtJ) \cdot ((\id^{\bt 2} \bt \Delta^{\bt 2})(\wtJ))^{(45)}  
    \cdot (p_+ \bt p_-)^{\bt 3} (\Phi^{135} \bt \Phi^{246}) \, .
\end{multline*}
The crucial observation now (compare 
~\cite[Lemma 2.3]{EtingofKazhdan:QuantizationI}) is 
that
\[
(p_+ \bt p_-)^{\bt 3} (\Phi^{135} \bt \Phi^{246}) = (1_+ \bt 1_-)^{\bt 3} \, .
\]
This follows from the fact that $\Phi$ is horizontal and the definitions of 
$\iota$ and $p_{\pm}$. We thus conclude that the following equality holds
in $(\VM_+\bt \VM_-)^{\bt 3}$:
\begin{multline}
\label{eqn:wtJbasic}
    \Delta^{\bt 3}(\Phi) \cdot (p_+ \bt p_-)^{\bt 3} ((\wtJ\bt 1^{\bt 2}) 
    \cdot ((\Delta^{\bt 2} \bt \id^{\bt 2})(\wtJ))^{(23)}) = \\
    (p_+ \bt p_-)^{\bt 3}
    ((1^{\bt 2}\bt \wtJ) \cdot ((\id^{\bt 2} \bt \Delta^{\bt 2})(\wtJ))^{(45)})
    \, .
\end{multline}

\smallskip
We introduce an auxiliary map 
$\psi: \VM_+ \bt \VM_- \ra (\VM_+ \bt \VM_-)^{\bt 2}$, which is defined
by
\[
    \psi(E) = (\wtJ^{(23)} \cdot (i_+ \bt i_-)(E))^{(23)} \, , \quad 
     E \in \VM_+ \bt \VM_- \, .
\]
Note that 
\begin{equation}
\label{eqn:psi1Plus1Minus}
    \psi(1_+ \bt 1_-)  = \wtJ \cdot (1_+ \bt 1_-)^{\bt 2} = \phi^{\bt 2}(J) \, .
\end{equation}
We claim that $\psi$ is an $\Aar$-map: Since $i_{\pm}$ are $\Aar$-maps,
\[
    \psi(\Delta(D) \cdot E) = 
	(\wtJ^{(23)} \cdot \Delta^{(4)}(D) \cdot (i_+ \bt i_-)(E))^{(23)} \, .
\]
for $D\in \Aar, E\in \VM_+ \bt \VM_-$. Since $\wtJ$ is an element coming 
from an $\A$-space, we may use Proposition ~\ref{prop:AarPwrInvar} and
obtain
\begin{align*}
    \psi(\Delta(D) \cdot E) & = 
	(\Delta^{(4)}(D) \cdot \wtJ^{23} \cdot (i_+ \bt i_-)(E))^{(23)} \\
	& = \Delta^{(4)}(D) \cdot (\wtJ^{23} \cdot (i_+ \bt i_-)(E))^{(23)}  \\
        & = \Delta^{(4)}(D) \cdot \psi(E) \, ,
\end{align*}
as required.

\smallskip
We next compute $(\psi\bt \id)\psi(1_+ \bt 1_-)$ and 
$(\id\bt\psi)\psi(1_+ \bt 1_-)$ in two ways. First, by the definition of 
$\psi$,
\begin{align*}
  (\id\bt\psi)\psi(1_+ \bt 1_-) & = 
    ((1^{\bt 2} \bt \wtJ^{(23)}) \cdot	
     (\id^{\bt 2} \bt i_+ \bt i_-)(\wtJ \cdot (1_+\bt 1_-)^{\bt 2}))^{(45)} \\
   & =  (p_+ \bt p_-)^{\bt 3}
   ((1^{\bt 2} \bt \wtJ) \cdot ((\id^{\bt 2} \bt \Delta^{\bt 2})(\wtJ))^{(45)})
   \, ,
\end{align*}
and similarly
\[
   (\psi\bt \id)\psi(1_+ \bt 1_-) = (p_+ \bt p_-)^{\bt 3} 
    ((\wtJ\bt 1^{\bt 2}) \cdot ((\Delta^{\bt 2} \bt \id^{\bt 2})(\wtJ))^{(23)})
   \, .
\]
Comparing with \eqref{eqn:wtJbasic}, we see that
\begin{equation}
\label{eqn:PsiCoassoc}
 \Delta^{\bt 3}(\Phi) \cdot (\psi\bt \id)\psi(1_+ \bt 1_-) = 
     (\id\bt\psi)\psi(1_+ \bt 1_-) \, .    
\end{equation}

Second, by \eqref{eqn:psi1Plus1Minus} we have
\begin{align*}
  (\psi\bt \id)\psi(1_+ \bt 1_-) & = (\psi\bt \id) \phi^{\bt 2}(J)  \\
	& = (\id^{\bt 2} \bt \phi) (\psi\phi \bt \id) (J) \, .
\end{align*}
Since both $\psi$ and $\phi$ are $\Aar$-maps,  
\[
  \psi\phi(D) = \Delta^{(4)}(D) \cdot \psi\phi(1) 
	=  \Delta^{(4)}(D) \cdot \phi^{\bt 2}(J) 
        = \phi^{\bt 2}(\Delta(D) \cdot J) \, ,  
\]
for all $D\in \Aar$.
Therefore, 
\begin{align*}
   (\id^{\bt 2} \bt \phi) (\psi\phi \bt \id) (J) & =
    (\id^{\bt 2} \bt \phi) (\phi^{\bt 2} \bt \id) 
    ((\Delta \bt \id)(J) \cdot J^{12}) \\
   & = \phi^{\bt 3} ((\Delta \bt \id)(J) \cdot J^{12}) 
    \, ,
\end{align*}
so that
\[
   (\psi\bt \id)\psi(1_+ \bt 1_-) 
   = \phi^{\bt 3} ((\Delta \bt \id)(J) \cdot J^{12})
   \, .
\]
A similar computation shows that
\[
   (\id\bt\psi)\psi(1_+ \bt 1_-) 
   = \phi^{\bt 3} ((\id \bt \Delta)(J) \cdot J^{23})
   \, .
\]
Together with \eqref{eqn:PsiCoassoc} we have
\[
   \Delta^{\bt 3}(\Phi) \cdot \phi^{\bt 3} ((\Delta \bt \id)(J) \cdot J^{12}) =
      \phi^{\bt 3} ((\id \bt \Delta)(J) \cdot J^{23}) \, ,  
\]
and using one more time that $\phi$ is an $\Aar$-map, we finally get
\[
    \phi^{\bt 3} (\Phi \cdot (\Delta \bt \id)(J) \cdot J^{12}) = 
       \phi^{\bt 3} ((\id \bt \Delta)(J) \cdot J^{23}) \, .
\]
The proof of \eqref{eqn:TwistToHopf} is now finished by applying the map
$(\phi^{-1})^{\bt 3} : (\VM_+ \bt \VM_-)^{\bt 3} \ra \lowupper{\Aar}{\bt 3}$.

\smallskip
The $R$-matrix of $\Aar_{\ek}$ is given, according to \eqref{eqn:TwistR},
by the formula
\[
   R_{\ek} = (J^{-1})^{21} \cdot 
               \exp(\frac{_1}{^2} 
  \printname{BCup}
  \setlength{\unitlength}{0.25\standardunitlength}
  \begin{array}{c}  \hspace{-1.7mm}
    \raisebox{-8pt}{\input draws/BCup.tex }
      \hspace{-1.9mm}
  \end{array}
) \cdot J \, .
\]    
The expansion \eqref{eqn:RekUptoDeg1} of $R_{\ek}$ is then an immediate 
consequence of \eqref{eqn:JuptoDeg1}.
\end{proof}

\subsection{The invariant $\Z_{\ek}$}
\label{_sect:Zek}

We denote by $\Z_{\ek}$ the invariant corresponding to $\Aar_{\ek}$ by 
the procedure of Section ~\ref{__sect:Z_H}. Since the invariant corresponding
to $\Aar_{\kz}$ is the directed Kontsevich integral $\iota\circ\ZK$, we 
deduce from Corollary ~\ref{corr:ZlinkIndepTwist} the following result.

\begin{theorem}
\lbl{thm:ZekEqualIotaZK}
For any framed link $L$,
\[
	\Z_{\ek}(L) = \iota\circ\ZK(L) \, .
\]
\end{theorem}

Note that if Conjecture ~\ref{conj:IotaIsInjective} (that $\iota$ is injective)
is true, then Theorem ~\ref{thm:ZekEqualIotaZK} implies that the framed
link invariant $\Z_{\ek}$ is {\em equivalent} to the Kontsevich integral. Since 
Conjecture ~\ref{conj:IotaIsInjective} is true at least at the level of
complex semisimple Lie algebras, it is safe to designate $\Z_{\ek}$
as an invariant `closely related' to the Kontsevich integral. Theorem 
~\ref{thm:ZekEqualIotaZK} thus shows that such an invariant can be
constructed from an $R$-matrix satisfying the QYBE, and this construction
captures the first nice property (RT1) of the Reshetikhin-Turaev invariants 
(see Section ~\ref{__sect:DgrmRTquest}).

\smallskip   
Admittedly, our construction of the $R$-matrix of $\Aar_{\ek}$ has little
success in mimicking the property (RT2), as the construction of the
diagrammatic Etingof-Kazhdan twist is far from being explicit. First,
it involves a specific value of the Kontsevich integral
(the element 
$\ZK(\raisebox{4pt}{$\hs{-2pt}
  \printname{Jbraid}
  \setlength{\unitlength}{0.2\standardunitlength}
  \begin{array}{c}  \hspace{-1.7mm}
    \raisebox{-8pt}{\input draws/Jbraid.tex }
      \hspace{-1.9mm}
  \end{array}
$\hs{-5pt}})\in \A^{\bt 4})$, 
which we can write down only by using the Drinfel'd associator $\Phi$. 
Second, our recipe for computing $\phi^{-1}: \VM_+ \bt \VM_- \ra \Aar$ is 
quite complicated. 

We can only assert, thanks to the non-constructive existence result of 
Theorem ~\ref{thm:RationalPhi}, that the invariant $\Z_{\ek}$ may 
be constructed based on an $R$-matrix having {\em rational} coefficients.

\smallskip
We turn to consider (RT3), namely the existence of a simple formula
for the invariant of the unknot. By \eqref{eqn:Zunknot}, the computation
of $\Z_{\ek}(\text{unknot})$ requires the knowledge of the $u$-element and
ribbon element of $\Aar_{\ek}$. (The $\alpha$ and $\beta$ elements are 
equal to $1$, since $\Aar_{\ek}$ is coassociative.) The latter element is 
easy, as by Proposition ~\ref{prop:RibbonTwist} it is equal to the ribbon
element of $\Aar_{\kz}$. Thus (see Proposition ~\ref{prop:AarKZ}), 
\begin{equation}
\label{eqn:EKribbon}
v_{\ek} = 
  \iota (\exp(-\frac{_1}{^2} \raisebox{4pt}{$
  \printname{HalfCirc}
  \setlength{\unitlength}{0.2\standardunitlength}
  \begin{array}{c}  \hspace{-1.7mm}
    \raisebox{-8pt}{\input draws/HalfCirc.tex }
      \hspace{-1.9mm}
  \end{array}
$})) =
  \exp(-\frac{_1}{^2}(\LHC + \RHC)) \, .  
\end{equation}
 
For the $u$-element, we make the following conjecture. Recall from 
Example ~\ref{exmpl:DgrmHalfSum} that the
diagrammatic half-sum of positive roots is the element
\begin{equation}
\label{eqn:HalfSumTwoForms}
  \varrho = \frac{_1}{^2}(\LHC - \RHC) = \frac{_1}{^2}(\TPD + \TPU) = 
	\frac{_1}{^2} \TP \in \Aar \, .
\end{equation}
Here we have introduced the convention that a diagram having both directed and 
undirected arcs denotes the sum of all directed Jacobi diagrams obtained by 
all possible assignments of directions to the undirected arcs.

\begin{conjecture}
\label{conj:Uek}
The $u$-element of $\Aar_{\ek}$ is given by
\begin{equation}
\label{eqn:Uek}
u_{\ek} = \exp(\varrho) v_{\ek} = 
   \exp(- \RHC) \, . 
\end{equation}
\end{conjecture}
(The fact that $v_{\ek}$ is the exponential of an element coming from $\A$,
which by Proposition ~\ref{prop:ImgIotaInvar} is central in $\Aar$, was
used when evaluating $\exp(\varrho) v_{\ek}$.)

\smallskip
The correctness of Conjecture ~\ref{conj:Uek} would imply, in view
of \eqref{eqn:Zunknot}, the following conjectural formulas for the
value in $\Aar(\bigcirc)$ of $\Z_{\ek}$ on the unknot:
\begin{conjecture}
\label{conj:ZekUnknot}
\begin{equation}
\label{eqn:ZekUnknot}
\Z_{\ek}(\text{unknot}) = \tr(\exp(\pm\frac{_1}{^2}\TP))  
\end{equation}
\end{conjecture} 

\begin{theorem}
\lbl{thm:ZekUnknotLieAlg}
Conjecture ~\ref{conj:ZekUnknot} is true at the level of complex semisimple
Lie algebras.
\end{theorem}
\begin{proof}
Let $\g$ be a complex semisimple Lie algebra. We have to prove that 
applying the map $\Targh:\Aar(\bigcirc) \ra U(\g)_{\g}\hb$ to both sides
of \eqref{eqn:ZekUnknot} yields an equality.

From Theorem ~\ref{thm:ZekEqualIotaZK} and \eqref{eqn:RTgEqualKZg}, we 
conclude that $\Targh(\Z_{\ek}(L)) = \RTg(L)$, for any framed link $L$. It now
follows from \eqref{eqn:RTgUnknot} and Example ~\ref{exmpl:DgrmHalfSum}  
that 
\begin{align*}
 \Targh(\Z_{\ek}(\text{unknot})) & = \RTg(\text{unknot})  \\
                                 & = \tr(e^{\pm \hbar \varrho_{\g}}) \\ 
                                 & = \tr(e^{\pm \Targh(\varrho)})  \\
                                 & = \Targh(e^{\pm \varrho}) \, , 
\end{align*}
as required.
\end{proof}

Since the maps $\Targh$ separate the points of $\A = \A(\bigcirc)$ up to 
degree 12 (see the discussion after Conjecture ~\ref{conj:IotaIsInjective}
in Section ~\ref{__sect:ImgKerIota}), it follows that Conjecture
~\ref{conj:ZekUnknot} holds up to that degree. It is nevertheless
instructive to furnish a direct diagrammatic proof for the degree 2 part 
(degrees 0 and 1 are trivial). 
Recall that the degree 2 part of $\ZK(\text{unknot})$ is
equal to $\upchi(\frac{_1}{^{48}}\text{\Large$\twowheel$})$. (This can be
deduced from the Wheels formula \eqref{eqn:BGRTformula}, but is a much more
elementary result.)

\begin{exercise}
Show that 
\[
   \tr(\frac{_1}{^2}\varrho^2) = 
    \tr(\iota\upchi(\frac{_1}{^{48}}\text{\Large$\twowheel$})) \, .
\]
\end{exercise}

\begin{remark}
Conjecture ~\ref{conj:Uek} also holds at the level of complex semisimple
Lie algebras (and this gives another proof of Theorem 
~\ref{thm:ZekUnknotLieAlg}). More precisely, the maps $\Targh$ intertwine 
$\Aar_{\ek}$ with $\Ugh_{\ek}$, which is the quasitriangular Hopf algebra 
obtained by twisting $\Ugh_{\kz}$ by the image of 
$J^{-1}\in \lowupper{\Aar}{\bt 2}$ 
under $\Targh$. Our claim then is that the $u$-element of $\Ugh_{\ek}$ is
equal to $e^{\hbar \varrho_{\g}} e^{-\hbar C_{\g}/2}$.

Indeed, by ~\cite[Theorem 3.7]{EtingofKazhdan:QuantizationVI}, $\Ugh_{\ek}$
is isomorphic to the Drinfel'd-Jimbo quantum group $\Uhg$. 
Moreover\footnote{We thank P.~Etingof for clarifying this point.}, the 
isomorphism from $\Uhg$ to $\Ugh_{\ek}$ has the property that
the elements $H_i$ are mapped to $h_i$, for $i=1,\ldots, n$.
Thus, the structure $\Ugh_{\ek}$ is an instance of the structure $\Ugh_{\qg}$. 
(Recall that the structure $\Ugh_{\qg}$ depends on the choice of an isomorphism
as in Theorem ~\ref{thm:Uhg-Iso-Ugh}). Our claim now follows from 
\eqref{eqn:uDJ}, which is a consequence of the known value of the $u$-element
of $\Uhg$. 
\end{remark}

\begin{remark}
An obvious corollary of Conjecture ~\ref{conj:ZekUnknot} would be that
$\tr(e^{\varrho}) = \tr(e^{-\varrho})$. This equality, which is a 
(conjectural) consequence of a simple topological fact
(reversing the orientation of the unknot does not change its isotopy class)
seems to be non-trivial in itself. We note that it is equivalent to the
series of equalities 
\[
	\tr(\varrho^{k}) = 0 \, , \quad k = 1, 3, 5, \ldots
\]
\end{remark}

\subsection{Summary of spaces, structures, and link invariants}
\label{_sect:SummaryQtrqH}

Some of the ribbon quasitriangular quasi-Hopf structures we have encountered,
together with the associated framed link invariants, are summarized in the 
diagram in Figure ~\ref{fig:SummaryQtrqH}. 
In this diagram, $\LF$ is the set of (isotopy classes of) framed links, and 
the target spaces of each invariant are indicated by its `base target space' 
(that is, its values on framed links with $n$ components are in the $n$th 
tensor power of the base target space). The wavy lines stand for twist 
equivalences. 

\begin{figure}[htpb]
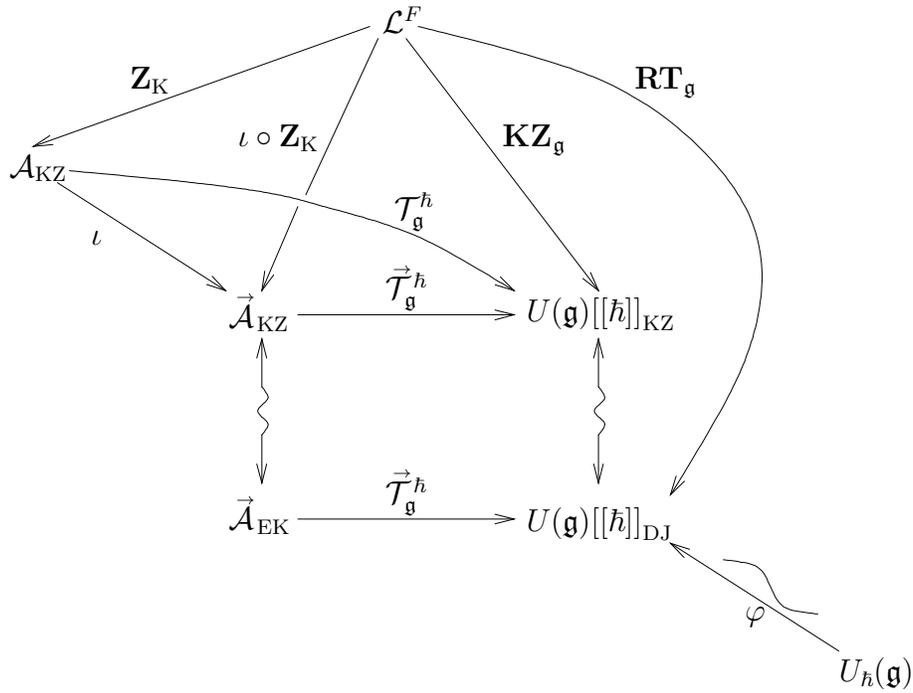

\[ 
	
  \printname{QtrqH}
  \setlength{\unitlength}{1.0\standardunitlength}
  \begin{array}{c}  \hspace{-1.7mm}
    \raisebox{-8pt}{\input draws/QtrqH.tex }
      \hspace{-1.9mm}
  \end{array}
 
\]
\caption{Ribbon quasitriangular quasi-Hopf structures 
	and associated link invariants.}
\lbl{fig:SummaryQtrqH}
\end{figure}

\subsection{The Polyak spaces}
\label{_sect:Arrow}

Let $\M$ be a $1$-dimensional skeleton. We denote by $\Dcar(\M)$ the collection
of $\M$-based diagrams whose body consists of directed chords.
\begin{definition}
The space $\Pcar(\M)$ is the quotient of the space formally spanned by 
$\Dcar(\M)$ (over the base field $\k$), by the {\em 6T relations} (compare
Exercise ~\ref{exer:DgrmCYBE}):
\begin{equation}
\label{eqn:6T}
  
  \printname{6T1}
  \setlength{\unitlength}{0.5\standardunitlength}
  \begin{array}{c}  \hspace{-1.7mm}
    \raisebox{-8pt}{\input draws/6T1.tex }
      \hspace{-1.9mm}
  \end{array}
 - 
  \printname{6T2}
  \setlength{\unitlength}{0.5\standardunitlength}
  \begin{array}{c}  \hspace{-1.7mm}
    \raisebox{-8pt}{\input draws/6T2.tex }
      \hspace{-1.9mm}
  \end{array}
 = 
  
  \printname{6T3}
  \setlength{\unitlength}{0.5\standardunitlength}
  \begin{array}{c}  \hspace{-1.7mm}
    \raisebox{-8pt}{\input draws/6T3.tex }
      \hspace{-1.9mm}
  \end{array}
 - 
  \printname{6T4}
  \setlength{\unitlength}{0.5\standardunitlength}
  \begin{array}{c}  \hspace{-1.7mm}
    \raisebox{-8pt}{\input draws/6T4.tex }
      \hspace{-1.9mm}
  \end{array}
 + 
  
  \printname{6T5}
  \setlength{\unitlength}{0.5\standardunitlength}
  \begin{array}{c}  \hspace{-1.7mm}
    \raisebox{-8pt}{\input draws/6T5.tex }
      \hspace{-1.9mm}
  \end{array}
 - 
  \printname{6T6}
  \setlength{\unitlength}{0.5\standardunitlength}
  \begin{array}{c}  \hspace{-1.7mm}
    \raisebox{-8pt}{\input draws/6T6.tex }
      \hspace{-1.9mm}
  \end{array}
 \, .
\end{equation}
\end{definition}
The spaces $\Pcar(\M)$ are graded by the number of chords in a diagram, and
are often interchanged with their graded completions. We put 
$\Pcar = \Pcar(\ua)$. The study of these spaces is motivated by the theory 
of finite-type invariants of {\em virtual knots} (see 
~\cite{GoussarovPolyakViro:VirtualKnots}).

\smallskip
Equation ~\ref{eqn:6T} is a diagrammatic form of the CYBE \eqref{eqn:CYBE}.
This leads, along the lines of Section ~\ref{__sect:DrctdJacobiToLie}, to
the construction of a map $\Tarr:\Pcar \ra U(\a)$, where $\a$ is a 
quasitriangular Lie bialgebra and $r$ is its $r$-matrix. (More generally,
there are maps $\Tarr:\Pcar(\M) \ra L(\M)$, where $L(\M)$ is the space
of $\M$-based $\a$-tensors.) 

Recall from Section ~\ref{__sect:ChordToJacobi} that $\A^{c}(\M)$ is the
chord diagram version of $\A(\M)$ (it is isomorphic to $\A'(\M)$. The maps 
$\iota: \A^{c}(\M) \ra \Pcar(\M)$ are defined by summing over all possible
assignments of directions to the chords of a diagram in $\A^{c}(\M)$. It is
not hard to verify that these maps are well-defined, namely that 4T relations
are mapped to zero. They are conjectured in ~\cite{Polyak:ArrowDgrmAlg} to be 
{\em injective} (compare Conjecture ~\ref{conj:IotaIsInjective}). If $\tg$ is 
the ambient Lie bialgebra of a Manin triple and $r$ is the natural $r$-matrix 
of $\tg$, then we have the following commutative diagram:
\[
\begin{CD}
\A^{c} @>{\iota}>> \Pcar \\
@V{\T_{\tg}}VV  @VV{\Tarr}V \\
U(\tg)^{\tg} @>{i}>> U(\tg) \, .
\end{CD}    
\]

\smallskip
Next comes the transition to uni-trivalent diagrams. Recall that a directed
graph is {\em acyclic} if it contains no directed cycles. A directed Jacobi
diagram is acyclic if its body graph is. We denote by $\Daar(\M)$ the 
collection of acyclic directed Jacobi diagrams based on a skeleton $\M$.
\begin{definition}
The space $\Par(\M)$ is the quotient of the formal $\k$-span of $\Daar(\M)$ by 
directed NS, AS, IHX, and STU relations (as in Definition ~\ref{def:Aar}), in 
which only acyclic diagrams participate.  
\end{definition} 
As usual, $\Par(\M)$ is graded by half the number of vertices in a diagram, and
we denote by $\Par$ the space $\Par(\ua)$.

For $\M$ a $1$-dimensional skeleton, The injection 
$i: \Dcar(\M) \ra \Daar(\M)$ 
descends to a well-defined map $i: \Pcar(\M) \ra \Par(\M)$, as a 6T relation 
is a consequence of directed STU relations. The following result (compare 
Proposition ~\ref{prop:Achord-eq-Ajacobi}) is due to Ohtsuki and Polyak.
\begin{theorem}
The map $i: \Pcar \ra \Par$ is an isomorphism.
\end{theorem}

We also have an injection $j:\Daar(\M) \ra \Dar(\M)$, which induces a 
map $j: \Par(\M) \ra \Aar(\M)$. This map is not onto (for example, the
degree $1$ diagram $\TPD$ is easily seen not to lie in the image of
$j:\Par \ra \Aar$.) I do not know whether it is injective or not.

\smallskip
Since the algebra\footnote{All algebraic structures on $\Par$-spaces are 
defined as for $\Aar$-spaces.} $\Par$ has (unlike $\Aar$) an equivalent
description in terms of directed chord diagrams,  
it would be interesting to carry out some of the constructions of Section 
~\ref{_sect:EKtwist} within $\Par$-spaces. A naive imitation is due to fail: 
If we try to define Verma modules by an acyclic version of Definition 
~\ref{def:Aar-pm}, then the map $\phi: \Par \ra \VM_+ \bt \VM_-$ (induced by
the coproduct) is not an isomorphism already in degree $1$.  

Nevertheless, the following question makes sense.
\begin{question}
Is there a solution of the QYBE in $\lowupper{\Par}{\bt 2}$ which is a 
quantization of the diagrammatic $r$-matrix 
$
  \printname{LeftArrow}
  \setlength{\unitlength}{0.25\standardunitlength}
  \begin{array}{c}  \hspace{-1.7mm}
    \raisebox{-8pt}{\input draws/LeftArrow.tex }
      \hspace{-1.9mm}
  \end{array}
$\hspace{-3pt}? More specifically, 
does the $R$-matrix $R_{\ek}$ of $\Aar_{\ek}$ belong to the image of 
$j: \lowupper{\Par}{\bt 2} \ra \lowupper{\Aar}{\bt 2}\,$?
\end{question}

\begin{remark}
The results of Etingof and Kazhdan in ~\cite{EtingofKazhdan:QuantizationII} 
(see also ~\cite[Lecture 20]{EtingofSchiffmann:LectQuantumGrps}) seem to hint 
at a positive answer, but I do not know whether their algebraic constructions 
and arguments can be translated into a diagrammatic framework.
\end{remark}

%% file: draws/A2hUnit.tex
%
\begingroup\makeatletter\ifx\SetFigFont\undefined%
\gdef\SetFigFont#1#2#3#4#5{%
  \reset@font\fontsize{#1}{#2pt}%
  \fontfamily{#3}\fontseries{#4}\fontshape{#5}%
  \selectfont}%
\fi\endgroup%
{\renewcommand{\dashlinestretch}{30}
\begin{picture}(1086,1431)(0,-10)
\thicklines
\path(93,33)(93,1383)
\blacken\path(153.000,1203.000)(93.000,1383.000)(33.000,1203.000)(153.000,1203.000)
\path(993,33)(993,1383)
\blacken\path(1053.000,1203.000)(993.000,1383.000)(933.000,1203.000)(1053.000,1203.000)
\end{picture}
}